\documentclass[10pt,reqno,letterpaper]{article}

\usepackage[
  letterpaper,
  margin=1in,
  includehead,includefoot
]{geometry}

\usepackage[utf8]{inputenc}
\usepackage[T1]{fontenc}
\usepackage{lmodern}
\usepackage{amsmath,amssymb,amsthm,mathtools,bm}

\usepackage{graphicx}
\usepackage{booktabs}
\usepackage[table]{xcolor}
\usepackage{makecell}
\usepackage{caption}
\usepackage{subcaption}
\usepackage{float}

\usepackage{tikz}
\usetikzlibrary{positioning,shapes.geometric,decorations.pathreplacing,arrows.meta,calc}
\usepackage{pgfplots}
\usepgfplotslibrary{groupplots}
\tikzset{font={\fontsize{7pt}{12}\selectfont}}

\usepackage[ruled,vlined]{algorithm2e}
\usepackage{algpseudocode}
\usepackage{matlab-prettifier}

\usepackage{microtype}
\usepackage{url}
\usepackage{hyperref}

\usepackage{enumitem}
\setlist{itemsep=2pt, topsep=4pt}

\usepackage{comment}
\excludecomment{draftnote}
\usepackage{soul}
\usepackage[normalem]{ulem}
\newtheorem{definition}{Definition}

\newtheorem{corollary}{Corollary}

\newtheorem{assumption}{Assumption}

\title{Stable Mesh-Free Variational Radial Basis Function Approximation\\
for Elliptic PDEs and Obstacle Problems}

\author{Phuong Dong Le, Giang Tran, and Hans De Sterck\\
Department of Applied Mathematics, University of Waterloo\\
Waterloo, ON N2L 3G1, Canada\\
\texttt{\{pdle,giang.tran,hans.desterck\}@uwaterloo.ca}
}
\date{}

\usetikzlibrary{positioning, decorations.pathreplacing, arrows.meta}
\usetikzlibrary{calc}

\usepackage{enumitem}
\usepackage{subcaption} 

\setlist{itemsep=2pt, topsep=4pt}

\usepackage{float}
\usepackage{comment}
\excludecomment{draftnote}
\usepackage{float}
\usepackage[font=small,labelfont=bf]{caption}
\begin{document}
\maketitle
\begin{abstract}
We present a comprehensive study of radial basis function (RBF) approximations for elliptic and obstacle-type boundary value problems under a variational formulation. Our focus is on practical accuracy, robustness and efficiency. To address ill-conditioning in dense systems, we apply truncated singular value decomposition (TSVD) and investigate its effect on stability and accuracy trade-offs. Numerical experiments report benchmarks on accuracy and show fast error decay. We investigate the trade-off between approximation and truncation errors for practical settings for the number of basis functions, the oversampling ratio and the truncation threshold. In comparison with other methods, RBF variational solvers deliver high accuracy at similar or lower cost for boundary value problems. 
\end{abstract}


\section{Introduction}\label{sec:introduction}

\hspace{5mm}Radial Basis Function (RBF) methods have been used as a tool in mesh-free approximation to solve the boundary value problem, offering high-order accuracy and flexibility \cite{buhmann2003radial, wendland2005scattered}. A common kernel choice in RBF methods is the Gaussian radial basis function because of its spectral accuracy \cite{buhmann2003radial}. However, this choice of basis function often leads to severe ill-conditioning of the system matrix when the shape parameter becomes small or when the number of centers is large. The system matrix then amplifies round-off errors which makes the numerical system unstable. This ill-conditioning has been studied in the literature \cite{schaback1995error, fasshauer2007meshfree}, and recent results have proposed regularization techniques to address it.

In \cite{adcock2024stable}, the authors demonstrate that accurate and stable least-squares approximations on a bounded domain can still be achieved from ill-conditioned RBF matrices in the linear scaling regime: when the scale parameter, \(b\), is chosen proportional to the number of basis functions, \(N\) as $(b\sim c \, N)$ for Gaussian RBFs, \cite{adcock2024stable,adcock2015tsvd, adcock2019frames} propose an overdetermined stabilized \emph{discrete least-squares} formulation on bounded domains, rather than interpolation. Their setting combines (i) placing a portion of RBF centers outside the problem domain, $\Omega$, to improve accuracy, and (ii) oversampling by using more least-squares evaluation points than centers, which improves numerical stability even when the associated RBF matrix is severely ill-conditioned. They further quantify the regime and provide explicit guidance for parameter choices, motivated by a redundancy/frame-based analysis. These results inspire our methodology in this work. Specifically, we are interested in expanding this technique to solve the variational problems, including the partial differential equations: elliptic and free-boundary value problems. 

We start by formulating elliptic boundary value problems as variational problems, which represent the solution as the minimizer of an energy functional. Similarly to \cite{yu2018deep}, we also incorporate the boundary conditions to the objective function by adding a suitable penalty term. To solve the corresponding optimization problem, we approximate the unknown solution using Gaussian radial basis functions, where some centers are outside of the domain. By approximating this minimizer using RBF networks, we project the infinite-dimensional variational problem into a finite-dimensional subspace spanned by radial basis functions. We examine various choices of the shape parameters, including optimal values suggested in \cite{adcock2024stable} for the one-dimensional problems, to obtain numerical solutions with high accuracy.

In addition to the linear elliptic problems, we extend our methodology to a class of free boundary value problems which impose an inequality constraint between solution and obstacle function. Such problems appear in many applications including contact mechanics and American financial options. The classical obstacle problem by \cite{friedman1982variational, glowinski1984numerical, ros2018obstacle} seeks a solution $u$ that satisfies the following constraints: 
\begin{equation}
u(x) \geq \psi(x), \quad - \Delta u(x) \geq 0, \quad \text{and} \quad \min \{ - \Delta u(x), u(x) - \psi(x) \} = 0 \quad \text{in} \quad \Omega, \quad u = g  \quad \text{on} \quad \partial \Omega 
\nonumber
\end{equation}
where $\psi(x)$ is the obstacle function. 

By formulating the obstacle problem as a constrained optimization problem and adopting an $\ell^{1}$-penalty approach based on \cite{tran20151, zosso2015fast}, we impose the constraint $u \geq \psi$ via a non-smooth penalty term, and for the boundary condition we apply the penalization via parameter $\beta$ such that: 
\begin{equation}
\min \limits_{u} \Big\{ \frac{1}{2} \int_{\Omega} \|  \nabla u \|_{2}^{2} \, d\vec{x} + \mu \int_{\Omega} (\psi -u)_{+}   \, d\vec{x}  + \frac{\beta}{2} \int_{\partial \Omega} (u - g)^2 \, d \vec{s} \Big\}
\nonumber
\end{equation}
where $\mu, \beta > 0$ are penalty parameters and $(\cdot)_{+} = \max\{ \cdot, 0 \}$. To solve the corresponding optimization problem, similar to the elliptic boundary value problem, we also approximate the unknown solution of the variational obstacle problem as a linear combination of Gaussian radial basis functions. Motivated by \cite{adcock2024stable}, we  sample the centers of those radial basis functions on an extended domain \(\overline{\Omega}_{T}\) where \(T > 0\) is some expansion domain factor. Our numerical experiments in Section~\ref{sec:numerical-experiments} indicate that the proposed method can deliver high accuracy solutions for free-boundary value problems. The discretized variational formulation leads to linear systems that are ill-conditioned particularly for small values of the shape parameter in the radial basis functions. For this reason, we combine an $\ell^1$ optimization method with truncated singular value decomposition (TSVD) to enforce the inequality constraint and stabilize the solution of the linear system. This formulation allows for robust numerical approximation and efficiency for the variational inequality. We use the alternating direction method of multipliers (ADMM) in \cite{boyd2011distributed} to solve the optimization problem. 

The main contribution of this paper is a mesh-free variational framework for elliptic and obstacle boundary value problems, combining stabilized RBF approximations with an $\ell^1$-based optimization method. Specifically:

\begin{itemize}
    \item We propose a mesh-free RBF approximation, and using TSVD \cite{adcock2019frames, adcock2024stable} and for $\ell_1$ optimization ADMM. 
    \item We use a penalty method for elliptic BVPs and obstacle problems as in \cite{tran20151, zosso2015fast} and \cite{zhang2011unified}. 
    \item We present numerical experiments validating the convergence, stability, and accuracy of our proposed method for elliptic and obstacle-type boundary value problems.
\end{itemize}

The remainder of our research paper is organized as follows.
Section~\ref{sec:problem-statement} states the problem setting and introduces the radial basis function networks representation, the assumptions on center distributions, and the variational equations for linear elliptic and obstacle-type boundary value problems. Section~\ref{sec:methodology} presents our methodology, including variational formulations for boundary value problems, the TSVD method for solving linear systems, and the ADMM framework used for non-smooth $\ell^1$ penalization. Section~\ref{sec:numerical-experiments} presents numerical experiments that assess. Finally, Section~\ref{sec:conclusion} summarizes the main findings and discusses limitations and directions for future work. 

\section{Problem Statement}\label{sec:problem-statement}

\hspace{5mm}In this section, we introduce the class of elliptic boundary value problems under consideration and a mesh-free approximation framework based on RBF approximation. The objective is to seek weak solutions of elliptic and obstacle-type partial differential equations on a bounded domain $\Omega \subseteq \mathbb{R}^{d}$ by a finite-dimensional approximation $\hat{u}$. The unknown coefficients of this representation will be determined through variational equations, which leads to stabilized linear least-squares systems or non-smooth constrained optimization problems.

\subsection{Variational Functionals for Boundary Value Problem}\label{sec:variational-function-bvp}

We are interested in solving elliptic and free-boundary partial differential equations with inequality constraints using a variational framework and penalty terms. Let $\Omega \subset \mathbb{R}^{d}$ be a bounded open domain  with a smooth boundary $\partial \Omega$. First, we start with the Poisson boundary value problem, which is defined as follows: 
\begin{equation}
\begin{cases}
- \Delta u = f( \vec{x} ), & \forall \vec{x} \in \Omega \\ 
u(\vec{x}) = g( \vec{x} ), & \vec{x} \in \partial \Omega. 
\end{cases}
\label{poisson-boundary-value-problem}
\end{equation}
The variational problem corresponding to Eq.~\eqref{poisson-boundary-value-problem} is given by: 
\begin{equation}
\frac{1}{2}\int_{\Omega} \| \nabla u \|_{2}^{2} d \, \vec{x} - \int_{\Omega} f \, u \, d \, \vec{x}, \quad \text{such that} \quad u(\vec{x}) = g(\vec{x}) \quad \forall \vec{x} \in \partial \Omega.  
\label{poisson-weak-problem-zeidler}
\end{equation}
Following \cite{yu2018deep}, we apply a penalty method to impose the boundary condition and consider minimizing a modified functional:
\[\min \limits_{u}\mathcal J_{\beta, \mathrm{poisson}}(u),\]
where 
\begin{equation}
\mathcal J_{\mathrm{poisson},\beta}(u) = \frac{1}{2} \int_{\Omega} \| \nabla u \|_{2}^{2} \, d\vec{x} - \int_{\Omega} f \, u \, d\vec{x} + \frac{\beta}{2} \int_{\partial \Omega} (u - g)^{2} \, d\vec{s},\quad (\beta > 0).
\label{weak-form-poisson-problem}
\end{equation}
Similarly, we consider the classical obstacle problem which is defined as follows: 
\begin{equation}
\begin{cases}
- \Delta u(\vec{x}) \geq 0, & \forall \vec{x} \in \Omega \\ 
u(\vec{x} ) = g(\vec{x}), & \forall \vec{x} \in \partial \Omega \\ 
u(\vec{x}) \geq \psi (\vec{x}), & \text{in} \, \Omega \\ 
(- \Delta u) \, (u - \psi) = 0, & \text{in}\, \Omega
\end{cases}
\label{obstacle-pde-bvp}
\end{equation}
where $\psi$ is a given obstacle function.
The following $\ell_{1}$-penalty approach provides a powerful alternative method to formulate the obstacle problem~\eqref{obstacle-pde-bvp} with \(g = 0\):
\[
\min\limits_u\frac{1}{2} \int_{\Omega} \| \nabla u \|_{2}^{2} \, d \vec{x} + \mu \int_{\Omega} (\psi - u)_{+} \, d \vec{x},
\]
where $\mu>0$ is the penalty coefficient. 

According to a result in \cite{tran20151}, for sufficiently large $\mu > 0$, the minimizer coincides with the solution of Eq.~\eqref{obstacle-pde-bvp}. Adding the boundary condition as in Eq.~\eqref{weak-form-poisson-problem}, we obtain the optimization problem: 
\[
\min\limits_u \mathcal{J}_{\mathrm{obstacle},\mu, \beta }(u),
\]
where
\begin{equation}
\mathcal{J}_{\mathrm{obstacle}, \mu, \beta}(u) = \frac{1}{2} \int_{\Omega} \| \nabla u \|_{2}^{2} \, d \vec{x} + \mu \int_{\Omega} ( \psi - u)_{+} d \vec{x}  + \frac{\beta}{2} \int_{\partial \Omega} (u - g)^{2}  \, d \vec{s}.
\label{weak-form-obstacle-problem}
\end{equation}
The boundary penalty term, $\int_{\partial \Omega} (u - g)^{2} d \vec{s}$, softly enforces the Dirichlet condition $u = g$ on $\partial \Omega$. Although it is possible to construct RBFs or modified ansatz functions that satisfy the boundary conditions exactly, such constructions often increase complexity and reduce flexibility for a generic function $g$. In contrast, enforcing the boundary via a penalty term retains the generality of the basis. 
To solve the optimization problems in Eqs.~\eqref{weak-form-poisson-problem} and \eqref{weak-form-obstacle-problem}, we will use a collocation method and approximate the solution using radial basis functions. In the following, we recall the radial basis function representation and related assumptions.

\subsection{Radial Basis Function Approximation}\label{sec:radial-basis-function-approximation}

\hspace{5mm}To formulate boundary value problems (BVPs), we consider a bounded open domain $\Omega \subset \mathbb{R}^{d}$ with a smooth boundary $\partial \Omega$. We are interested in finding a function $u: \overline{\Omega} \to \mathbb{R}$ in some Sobolev space satisfying a prescribed differential equation in the domain and adhering to the boundary conditions. In this setting, we seek an approximate solution constructed by a linear combination of \(N\) radial basis functions centered at points $\{\vec{c}_{i} \}_{i = 1}^{N}$ where each $\vec{c}_{i} \in \mathbb{R}^{d}$. For the sake of completeness, we recall the definition and some important properties of radial basis functions (see also \cite{powell1987radial, micchelli1984interpolation, fasshauer2007meshfree}).

\begin{definition} \label{definition-radial-basis-kernel}
A function $\Phi : \mathbb{R}^{d} \to \mathbb{R}$ is called \emph{radial} when it depends on the distance of the argument from a center $\vec{c}\in\mathbb{R}^d$; equivalently, there exists a univariate mapping $\phi: [ 0, \infty) \to \mathbb{R}$ such that
\[
\Phi(\vec{x}) = \phi(r),\quad \text{where}\quad r = \| \vec{x} - \vec{c} \|  
\]
and $\| \cdot \|$ denotes a norm on $\mathbb{R}^{d}$. 
\end{definition}

Let $\phi: [0, \infty) \to \mathbb{R}$ be a well-defined basis function. Given a set of $N$ center vectors $\{\vec{c}_i\}_{i=1}^N$ in $\mathbb{R}^d$, a radial basis function approximation $\hat{u}: \overline{\Omega} \to \mathbb{R}$ of a function $u$ is defined as:
\begin{equation}
\hat{u}(\vec{x}) \coloneqq \frac{1}{\sqrt{N}} \sum_{i = 1}^{N} w_{i} \phi( \| \vec{x} - \vec{c}_{i} \|_{2} ),
\label{eq:main-approximation-form-rbf}
\end{equation}
where $\vec{x} \in \overline{\Omega}$ and the weight coefficient vector is given by $\vec{w}= (w_1,\ldots,w_N)^{\top}\in\mathbb{R}^N$. The scaling factor $\frac{1}{\sqrt{N}}$ normalizes the formulation with respect to the number of centers. 

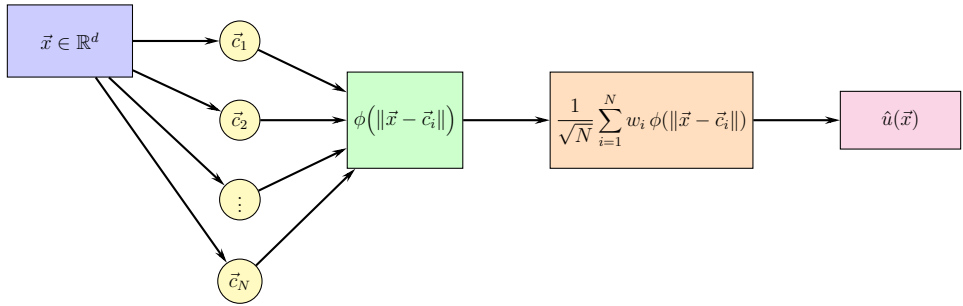
\begin{figure}[H]
\centering
    \begin{tikzpicture}[
  scale=0.64, every node/.style={transform shape,font=\large},
  node distance = 0.8cm and 1.0cm,
  input/.style  ={rectangle, draw, fill=blue!20,
                  minimum width=26mm, minimum height=15mm, align=center},
  center/.style ={circle,    draw, fill=yellow!30,
                  minimum size=8mm,  align=center},
  phi/.style    ={rectangle, draw, fill=green!20,
                  minimum width=24mm, minimum height=20mm, align=center},
  weight/.style ={rectangle, draw, fill=orange!25,
                  minimum width=30mm, minimum height=20mm, align=center},
  output/.style ={rectangle, draw, fill=magenta!20,
                  minimum width=25mm, minimum height=12mm, align=center},
  arrow/.style  ={-{Stealth[length=2mm,width=1mm]}, thick}
]

\node[input] (x) {$\vec{x}\in\mathbb{R}^d$};
\node[center, right=1.8cm of x] (c1) {$\vec{c}_1$};
\node[center, below=of c1]      (c2) {$\vec{c}_2$};
\node[center, below=of c2]      (vd) {\vdots};
\node[center, below=of vd]      (cN) {$\vec{c}_N$};

\node[phi, right=1.8cm of c2] (phi)
      {$\phi \bigl(\|\vec{x}-\vec{c}_i\|\bigr)$};

\node[weight, right=1.8cm of phi] (sum)
      {$\displaystyle\frac{1}{\sqrt{N}}\sum_{i=1}^N w_i\,\phi(\| \vec{x} - \vec{c}_{i} \| )$};

\node[output, right=1.8cm of sum] (u) {$\hat{u}(\vec{x})$};

\foreach \c in {c1,c2,vd,cN} {
    \draw[arrow] (x) -- (\c);
    \draw[arrow] (\c) -- (phi);
}
\draw[arrow] (phi) -- (sum);
\draw[arrow] (sum) -- (u);

\end{tikzpicture}
    \caption{Illustration of a RBF network. A single input vector $\vec{x} \in \mathbb{R}^d$ is combined with a set of center points $\vec{c}_i \in \mathbb{R}^d$, $i = 1, \ldots, N$, to evaluate radial kernels $\phi (\| \vec{x} - \vec{c}_i\|)$. The resulting basis function values $\phi$'s are combined linearly using weights $w_i$ to obtain the approximation $\hat{u}(\vec{x}) = \frac{1}{\sqrt{N}}\sum_{i=1}^N w_i \phi(\| \vec{x} - \vec{c}_{i} \|)$.}
    \label{fig:rbf-architecture}
\end{figure}

The universal approximation theorem in \cite{park1991universal} establishes that a single hidden layer with fixed-shape radial basis units can uniformly approximate any integrable function on a compact subset of $\mathbb{R}^{d}$. Mesh-free approximation frameworks use the smoothness and locality properties of radial basis functions to approximate solutions \cite{wendland2005scattered, fasshauer2007meshfree}. The coefficients in the representation are determined by considering residuals of the differential operator and boundary conditions using collocation or variational methods. A one-layer RBF network is illustrated in Figure~\ref{fig:rbf-architecture}.
 
Table~\ref{tab:rbf-types} summarizes some of the most popular radial basis kernels, including global functions and compactly supported functions. 

\begin{table}[H]
\centering
\renewcommand{\arraystretch}{1.4}
\begin{tabular}{| c |c c c | c }
\hline
\textbf{Name} & \textbf{Formula $\phi(r)$} & \textbf{Smoothness} & \textbf{Locality} \\
\hline
Gaussian       & $e^{-(b r)^2}$                    & $C^\infty$         & Global \\
Multiquadric   & $\sqrt{1 + (b r)^2}$              & $C^\infty$         & Global \\
Inverse Quad  & $\frac{1}{1 + (b r)^2}$           & $C^\infty$         & Global \\
Inverse MQ     & $\frac{1}{\sqrt{1 + (b r)^2}}$    & $C^\infty$         & Global \\
Wendland $C^2$ & $\chi_{[0,1]}(br) (1 - br)^4 (4br + 1)$        & $C^2$      & Compact \\
Wendland $C^4$ & $\chi_{[0,1]}(br) (1 - br)^6 (35(br)^2 + 18br + 3)$ & $C^4$ & Compact \\
\hline
\end{tabular}
\caption{Examples of radial basis functions with the shape parameter \(b > 0\). Here, the indicator function, $\chi_{[0,1]}(br)$, equals 1 when $br \leq 1$, and 0 otherwise.}
\label{tab:rbf-types}
\end{table}

Each RBF in the network is activated when the input $\vec{x}$ lies near its corresponding center $\vec{c}_{i}$. We now outline the assumptions that we make when using RBFs to approximate solutions of PDEs in this paper. 

\begin{assumption}
In numerical experiments, we normally use the Gaussian kernel defined in Table~\ref{tab:rbf-types}, 

\[
\phi(r) = e^{-(br)^{2}}.
\]
In some numerical tests, we state the use of different kernels. 
\label{assumption-1-core}
\end{assumption}

\begin{assumption}[Fill distance, see \cite{fasshauer2007meshfree,narcowich2006sobolev}]
The center points are sampled from a distribution such that $\{\vec{c}_{1}, \dots, \vec{c}_{N} \} \subset \overline{\Omega}$. Then, the fill distance property states that
\[
h_{\Omega} \coloneq \sup \limits_{ \vec{x} \in \Omega } \min \limits_{1 \leq i \leq N} \, | \vec{x} - \vec{c}_{i} | \rightarrow 0, \quad \text{as}\quad N \rightarrow \infty. 
\]
\label{assumption-2-core}
\end{assumption}
\begin{assumption}[Separation distance, see \cite{fasshauer2007meshfree, narcowich2006sobolev}] For a set of center points $\{ \vec{c}_{1}, \dots, \vec{c}_{N}  \} \subset \overline{\Omega} $, the separation distance $q_{\Omega}$ satisfies 
\[
q_{\Omega} \coloneq \frac{1}{2} \, \min \limits_{\substack{i, j = 1,\dots, N  \\ i \neq j   }} \| \vec{c}_{i} - \vec{c}_{j} \|_{2} \rightarrow 0 , \quad \text{as}\quad N \rightarrow \infty.
\]
\label{assumption-3-core}
\end{assumption}

\begin{corollary}[See Corollary 10.13 in \cite{wendland2005scattered} and Section 2.1 in \cite{narcowich2006sobolev}]
\label{assumption-4-frame}
Let $\phi(r)=e^{-(br)^2}$ as in Assumption~\ref{assumption-1-core}, and choose center points $C=\{\vec c_i\}_{i=1}^N\subset\overline{\Omega}$ uniformly distributed as in Assumptions~\ref{assumption-2-core}–\ref{assumption-3-core}. Then there exist constants $B_{1}, B_{2} >0$, independent of $N$, such that for all $ \vec{\omega} \in\mathbb{R}^N$,
\[
B_{1} \sum_{i=1}^N |\omega_i|^2
\;\le\;
\Big\|\sum_{i=1}^N \omega_i \,\phi(\| \vec{x}-\vec c_i\|_{2} )\Big\|_{W^{s,2}(\Omega)}^{2}
\;\le\;
B_{2} \sum_{i=1}^N |\omega_i|^2.
\]
Thus $\{\phi(\| \vec{x}-\vec c_i \|_{2})\}_{i=1}^N$ forms a frame for its closed span in $W^{s,2}(\Omega)$.
\end{corollary}

In practice, we place centers in an expanded domain, see Section~\ref{sec:numerical-experiments}. The first assumption states we use Gaussian kernels, which meet the requirements of the universal approximation property of RBFs. The Assumptions ~\ref{assumption-2-core} and ~\ref{assumption-3-core} guide our construction of center points. While the fill distance $h_{\Omega}$ measures the largest distance from any point in $\Omega$ to the nearest center, the separation distance $q_{\Omega}$ measures the minimal spacing between centers.
Corollary~\ref{assumption-4-frame} states that using our radial kernel from Assumption~\ref{assumption-1-core} results in a frame sequence for the closed span $\{ \phi(\| \vec{x} - \vec{c}_{i} \|_{2}) \}_{i = 1}^{N}$, giving uniform equivalence between the coefficient and function norms. In \cite{wendland2005scattered}, the author expands this from the Gaussian RBF to any compactly supported function $\phi \in L^{1}(\mathbb{R}^{d}) \cap C(\mathbb{R}^{d})$. This would guaranty that the kernels in Table~\ref{tab:rbf-types}, including both global and compactly supported functions, to satisfy the property in Corollary~\ref{assumption-4-frame}. The property in Corollary~\ref{assumption-4-frame} is to allow the application of TSVD as in \cite{adcock2015tsvd, adcock2019frames}. 

In this work, we approximate the solution by the representation $\hat{u}$ in Eq.~\eqref{eq:main-approximation-form-rbf}. The weighted coefficients are optimized using TSVD for linear elliptic boundary value problems, and the ADMM is used for obstacle-type problems.

\section{Methodology}\label{sec:methodology}

\hspace{5mm}We now describe the numerical framework used to approximate the linear elliptic and free-boundary value problems. Specifically, we aim to find the minimizers of the corresponding variational PDEs in Eq. \eqref{weak-form-poisson-problem} and Eq.\eqref{weak-form-obstacle-problem}. We approximate the minimizers by radial basis functions as in Eq.\eqref{eq:main-approximation-form-rbf}. To approximate the integrals in the variational PDEs, we select $m$ collocation points from the domain $\overline{\Omega} = \Omega \cup \partial \Omega$, including $m_{\Omega}$ interior points $\vec{x}^{(\Omega)}_{1},\ldots,\vec{x}^{(\Omega)}_{m_\Omega}$  and $m_{\partial \Omega}$ boundary points $\vec{x}^{(\partial \Omega)}_{1},\ldots,\vec{x}^{(\partial \Omega)}_{m_{\partial \Omega}}$. The corresponding discrete form of Eq.~\eqref{weak-form-poisson-problem} for Poisson's equation is as follows:

\begin{align}
\widehat{\mathcal{J}}_{\mathrm{poisson},\beta}(\vec{x}) 
&= \frac{| \Omega |}{m_{\Omega}}\sum_{j = 1}^{m_{\Omega}} \left( 
    \frac{1}{2} \| \nabla  \hat{u}(\vec{x}^{(\Omega)}_{j}) \|_{2}^{2} 
    - f(\vec{x}^{(\Omega)}_{j}) \hat{u}(\vec{x}^{(\Omega)}_{j}) 
\right) + \frac{\beta \, | \partial \Omega| }{2 \, m_{\partial \Omega}} 
\sum_{n = 1}^{m_{\partial \Omega}} \left( 
    \hat{u}(\vec{x}^{(\partial \Omega)}_{n}) - g(\vec{x}^{(\partial \Omega)}_{n}) 
\right)^{2}.
\label{eq:discretized-weak-form-poisson} 
\end{align}

In the same manner, the discrete form of Eq.~\eqref{weak-form-obstacle-problem} for the obstacle problem is defined as follows: 
\begin{align}
\widehat{\mathcal{J}}_{\mathrm{obstacle},\mu, \beta}(\vec{x}) 
&= \frac{| \Omega |}{m_{\Omega}} \sum_{j = 1}^{m_{\Omega}} \left( 
    \frac{1}{2} \| \nabla  \hat{u}(\vec{x}^{(\Omega)}_{j}) \|_{2}^{2} 
  +   \mu 
    \max\left( \psi(\vec{x}^{(\Omega)}_{j}) - \hat{u}(\vec{x}^{(\Omega)}_{j}), 0 \right) \right)+ \frac{\beta \, | \partial \Omega| }{2 \, m_{\partial \Omega}}
\sum_{n = 1}^{m_{\partial \Omega}} \left( 
    \hat{u}(\vec{x}^{(\partial \Omega)}_{n}) - g(\vec{x}^{(\partial \Omega)}_{n}) 
\right)^{2}.
\label{eq:discretized-weak-form-obstacle}
\end{align}
Here, $|\Omega|$ and $|\partial \Omega|$ denote the Lebesgue measure of the computational domains $\Omega, \partial \Omega$, respectively.

It remains for us to compute the explicit derivatives for the discretized variational formulation in each boundary value problem. Under Assumption~\ref{assumption-1-core} and the RBF approximation in Eq.~\eqref{eq:main-approximation-form-rbf}, we use the Gaussian basis function such that $\phi(r) = \exp(- (b \, r)^{2} ) $. In particular,  we get $\ \phi'(r) = - 2 b^{2} \, r \, \exp(- (b \, r)^{2} )  $. Thus, for $\vec{c}_{i} \in \mathbb{R}^{d}$, we compute the gradient $\nabla \hat{u} (\vec{x}_{j}^{(\Omega)} )$ as follows:  
\begin{equation}
\nabla \hat{u} (\vec{x}_{j}^{(\Omega)} )
= \frac{1}{\sqrt{N}}  \sum_{i = 1}^{N}  w_{i} \nabla \phi (\| \vec{x}_{j}^{(\Omega)} - \vec{c}_{i} \|_{2}), \quad 1 \leq j \leq m_{\Omega}. 
\label{eq:nabla-derivative-for-u}
\end{equation}

Hence, we deduce that

\begin{equation}
\| \nabla \hat{u}(\vec{x}_{j}^{(\Omega)} ) \|_{2}^{2} = \frac{1}{N} \sum_{i = 1}^{N} \sum_{l = 1}^{N} w_{i} \, w_{\ell} \big( \nabla \phi (\| \vec{x}_{j}^{(\Omega)} - \vec{c}_{i} \|_{2} ) \big)^{\top} \, \nabla \phi (\| \vec{x}_{j}^{(\Omega)} - \vec{c}_{\ell} \|_{2}), \quad 1 \leq j \leq m_{\Omega}. 
\label{eq:nabla-energy-form-for-u}
\end{equation}

The remainder of this section is divided into two parts: (i) an unconstrained optimization framework for linear elliptic boundary value problems and (ii) a constrained optimization framework for obstacle-type problems.

\subsection{Unconstrained Optimization for Elliptic Boundary Value Problems}\label{sec:unconstrained-optimization}

In this section, we derive a numerical scheme to solve the Poisson boundary value problem by minimizing the convex functional of Eq.~\eqref{eq:discretized-weak-form-poisson}. The details of the derivation are given by Appendix~\ref{appendix:subsection-derivation-optimization-elliptic}. Using the basis kernel's differentiability, we need to solve the following: 
\begin{align}
 & 
\min_{\vec{w}\in\mathbb{R}^{N}}
\widehat{\mathcal{J}}_{\text{poisson},\beta}(\vec{w}),
\nonumber
\\[4pt]
\, \text{where}\quad \widehat{\mathcal{J}}_{\text{poisson},\beta}(\vec{w}) &=
\frac{1}{2}\,\vec{w}^{\top}A^{\small(\mathrm{poisson})}_{1}\vec{w}
-\vec{f}^{\top}A_{2}\vec{w}
+\frac{\beta}{2}\,\lVert A_{3}\vec{w}-\vec{g}\rVert_{2}^{2},
\label{eq:poisson_loss}
\end{align}
and the matrices $A_1,A_2$ and $A_3$ are given by \\

\begin{equation}
(A_{1})_{i,l} = \frac{|\Omega|}{m_{\Omega} \, N } \sum_{j = 1}^{m_{\Omega}} \nabla \phi^{\top}(\| \vec{x}^{(\Omega)}_{j} - \vec{c}_{i} \|_{2} ) \, \nabla\phi(\| \vec{x}^{(\Omega)}_{j} - \vec{c}_{\ell} \|_{2}), \quad \text{for} \, \, 1 \leq i, \ell \leq N.   
\label{eq:A1-matrix-term-entry-wise-poisson}
\end{equation}

\begin{align}
A_{2} &\in \mathbb{R}^{m_{\Omega}\times N},      &\quad (A_{2})_{j,i} &= \frac{|\Omega|}{m_{\Omega}\sqrt{N}}\phi(\| \vec{x}^{(\Omega)}_{j} - \vec{c}_{i} \|_{2} ),
&\; j &= 1,\dots,m_{\Omega},\; i = 1,\dots,N,
\label{eq:A2-matrix-term-entry-wise}\\[6pt]
A_{3} &\in \mathbb{R}^{m_{\partial\Omega}\times N}, & (A_{3})_{n,i} &= \sqrt{\frac{|\partial \Omega|}{ m_{\partial \Omega}\, N} }\phi( \| \vec{x}^{(\partial \Omega)}_{n} - \vec{c}_{i} \|_{2}),
&\; n &= 1,\dots,m_{\partial\Omega}, \; i = 1,\dots,N.\label{eq:A3-matrix-term-entry-wise}
\end{align}

The vectors of the source function and boundary data are given by 

\begin{equation}
\vec{f} \coloneq \begin{bmatrix}
    f \left( \vec{x}^{(\Omega)}_{1} \right) \\ 
    \vdots \\
    f \left( \vec{x}^{(\Omega)}_{m_{\Omega}   } \right) 
\end{bmatrix} \in \mathbb{R}^{m_{\Omega}}, \quad \vec{g} \coloneq \sqrt{\frac{|\partial \Omega|}{m_{\partial \Omega}}} \, \begin{bmatrix}
   g \left( \vec{x}^{(\partial \Omega)}_{1} \right) \\ 
   \vdots \\
   g \left( \vec{x}^{(\partial \Omega)}_{m_{\partial \Omega}} \right)
\end{bmatrix} \in \mathbb{R}^{m_{\partial \Omega}}. 
\label{eq:source-boundary-data-matrix-term}
\end{equation}

The minimizer $\vec{w}^{*} $ satisfies the following equation:
\begin{equation}
 A_{\mathrm{elliptic}}\vec{w}^{*} = \vec{b}_{\mathrm{elliptic}},
 \label{eqn:ellipticPDE}
\end{equation}
where $A_{\mathrm{elliptic}} \coloneq A_{1} + \beta A_{3}^{\top} A_{3} \in \mathbb{R}^{N \times N}$ and $\vec{b}_{\mathrm{elliptic}} \coloneq A_{2}^{\top} \, \vec{f} + \beta A_{3}^{\top} \vec{g} \in \mathbb{R}^{N}$. However, the choice of the shape parameter $b$ in the kernel function $\phi $ affects the condition number of the stiffness matrix in Eq.~\eqref{eq:poisson_loss}. To address this, we compute the TSVD solution:
\begin{equation}
\vec{w}_{\tau} \coloneq A^{\dagger}_{\mathrm{elliptic},\tau} \vec{b}_{\mathrm{elliptic}} = V_{p} \,  \Sigma_{\tau}^{-1} U^{\top}_{p}  \, \vec{b}_{\mathrm{elliptic}},
\label{eq:pseudo-inverse-elliptic-tau-tsvd}
\end{equation}
where $\Sigma_{\tau}^{-1} = \mathrm{diag}(1/\sigma_{1},\cdots, 1/\sigma_{p})$ with $\sigma_{p} > \sigma_{1} \, \tau\geq \sigma_{p+1}$, $A_{\mathrm{elliptic}} = U \Sigma  V^{\top}$ is the SVD of $A_{\mathrm{elliptic}}$, and \(U_{p} \in \mathbb{R}^{N \times p}, V_{p} \in \mathbb{R}^{N \times p}\) are matrices of the first \(p\) columns of \(U\) and \(V\), respectively. It is proved in \cite{adcock2024stable} that, with a suitable choice of parameters for RBFs, the TSVD solution yields an accurate function approximation, 
\begin{equation}
\hat{u}_{\tau} (\vec{x} ) = \frac{1}{\sqrt{N}} \sum_{i = 1}^{N} w_{\tau, i} \phi (\| \vec{x} - \vec{c}_{i} \|_{2} ).  
\label{eq:approximation-solution-u-hat-tau-elliptic}
\end{equation}
Hence, the following algorithm~\ref{algorithm1-rbf-solver-for-elliptic-pde} summarizes our schematic solver for linear elliptic partial differential equations and establishes the foundation for the subsequent algorithmic scheme addressing obstacle-type problems.

\begin{algorithm}[H]\label{algorithm1-rbf-solver-for-elliptic-pde}
\small 
\caption{Stable radial basis function weak‑form solver for variational elliptic PDEs using TSVD solutions.}
\SetAlgoLined
\KwIn{$\phi(\cdot)$ (kernel function); $\{ \vec{c}_{i} \}_{i = 1}^{N} \subseteq \overline{\Omega}_{T}$ (sampled center nodes from the expansion domain);
      $\{\vec{x}^{(\Omega)}_j\}_{j=1}^{m_\Omega}$ (interior nodes); 
      $\{\vec{x}^{(\partial \Omega)} _n\}_{n=1}^{m_{\partial\Omega}}$ (boundary nodes); 
      $f(\cdot),\,g(\cdot)$ (source/boundary data); 
      $\beta>0$ (boundary penalty parameter); 
      $\tau>0$ (truncaction tolerance)}
 
\KwOut{Weight vector $\vec{w}_{\tau}\in\mathbb{R}^{N}$; Approximation $\hat{u}_{\tau}(\vec{x}) \in \mathbb{R}$.}

\BlankLine
\textbf{Assemble discrete operators}\;
\Indp
Compute $A_1\in\mathbb{R}^{N\times N}$ using Eq.~\eqref{eq:A1-matrix-term-entry-wise-poisson}\; 
Compute $A_2\in\mathbb{R}^{m_\Omega\times N}$ and $A_3\in\mathbb{R}^{m_{\partial\Omega}\times N}$ using
Eqs.~\eqref{eq:A2-matrix-term-entry-wise}–\eqref{eq:A3-matrix-term-entry-wise}\;
Form vectors
$\vec{f}$, $\vec{g}$ as in Eq.~\eqref{eq:source-boundary-data-matrix-term} \;
\Indm
\BlankLine
\textbf{Build normal‑equation system according to Eq.~\eqref{eqn:ellipticPDE}}\;
\Indp
$A_{\mathrm{elliptic}}\gets A_{1}+\beta A_3^{\top}A_3$\;
$\vec{b}_{\mathrm{elliptic}}\gets A_2^{\top}\vec{f} + \beta A_3^{\top}\vec{g}$\;
\Indm
\BlankLine
\textbf{Stabilize via TSVD}\;
\Indp
Compute TSVD: $A_{\mathrm{elliptic}} = U \Sigma V^{\top}$\;
Identify index set $\mathcal{T}_{\tau}=\{k : \sigma_k > \tau \,  \sigma_{1}  \}$\;
Form truncated pseudo‑inverse
$\Sigma_{\tau}^{-1}$ where $(\Sigma_{\tau}^{-1})_{kk}=1/\sigma_k$ if $k\in\mathcal{T}_{\tau}$ and $0$ otherwise\;
Form the pseudo-inverse matrix $A_{\mathrm{elliptic}, \tau}^{\dagger}\gets V_{p} \, \big( \Sigma_{\tau}^{-1} \big) \, U_{p}^{\top}$\;
\Indm
\BlankLine
\textbf{Compute the weights}\;
\Indp
$\vec{w}_{\tau}\gets A_{\mathrm{elliptic}, \tau}^{\dagger}\,\vec{b}_{\mathrm{elliptic}}$ using Eq.~\eqref{eq:pseudo-inverse-elliptic-tau-tsvd}  \;
\Indm
\BlankLine
\Return $\vec{w}_{\tau} \in \mathbb{R}^{N}$, \, $\hat{u}_{\tau}(\vec{x}) = \frac{1}{\sqrt{N}}\sum_{i = 1}^{N} w_{(\tau),i} \, \phi(\| \vec{x} - \vec{c}_{i} \|_{2} )$ as in Eq.~\eqref{eq:approximation-solution-u-hat-tau-elliptic}. 
\end{algorithm}

\subsection{Constrained Optimization for Obstacle-Type Problems}\label{sec:constrained-optimization}

Using the radial basis function representation for $\hat{u}$ (see Eq.~\eqref{eq:main-approximation-form-rbf})and the formulation of the gradient of $\hat{u}$ in terms of the weight vector (see Eq.~\eqref{eq:nabla-derivative-for-u}), the discrete form for the obstacle problem in Eq.~\eqref{eq:discretized-weak-form-obstacle} becomes:
\[
\widehat{\mathcal{J}}_{\mathrm{obstacle}, \mu, \beta}(\vec{w}) = \frac{1}{2} \vec{w}^{\top} A_{1} \vec{w} + \mu \| (\vec{\psi} - A_{2} \vec{w} )_{+} \|_{1} + \frac{\beta}{2} \| A_{3} \vec{w} - \vec{g} \|_{2}^{2} 
\]
where the matrices $A_{1}$, $A_2$, and $A_{3}$ are the same as those in Eqs.~\eqref{eq:A1-matrix-term-entry-wise-poisson}-\eqref{eq:A3-matrix-term-entry-wise}. The vector of the obstacle function is defined as follows: 
\begin{equation}
\vec{\psi} \coloneq \frac{|\Omega|}{m_{\Omega}} \, \begin{bmatrix}
    \psi(\vec{x}^{(\Omega)}_{1} ) \\ 
    \vdots \\ 
    \psi(\vec{x}^{(\Omega)}_{m_{\Omega}} )
\end{bmatrix} \in \mathbb{R}^{m_{\Omega}}.
\label{eq:psi-matrix-entry-wise-obstacle}
\end{equation}

The unconstrained optimization problem which involves the $\ell^{1}$-norm is non-smooth and convex. We solve the corresponding optimization problem by the ADMM method. Specifically, we introduce a slack variable $ \vec{v} = \vec{\psi} - A_{2} \vec{w}$. This results in a constrained optimization problem:

\begin{equation}
\begin{cases}
\min\limits_{\vec{w}, \vec{v} } \widehat{\mathcal{J}}_{\mathrm{obstacle},\mu, \beta}(\vec{w}, \vec{v}) \\
\widehat{\mathcal{J}}_{\mathrm{obstacle}, \mu, \beta}(\vec{w}, \vec{v}) = \frac{1}{2} \vec{w}^{\top} A_{1} \vec{w} + \mu \| \vec{v}_{+} \|_{1} + \frac{\beta}{2} \| A_{3} \vec{w} - \vec{g} \|_{2}^{2} \\
\text{s.t.}\quad \vec{v} - \vec{\psi} + A_{2} \vec{w} = 0. 
\end{cases}
\label{eq:constrained-optimization-obstacle-problem-equation}
\end{equation}

The corresponding augmented Lagrangian functional with a penalty parameter $\rho > 0$ and a scaled dual variable $\vec{z}$ is 

\begin{equation}
\widehat{L}_{\rho}(\vec{w}, \vec{v}, \vec{z}) = \frac{1}{2} \vec{w}^{\top} A_{1} \vec{w} + \frac{\beta}{2} \| A_{3} \vec{w} - \vec{g} \|_{2}^{2} + \mu \, \| \vec{v}_{+} \|_{1} + \frac{\rho}{2} \| \vec{v} - \vec{\psi} + A_{2}\vec{w} + \vec{z}  \|_{2}^{2} - \frac{\rho}{2} \| \vec{z} \|_{2}^{2}.
\label{eq:augmented-lagrangian-function-obstacle-problem}
\end{equation}

The ADMM algorithm iteratively updates the primal variables \(\vec{w}, \vec{v}\) and the dual variable \(\vec{z}\) as follows: 
\begin{align}
\vec{w}^{(i+1)} & \coloneq \arg \min \limits_{\vec{w}} \Big\{ \frac{1}{2} \vec{w}^{\top}  A_{1} \vec{w} + \frac{\beta}{2} \| A_{3} \vec{w} - \vec{g} \|_{2}^{2} + \frac{\rho}{2} \| \vec{v}^{(i)} - \vec{\psi} + A_{2} \vec{w} + \vec{z}^{(i)} \|_{2}^{2}   \Big\}, 
\\ 
\vec{v}^{(i+1)} & \coloneq \arg \min \limits_{\vec{v}} \Big\{ \mu \| \vec{v}_{+} \|_{1} + \frac{\rho}{2} \, \| \vec{v} - \vec{\psi} + A_{2} \vec{w}^{(i+1)} +  \vec{z}^{(i)} \|_{2}^{2} \Big\},  
\\
\vec{z}^{(i+1)} & \coloneq \vec{z}^{(i)} + \big( \vec{v}^{(i+1)} - \vec{\psi} + A_{2} \vec{w}^{(i+1)}    \big).
\end{align}

The solution \(\vec{w}^{(i+1)}\) to the \(\vec{w}\)-sub-minimization problem satisfies the following: 

\[
\Big( A_{1} + \beta A_{3}^{\top} A_{3} + \rho A_{2}^{\top} A_{2}  \Big) \, \vec{w}^{(i+1)} = \beta A_{3}^{\top} \vec{g} + \rho A_{2}^{\top} \Big( \vec{\psi} - \vec{v}^{(i)} - \vec{z}^{(i)}   \Big). 
\]

Since the total stiffness matrix is ill-conditioned in the $\vec{w}$-update, we stabilize the symmetric positive-definite linear system by computing the TSVD solution with the parameter \(\tau\): 

\begin{equation}
\vec{w}_{\tau}^{(i+1)} = A_{\mathrm{obstacle}, \tau}^{\dagger} \, \vec{b}^{(i)}_{\mathrm{obstacle}}, 
\label{eq:w-minimization-update-obstacle-problem}
\end{equation}
where $A_{\mathrm{obstacle}} \coloneq A_{1} + \beta A_{3}^{\top} A_{3} + \rho A_{2}^{\top} A_{2} \in \mathbb{R}^{N \times N}$ and $\vec{b}^{(i)}_{\mathrm{obstacle}} \coloneq \beta A_{3}^{\top} \vec{g} + \rho A_{2}^{\top} (\vec{\psi} - \vec{v}^{(i)} - \vec{z}^{(i)} ) \in \mathbb{R}^{N}$. The matrix \(A^{\dagger}_{\mathrm{obstacle},\tau}\) is computed in the same way as the matrix \(A^{\dagger}_{\mathrm{elliptic},\tau}\) in Eq.~\eqref{eq:pseudo-inverse-elliptic-tau-tsvd}. 

Let $\vec{t}^{(i)} \coloneq \vec{\psi} - A_{2} \vec{w}^{(i+1)} - \vec{z}^{(i)} \in \mathbb{R}^{m_{\Omega}}$. For each component $j = 1, \dots, m_{\Omega}$, the $\vec{v}$-update is given by
\begin{equation}
v_{j}^{(i+1)} = \begin{cases}
    t_{j}^{(i)}, & t_{j}^{(i)} < 0, \\ 
    0, & 0 \leq t_{j}^{(i)} \leq \frac{\mu}{\rho}, \\ 
    t_{j}^{(i)} - \frac{\mu}{\rho}, & t_{j}^{(i)} > \frac{\mu}{\rho}. 
\end{cases}
\label{eq:v-minimization-update-obstacle-problem}
\end{equation}
The primal and dual residuals in the iteration $(i+1)$ are defined as follows: 
\[
\vec{r}^{(i+1)} \coloneq \vec{v}^{(i+1)} - \vec{\psi} + A_{2} \vec{w}^{(i+1)},\quad \vec{s}^{(i+1)} \coloneq \rho A_{2}^{\top} \, (\vec{v}^{(i+1)} - \vec{v}^{(i)} ). 
\]
Following the suggestion in \cite{boyd2011distributed}, our algorithm terminates when both the primal and dual residuals, \(\vec{r}^{(i+1)}\) and \(\vec{s}^{(i+1)}\), respectively, meet the relative tolerance thresholds \(\varepsilon_{\mathrm{primal}}, \varepsilon_{\mathrm{dual}} > 0\): 
\begin{align}
\frac{\|  \vec{r}^{(i+1)} \|_{2} }{ \max \Big( \| A_{2} \vec{w}^{(i+1)} \|_{2}, \|  \vec{v}^{(i+1)} \|_{2}, \| \vec{\psi} \|_{2} \Big) } &\leq \epsilon_{\text{primal}}, 
\label{eq:relative-primal-residuals-obstacle-problem}
\\
\frac{ \| \vec{s}^{(i+1)} \|_{2} }{ \| \rho A_{2}^{\top} \,  \vec{z}^{(i+1)} \|_{2} } &\leq \epsilon_{\text{dual}}.
\label{eq:relative-dual-residuals-obstacle-problem}
\end{align}
The proposed algorithm for the obstacle problem is given in Algorithm~\ref{algorithm2-rbf-obstacle-admm}. In Section~\ref{sec:numerical-experiments}, we employ numerical experiments for linear elliptic and obstacle-type boundary value problems to benchmark the accuracy of Algorithms~\ref{algorithm1-rbf-solver-for-elliptic-pde}-\ref{algorithm2-rbf-obstacle-admm}. 

\begin{algorithm}[H]\label{algorithm2-rbf-obstacle-admm}
\small 
\caption{ADMM-stabilized RBF solver for the classical obstacle problem}

\KwIn{$\phi(\cdot)$ (kernel function);
       $\{\vec{c}_{\ell}\}_{\ell = 1}^{N} \subseteq \overline{\Omega}_{T}$ (sampled center nodes from the expanded domain);
      $\{\vec{x}^{(\Omega)}_j\}_{j=1}^{m_\Omega}$ (interior nodes);
      $\{\vec{x}^{(\partial \Omega)}_n\}_{n=1}^{m_{\partial\Omega}}$ (boundary nodes);
      $f(\cdot),\,g(\cdot)$ (source/boundary data);
      $\psi(\cdot)$ (obstacle data);
      $\beta>0$ (boundary penalty);
      $\mu>0$ (obstacle weight);
      $\rho>0$ (ADMM penalty);
      $\tau>0$ (truncation tolerance);
      $\mathrm{maxIter}$ (maximum number of iterations);
$\varepsilon_{\mathrm{pri}},\varepsilon_{\mathrm{dual}}$ (tolerances for primal and dual residuals)}

\KwOut{Weights $\vec{w}_{\tau} \in\mathbb{R}^{N}$; approximation $\hat u_{\tau}(\vec x) \in \mathbb{R}^{d}$}

\BlankLine
\textbf{Assemble operators and vectors}\;
\Indp
Compute $A_1, A_2, A_3$ as stated in Eqs.~\eqref{eq:A1-matrix-term-entry-wise-poisson}, ~\eqref{eq:A2-matrix-term-entry-wise} and ~\eqref{eq:A3-matrix-term-entry-wise}, respectively    \;

Form $A_{\mathrm{obs}}
    \gets A_1+\beta A_3^{\top}A_3+\rho A_2^{\top}A_2$\;
    
Form $\vec{f},\;\vec{g},\;\vec{\psi}$ using Eqs.~\eqref{eq:source-boundary-data-matrix-term} and ~\eqref{eq:psi-matrix-entry-wise-obstacle}  \;

\Indm\BlankLine

\textbf{Initialize} $\vec{w}^{(0)},\;\vec{v}^{(0)},\;\vec{z}^{(0)}$\;

\For{$i=0,\dots,\mathrm{maxIter}$}{

\tcp{\small1. $\vec{w}$-update (TSVD solve)}
  $\vec{b}^{(i)}_{\mathrm{obstacle}}
      \leftarrow
      \beta A_3^{\top}\vec{g}
      +\rho A_2^{\top}(\vec{\psi}-\vec{v}^{(i)}-\vec{z}^{(i)})$\;
  $\vec{w}_{\tau}^{(i+1)}\leftarrow A_{\mathrm{obstacle}, \tau}^{\dagger}\,\vec{b}^{(i)}_{\mathrm{obstacle}}$ using Eq.~\eqref{eq:w-minimization-update-obstacle-problem}\;\BlankLine

  \tcp{\small2. $\vec{v}$-update (positive hard‑shrinkage)}
  $\vec{t}^{(i)}
      \leftarrow \vec{\psi}-A_2\vec{w}^{(i+1)}-\vec{z}^{(i)}$\;
  \ForEach{$j=1,\dots,m_\Omega$}{
    \uIf{$t_j^{(i)}<0$}{$v_j^{(i+1)}\leftarrow t_j^{(i)}$}
    \uElseIf{$0\le t_j^{(i)}\le \mu/\rho$}{$v_j^{(i+1)}\leftarrow 0$}
    \Else{$v_j^{(i+1)}\leftarrow t_j^{(i)}-\mu/\rho$}
  }\BlankLine

    \tcp{\small3. Dual update}
  $\vec{z}^{(i+1)}
      \leftarrow
      \vec{z}^{(i)}+\vec{v}^{(i+1)}-\vec{\psi}+A_2\vec{w}^{(i+1)}$\;\BlankLine
  
  \tcp{\small4. Primal and dual residuals; stopping criteria}
  $\vec{r}^{(i+1)}
      \leftarrow \vec{v}^{(i+1)}-\vec{\psi}+A_2\vec{w}^{(i+1)}$\;
  $\vec{s}^{(i+1)}
      \leftarrow \rho A_2^{\top}(\vec{v}^{(i+1)}-\vec{v}^{(i)})$\;
  \If{$  \frac{\|\vec{r}^{(i+1)}\|_2}
         {\max\{\|A_2\vec{w}^{(i+1)}\|_2 \, \|\vec{v}^{(i+1)}\|_2 \|\vec{\psi}\|_2\} }
        \le \varepsilon_{\mathrm{pri}}
        \;\wedge\;
        \frac{\|\vec{s}^{(i+1)}\|_2}
        {\|\rho A_2^{\top}\vec{z}^{(i+1)}\|_2}
        \le \varepsilon_{\mathrm{dual}}$}{
        \textbf{break}}
}

\Return{$\vec{w}_{\tau}^{(i+1)} \in \mathbb{R}^{N}$, $\displaystyle \hat u_{\tau}(\vec x)= \frac{1}{\sqrt{N}} \sum_{l=1}^{N}w^{(i+1)}_{(\tau), l}  \, \phi (\|\vec x-\vec c_l\|_2)$}\; 
\end{algorithm}

\pagebreak

\section{Numerical Experiments}\label{sec:numerical-experiments}

\hspace{5mm}In the following experiments, we assess the accuracy, stability, and efficiency of the RBF-based variational method in Eq.~\eqref{eq:main-approximation-form-rbf} across different elliptic partial differential equations, including Poisson and obstacle-type boundary value problems. The tests are carried out with non-homogeneous Dirichlet boundary conditions. The Gaussian kernel in Table~\ref{tab:rbf-types} is employed with a shape parameter $b > 0$ in all experiments unless we compare performance across different kernels. We use dimensional scaling for the kernel shape: in the one-dimensional case, we take $b = c \, N$, while in two-dimensional experiments, we choose $b = c \, \sqrt{N}$, with $c = c(T, \tau)$ selected based on \cite{adcock2024stable}. The formula for the optimal value \(c_{\mathrm{opt}} = c(T,\tau)\) is given by: 
\begin{equation}
c(T,\tau) = \displaystyle \min \left\{1, \, \frac{\pi}{T \, \sqrt{\log(1+\tau^{-2})}}    \right\},
\label{eq:optimal-value-c-formula-from-adcock}
\end{equation}
where \(T \geq 1\) and \(\tau \in (0,1)\) are the domain expansion factor and the truncation threshold value, respectively. 
The multiplier \(c(T,\tau)\), as stated in Eq.~\eqref{eq:optimal-value-c-formula-from-adcock}, governs the trade-off between approximation accuracy and numerical stability through the shape parameter of the RBFs. A larger \(c\) yields more localized basis functions, which generally
improve conditioning but may reduce approximation accuracy. In contrast, a smaller \(c\)
produces flatter and more overlapping basis functions, which can enhance approximation
accuracy for smooth solutions but may also lead to ill-conditioning and greater sensitivity
to truncation errors. Expression~\eqref{eq:optimal-value-c-formula-from-adcock} was derived for the case of least-squares function approximation in one dimension. While it is not necessarily optimal in our case of variational PDE solvers, we use it as a guideline and investigate its effectiveness in our numerical experiments. Following the approach of \cite{adcock2024stable} in 1D, the center points are not strictly restricted to $\overline{\Omega} = [b_{1}, b_{2}]$, but are drawn uniformly from an expanded domain $[q_{\mathrm{mid}} - T \, L_{\mathrm{half}}, q_{\mathrm{mid}} + T \, L_{\mathrm{half}}]$ where $q_{\mathrm{mid}} = \frac{b_{1} + b_{2}}{2}$ is the midpoint of the domain, $L_{\mathrm{half}} = \frac{b_{2} - b_{1}}{2}$ is the half-length, and $T \geq 1$ is the domain expansion factor, as stated in Eq.~\eqref{eq:optimal-value-c-formula-from-adcock}. This oversampling procedure, which operates outside the physical domain, reduces boundary effects on approximation accuracy and improves numerical stability within the collocation system. 
We compare the proposed RBF-based variational solver with finite-difference ADMM, neural network and Galerkin methods in the linear elliptic boundary value problems. For free-boundary value problems, we compare the proposed RBF-based variational method, finite-difference solvers with ADMM, proximal-gradient descent (PGD), and neural networks. For the finite-difference benchmark, we use the method of ADMM\footnote{The implementation of Finite-difference (ADMM) is at \url{https://github.com/jwcalder/MinimalSurfaces}} in \cite{tran20151} to solve the variational form of the Poisson problem under the same Dirichlet boundary conditions. In the two-dimensional setting, we compare with finite-difference and solve the optimization problems with PGD\footnote{The implementation of \(\texttt{MGProx}\) can be accessed at \url{https://github.com/lumeilevel/MGProx}} as in \cite{Ang_2024}. For the neural-network benchmark, we use physics-informed neural networks (PINNs)\footnote{The package of PINNs is accessed at \url{https://github.com/maziarraissi/PINNs}} from \cite{raissi2017physics,raissi2017physicsinformeddeeplearning, raissi2018deep}, in which the network parameters are trained by minimizing the corresponding variational loss from Section~\ref{sec:methodology}. In one-dimensional experiments, we set the fully-connected neural network to consist of one hidden layer with hyperbolic tangent activation function using \(N\) neurons, followed by a linear output layer. In two-dimensional experiments, we consider two architectures: the same one-hidden-layer as in the one-dimensional case, and a deeper network with three connected layers which have \(N\) neurons each, again using \(\texttt{tanh}\) activations. In the Galerkin method, we solve the weak form of the boundary value problems using finite-element discretization; a common open-source implementation is \(\texttt{FEniCSx}\)\footnote{The package of \(\texttt{FEniCSx}\) is implemented at \url{https://github.com/METHODS-Group/ProximalGalerkin}} in \cite{dokken2025latent}. We use \(N\) as the mesh resolution in the Galerkin method. For reproducibility, we report the main implementation parameters for each method. In the neural network, we use the same parameter settings as in \cite{raissi2017physics} with a learning rate of \( 10^{-2}\), the Adam optimizer, and a hyperbolic tangent activation function. For the finite-difference and Galerkin methods, the stopping criterion is chosen so that the relative error is below a prescribed tolerance. 
For each problem, we sample $m=  \zeta \times N$ collocation points, including $m_{\Omega}$ interior points and $m_{\partial\Omega}$ boundary points, where $\zeta$ is an oversampling parameter, and $N$ is the number of RBF centers. In the one-dimensional setting, the boundary consists solely of two endpoints \(b_{1}\) and \(b_{2}\). Thus, the number of interior points is \(m_{\Omega} = m - 2\). In the two-dimensional setting, we uniformly randomly sample \(\sqrt{\zeta \, N}\) points on each side of the boundary domain \(\partial\Omega\), yielding \(m_{\partial \Omega} = 4 \sqrt{\zeta\, N}\). For interior collocation points, we uniformly randomly sample \(m_{\Omega} = \zeta \, N - 4 \sqrt{\zeta \, N}\) points. 
For the comparison between methods, we compute the relative $\ell^2$ error,
        \begin{equation}
        E(N) \coloneq \frac{\| \hat{u}_{N}  - u_{exact} \|_{2}  }{ \| u_{exact} \|_{2}  }, 
        \label{relative-l2-error-norm}
        \end{equation}
where $\hat{u}_{N}$ is the numerical approximation of the optimization problems (see  Eq.~\eqref{eq:approximation-solution-u-hat-tau-elliptic}) and the 2-norms are computed by sums over the collocation points. 

Empirically, we also calculate the local estimation of the convergence order $p_{i}$ and their average order $\hat{p}$: 
    \begin{equation}
    p \coloneq \log \left( \frac{E(N/2)  }{ E(N)}\right).
    \end{equation}
Those experiments quantify the order of accuracy of our method for Poisson and free-boundary value problems. In all experiments, we set hyper-parameters according to Table~\ref{tab:params-numerical-experiments}.

\begin{table}[H]
\centering
\small
\caption{Parameters for the RBF-algorithm solver in numerical experiments}
\label{tab:params-numerical-experiments}
\renewcommand{\arraystretch}{1.2}
\setlength{\tabcolsep}{4pt}
\begin{tabular}{|p{5.0cm}|c|p{8.5cm}|}
\hline
\textbf{Parameter} & \textbf{Symbol} & \textbf{Value/Formula} \\
\hline
Number of basis functions & \(N\) & \texttt{problem-dependent} \\
\hline
Oversampling factor & \(\zeta\) & \(2\) (1D); \(4\) (2D) \\
\hline
Total number of all collocation points & \(m\) & \(\zeta \times N\) \\
\hline
Half-length of interval \([b_{1}, b_{2}]\) & \(L_{\mathrm{half}}\) & \(L_{\mathrm{half}} = \dfrac{b_{2} - b_{1}}{2}\) \\
\hline
Midpoint of the interval \([b_{1}, b_{2}]\) & \(q_{\mathrm{mid}}\) & \(q_{\mathrm{mid}} = \dfrac{b_{1} + b_{2}}{2}\) \\
\hline
Domain expansion factor & \(T\) & Centers are sampled in 1D from 
\[
\bigl[q_{\mathrm{mid}} - T L_{\mathrm{half}},\;
q_{\mathrm{mid}} + T L_{\mathrm{half}}\bigr]
\] \\
\hline
Truncation threshold value & \(\tau\) & \(10^{-15}\) \\
\hline
Shape-factor multiplier & \(c\) & \(\displaystyle \min\!\left\{1,\;
\frac{\pi}{T\sqrt{2\log(1+\tau^{-2})}}\right\}\), as in \cite{adcock2024stable} \\
\hline
Kernel shape parameter & \(b\) & \(cN\) (1D); \(c\sqrt{N}\) (2D) \\
\hline
Obstacle penalty weight & \(\mu\) & \texttt{problem-dependent} \\
\hline
Boundary-condition penalty term & \(\beta\) & \texttt{problem-dependent} \\
\hline
ADMM penalty parameter & \(\rho\) & \texttt{problem-dependent} \\
\hline
Maximum number of iterations in ADMM & \texttt{maxIter} & \(50{,}000\) \\
\hline
\end{tabular}
\end{table}

\subsection{Elliptic Boundary Value Problem}\label{sec:elliptic-boundary-value-problem-numerical-experiment}

\hspace{5mm}In this subsection, we consider a class of elliptic boundary value problems defined on a bounded closed domain $\overline{\Omega} \subset \mathbb{R}$. The boundary value problem is the following Poisson equation subject to inhomogeneous Dirichlet boundary conditions:
\begin{equation}
\begin{cases}
    - \Delta u = f(x), & \forall x \in \Omega = (0,1) \\ 
    u(0) = a_{1}, \\ 
    u(1) = a_{2}.
\end{cases}
\label{numerical-experiment-poisson-bvp}
\end{equation}
We use the following problem to test the variational form of the Poisson equation. The solution of the Poisson boundary value problem  is: 
\begin{equation}
\begin{cases}
u(x) = \sin(\pi x) -  (a_{1} - a_{2}) x + a_{1},\quad x \in \Omega = (0,1) \\ 
u(0) = a_{1}, u(1) = a_{2}, 
\end{cases}
\label{poisson-analytical-solution}
\end{equation}
where the external forcing term $f(x) = - \pi^{2} \sin(\pi x)$, and the boundary values are set at $a_{1} = -1.0$ and $a_{2} = -\frac{3}{2}$.

\begin{figure}[H]
   \centering
    \includegraphics[width = 0.75\linewidth]{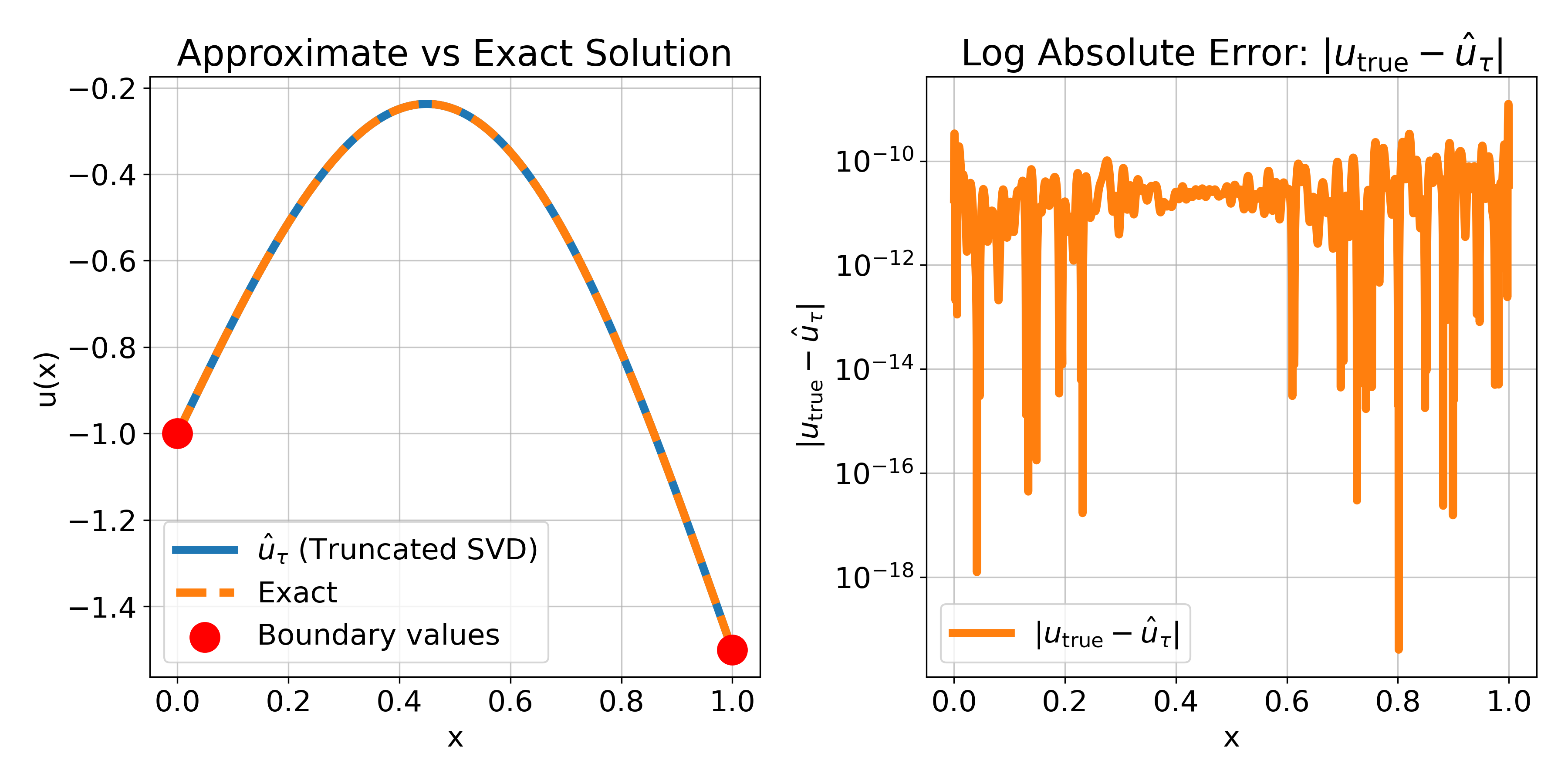}
    \caption{Numerical approximation with RBF for Poisson equation $- \Delta u = f(x)$. The relative error, $E(N) = 7.419 \times 10^{-11}$. The resolution parameters are $(N, m) = (4096, 8192)$. Center points are sampled uniformly with $\zeta = 2, T = 8$ and settings are $\tau =  10^{-15}, \beta = 3 \times 10^{5}$.}
    \label{fig1_poisson_pde}
\end{figure}

Our method uses an RBF representation according to Eq.~\eqref{eq:main-approximation-form-rbf} regularized by truncated SVD at the cutoff $\tau$, together with a variational penalty parameter $\beta$ on the Dirichlet boundary conditions. In Fig.~\ref{fig1_poisson_pde}, the top-left panel shows accurate agreement between the approximation and the solution in Eq.~\eqref{poisson-analytical-solution} with relative $\ell^2$ error $ 7.419 \times 10^{-11}$ where $N = 4096$. At the boundary, we penalize Dirichlet conditions via parameter $\beta = 3 \times 10^{5}$. The pointwise error remains uniformly small according to the right panel in Fig.~\ref{fig1_poisson_pde}. In particular, $| u_{exact}(x) - \hat{u}_{\tau}(x)|$ oscillates between $10^{-9}$ and $10^{-12}$ at the interior collocation points, dipping to around $10^{-14}$ in the boundary regions.    

\begin{figure}[H]
    \centering
    \begin{subfigure}[t]{0.5\textwidth}
    \centering\includegraphics[width = 1.0\linewidth]{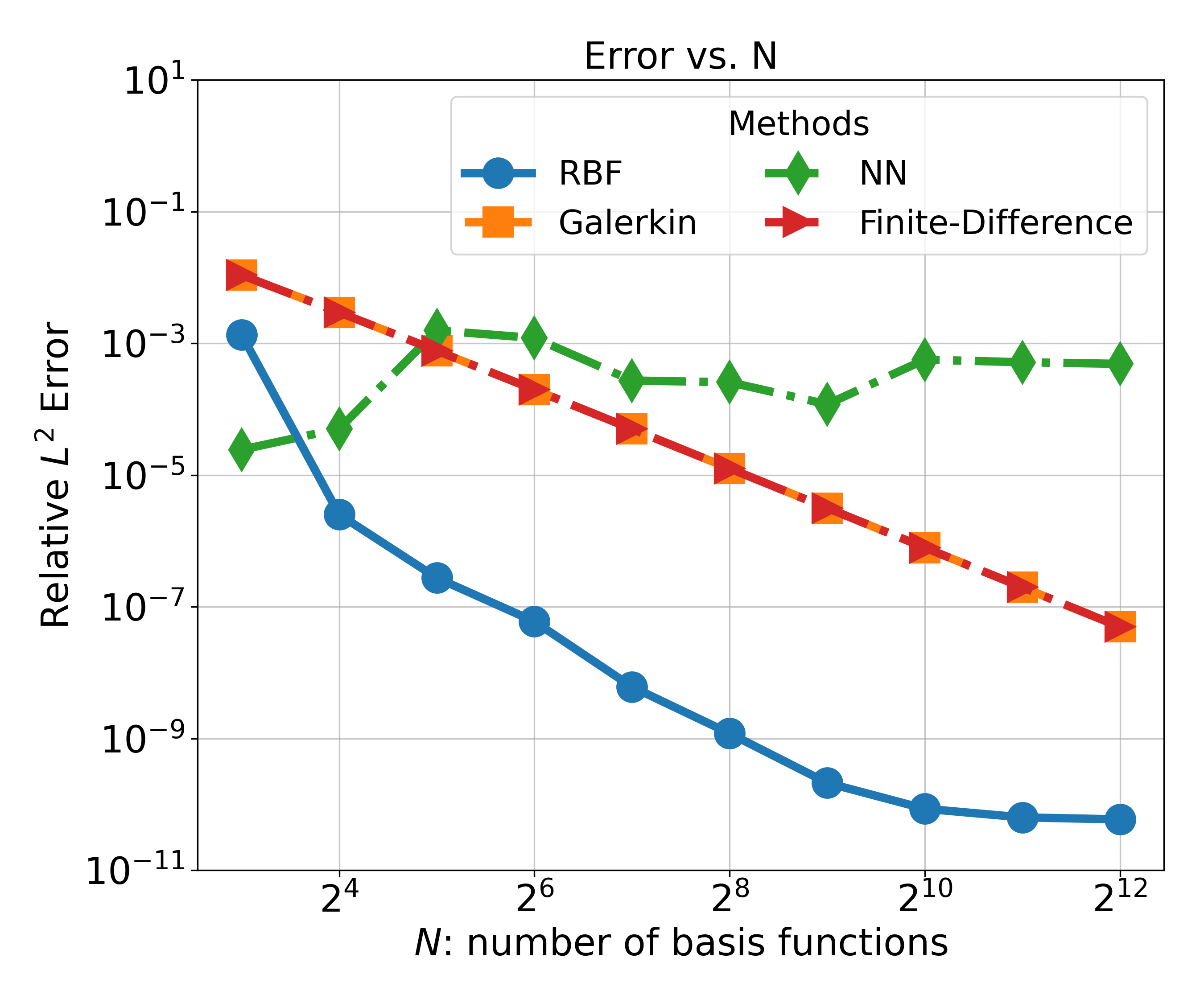}
    \caption{Convergence: relative $\ell^2$ error vs.\ $N$}
    \end{subfigure}\hfill 
    \begin{subfigure}[t]{0.5\textwidth}
    \centering\includegraphics[width = 1.0\linewidth]{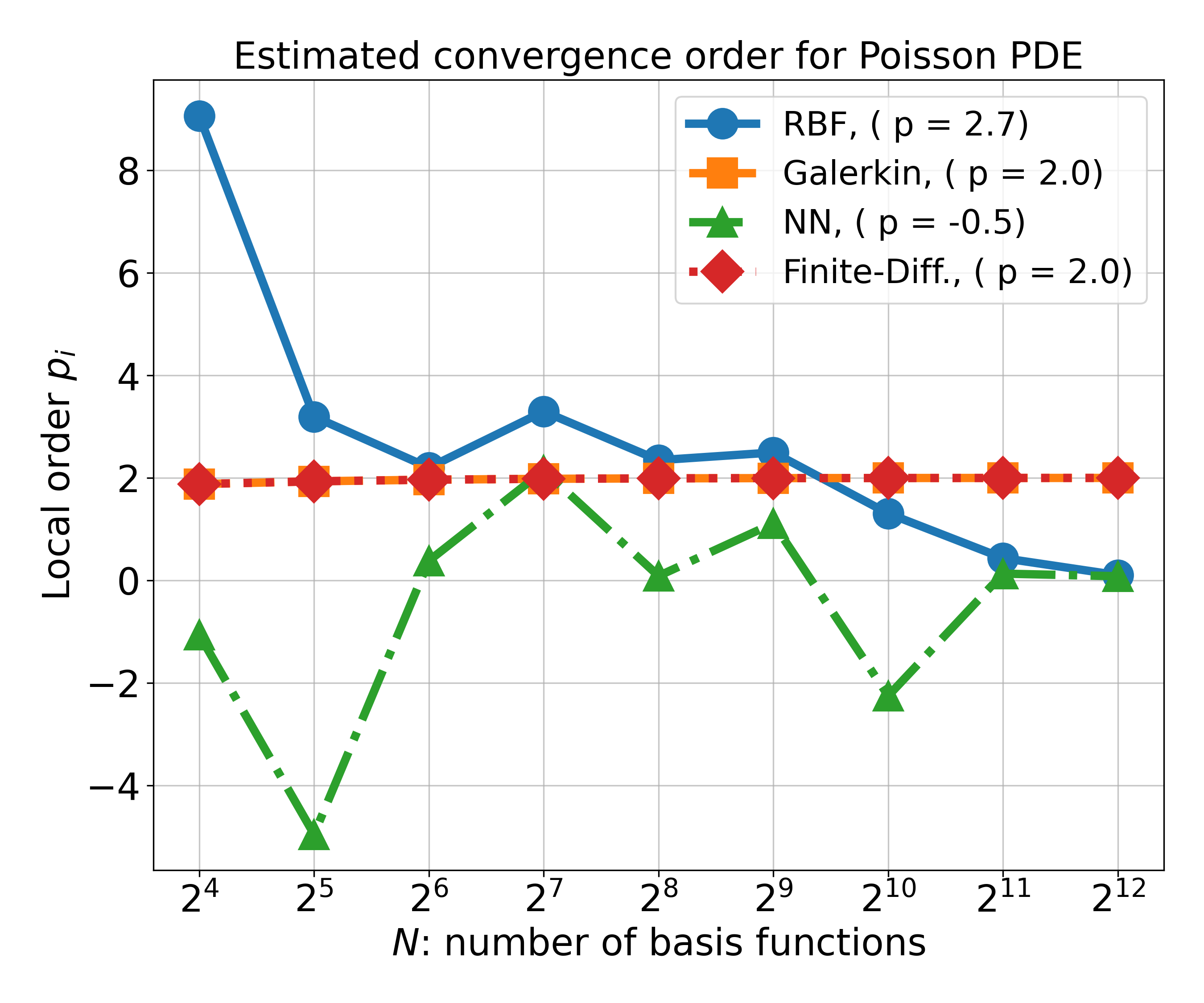}
    \caption{Local estimated convergence order $p$.}
    \end{subfigure}
    \caption{Poisson PDE $-\Delta  u = f(x)$: accuracy and convergence order. In this experiment, we set the following hyperparameters for the Poisson problem: \(T = 8, \,\beta = 3\times 10^{5}, \, \tau = 10^{-15}\)}
    \label{fig2_poisson_pde}
\end{figure}

Fig.~\ref{fig2_poisson_pde} illustrates the comparative performance of different numerical methods for the Poisson equation, including the proposed RBF method, as well as neural-network, Galerkin, and finite difference methods. 

For the benchmark tests in Figs.~\ref{fig1_poisson_pde}--\ref{fig2_poisson_pde}, the compared methods are run with the following hyperparameter settings. The proposed RBF method uses Gaussian basis functions with \(\zeta=2\), \(T=8\), \(\tau=10^{-15}\), and \(\beta=3\times 10^{5}\). The neural-network baseline uses one hidden layer with \(\tanh\) activation and a linear output layer, with width chosen so that the number of trainable parameters is comparable to \(N\). The Galerkin method uses the standard weak formulation with the same nominal number of basis functions, and the finite-difference method is computed on a uniform mesh with comparable resolution. All stopping criteria are kept fixed across runs.

Panel Fig.~\ref{fig2_poisson_pde}(a) demonstrates that the RBF method achieves a rapid decrease in relative $\ell^2$ error as the number of basis functions $N$ increases, exhibiting fast convergence. The Galerkin and finite-difference methods display the expected second-order accuracy, while the neural network (1 hidden layer) converges slowly and stagnates at higher resolutions. Panel Fig.~\ref{fig2_poisson_pde}(b) confirms the trends by showing the estimated local convergence order $p$. The RBF method attains close to $2.7$ order on average, surpassing the second-order rates of the Galerkin and finite-difference schemes.

\subsection{Obstacle-Type Boundary Value Problem}

In the experiments for obstacle-type boundary value problems, we seek a minimum of Eq.~\eqref{weak-form-obstacle-problem} to solve the boundary value problem in Eq.~\eqref{obstacle-pde-bvp}. The hyper-parameters in the variational form~\ref{weak-form-obstacle-problem} are set according to Table~\ref{tab:params-numerical-experiments}. The pair parameters $(\mu, \rho)$ in Eq.~\eqref{eq:augmented-lagrangian-function-obstacle-problem} are tuned according to \cite{tran20151}. 

We vary the $N$ number of basis functions to assess accuracy, tracking the relative $\ell^2$ error, and monitor ADMM convergence through the relative primal and dual residuals in Eqs.~\eqref{eq:relative-primal-residuals-obstacle-problem}-\eqref{eq:relative-dual-residuals-obstacle-problem}, respectively. Generally, we set the tolerances as $\epsilon_{\mathrm{tol.}} = \epsilon_{primal} = \epsilon_{dual}$, depending on the desired accuracy for the experiment. Finally, we record the peak memory usage and computational time in seconds to characterize the scalability.

\subsubsection{One-Dimensional Obstacle Problem}

We solve the one-dimensional obstacle problem using an ADMM-based penalized formulation according to Section~\ref{sec:constrained-optimization} with piecewise-defined obstacles $\psi_1$ (one bump) and $\psi_{2}$ (two bumps) given by: 
\begin{align}
\psi_{1}(x) &=
\begin{cases}
100\,x^{2}, & x\in[0,0.25],\\
100\,x(1-x)-12.5, & x\in(0.25,0.5],\\
\psi_{1}(1-x), & x\in(0.5,1],
\end{cases}
\label{psi1-function}
\end{align}
and
\begin{align}
\psi_{2}(x) &=
\begin{cases}
10\sin(2\pi x), & x\in[0,0.25],\\
5\cos\bigl(\pi(4x-1)\bigr)+5, & x\in(0.25,0.5],\\
\psi_{2}(1-x), & x\in(0.5,1].
\end{cases}
\label{psi2-function}
\end{align}
For each obstacle problem, the exact solutions are given as follows: 
\begin{align}
u_{1,\mathrm{exact}}(x) &=
\begin{cases}
(100 - 50\sqrt{2})\,x, & x \in [0,\tfrac{1}{2\sqrt{2}}],\\
100\,x(1-x) - 12.5,     & x \in (\tfrac{1}{2\sqrt{2}},0.5],\\
u_{1,\mathrm{exact}}(1-x), & x \in (0.5,1],
\end{cases}
\label{u1_exact_obstacle-problem-1d}
\end{align}
and
\begin{align}
u_{2,\mathrm{exact}}(x) &=
\begin{cases}
10\sin(2\pi x), & x \in [0,0.25],\\
10,             & x \in (0.25,0.5],\\
u_{2,\mathrm{exact}}(1-x), & x \in (0.5,1].
\end{cases}
\label{u2_exact_obstacle-problem-1d}
\end{align}
We test our ADMM solver on the unit domain $[0,1]$ for obstacles~\eqref{psi1-function} and ~\eqref{psi2-function}. For each obstacle problem, we let \(\beta = 10^{6}\) in the one bump and \(\beta = 10^{8}\) in the two bumps problem. The discussion for selecting the boundary penalty parameter is in Appendix~\ref{appendix:experiment_boundary_penalty_parameter}. Across all panels, the error decreases sharply as $\beta$ increases. In one-bump problem, when $\beta$ is beyond $10^{6}$ the error plateaus due to numerical stiffness and loss of conditioning. For the two-bump problem, the error plateaus up to $\beta \geq  10^{8}$. We decide to choose $\beta = 10^{6}$ for the one-bump problem and $\beta = 10^{8}$ for the two-bump problem. Section~\ref{sec:comparison-expansion-domain-factor-T} will justify the choice of domain expansion factor \(T\) for each obstacle problem. Specifically, we choose \(T = 2\) for the one-bump obstacle problem while \(T = 3\) is chosen for the two-bump problem. In ADMM, we let the parameters \(\mu = 300, \, \rho = 45\) for 1D one-bump obstacle problem, while \(\mu = 2.5 \times 10^{4}, \rho = 250\) in the two-bump problem. Figs.~\ref{fig1_one_bump_obstacle1D}-~\ref{fig1_two_bump_obstacle1D} present the RBF approximations of the 1D obstacle problems with a specific choice of hyperparameters. 

\begin{figure}[h!]
\centering
\includegraphics[width=\linewidth]{Result/ObstaclePDE/1Dcase/Problem_1/fig_obstacle1D_single_bump.png}
    \caption{\textbf{RBF-Approximation of the 1-D obstacle problem with a smooth single-bump obstacle function $\psi(x)$ on the domain $\Omega = [0,1]$ with $(N, m) = (2048, 4096)$ and parameters $(\mu, \rho) = (300, 45)$.}
    (a) Exact solution $u(x)$, approximation $\hat u(x)$, and obstacle function $\psi(x)$.
    (b) Absolute error $\lvert \hat u -u \lvert$ on the log-scale.
    (c) ADMM convergence demonstrates rapid decay of relative primal and dual residuals with tolerance $\epsilon_{tol.} =  10^{-6}$. The method converges in 7\,234 iterations. 
    (d) Relative $\ell^2$ error is reaches $1.904 \times 10^{-12}$.}
    \label{fig1_one_bump_obstacle1D}
\end{figure}

\pagebreak

\begin{figure}[h!]
\centering
\includegraphics[width=\linewidth]{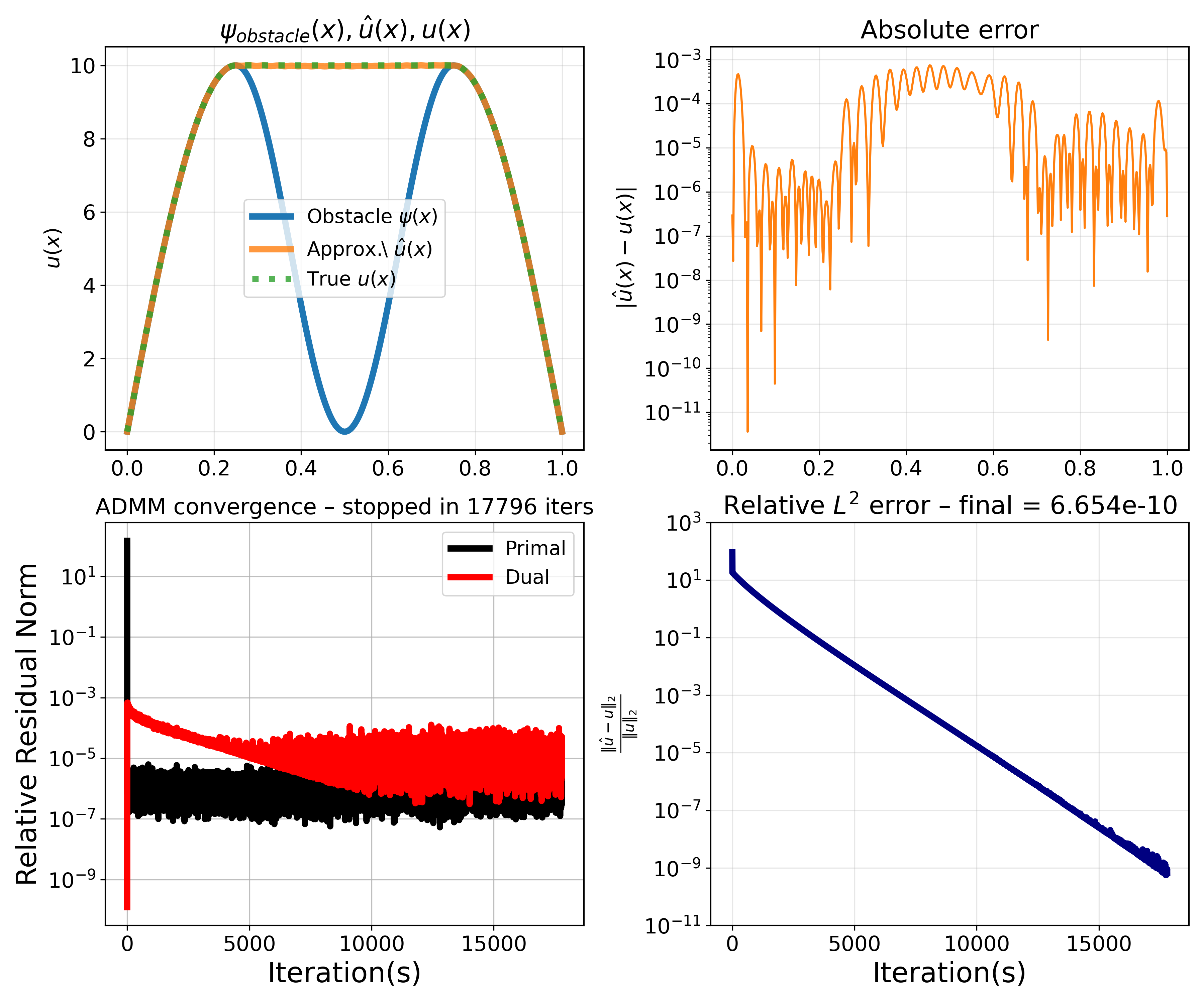}
    \caption{\textbf{Result of the 1-D obstacle problem with a piecewise two-bump $\psi(x)$ over the domain $\overline{\Omega} = [0,1]$ using $(N, m) = (1024, 2048)$ and hyper-parameters $(\mu, \rho) = (2.5 \times 10^{4}, 250)$.}
    (a) Exact solution $u(x)$, RBF-approximation $\hat u(x)$, and obstacle function $\psi(x)$. 
    (b) Log-scale plot of the absolute error $\lvert \hat u - u \lvert$. 
    (c) ADMM convergence showing steady decay of relative primal and dual residuals over 17\,796 iterations with tolerance $\epsilon_{tol.} = 1\times 10^{-6}$. 
    (d) Relative $\ell^2$ error reaches $6.654 \times 10^{-10}$. 
    }
    \label{fig1_two_bump_obstacle1D}
\end{figure}

\pagebreak

\begin{figure}[h!]
    \centering
    \begin{subfigure}[t]{0.5\textwidth}  \centering\includegraphics[width=1.0\linewidth]{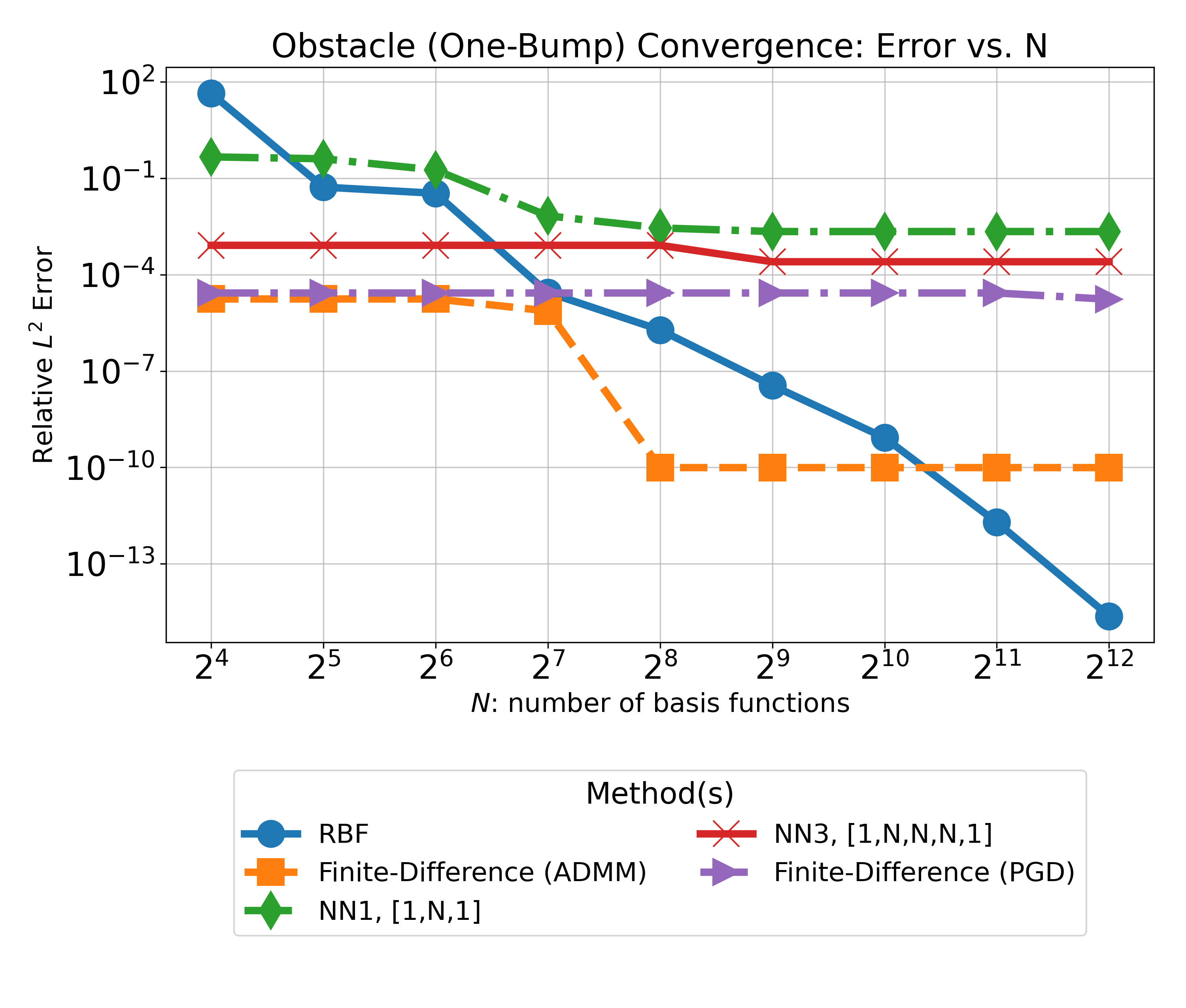}
    \caption{$T = 2.0, \zeta = 2, \tau = 10^{-15}, \beta = 10^{6}$}
    \end{subfigure}\hfill
    \begin{subfigure}[t]{0.5\textwidth}
    \centering\includegraphics[width=1.0\linewidth]{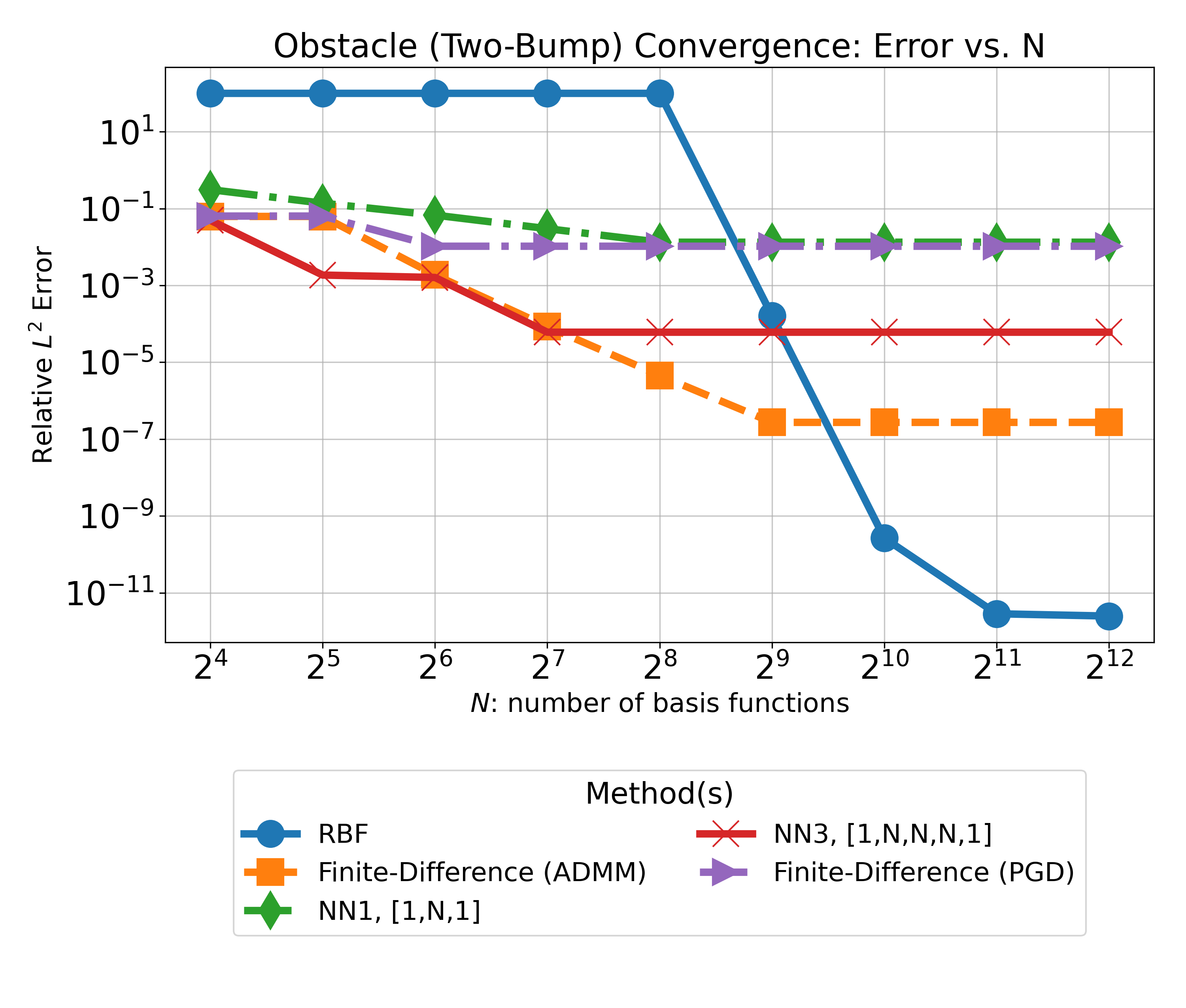}
    \caption{$T = 3.0, \zeta = 2, \tau = 10^{-15}, \beta = 10^{8}$}
  \end{subfigure}
 \caption{Convergence for 1D obstacle problem for both problems: relative $\ell^2$ error vs.\ $N$.}
 \label{fig2_1d_obstacle_problem_convergence_benchmark}
\end{figure}

\begin{figure}[h!]
    \centering
    \begin{subfigure}[t]{0.5\textwidth}  \centering\includegraphics[width=1.0\linewidth]{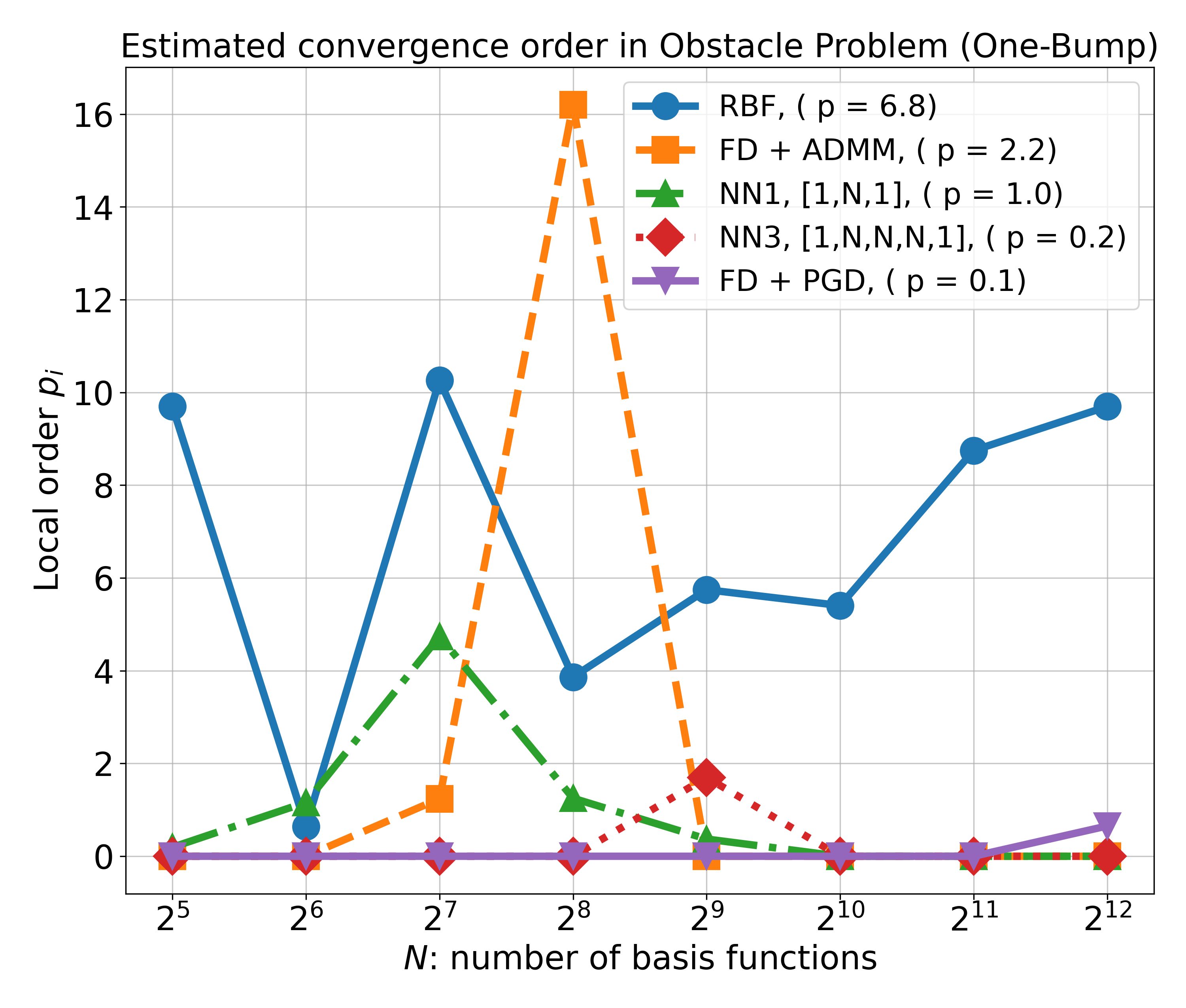}
    \caption{$T = 2.0, \zeta = 2, \tau = 10^{-15} $}
    \end{subfigure}\hfill
    \begin{subfigure}[t]{0.5\textwidth}
    \centering\includegraphics[width=1.0\linewidth]{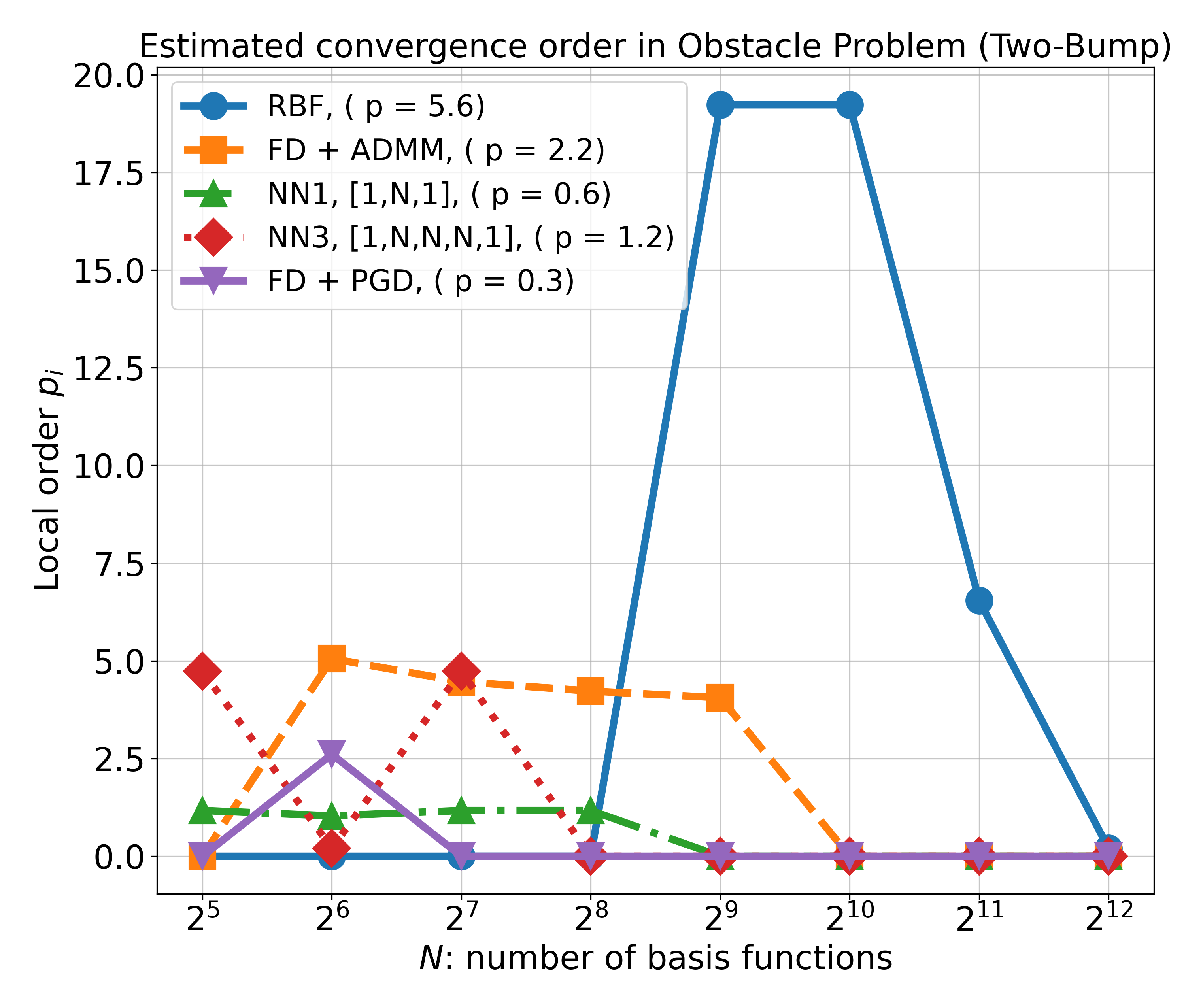}
    \caption{$T = 3.0, \zeta = 2, \tau = 10^{-15}$}
  \end{subfigure}
 \caption{Similar to Figure~\ref{fig2_1d_obstacle_problem_convergence_benchmark}. We estimate the local convergence order $p$ for each $N$ and the corresponding average order $\hat{p}$ for different solution methods.}
 \label{fig2_1d_obstacle_problem_convergence_order_benchmark}
\end{figure}

Figs.~\ref{fig2_1d_obstacle_problem_convergence_benchmark}-~\ref{fig2_1d_obstacle_problem_convergence_order_benchmark} present comparisons of different methods and estimates of the convergence order. We see that the Gaussian RBF method demonstrates fast convergence, with the relative $\ell^2$ error reaching machine-precision as \(N\) increases for both problems. In contrast, finite-difference (ADMM) performs consistently but flattens out for large \(N\). Finite-difference (PGD) shows much slower or stagnant convergence, typically limited to error plateaus between $10^{-4}$ and $10^{-6}$. Meanwhile, NN-based methods exhibit mild fluctuations and low accuracy due to the difficulty in capturing contact regions in free-boundary problems. Fig.~\ref{fig2_1d_obstacle_problem_convergence_order_benchmark} quantifies the corresponding local convergence orders $p$. The panels confirm the accuracy of the RBF method. Finite-difference methods maintain the expected second-order behavior. The NN-based methods display poor convergence. 

\pagebreak

\pagebreak

\subsubsection{Two-Dimensional Dome Obstacle Problem}

\begin{figure}[h!]
    \centering
    \begin{subfigure}[t]{0.5\textwidth}
    \centering
    \includegraphics[width = 0.7\linewidth]{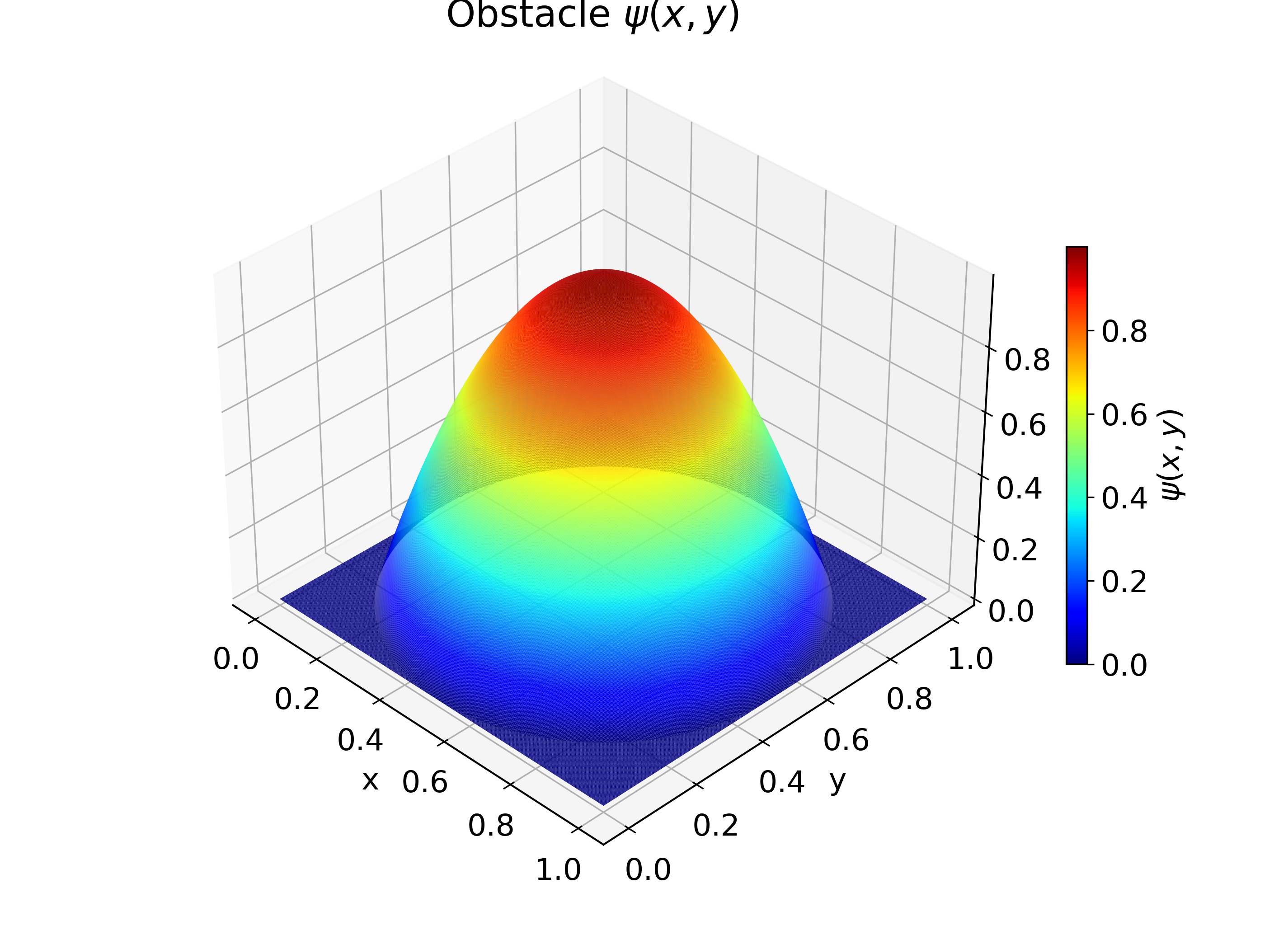}
    \caption{The obstacles $\psi(x,y)$}
    \end{subfigure}\hfill
    \begin{subfigure}[t]{0.5\textwidth}
    \centering
    \includegraphics[width = 0.7\linewidth]{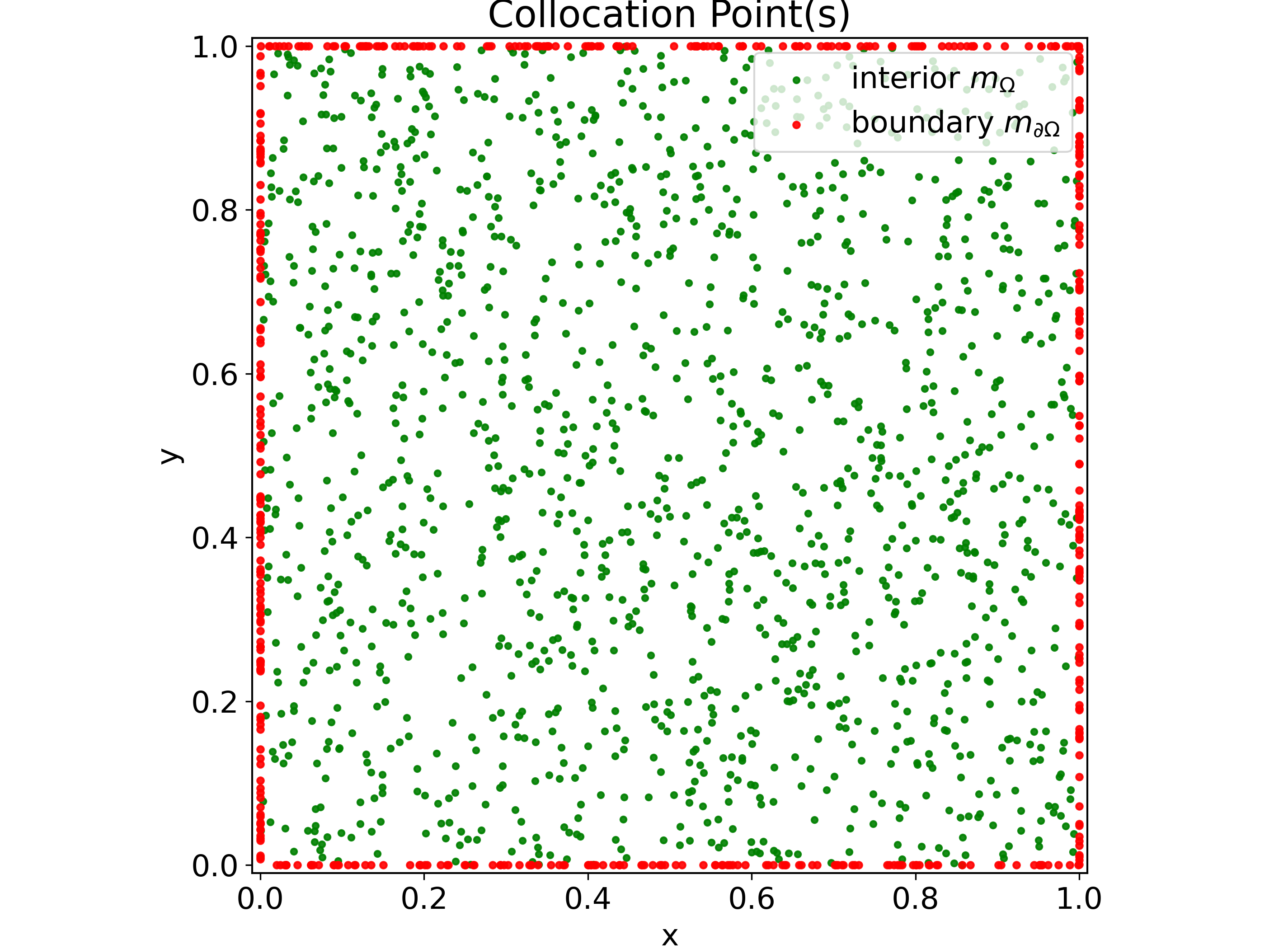}
    \caption{Collocation points}
\end{subfigure}
\caption{(a): Obstacle function $\psi(x,y)$ on the domain $\overline{\Omega} = [0,1]^{2}$. (b). Collocation points used by the variational method: green dots denote interior points and red dots are boundaries points.}
\label{fig1_obstacle2d_psi_collocation_points}
\end{figure}

In this experiment, we test the performance of our proposed method using the dome obstacle problem from \cite{tran20151, zosso2015fast}. The interior domain is $\Omega = (0,1)^2$. First, we let the obstacle be a single dome on the unit square centered at $(\frac{1}{2}, \frac{1}{2})$ as follows: 
\begin{equation}
\begin{aligned}
\psi(x,y) &=
\begin{cases}
\displaystyle 1 - \dfrac{r(x,y)^2}{r_{c}^2}, & \text{for }  r(x,y) \le r_{c},\\[10pt]
0, & \text{for } r(x,y) > r_{c},
\end{cases}
\end{aligned}
\label{eq:2d_psi_func_dome}
\end{equation}
where the inner radius is $r_{c} \coloneq \frac{1}{2}$, and $r (x,y) \coloneq \sqrt{(x - \frac{1}{2 })^{2} + (y - \frac{1}{2} )^{2}}$. The dome obstacle function $\psi(x,y)$ is shown on the left in Fig.~\ref{fig1_obstacle2d_psi_collocation_points}. 

The linear obstacle problem admits the following radially-symmetrical analytical solution:
\begin{equation}
u_{\mathrm{exact}}(x,y) = \begin{cases}
1 - \frac{r(x,y)^{2}}{r_{c}^{2} }     , & r(x,y) \leq r_{*} \\ 
1 - \frac{r_{*}^{2} }{r_{c}^{2}} \Big( 1 + 2 \log \left( \frac{r(x,y)}{r_{*}}  \right)\Big),& r(x,y)> r_{*}, \\
\end{cases} 
\label{eq:exact_solution_2d_obstacle_problem_dome}
\end{equation}
where the contact radius $r_{*} \approx 0.260197$ is determined by solving the following equation: 
\begin{equation}
\frac{r_{*}^{2}}{r_{c}^{2}} \, \Big( 1 + 2 \log \left( \frac{1}{r_{*} }\right) \Big) = 1, \quad \text{where } 0 < r_{*} < r_{c} 
\label{eq:roots_equation_2d_obstacle_problem_dome}
\end{equation}
For the boundary conditions, we set $g(x,y) = u_{\mathrm{exact}}(x,y)$ for any $(x,y) \in \partial \Omega$. We sample $m$ random points from the grid $(x,y)$. As illustrated in Fig.~\ref{fig1_obstacle2d_psi_collocation_points}, we use $N = 400$ total number of basis functions, which corresponds to \(m = 1,600\) collocation points over the domain $\overline{\Omega}$ since \(\zeta = 4\). Significantly, we use the linear kernel scaling via $b = c \, \sqrt{N}$, where $c = c(T,\tau)$ is determined as in Table~\ref{tab:params-numerical-experiments}. The choice of kernel shape parameter therefore keeps the length proportional to the fill distance $h_{\Omega}$ in Assumption~\ref{assumption-3-core} as the number of basis function $N$ increases.

\begin{figure}[h!]
    \centering

    \begin{subfigure}[t]{0.48\linewidth}
        \centering
        \includegraphics[width=\linewidth]{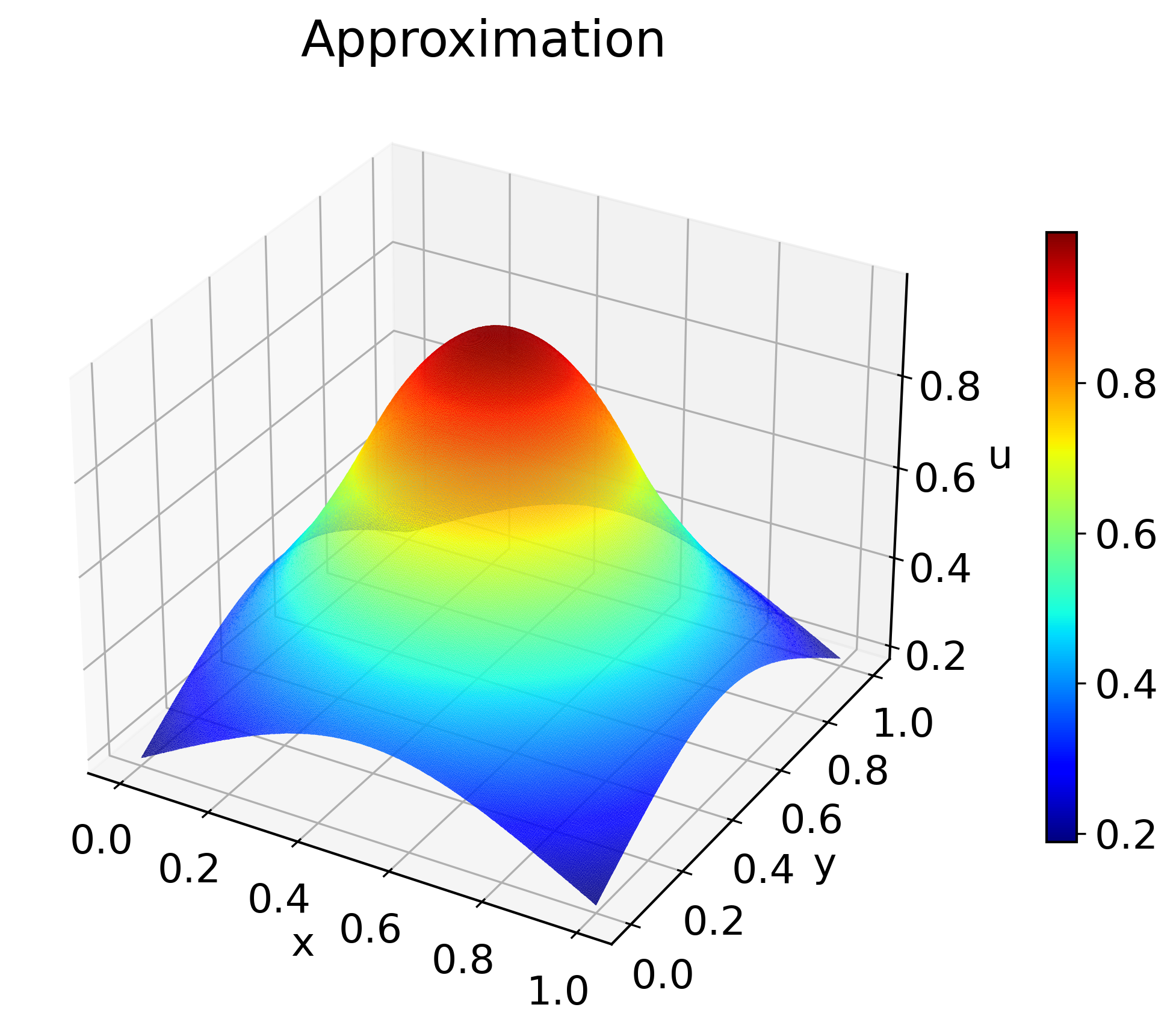}
        \caption{Approximation $\hat{u}(x,y)$ obtained with ADMM optimization using hyper-parameters $\mu = 10/h^{2}$, $\rho = 20.3$, and $\beta = 10^{6}$.}
    \end{subfigure}
    \hfill
    \begin{subfigure}[t]{0.48\linewidth}
        \centering
        \includegraphics[width=\linewidth]{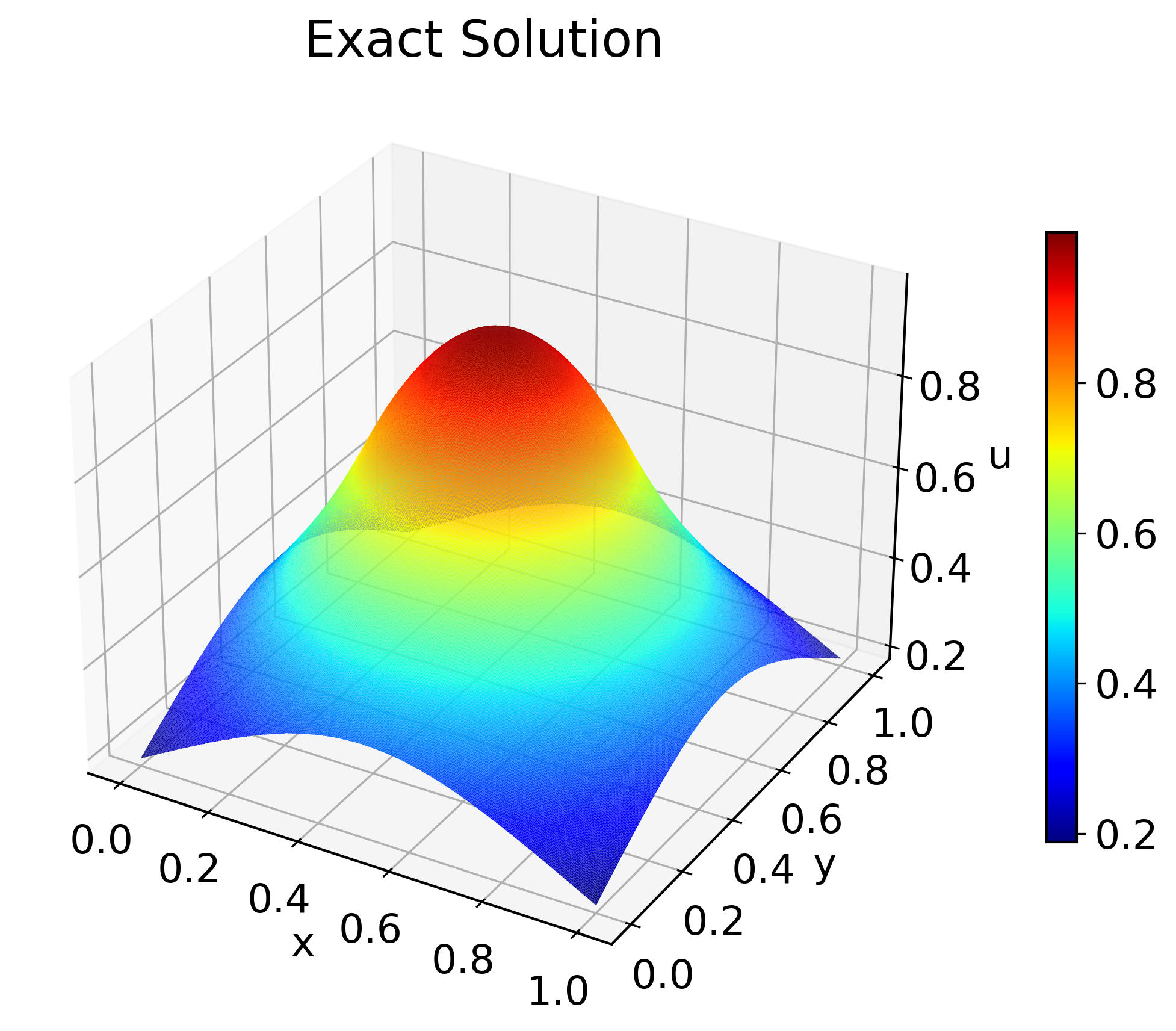}
        \caption{Exact solution $u_{\mathrm{exact}}(x,y)$.}
    \end{subfigure}

    \vspace{0.5em}

    \begin{subfigure}[t]{0.48\linewidth}
        \centering
        \includegraphics[width=\linewidth]{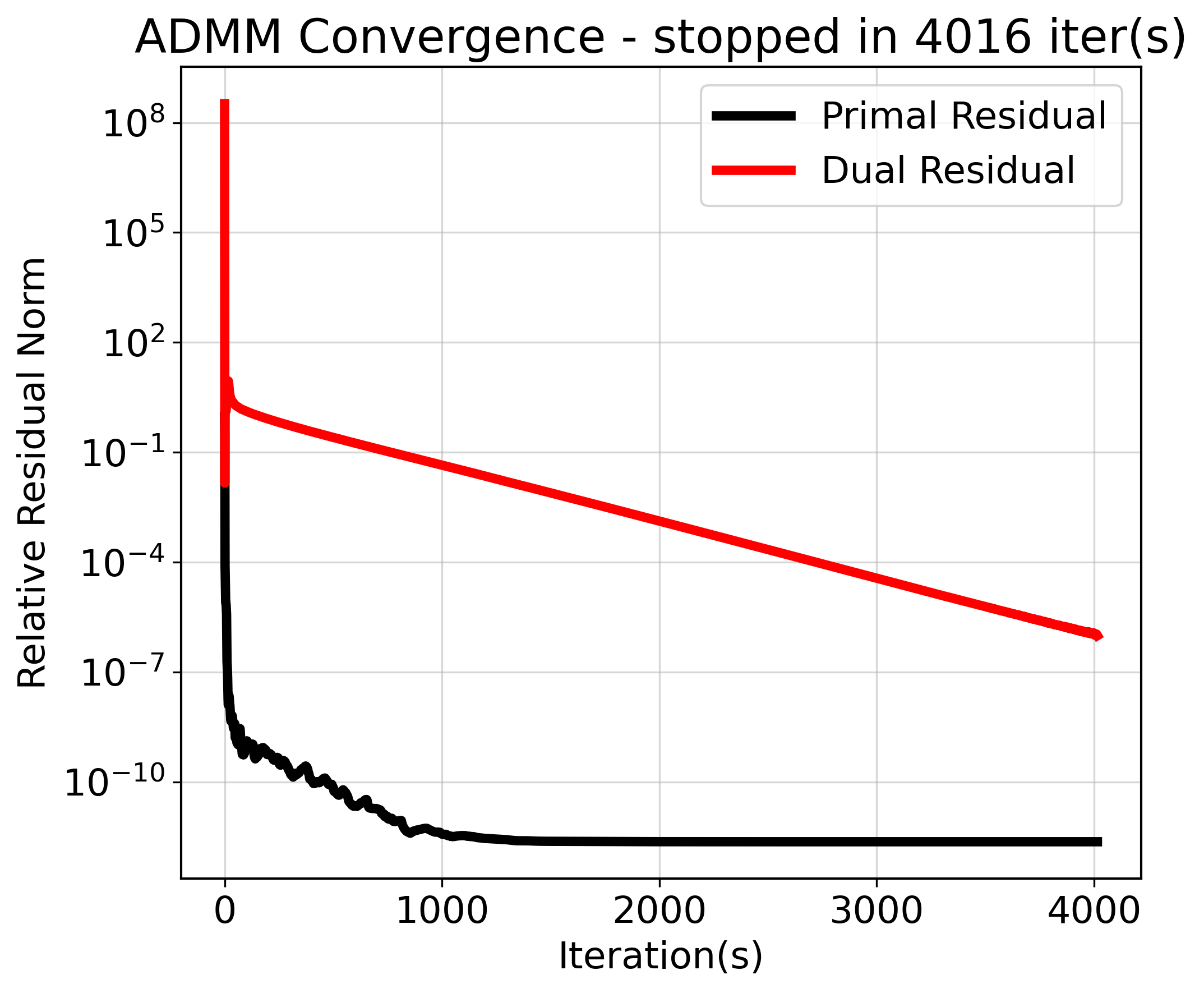}
        \caption{Relative primal and dual residuals showing convergence to the tolerance $\epsilon_{\mathrm{tol}} =  10^{-6}$ in 4,016 iterations.}
    \end{subfigure}
    \hfill
    \begin{subfigure}[t]{0.48\linewidth}
        \centering
        \includegraphics[width=\linewidth]{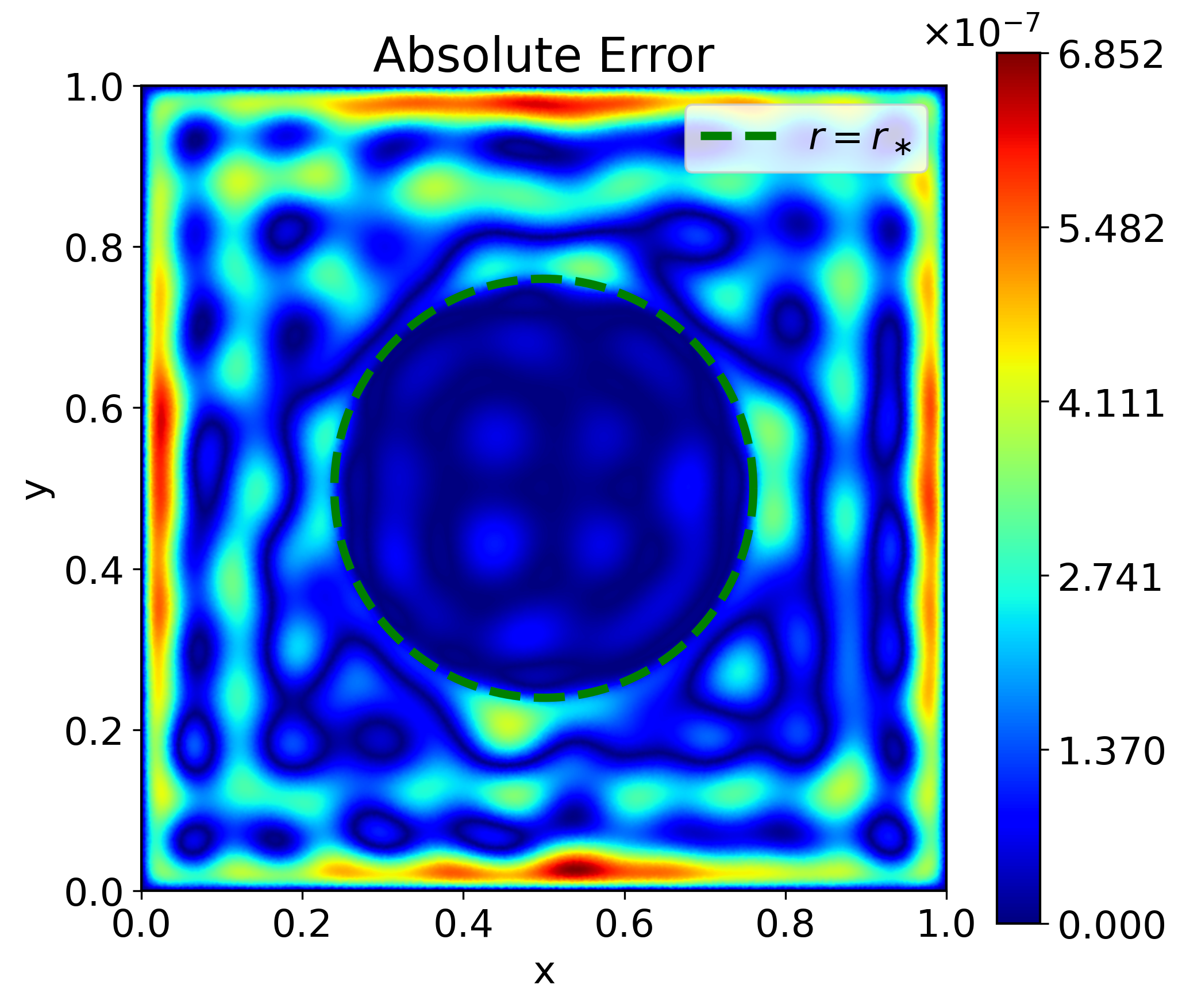}
        \caption{Pointwise error $\lvert \hat{u}-u\rvert$.}
    \end{subfigure}

\caption{\textbf{Meshfree RBF approximation with ADMM solver for a radially symmetric obstacle problem on $\overline{\Omega}=[0,1]^2$ with $N=400$ basis functions and $m = 1,600$.} The top-left panel shows the approximation $\hat{u}(x,y)$ obtained by ADMM optimization using the hyper-parameters $\mu = 10/h^2$, $\rho = 20.3$, and $\beta = 10^6$. The top-right panel shows the exact solution $u_{\mathrm{exact}}(x,y)$. The bottom-left panel shows the relative primal and dual residuals, indicating convergence to the tolerance $\epsilon_{\mathrm{tol}} = 10^{-6}$ in 4,016 iterations. The bottom-right panel shows the pointwise error $|\hat{u}-u|$; the relative $\ell^2$ error between the approximation and the analytical solution is $E(N)=4.681\times 10^{-11}$.}
    \label{fig1_obstacle2d_dome}
\end{figure}

\begin{figure}[H]
\centering
\begin{subfigure}[t]{0.5\textwidth}
\centering\includegraphics[width = 1.0 \linewidth]{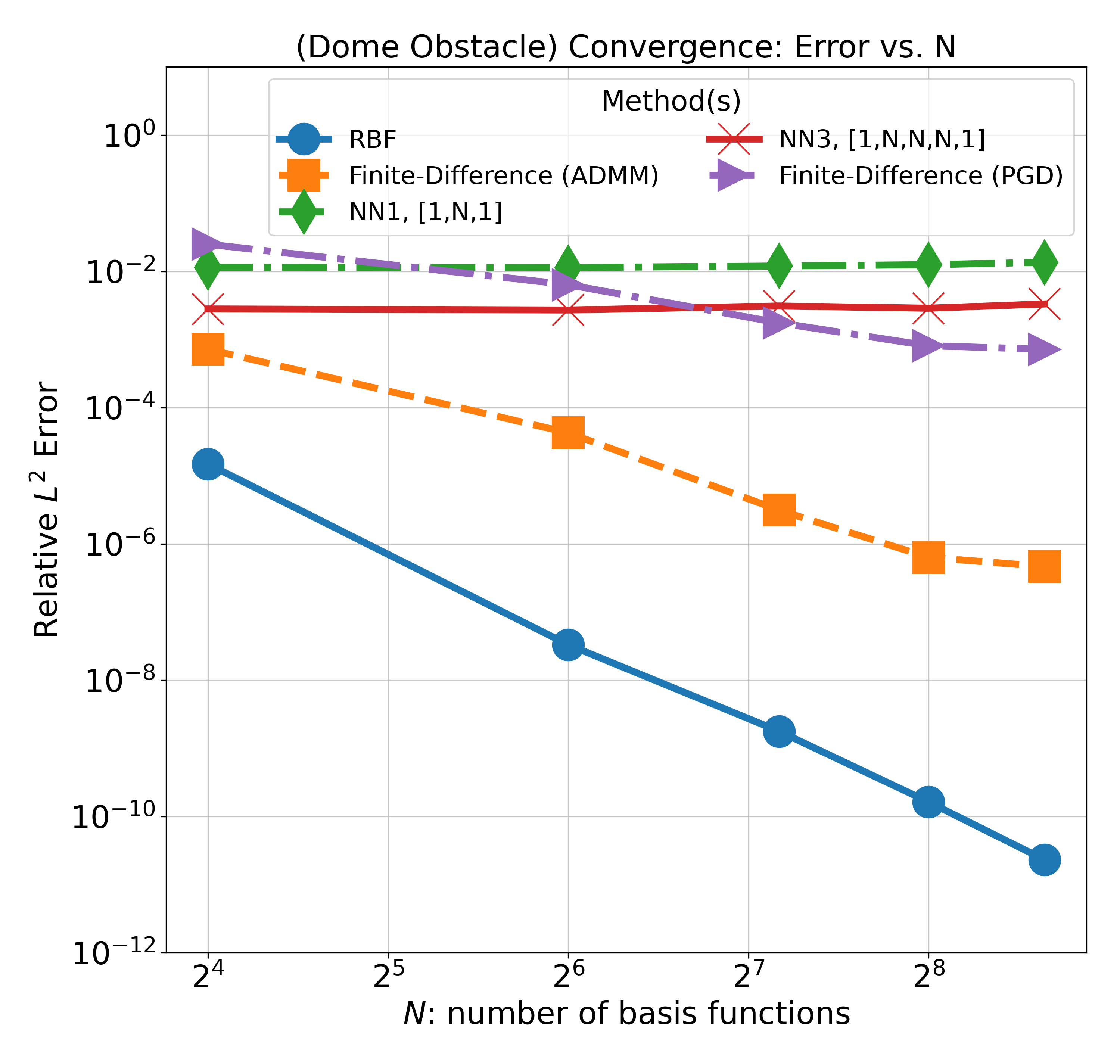}
\caption{Convergence: relative \(\ell^2\) error vs. \(N\)}
\end{subfigure}\hfill 
\begin{subfigure}[t]{0.5\textwidth}
\centering\includegraphics[width = 1.0 \linewidth]{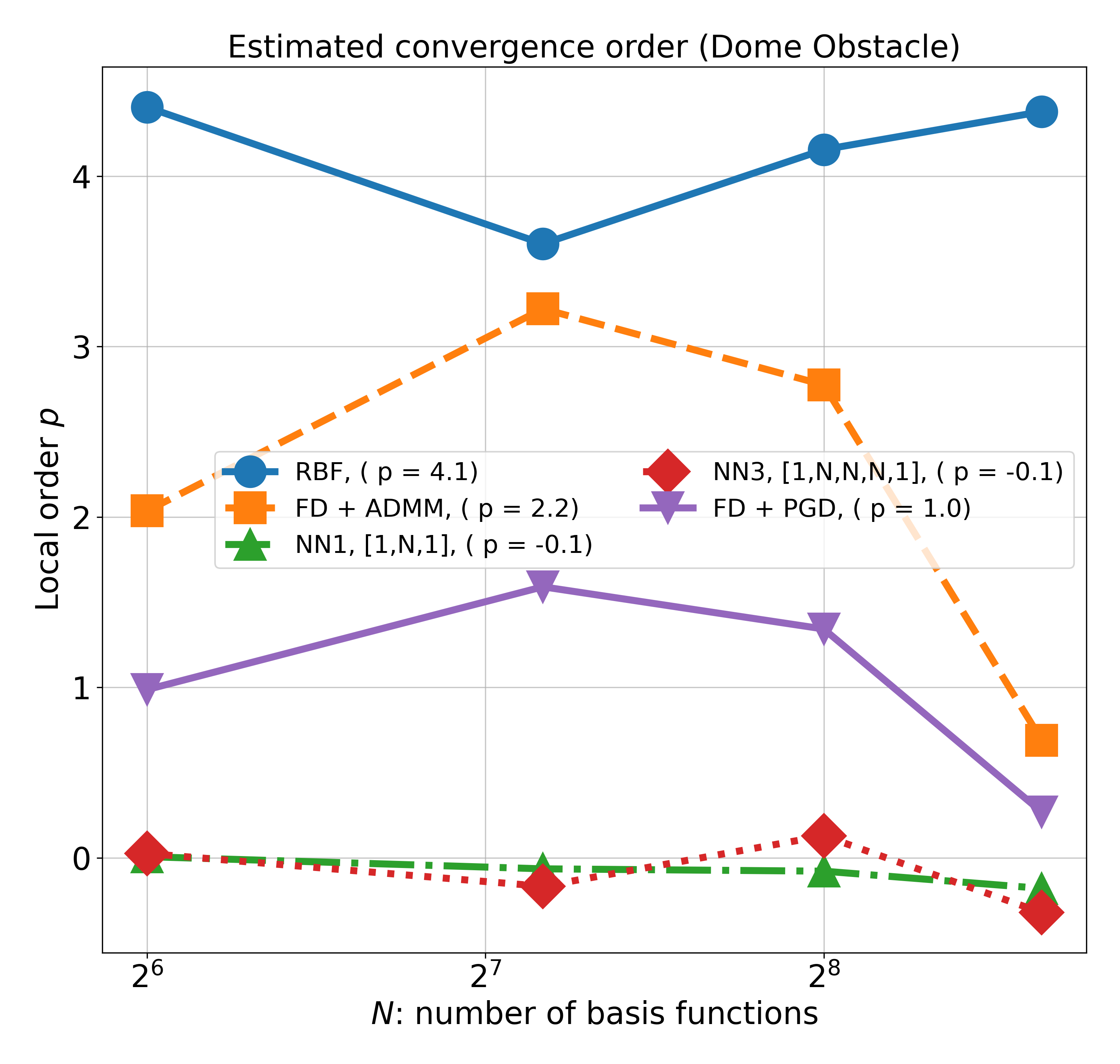}
\caption{Local estimated convergence order \(p\)}
\end{subfigure}
\caption{Convergence comparison for the dome obstacle problem: (a) Relative $\ell^2$ error versus the number of basis functions $N$ across different methods. (b) Local estimated convergence order $p$ for the corresponding methods.}
\label{fig:2d_obstacle_problem_dome_benchmark_method_solvers}
\end{figure}

As shown in Fig.~\ref{fig1_obstacle2d_dome}, the RBF-ADMM from Section~\ref{sec:methodology} reproduces the solution with high fidelity after $4,016$ iterations. The relative $\ell^2$ error is $4.681 \times 10^{-11}$ confirming the recovery of the exact solution $u_{exact}(x,y)$ in Eq.~\eqref{eq:exact_solution_2d_obstacle_problem_dome}. Both the relative primal and dual residuals meet the tolerance $\epsilon_{\mathrm{tol.}} = 10^{-6}$ showing convergence of the ADMM. The absolute point-wise error stays below $10^{-6}$ in the interior and near the boundary. In this problem, we choose \(T = 3\), which will be justified in Section~\ref{sec:comparison-expansion-domain-factor-T}. The boundary penalty weight is set at \(\beta = 10^{6}\), which is also discussed later in Appendix~\ref{appendix:experiment_boundary_penalty_parameter}. Fig.~\ref{fig:2d_obstacle_problem_dome_benchmark_method_solvers} shows that the variational RBF method performs substantially better than finite-difference and neural-network methods, achieving fast convergence and low computational errors. 

\pagebreak

\subsection{Sensitivity Analysis of Expansion Domain Factor and Radial Basis Functions}\label{sec:sensitivity-analysis-T-rbf-kernel}

\hspace{5mm}In this section, we investigate the sensitivity of the proposed RBF method with respect to some key parameter choices: the expansion domain factor \(T\) and the choice of radial basis function \(\phi\) in Table~\ref{tab:rbf-types}. Both parameters have a direct influence on the numerical stability and approximation quality of the method. The expansion factor \(T\) determines how far the RBF centers are sampled beyond the physical domain, and hence controls the trade-off between reducing boundary effects and preserving local resolution near the domain of interest. The kernel function, on the other hand, determines the smoothness, support, and approximation behavior. Therefore, a systematic comparison of these parameters is necessary to identify robust configurations for the Poisson and obstacle problems considered in this work. For other parameters including penalty weight \(\beta\) and truncation threshold value \(\tau\), we justify choices in Appendix~\ref{appendix:sect_experiments_parameter_settings}. In particular, for the 1D Poisson equation we use the penalty weight for boundary condition \(\beta = 3 \times 10^{5}\). For the 1D obstacle problems with one or two bumps, we choose \(\beta =  10^{6}, \,  10^{8}\), respectively. Lastly, in the 2D dome obstacle problem we let \(\beta  =  10^{6}\). For all experiments, we choose the truncation threshold value at \(\tau =  10^{-15}\). 

\subsubsection{Comparison of Domain Expansion Factor \(T\)}
\label{sec:comparison-expansion-domain-factor-T}

\begin{figure}[h!]
    \centering
    \begin{subfigure}[t]{0.48\textwidth}
        \centering
        \includegraphics[width=\linewidth]{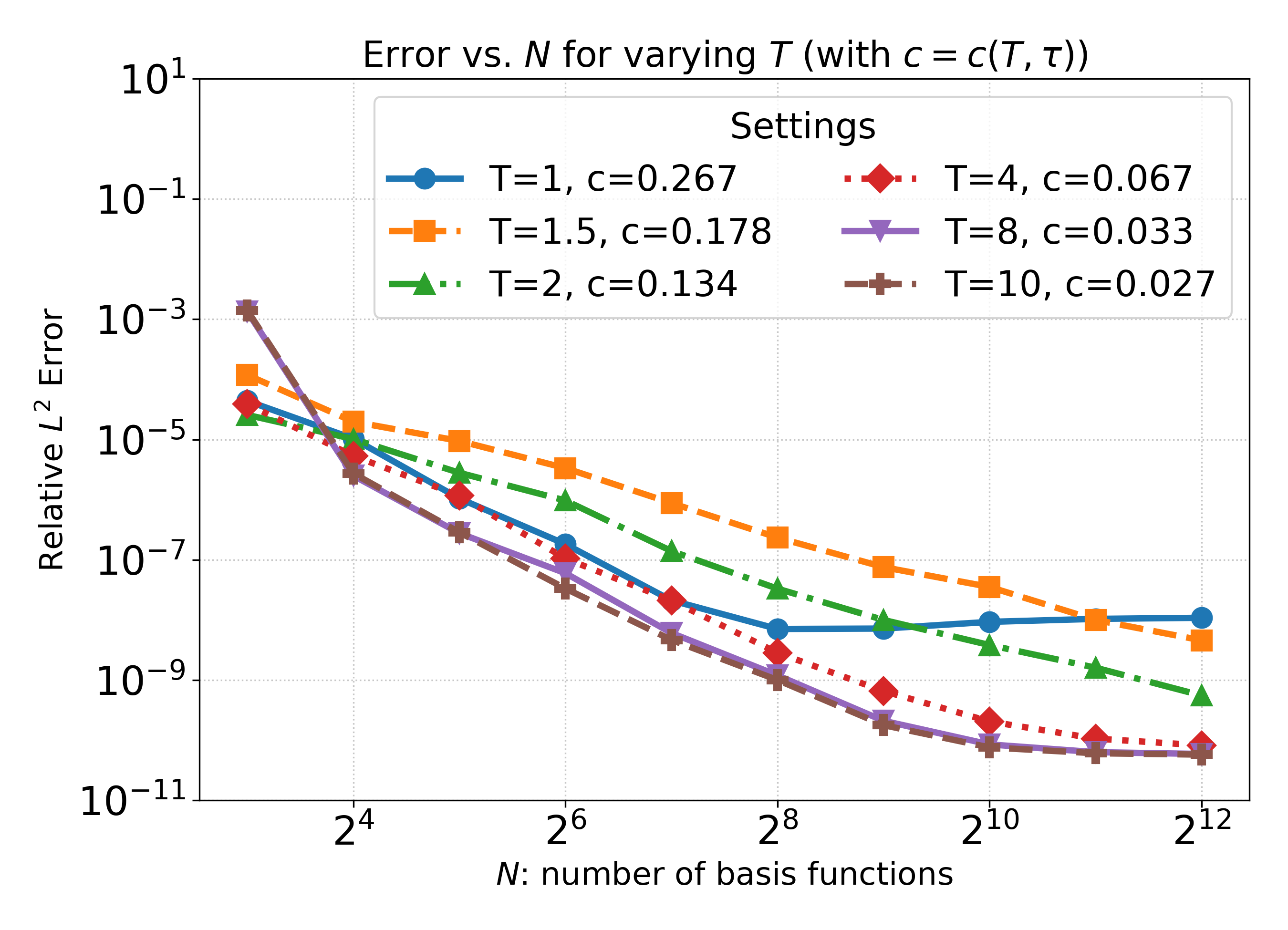}
        \caption{1D Poisson problem, \(\beta = 3 \times 10^{5}\).}
    \end{subfigure}\hfill
    \begin{subfigure}[t]{0.48\textwidth}
        \centering
        \includegraphics[width=\linewidth]{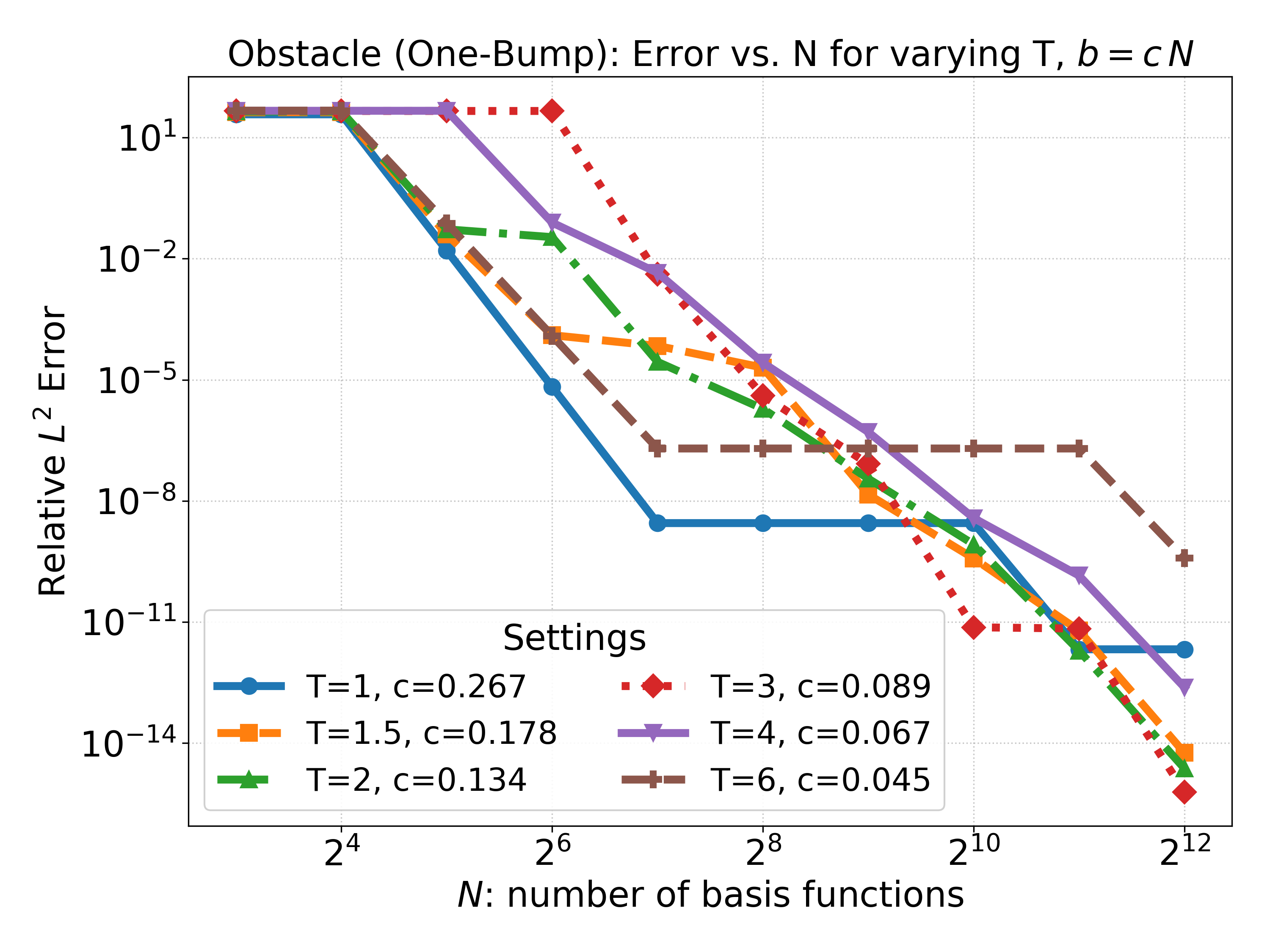}
        \caption{1D one-bump obstacle problem, \(\beta = 10^{6}\).}
    \end{subfigure}

    \vspace{0.5em}

    \begin{subfigure}[t]{0.48\textwidth}
        \centering
        \includegraphics[width=\linewidth]{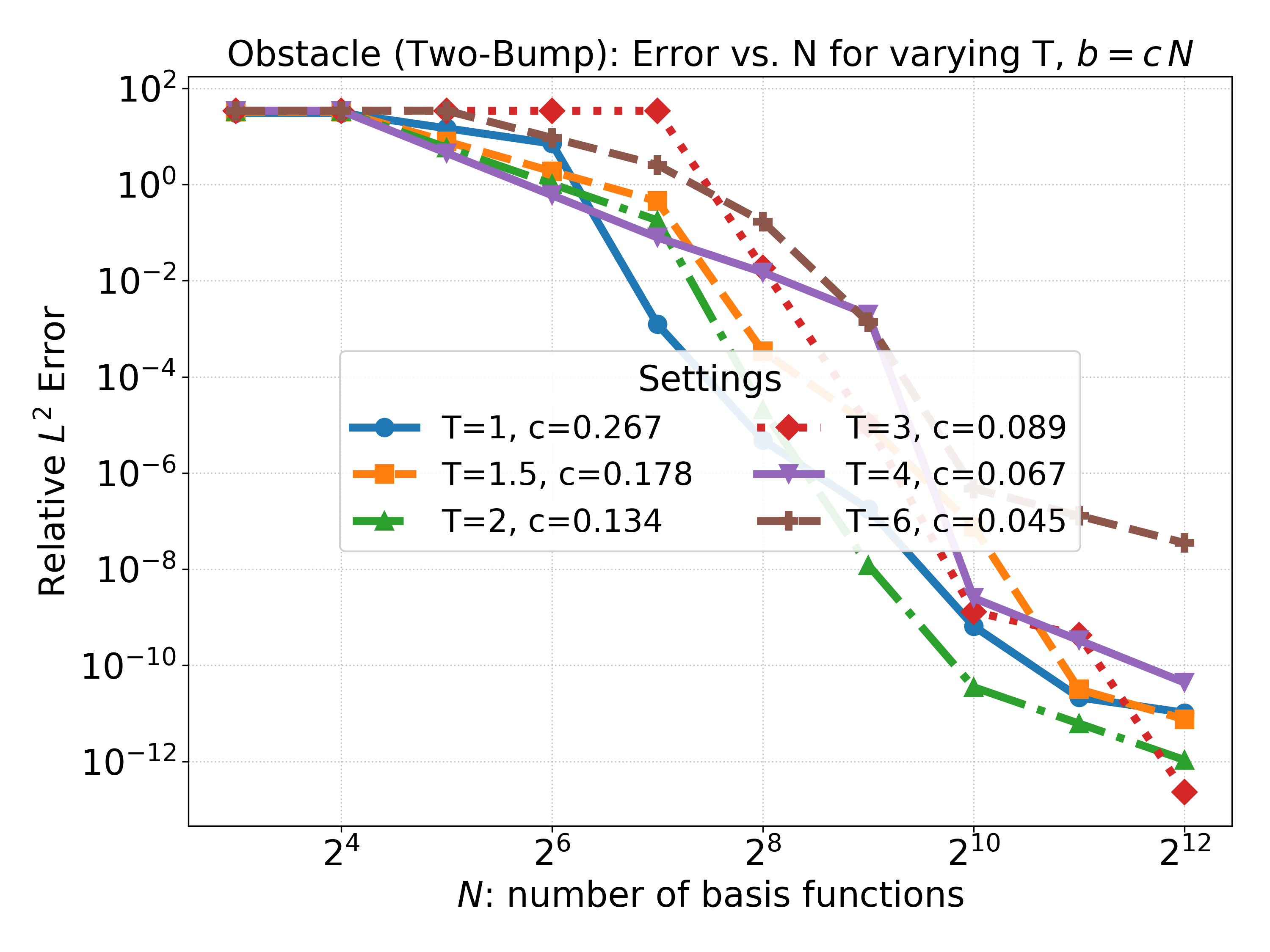}
        \caption{1D two-bump obstacle problem, \(\beta =  10^{8}\).}
    \end{subfigure}\hfill
    \begin{subfigure}[t]{0.48\textwidth}
        \centering
        \includegraphics[width=\linewidth]{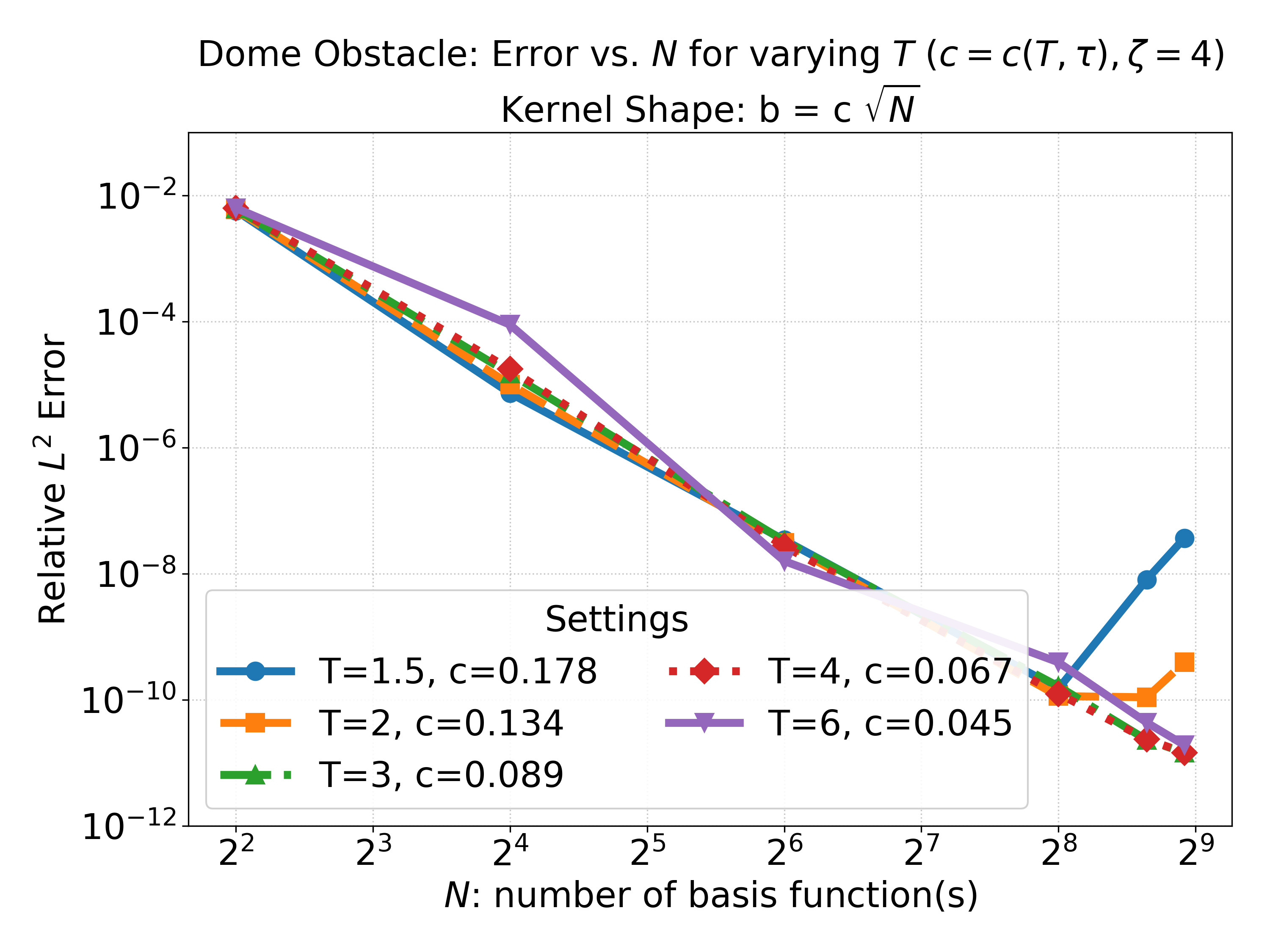}
        \caption{2D dome obstacle problem, \(\beta =  10^{6}\).}
    \end{subfigure}

    \caption{\textbf{Comparison of the domain expansion factor \(T\).}
    Relative \(\ell^{2}\) error versus \(N\) for several choices of the expansion factor \(T\) across different test problems. In all experiments, we calculate \(c = c(T, \tau)\) using the formula provided in Table~\ref{tab:params-numerical-experiments}. 
    In (a), the 1D Poisson problem is shown with \(\zeta = 2,\, \tau = 10^{-15}, \beta = 3 \times 10^{5}\). 
    In (b) and (c), the 1D obstacle problems are shown with \(\zeta = 2\), \(\tau = 10^{-15}\), and \(\beta\) is chosen \( 10^{6}\) and \( 10^{8}\), respectively. 
    In (d), the 2D dome obstacle problem is shown with \(\zeta = 4\), \(\tau = 10^{-15}\) and \(\beta =  10^{6}\).}
    \label{fig:sensitivity_T_2x2}
\end{figure}

\hspace{5mm}Figure~\ref{fig:sensitivity_T_2x2} suggests that the preferred value of \(T\) depends on the underlying problem. For the 1D Poisson problem, we use \(T=8\), since it gives a more stable error decay than the smaller expansion factors and remains among the best-performing choices over the tested range of \(N\). For the 1D one-bump and two-bump obstacle problems, we use \(T=2, \, 3\) respectively, because this value captures the main improvement from domain expansion while maintaining good approximation quality for the localized obstacle structure. Increasing \(T\) further does not produce a uniformly better behavior in these two examples. For the 2D dome obstacle problem, we use \(T=3\), which gives one of the most stable error profiles among the tested values and provides a good compromise between reducing boundary effects and preserving effective resolution near the physical domain. 

\subsubsection{Comparison of Kernel Functions \(\phi\)}
\label{sec:comparison-kernel-rbf-functions}
\begin{figure}[h!]
    \centering
    \begin{subfigure}[t]{0.48\textwidth}
        \centering
        \includegraphics[width=\linewidth]{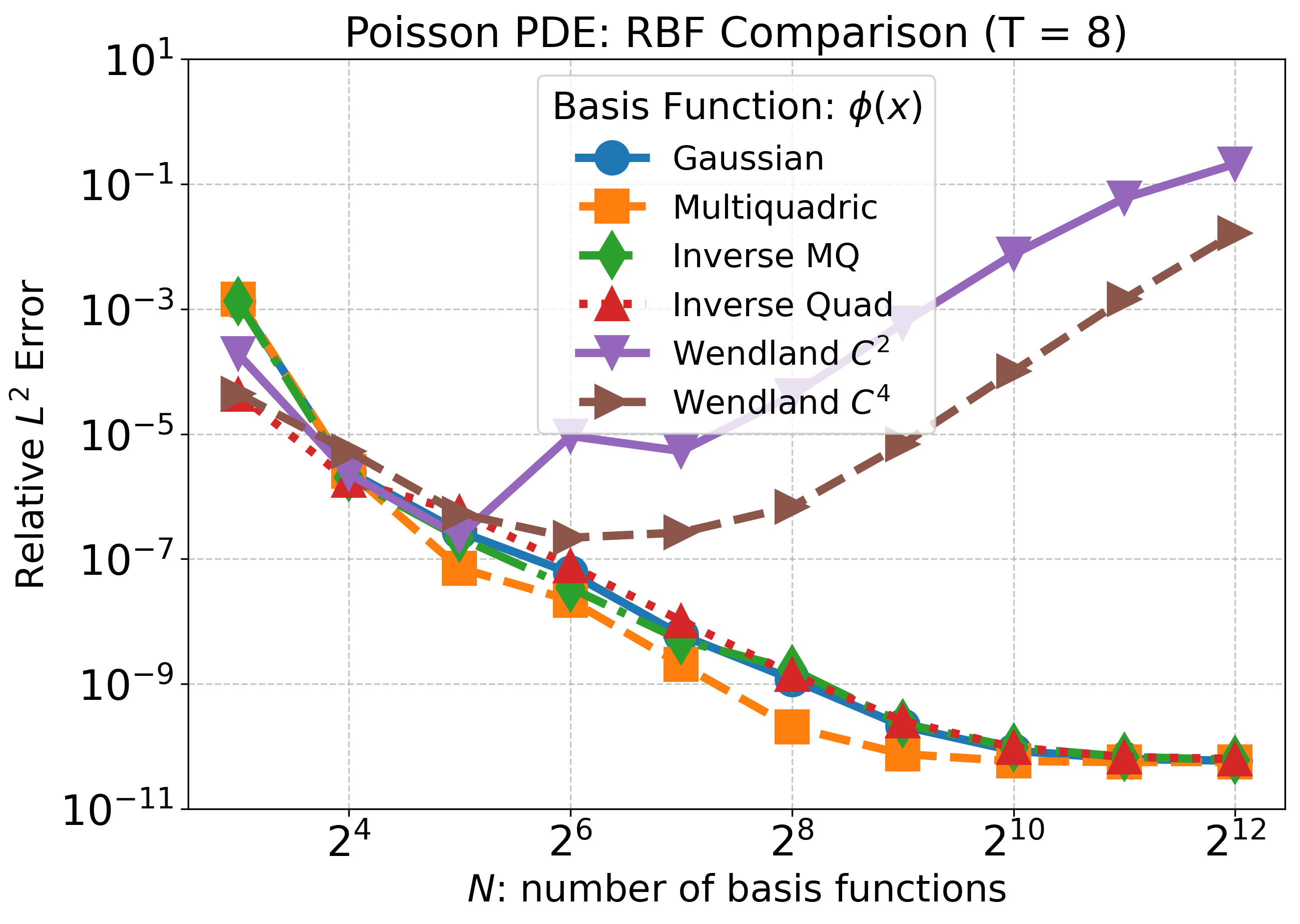}
        \caption{1D Poisson problem, \(\beta =  10^{6}\)}
    \end{subfigure}\hfill
    \begin{subfigure}[t]{0.48\textwidth}
        \centering
        \includegraphics[width=\linewidth]{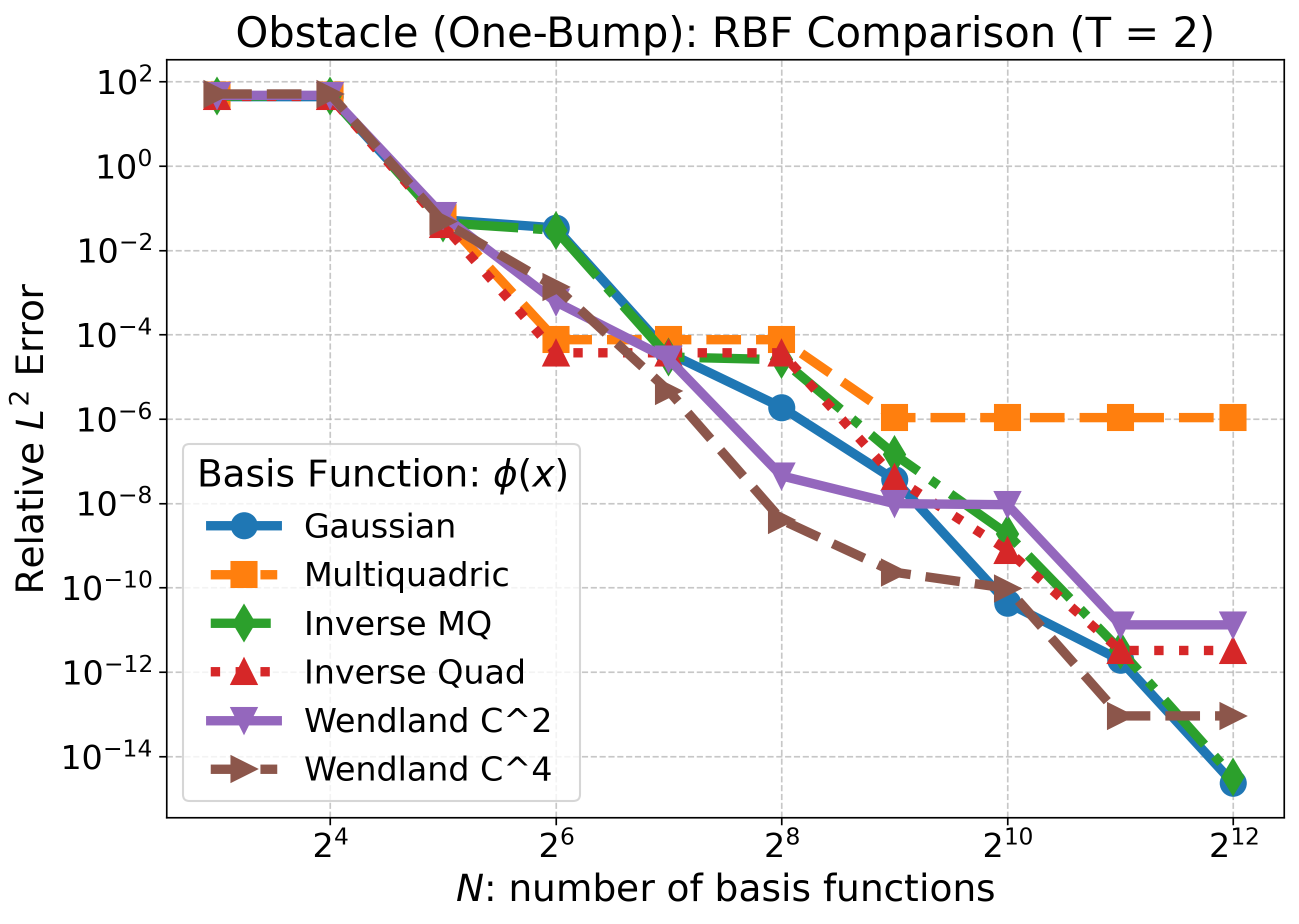}
        \caption{1D one-bump obstacle problem, \(\beta = 10^{6}\)}
    \end{subfigure}

    \vspace{0.5em}

    \begin{subfigure}[t]{0.48\textwidth}
        \centering
        \includegraphics[width=\linewidth]{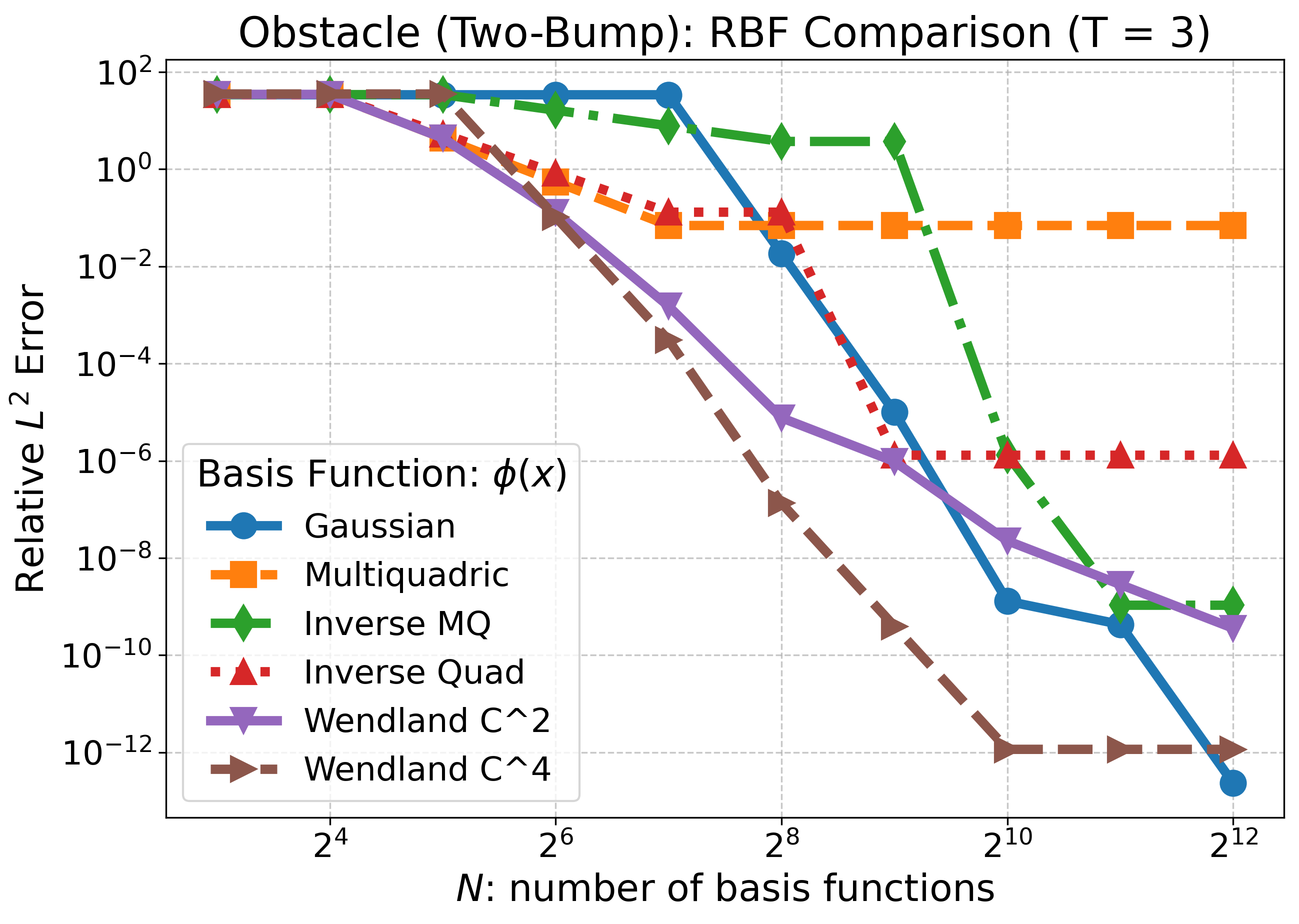}
        \caption{1D two-bump obstacle problem, \(\beta = 10^{8}\)}
    \end{subfigure}\hfill
    \begin{subfigure}[t]{0.48\textwidth}
        \centering
        \includegraphics[width=\linewidth]{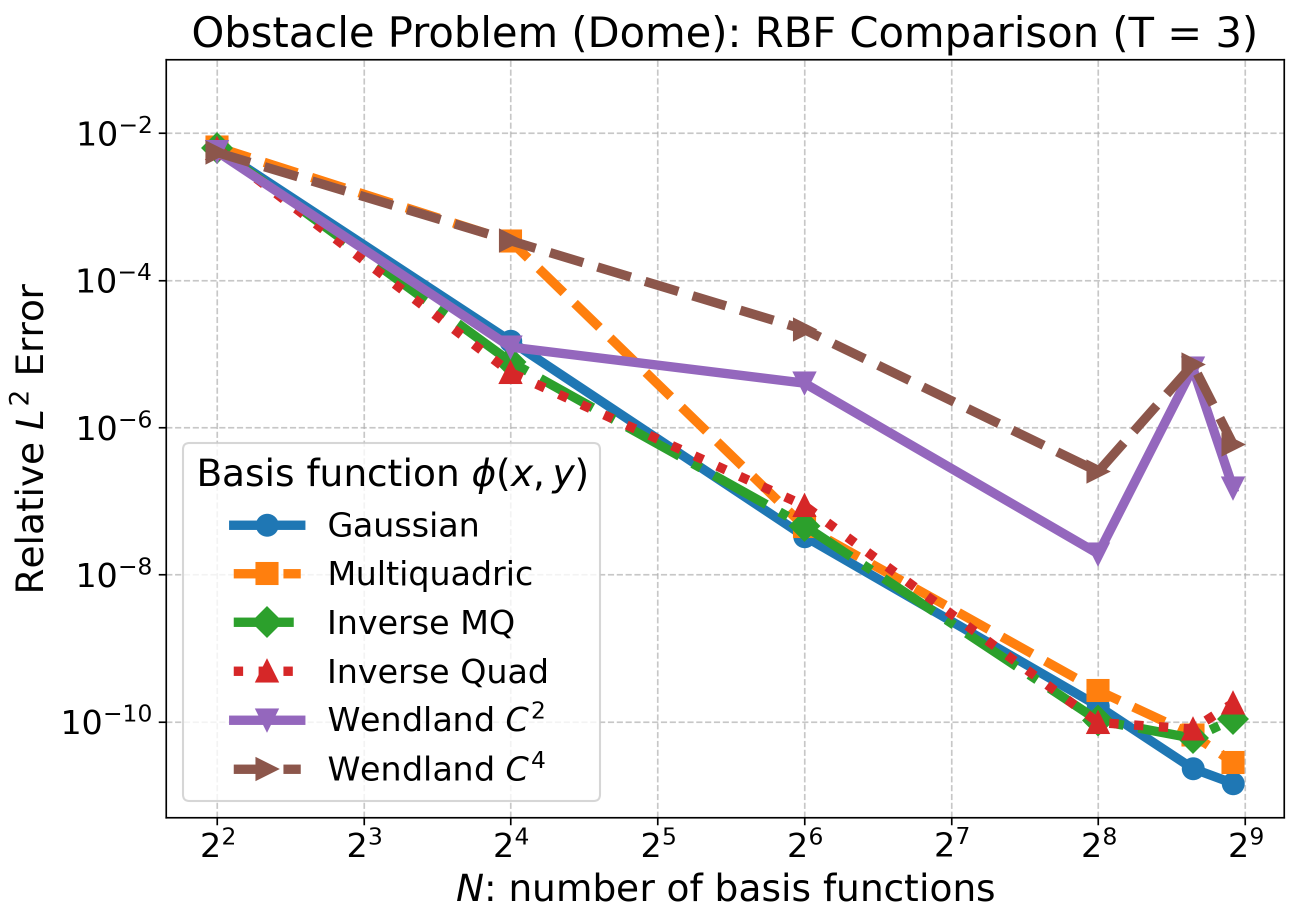}
        \caption{2D dome obstacle problem, \(\beta =  10^{6}\)}
    \end{subfigure}

    \caption{\textbf{Comparison of kernel functions.}
    Relative \(\ell^{2}\) error versus \(N\) for the radial basis functions listed in Table~\ref{tab:rbf-types}. 
    In (a), the 1D Poisson problem is shown with \(\zeta = 2\), \(T = 8\), and \(\tau = 10^{-15}\). 
    In (b) and (c), the 1D obstacle problems are shown with \(T=2.0\) and \(3.0\), \(\zeta = 2\), and \(\tau = 10^{-15}\), using \(\beta = 10^{6}\) and \(\beta = 10^{8}\), respectively. 
    In (d), the 2D dome obstacle problem is shown with \(T=3\), \(\tau = 10^{-15}\), \(\beta = 10^{6}\), and \(\zeta = 4\).}
    \label{fig:sensitivity_kernel_2x2}
\end{figure}

\hspace{5mm}The kernel comparison in Figure~\ref{fig:sensitivity_kernel_2x2} justifies the use of the Gaussian radial basis function in our experiments. Among all kernels in Table~\ref{tab:rbf-types}, the Gaussian kernel delivers among the most robust error decay and the best overall approximation accuracy across both the Poisson and obstacle problems. The Gaussian kernel is infinitely differentiable and globally supported, which makes it particularly effective for representing smooth solutions and for capturing global features of the variational solution space. In contrast, the compactly supported Wendland kernels are less smooth and more localized, so although they may offer computational advantages in other contexts, they are less competitive here in terms of approximation accuracy. Consequently, the Gaussian kernel is justified as the default kernel in this work.

\section{Conclusion}\label{sec:conclusion}

We have presented a mesh-free variational framework to solve the elliptic boundary value problem and obstacle-type variational problems using radial basis function (RBF) network approximation according to Eq.~\eqref{eq:main-approximation-form-rbf}. To address ill-conditioning inherent in global RBF discretizations, we incorporated truncated singular value decomposition (TSVD) as a stabilization mechanism in the least-squares solver.

We extend the methodology to solve classical obstacle problems according to Eq.~\eqref{obstacle-pde-bvp} by introducing a variational formulation with $ \ell^1$-penalty on the contact constraint. The resulting non-smooth optimization problem is efficiently solved using the alternating direction method of multipliers (ADMM), as described in Section~\ref{sec:constrained-optimization}. In our methodology of Section~\ref{sec:methodology}, we enlarge the effective basis support through the expansion factor $T$. The shape parameters is chosen using linear scaling $b = c(T, \tau) \, N$ in one dimension and $b = c(T, \tau) \, \sqrt{N} $ in two dimensions, which leads to stable convergence and high accuracy. Extensive numerical experiments in Section~\ref{sec:numerical-experiments} confirm that our RBF solution method achieves fast convergence for elliptic and obstacle problems. The results offer pragmatic guidance for selecting the $(N, T, \tau)$ parameters to balance accuracy, condition number, and computational stability in mesh-free variational solvers. 

Future research directions include: 

\begin{itemize}
  \item Generalizing the framework to high‐dimensional problems including adaptive center distributions.
  \item Investigating alternative penalty formulations (e.g. Huber or augmented Lagrangian methods) in order to improve conditioning. 
  \item Developing a systematic theoretical characterization of the trade-offs between hyper-parameters $(N, T, \tau)$, in terms of accuracy and conditioning, by defining a relation like Eq.~\eqref{eq:optimal-value-c-formula-from-adcock} specifically tailored to our variational RBF method in multiple dimensions.  
  \item Extending the framework to other nonlinear PDEs and time-dependent variational problems. 
\end{itemize}

\paragraph*{Data and Code Availability}
All coding scripts used to produce the numerical experiments are openly available at
\url{https://github.com/phuongdongle/Variational-Formulation-RBFNets.git}.

\bibliography{references} 

@book{buhmann2003radial,
  title={Radial Basis Functions: Theory and Implementations},
  author={Buhmann, Martin},
  year={2003},
  publisher={Cambridge University Press}
}

@book{wendland2005scattered,
  title={Scattered Data Approximation},
  author={Wendland, Holger},
  year={2005},
  publisher={Cambridge University Press}
}

@incollection{micchelli1984interpolation,
  title={Interpolation of scattered data: distance matrices and conditionally positive definite functions},
  author={Micchelli, Charles},
  booktitle={Approximation theory and spline functions},
  pages={143--145},
  year={1984},
  publisher={Springer}
}

@article{powell1987radial,
  title={Radial basis functions for multivariable interpolation: a review},
  author={Powell, Michael James David},
  journal={Algorithms for approximation},
  pages={143--167},
  year={1987}
}

@book{fasshauer2007meshfree,
  title={Meshfree Approximation Methods with MATLAB},
  author={Fasshauer, Gregory Eric},
  year={2007},
  publisher={World Scientific}
}

@article{schaback1995error,
  title={Error estimates and condition numbers for radial basis function interpolation},
  author={Schaback, Robert},
  journal={Advances in Computational Mathematics},
  volume={3},
  number={3},
  pages={251--264},
  year={1995},
  publisher={Springer}
}

@book{friedman1982variational,
  title={Variational Principles and Free-Boundary Problems},
  author={Friedman, Avner},
  year={1982},
  publisher={Wiley}
}

@article{adcock2019frames,
  title={Frames and numerical approximation II: Generalized sampling},
  author={Adcock, Ben and Huybrechs, Daan},
  journal={Journal of Fourier Analysis and Applications},
  volume={26},
  number={6},
  pages={87},
  year={2020},
  publisher={Springer}
}

@article{adcock2024stable,
  title={Stable and accurate least squares radial basis function approximations on bounded domains},
  author={Adcock, Ben and Huybrechs, Daan and Piret, {Cécile}},
  journal={SIAM Journal on Numerical Analysis},
  volume={62},
  number={6},
  pages={2698--2718},
  year={2024},
  publisher={SIAM}
}

@article{tran20151,
  title={An {L\^{}1} penalty method for general obstacle problems},
  author={Tran, Giang and Schaeffer, Hayden and Feldman, William M and Osher, Stanley},
  journal={SIAM Journal on Applied Mathematics},
  volume={75},
  number={4},
  pages={1424--1444},
  year={2015}
}

@article{zhang2011unified,
  title={A unified primal-dual algorithm framework based on {B}regman iteration},
  author={Zhang, Xiaoqun and Burger, Martin and Osher, Stanley},
  journal={Journal of Scientific Computing},
  volume={46},
  pages={20--46},
  year={2011}
}

@article{boyd2011distributed,
  title={Distributed optimization and statistical learning via the alternating direction method of multipliers},
  author={Boyd, Stephen and Parikh, Neal and Chu, Eric and Peleato, Borja and Eckstein, Jonathan},
  journal={Foundations and Trends in Machine learning},
  volume={3},
  number={1},
  pages={1--122},
  year={2011}
}

@article{ros2018obstacle,
  title={Obstacle problems and free boundaries: an overview},
  author={Ros-Oton, Xavier},
  journal={SeMA Journal},
  volume={75},
  number={3},
  pages={399--419},
  year={2018}
}

@article{adcock2015tsvd,
  title={Frames and numerical approximation},
 author={Adcock, Ben and Huybrechs, Daan},
  journal={SIAM Review},
  volume={61},
  number={3},
  pages={443--473},
  year={2019},
  publisher={SIAM}
}

@book{glowinski1984numerical,
  title={Numerical Methods for Nonlinear Variational Problems},
  author={Glowinski, Roland},
  year={1984},
  publisher={Springer-Verlag}
}

@article{yu2018deep,
  title={The {Deep Ritz} method: a deep learning-based numerical algorithm for solving variational problems},
  author={Bing Yu and Weinan E},
  journal={Communications in Mathematics and Statistics},
  volume={6},
  number={1},
  pages={1--12},
  year={2018},
  publisher={Springer}
}

@article{park1991universal,
  title={Universal approximation using radial-basis-function networks},
  author={Park, Jooyoung and Sandberg, Irwin W},
  journal={Neural computation},
  volume={3},
  number={2},
  pages={246--257},
  year={1991},
  publisher={MIT Press}
}

@article{narcowich2006sobolev,
  title={Sobolev error estimates and a {B}ernstein inequality for scattered data interpolation via radial basis functions},
  author={Narcowich, Francis J and Ward, Joseph D and Wendland, Holger},
  journal={Constructive Approximation},
  volume={24},
  number={2},
  pages={175--186},
  year={2006},
  publisher={Springer}
}

@article{zosso2015fast,
  title={An efficient primal-dual method for the obstacle problem},
  author={Zosso, Dominique and Osting, Braxton and Xia, Mandy and Osher, Stanley J},
  journal={Journal of Scientific Computing},
  publisher = {Springer}, 
  pages={416-437},
  volume = {73}, 
  year={2017}
}

@article{raissi2017physics,
  title={Physics informed deep learning (part i): Data-driven solutions of nonlinear partial differential equations},
  author={Raissi, Maziar and Perdikaris, Paris and Karniadakis, George Em},
  journal={arXiv preprint arXiv:1711.10561},
  year={2017}
}

@misc{raissi2017physicsinformeddeeplearning,
      title={Physics Informed Deep Learning (Part II): Data-driven Discovery of Nonlinear Partial Differential Equations}, 
      author={Maziar Raissi and Paris Perdikaris and George Em Karniadakis},
      year={2017},
      eprint={1711.10566},
      archivePrefix={arXiv},
      primaryClass={cs.AI},
      url={https://arxiv.org/abs/1711.10566}, 
}

@article{raissi2018deep,
  title={Deep hidden physics models: Deep learning of nonlinear partial differential equations},
  author={Raissi, Maziar},
  journal={Journal of Machine Learning Research},
  volume={19},
  number={25},
  pages={1--24},
  year={2018}
}

@misc{dokken2025latent,
      title={The latent variable proximal point algorithm for variational problems with constraints},
      author={Dokken, {J\o rgen} S. and Farrell, Patrick~E. and Keith, Brendan and Papadopoulos, Ioannis~P.A. and Surowiec, Thomas~M.},
      year={2025},
}

@article{Ang_2024,
   title={{MGProx}: A Nonsmooth Multigrid Proximal Gradient Method with Adaptive Restriction for Strongly Convex Optimization},
   volume={34},
   ISSN={1095-7189},
   url={http://dx.doi.org/10.1137/23M1552140},
   DOI={10.1137/23m1552140},
   number={3},
   journal={SIAM Journal on Optimization},
   publisher={Society for Industrial & Applied Mathematics (SIAM)},
   author={Ang, Andersen and De Sterck, Hans and Vavasis, Stephen},
   year={2024},
   month=aug, pages={2788–2820} }

\pagebreak

\appendix

\section{Table of Notations}
\subsection{General Mathematical Notation}

\begin{table}[H]
\centering
\renewcommand{\arraystretch}{1.5}
\begin{tabular}{@{}ll@{}}
\hline
\textbf{Symbol} & \textbf{Description} \\
\hline
$\Omega$                & Open bounded domain in $\mathbb{R}^d$ \\
$\partial \Omega$       & Boundary of the domain $\Omega$ \\
$\overline{\Omega}$          & Closure of $\Omega$, i.e., $\Omega \cup \partial \Omega$ \\
$[b_{1}, b_{2} ] \subseteq \mathbb{R}$   & Closed domain in \(\mathbb{R}\) \\
$\overline{\Omega}_{T}$ & Expanded domain with expansion factor \(T \geq 1\) \\
$ \vec{x} \in \mathbb{R}^d$     & Point in Euclidean space \\
$u(\vec{x})$                  & Candidate or exact solution function \\
$f(\vec{x})$                  & Source term in PDE \\
$g(\vec{x})$                  & Boundary data function \\
$\nabla u$              & Gradient of the function $u$ \\
$\Delta u$              & Laplacian of $u$, i.e., $\nabla \cdot \nabla u$ \\
\hline
\end{tabular}
\caption{General notation for domain, solution properties, and differential operators.}
\label{tab:notation-general}
\end{table}

\subsection{Problem-Specific Notation}

\begin{table}[H]
\centering
\renewcommand{\arraystretch}{1.5}
\begin{tabular}{@{}ll@{}}
\hline
\textbf{Symbol} & \textbf{Description} \\
\hline
$\mathcal{J}(u)$        & Variational functional associated with the PDE \\
$\beta$                 & Penalty parameter for boundary condition enforcement \\
$\mu$                   & Penalty parameter for obstacle constraints \\
$\psi( \vec{x})$              & Obstacle function \\
$(\psi - u)_+$          & Contact gap for obstacle problem: $\max(0, \psi - u)$ \\

$\phi(r)$               & Radial basis function with shape parameter $b$ \\
$\vec{c}_{i} $                   & Radial basis function center points \\

$h_{\Omega}$ & The maximum distance from any point in $\Omega$ to its nearest center $\vec{c}_{j}$ \\

$q_{\Omega}$ & Half the minimal pairwise spacing of centers \\

$L_{\mathrm{half}}$ & Half-length of the closed domain \([b_{1}, b_{2}]\), i.e., \(L_{\mathrm{half}} = \frac{b_{2} - b_{1}}{2}\) \\

$q_{\mathrm{mid}}$   & Midpoint of the closed domain \([b_{1}, b_{2} ]\), i.e., \(q_{\mathrm{mid}} = \frac{b_{1} + b_{2} }{2}\)\\

$\hat{u}(\vec{x})$            & Radial basis function approximation to $u(x)$ \\

$\tau$ & Threshold for the singular values in the truncated SVD\\ 
$ \mathcal{I}_{\tau} $ & Index set of discarded singular values $\mathcal{I}_{\tau} \coloneq \{ i : \sigma_{i} \leq \sigma_{1} \, \tau \} $ \\ 
$k$ & Cardinality of the discarded index set $k = | \mathcal{I}_{\tau} |$ \\ 
$\kappa$ & Absolute condition number, satisfying $\kappa \leq \min \{ 1/\sigma_{\mathrm{min}}(\mathcal{I}_{\tau} ), 1/\tau \}$ \\ 
$\rho$                  & Augmented Lagrangian penalty parameter (ADMM) \\
$\mathrm{maxIter}$                      & Maximum number of iterations (ADMM)\\  
$\vec{v} \in \mathbb{R}^N$ & Auxiliary variable in optimization (e.g., for constraints) \\
$\vec{z} \in \mathbb{R}^N$ & Projection variable in operator splitting methods \\
\hline
\end{tabular}
\caption{Notation related to variational formulations, RBF approximation, and constraints.}
\label{tab:notation-problem}
\end{table}

\section{Derivation of Optimization Problem}\label{appendix:derivation-optimization-form}

We outline the steps of deriving the minimizations in Section~\ref{sec:methodology} in the following.

\subsection{Linear Elliptic Boundary Value Problem}\label{appendix:subsection-derivation-optimization-elliptic}

We recall the weak-form variational functionals that correspond to the Poisson and reaction-diffusion PDEs, see Eq.~\eqref{weak-form-poisson-problem}. Let $m$ be the total number of collocation points for domain $\overline{\Omega}$ including \(m_{\Omega}\) interior and $m_{\partial \Omega}$ boundary points. We take a set of \(N\) uniformly spaced centers for the basis kernels $\phi_{i}$ for all $1 \leq i \leq N$. By Eq.~\eqref{eq:main-approximation-form-rbf}, we have the following: 

\[
\nabla \hat{u}(\vec{x}) = \frac{1}{\sqrt{N}} \sum\limits_{i=1}^N w_{i} \nabla \phi ( \| \vec{x} - \vec{c}_{i} \|_{2}) \in \mathbb{R}^{d}. 
\]

Defining the center matrix \(C\) as

\[
C \coloneq \begin{bmatrix}
\quad \vec{c}_{1} \quad | &   \vec c_{2} \quad | & \cdots \quad  | &  \vec c_{N-1} \quad | & \vec c_{N} \quad 
\end{bmatrix} \in \mathbb{R}^{d \times N}, 
\]

we further define 

\[\Phi(\vec{x},C) \coloneq \begin{bmatrix}
\phi(\| \vec{x} - \vec{c}_{1} \|_{2}) \\ 
\vdots \\
\phi(\| \vec{x} - \vec{c}_{N} \|_{2})   
\end{bmatrix}\in \mathbb{R}^{N}, \quad \nabla\Phi(\vec{x}, C) \coloneq \begin{bmatrix}
\nabla \phi^{\top}(\| \vec{x} - \vec{c}_{1} \|_{2} ) \\ 
\vdots \\
\nabla \phi^{\top} (\|\vec{x} - \vec{c}_{N} \|_{2})  
\end{bmatrix}\in \mathbb{R}^{N\times d}. 
\]

We denote the matrix of sampled points. This works for any \(d \geq 1\) as follows: 

\[
X_{m_{\Omega}} = 
\begin{bmatrix}
\left( \vec{x}_{1}^{(\Omega)} \right)^{\top} \\
\vdots \\
\left(  \vec{x}_{m_{\Omega}}^{(\Omega)}   \right)^{\top} 
\end{bmatrix} \in \mathbb{R}^{m_{\Omega} \times d}. 
\]

Defining \(\phi_{i} (\vec{x}) = \phi(\| \vec{x} - \vec{c}_{i} \|_{2}) \), we also use the notation. 
\[ 
\Phi(X_{m_\Omega},C): = \begin{pmatrix}
\phi_1(\vec{x}_1^{(\Omega)}) & \ldots & \phi_N(\vec{x}_1^{(\Omega)}) \\
\phi_1(\vec{x}_2^{(\Omega)}) & \ldots & \phi_N(\vec{x}_2^{(\Omega)}) \\
\vdots& &\vdots \\
\phi_1(\vec{x}_{m_{\Omega}}^{(\Omega)}) & \ldots & \phi_N(\vec{x}_{m_{\Omega}}^{(\Omega)}) \\
\end{pmatrix} \in \mathbb{R}^{m_{\Omega} \times N}.\]

For \(d = 1\), we also use the notation

\[
\nabla \Phi(X_{m_\Omega},C): = \begin{pmatrix}
\phi_{1}'(\vec{x}_1^{(\Omega)}) & \ldots & \phi_{N}'(\vec{x}_1^{(\Omega)}) \\
\phi_{1}'(\vec{x}_2^{(\Omega)}) & \ldots & \phi_{N}'(\vec{x}_2^{(\Omega)}) \\
\vdots& &\vdots \\
\phi_{1}'(\vec{x}_{m_{\Omega}}^{(\Omega)}) & \ldots & \phi_{N}'(\vec{x}_{m_{\Omega}}^{(\Omega)}) \\
\end{pmatrix} \in \mathbb{R}^{m_{\Omega} \times N}.
\]

This yields that: 
\begin{equation}\label{eqn:LapMat_entrywise}
\| \nabla \hat{u}(\vec{x}) \|_{2}^{2} = \frac{1}{N} \sum_{i, l} w_{i} w_{l} \nabla \phi^{\top} (\| \vec{x} - \vec{c}_{i} \|_{2} ) \, \nabla \phi(\| \vec{x} - \vec{c}_{l} \|_{2})  = \frac{1}{N} \, \vec{w}^{\top} \nabla \Phi(\vec{x}, C) \, \nabla \Phi^{\top}(\vec{x}, C) \, \vec{w},\quad \text{where}\, 1 \leq i,\ell \leq N. 
\end{equation} 

We deduce that for $1 \leq i, \ell \leq N$: 

\begin{equation}\label{eqn:LapMat}
\begin{aligned}
\int_{\Omega} \frac{1}{2} \| \nabla \hat{u} (\vec{x}) \|_{2}^{2}  \, dx & \approx \frac{|\Omega|}{2 \, m_{\Omega}} \sum_{j = 1}^{m_{\Omega}} \| \nabla \hat{u}(\vec{x}^{(\Omega)}_{j}) \|_{2}^{2}  \\
&\approx \frac{|\Omega|}{2N  \, m_{\Omega}} \sum_{j = 1}^{m_{\Omega}}\sum_{i, l} w_{i} w_{l} \nabla \phi^{\top} (\| \vec{x}^{(\Omega)}_{j} - \vec{c}_{i} \|_{2} ) \, \nabla \phi(\| \vec{x}^{(\Omega)}_{j} - \vec{c}_{l} \|_{2})  \\ 
&\approx  \frac{|\Omega|}{2 \, N \, m_{\Omega}}  \vec{w}^{\top} \left( \sum_{j = 1}^{m_{\Omega}}  \nabla \Phi(\vec{x}^{(\Omega)}_{j}, C) \nabla \Phi^{\top}(\vec{x}^{(\Omega)}_{j}, C)   \right) \, \vec{w}\\
&\approx \frac{1}{2}\,\vec{w}^{\top}A^{\small(\mathrm{poisson})}_{1}\vec{w},
\end{aligned}
\end{equation}
where 
\begin{equation}\label{eqn:A1}
A_1^{\small(\mathrm{poisson})} =  \frac{| \Omega | }{ N m_{\Omega}} \left(\sum_{j = 1}^{m_{\Omega}}  \nabla \Phi(\vec{x}^{(\Omega)}_{j}, C) \nabla \Phi^{\top}(\vec{x}^{(\Omega)}_{j}, C) \right).
\end{equation}

Note that, for $d=1$, we can rewrite the matrix $A_1$ as follows:
\[
A_1^{\small(\mathrm{poisson})}=  \frac{| \Omega | }{ N m_{\Omega}} \left(\nabla \Phi(X_{m_\Omega},C) \right)^{\top}\nabla \Phi(X_{m_\Omega},C).
\]

In the following steps, we derive the load matrices and vectors for interior and boundary collocation points. We have
\[
\int_{\Omega} f \, \hat{u} \, dx \approx \frac{|\Omega|}{m_{\Omega}} \sum_{j} f(\vec{x}^{(\Omega)}_{j}) \, \hat{u}(\vec{x}^{(\Omega)}_{j}) \approx \frac{|\Omega|}{\sqrt{N} \, m_{\Omega}}  \vec{w}^{\top} \sum_{j = 1}^{m_{\Omega}} f(\vec{x}^{(\Omega)}_{j} ) \Phi(\vec{x}^{(\Omega)}_{j}, C) = [f(\vec{x}^{(\Omega)}_{1}),..., f(\vec{x}^{(\Omega)}_{m_{\Omega}})]^{T} A_{2} \vec{w} = \vec{f}^{\top} A_{2} \vec{w},
\]

 where 
\begin{equation}\label{eqn:A2}
A_2= \frac{|\Omega|}{\sqrt{N} \, m_{\Omega}} \begin{bmatrix}
\Phi( \vec{x}_{1}^{(\Omega)}, C)^{\top} \\ 
\vdots \\ 
\Phi(\vec{x}_{m_{\Omega}}^{(\Omega)}, C)^{\top}  
\end{bmatrix}  = \frac{| \Omega | }{\sqrt{N} \, m_{\Omega}}\, \Phi(X_{m_{\Omega}},C)\in\mathbb{R}^{m_\Omega\times N}\quad \text{and}\quad \vec{f} = \begin{bmatrix}
        f(\vec{x}^{(\Omega)}_{1}) \\ 
        \vdots \\ 
        f(\vec{x}^{(\Omega)}_{m_{\Omega}})
        \end{bmatrix} \in \mathbb{R}^{m_{\Omega}}.
\end{equation}

For $d = 1$, we get that: 

\begin{equation}
A_2 = \frac{| \Omega | }{\sqrt{N} \, m_{\Omega}} \, \Phi(X_{m_{\Omega}}, C) \in \mathbb{R}^{m_{\Omega} \times N} 
\end{equation}

Finally, the enforce of boundary condition yields that
\begin{equation}\label{eqn:bdry}
\begin{aligned}
\frac{\beta}{2} \int_{\partial \Omega} (\hat{u} - g)^{2} \, dx & \approx \frac{\beta \, |\partial \Omega|}{2 \, m_{\partial \Omega}} \sum_{l = 1}^{m_{\partial \Omega}} ( \hat{u}(\vec{x}^{(\partial\Omega)}_{l}) - g(\vec{x}^{(\partial\Omega)}_{l}) )^2\\
&\approx
\frac{\beta \, |\partial \Omega|}{2 \, m_{\partial \Omega}}  \, \sum\limits_{l=1}^{m_{\partial \Omega}} \Big( \frac{1}{\sqrt{N}} \vec{w}^{\top} \Phi(\vec{x}^{(\partial\Omega)}_{l}, C) - g(\vec{x}^{(\partial\Omega)}_{l}) \Big)^2\\
&\approx
\frac{\beta\,|\partial\Omega|}{2\,m_{\partial\Omega}}
\left\lVert
\begin{bmatrix}
\tfrac{1}{\sqrt{N}}\vec{w}^{\top}\Phi(\vec{x}^{(\partial\Omega)}_{1},C)-g(\vec{x}^{(\partial\Omega)}_{1})\\
\vdots\\
\tfrac{1}{\sqrt{N}}\vec{w}^{\top}\Phi(\vec{x}^{(\partial\Omega)}_{m_{\partial\Omega}},C)-g(\vec{x}^{(\partial\Omega)}_{m_{\partial\Omega}})
\end{bmatrix}
\right\rVert_{2}^{2}\\[6pt]
&\approx
\frac{\beta\,|\partial\Omega|}{2\,m_{\partial\Omega}}
\left\lVert
\tfrac{1}{\sqrt{N}}
\begin{bmatrix}
\Phi(\vec{x}^{(\partial \Omega)}_{1},C)^{\top}\\
\vdots\\
\Phi(\vec{x}^{(\partial \Omega)}_{m_{\partial\Omega}},C)^{\top}
\end{bmatrix}
\vec{w}
-
\begin{bmatrix}
g(\vec{x}^{(\partial \Omega)}_{1})\\
\vdots\\
g(\vec{x}^{(\partial \Omega)}_{m_{\partial\Omega}})
\end{bmatrix}
\right\rVert_{2}^{2}\\
&\approx
\frac{\beta}{2}\,\bigl\lVert A_{3}\vec{w}-\vec{r}_{3}\bigr\rVert_{2}^{2},
\end{aligned}
\end{equation}
where 
\begin{equation}\label{eqn:A3}
A_3= \sqrt{\frac{|\partial\Omega|}{N m_{\partial\Omega}}} \begin{bmatrix}
\Phi(\vec{x}^{(\partial \Omega)}_{1},C)^{\top}\\
\vdots\\
\Phi(\vec{x}^{(\partial \Omega)}_{m_{\partial\Omega}},C)^{\top}
\end{bmatrix} = \sqrt{\frac{| \partial \Omega |}{N \, m_{\partial \Omega}}} \, \Phi\left( X_{m_{\partial \Omega}}, C   \right)     \in\mathbb{R}^{m_{\partial\Omega} \times N}\quad\text{and}\quad   \vec{g} = \sqrt{ \frac{|\partial \Omega|}{m_{\partial \Omega}}} \, \begin{bmatrix}
        g(\vec{x}^{(\partial\Omega)}_{1}) \\ 
        \vdots \\
        g(\vec{x}^{(\partial\Omega)}_{m_{\partial \Omega}})  
    \end{bmatrix} \in \mathbb{R}^{m_{\partial \Omega}}.
\end{equation}

For reaction-diffusion, the terms are similar, except for stiffness matrix, which is given by 
    
    \begin{equation}\label{eqn:RDmat}
    A_{1} = \frac{|\Omega|}{N \, m_{\Omega}} \sum_{j = 1}^{m_{\Omega}} \Big( \nabla \Phi(\vec{x}^{(\Omega)}_{j}, C) \, \nabla \Phi^{\top}(\vec{x}^{(\Omega)}_{j}, C) + \Phi^{\top}(\vec{x}^{(\Omega)}_{j}, C) \, \Phi(\vec{x}^{(\Omega)}_{j}, C) \Big)  \in \mathbb{R}^{N \times N}. 
    \end{equation}

For the Poisson and reaction-diffusion equations, we combine those terms, which yields the discrete optimization problem
\begin{equation}
\min \limits_{\vec{w} \in \mathbb{R}^{N}} \, \Big\{ \frac{1}{2} \vec{w}^{\top} A_{1} \, \vec{w} - \vec{f}^{\top} A_{2} \, \vec{w} + \frac{\beta}{2} \| A_{3} \, \vec{w}  - \vec{g} \|_{2}^{2}\Big\}
\nonumber
\end{equation}

To solve this problem, we take the derivative with respect to $\vec{w}$, 

\[
\nabla \mathcal{\widehat{J}} = A_{1} \vec{w} - A_{2}^{\top} \vec{f} + \beta A_{3}^{\top} A_{3} \vec{w} - \beta A_{3}^{\top} \vec{g}, 
\]

and set $\nabla \mathcal{\widehat{J}} = 0$:  

\[
\Big( A_{1} + \beta A_{3}^{\top} A_{3} ) \vec{w} = \Big(\beta A_{3}^{\top} \vec{g} + A_{2}^{\top} \vec{f}\Big). 
\]

\subsection{Obstacle-Type Boundary Value Problem}\label{appendix:subsection-derivation-optimization-obstacle}

In the classical obstacle problem, it remains to consider the $\ell^1$-norm penalty in Eq.~\eqref{weak-form-obstacle-problem}:  

\[
\mu \int_{\Omega} (\psi - \hat{u})_{+} \, dx \approx \frac{\mu |\Omega|}{m_{\Omega}} \sum_{j} \max( \psi(\vec{x}^{(\Omega)}_{j} ) - \hat{u}(\vec{x}^{(\Omega)}_{j}), 0) \approx \frac{\mu |\Omega|}{m_{\Omega}} \sum_{j} \max( \vec{\psi}_{j} - (A_{2} \vec{w} )_{j}, 0) = \mu \| (\vec{\psi} - A_{2} \vec{w} )_{+} \|_{1}
\]
where $(\cdot)_{+} = \max(\cdot, 0)$. 

This yields
\[
A_{2} = \frac{|\Omega|}{m_{\Omega}} \begin{bmatrix}
\Phi^{\top}(\vec{x}^{(\Omega)}_{1}, C) \\ 
\vdots \\
\Phi^{\top}(\vec{x}^{(\Omega)}_{m_{\Omega}}, C)
\end{bmatrix} \in \mathbb{R}^{m_{\Omega} \times N},\quad 
\vec{\psi} = \frac{|\Omega|}{m_{\Omega}}
\begin{bmatrix}
\psi(\vec{x}^{(\Omega)}_{1}) \\
\vdots \\
\psi(\vec{x}^{(\Omega)}_{m_{\Omega}})  
\end{bmatrix} \in \mathbb{R}^{m_{\Omega}}. 
\]
Therefore, the optimization problem is
\[
\min \limits_{\vec{w} \in \mathbb{R}^{N}} \Big\{\frac{1}{2} \vec{w}^{\top} A_{1} \vec{w} + \frac{\beta}{2} \|A_{3} \vec{w} - \vec{g} \|_{2}^{2} + \mu \| (\vec{\psi} - A_{2} \vec{w})_{+} \|_{1}       \Big\}. 
\]
We introduce the slack variable $\vec{v} = \vec{\psi} - A_{2} \vec{w}$, and consider the augmented Lagrangian functional with a scaling penalty parameter $\rho > 0$

\[
L_{\rho}(\vec{w}, \vec{v}, \vec{z}) = \frac{1}{2} \vec{w}^{\top} A_{1} \vec{w} + \frac{\beta}{2} \|A_{3} \vec{w} - \vec{g} \|_{2}^{2} + \mu \| (\vec{v})_{+} \|_{1} + \frac{\rho}{2} \| \vec{v} - \vec{\psi} + A_{2} \vec{w} + \vec{z} \|_{2}^{2}  - \frac{\rho}{2} \| \vec{z} \|_{2}^{2}. 
\]

In the $\vec{w}$-minimization problem, we have 

\[
\vec{w}^{(i+1)} = \min \limits_{\vec{w}} \frac{1}{2} \vec{w}^{\top} A_{1} \vec{w} + \frac{\beta}{2} \|A_{3} \vec{w} - \vec{g} \|_{2}^{2} + \frac{\rho}{2} \| \vec{v}^{(i)} - \vec{\psi} + A_{2} \vec{w} + \vec{z}^{(i)} \|_{2}^{2} .  
\]

This is a convex quadratic with closed-form solution: 

\[
(A_{1} + \beta A_{3}^{\top} A_{3} + \rho A_{2}^{\top} A_{2} ) \vec{w}^{(i+1)} = \beta A_{3}^{\top} \vec{g} + \rho A_{2}^{\top} (\vec{\psi} - \vec{v}^{(i)} - \vec{z}^{(i)} ). 
\]

In the $\vec{v}$-minimization problem, we let
\[
\vec{t}^{(i)} \coloneq \vec{\psi} - A_{2} \vec{w}^{(i+1)} - \vec{z}^{(i)}. 
\]
Then
\[
\vec{v}^{(i+1)} = \min \limits_{\vec{v}} \mu \| (\vec{v})_{+} \|_{1} + \frac{\rho}{2} \| \vec{v} - \vec{t}^{(i)} \|_{2}^{2}. 
\]
This decouples over coordinates. For each component $j$: 
\[
v_{j}^{(i+1)} = \begin{cases}
    t_{j}^{(i)}, & t_{j}^{(i)} \leq 0, \\
    \max(t_{j}^{(i)} - \frac{\mu}{\rho}, 0), & t_{j}^{(i)} > 0, 
\end{cases} = \min \left\{t_{j}^{(i)}, \, \left( t_{j}^{(i)} - \frac{\mu}{\rho}\right)_{+} \right\}. 
\]
We use a proximal operator for the positive-part $\ell^1$-norm which shrinks positive entries and leaves negative entries untouched. 
Finally, the dual variable is updated as follows: 
\[
\vec{z}^{(i+1)} = \vec{z}^{(i)} + (\vec{v}^{(i+1)} - \vec{\psi} + A_{2} \vec{w}^{(i+1)} ).
\]
\pagebreak

\section{Experiments on Parameter Settings}\label{appendix:sect_experiments_parameter_settings}

This appendix collects supplementary numerical experiments used to examine the sensitivity of the proposed method with respect to several key hyperparameters. In particular, we study the influence of the boundary penalty parameter \(\beta\), the truncation threshold \(\tau\) used in the stabilized linear algebra procedure, and the parameter \(c\) appearing in the scaling of the shape parameter. These experiments are intended to support the parameter choices adopted in the main text and to provide additional empirical evidence on the robustness and stability of the method.

More specifically, Appendix~\ref{appendix:experiment_boundary_penalty_parameter} investigates the role of the boundary penalty parameter \(\beta\) in balancing accurate boundary enforcement and overall accuracy across different boundary value problems. Appendix~\ref{appendix:experiment_optimal_value_c} provides numerical justification for the choice of the constant \(c=c(T,\tau)\), which governs the scaling of the shape parameter and plays an important role in the overall accuracy of the approximation.

Finally, Appendix~\ref{appendix:experiment_truncation_threshold} studies the effect of the truncation threshold \(\tau\) on stability and approximation accuracy when truncated singular value decomposition is employed. 

Unless otherwise stated, all other numerical settings remain the same as those used in the corresponding experiments in the main text so that the effect of each parameter can be isolated clearly.

\subsection{Experiments on boundary penalty parameter \(\beta\)}
\label{appendix:experiment_boundary_penalty_parameter}

We examine the effect of the boundary penalty parameter \(\beta\) on the accuracy of the proposed variational Gaussian RBF method. Since \(\beta\) controls the strength of boundary-condition enforcement, its value must balance boundary fidelity and overall accuracy.

\begin{figure}[H]
    \centering
    \includegraphics[width=0.6\linewidth]{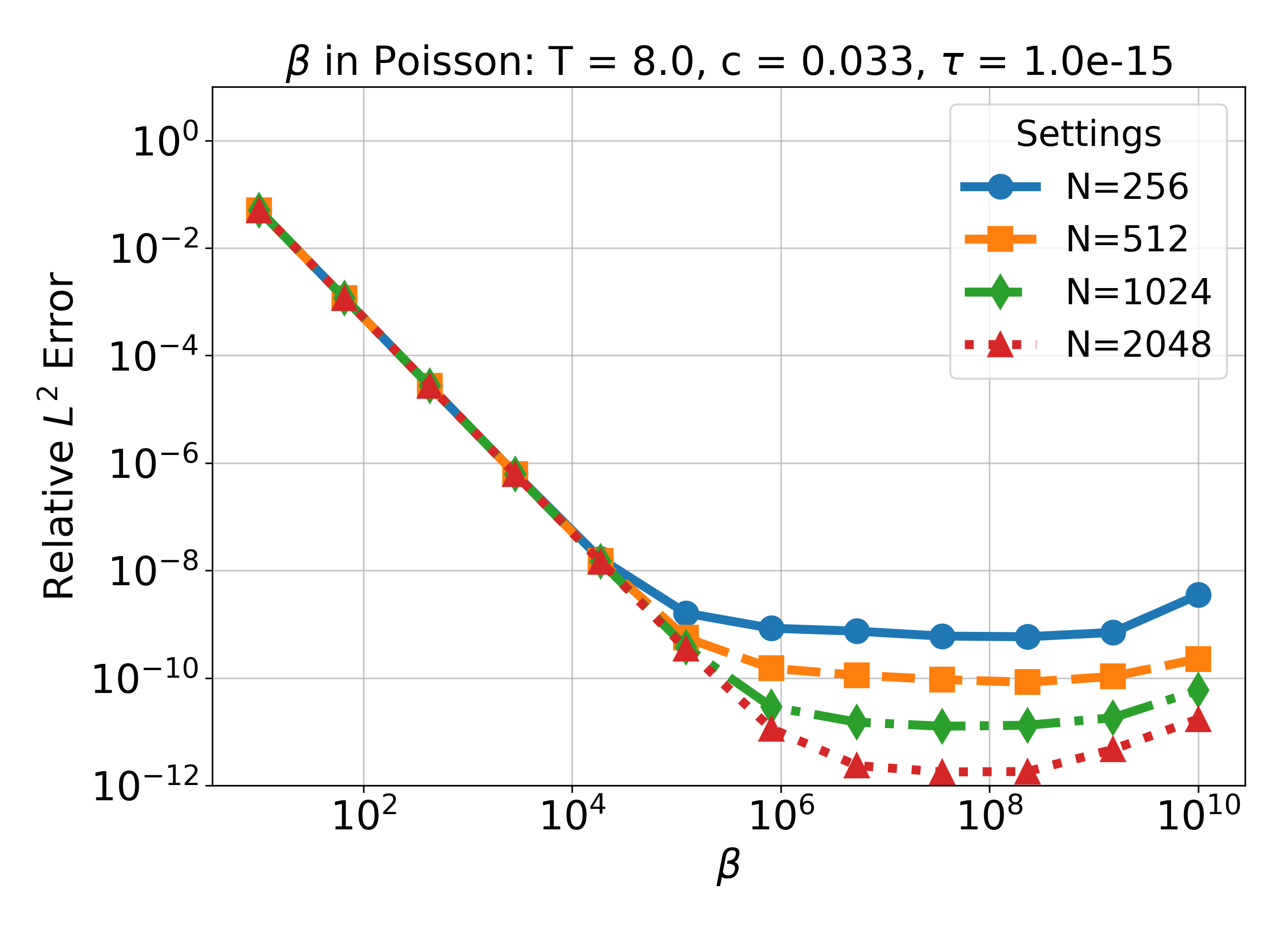}
    \caption{\textbf{Effect of the penalty parameter \(\beta\) for the 1D Poisson problem.} Relative \(\ell^{2}\) error versus \(\beta\) for \(N=256,512,1024,\) and \(2048\), with fixed \(T=8.0\), \(c=c(T, \tau) = 0.033\), and \(\tau=10^{-15}\). The error decreases as \(\beta\) increases and stabilizes for sufficiently large \(\beta\), while larger \(N\) yields higher accuracy.}
    \label{fig:1d_poisson_problem_benchmark_beta_boundary}
\end{figure}

Figure~\ref{fig:1d_poisson_problem_benchmark_beta_boundary} shows the relative \(\ell^2\) error versus \(\beta\) for the 1D Poisson problem with \(N=256,512,1024,\) and \(2048\), using fixed \(T=8.0\), \(c=c(T,\tau)=0.033\), and \(\tau=10^{-15}\). For all values of \(N\), the error decreases rapidly as \(\beta\) increases from small values, indicating that weak penalties lead to insufficient boundary enforcement. Once \(\beta\) is sufficiently large, the error reaches a plateau, showing that further increases in \(\beta\) have little effect. Larger \(N\) consistently yields smaller errors. Overall, the Poisson results indicate that the method is robust with respect to \(\beta\). Therefore, in this experiment, we decide to choose \(\beta = 3 \times 10^{5}\), which is enough for convergence. 

\begin{figure}[H]
       \centering
    \begin{subfigure}[t]{0.5\textwidth}
    \centering\includegraphics[width=1.0\linewidth]{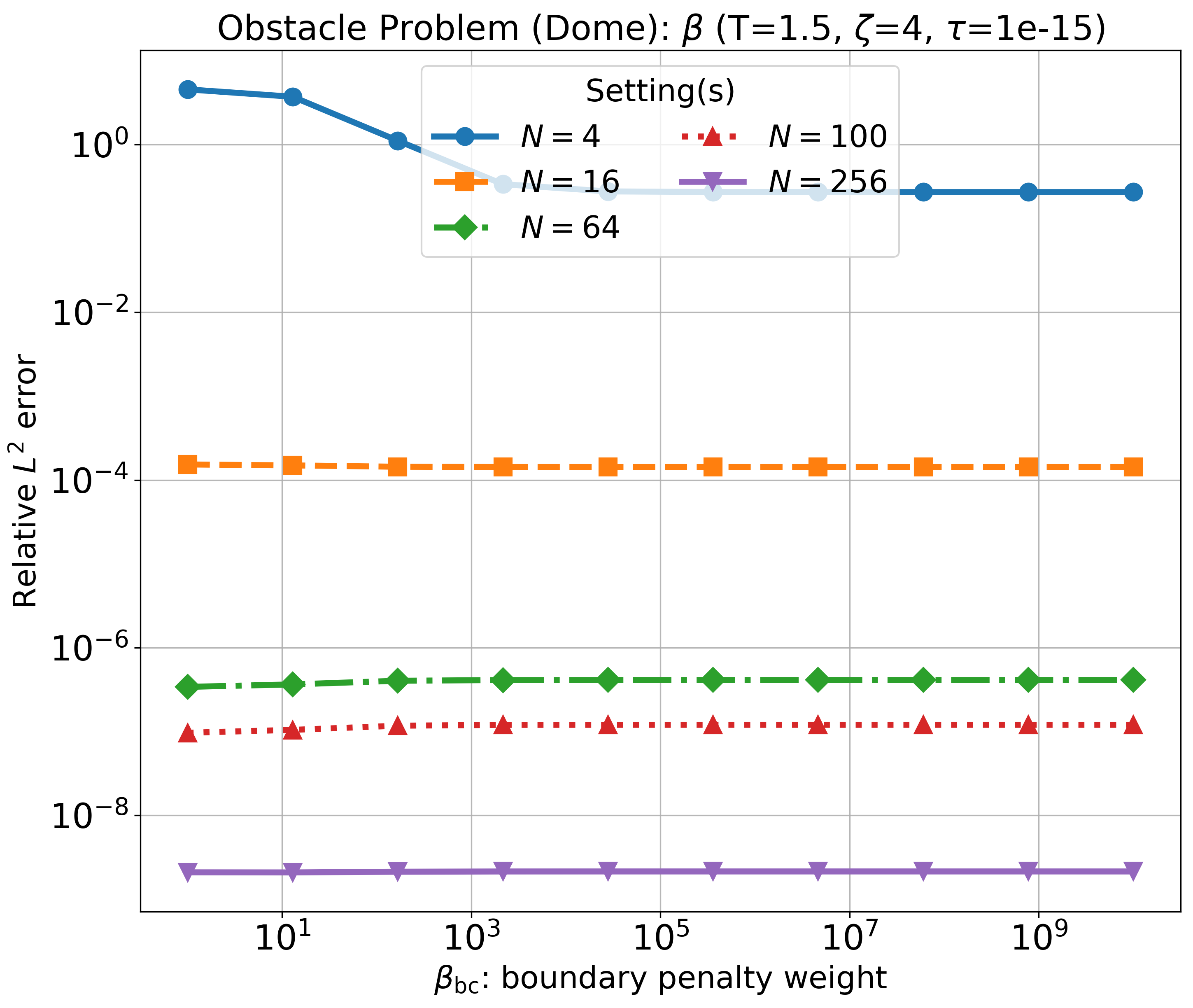}
    \caption{Settings: $T = 1.5, \zeta = 4,\tau = 10^{-15}$}
    \end{subfigure}\hfill
    \begin{subfigure}[t]{0.5\textwidth}
    \centering\includegraphics[width=1.0\linewidth]{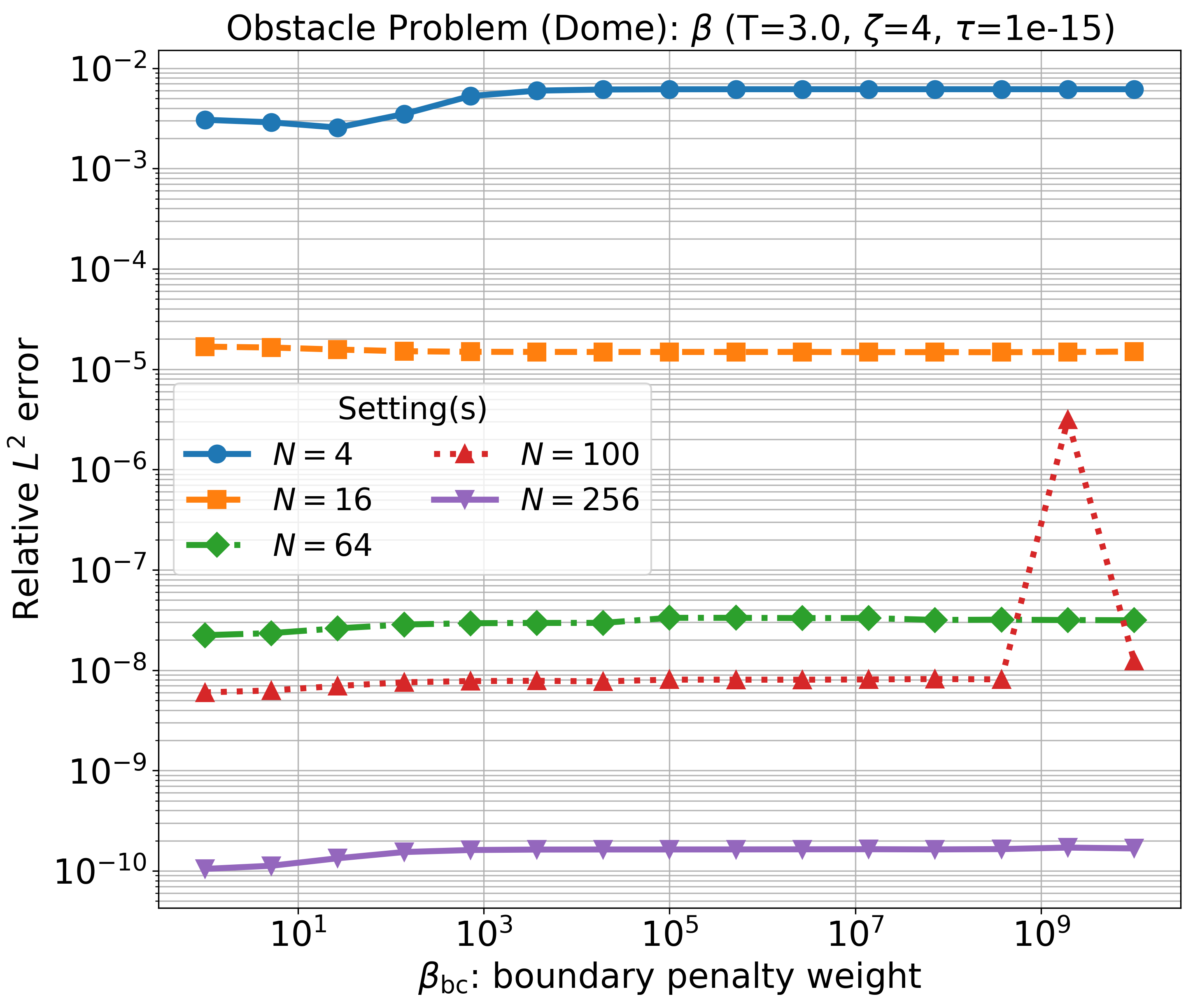}
    \caption{Settings: $T = 3, \zeta = 4, \tau = 10^{-15}$}
  \end{subfigure}
  \caption{Relative $\ell^2$ error vs.\ boundary penalty weight $\beta$ for the dome obstacle. 
For $T=1.5$ (a), errors plateau once $\beta\gtrsim 10^3$ and improve with $N$. 
For $T=3$ (b), overly large $\beta$ degrades accuracy; moderate values of \(\beta\) between \(10^{4}\) and \(10^{6}\) is stable.}

\label{fig:2d_obstacle_problem_dome_benchmark_beta_boundary}
\end{figure}

In Figure~\ref{fig:2d_obstacle_problem_dome_benchmark_beta_boundary} for the dome obstacle problem, \(\beta=10^6\) also falls in the range where the error is already stable for \(T=1.5\), and remains within a reliable regime for \(T=3\). These results support the use of \(\beta=10^6\) as a practical and robust choice in the experiments.

\begin{figure}[p]
    \centering
    \begin{subfigure}[t]{0.5\textwidth}
\centering\includegraphics[width=1\linewidth]{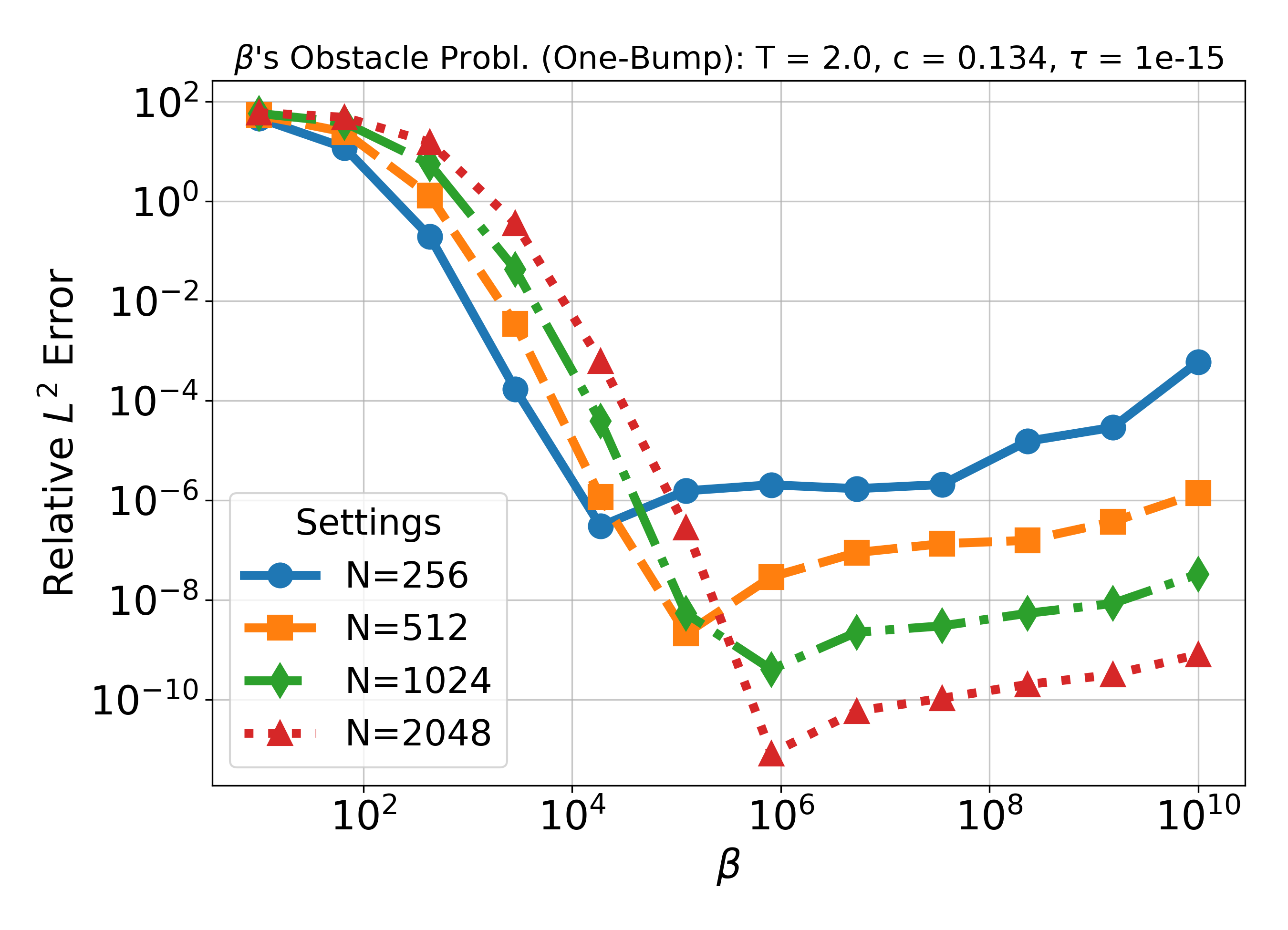}
    \caption{$T = 2.0, \tau = 10^{-15}, c = 0.134$}
    \end{subfigure}\hfill
    \begin{subfigure}[t]{0.5\textwidth}
\centering\includegraphics[width=1\linewidth]{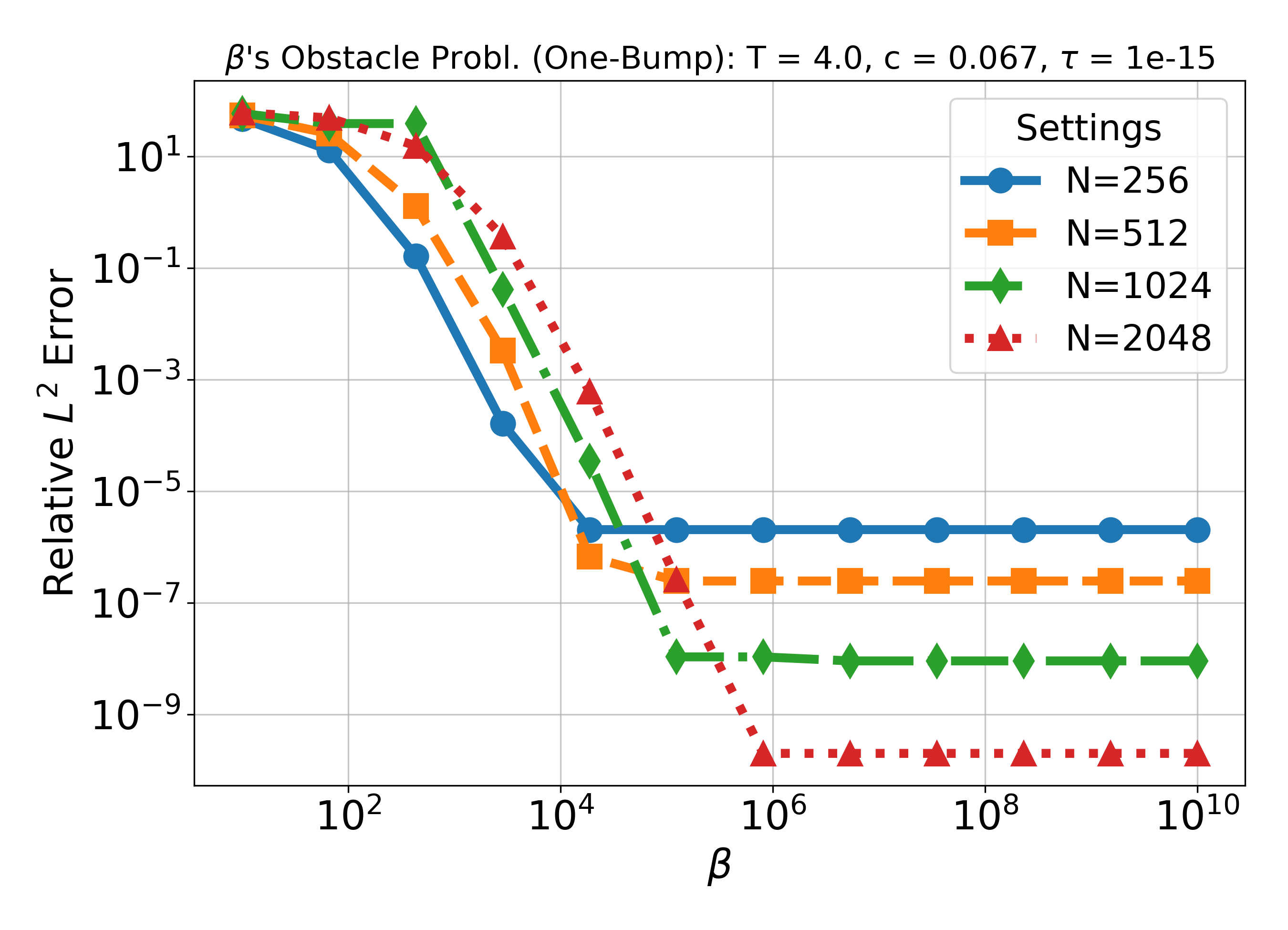}
\caption{$T = 4.0, \tau = 10^{-15}, c = 0.067$}
    \end{subfigure}
\vspace{0.1ex}
    \begin{subfigure}[t]{0.5\textwidth}
    \centering\includegraphics[width=1\linewidth]{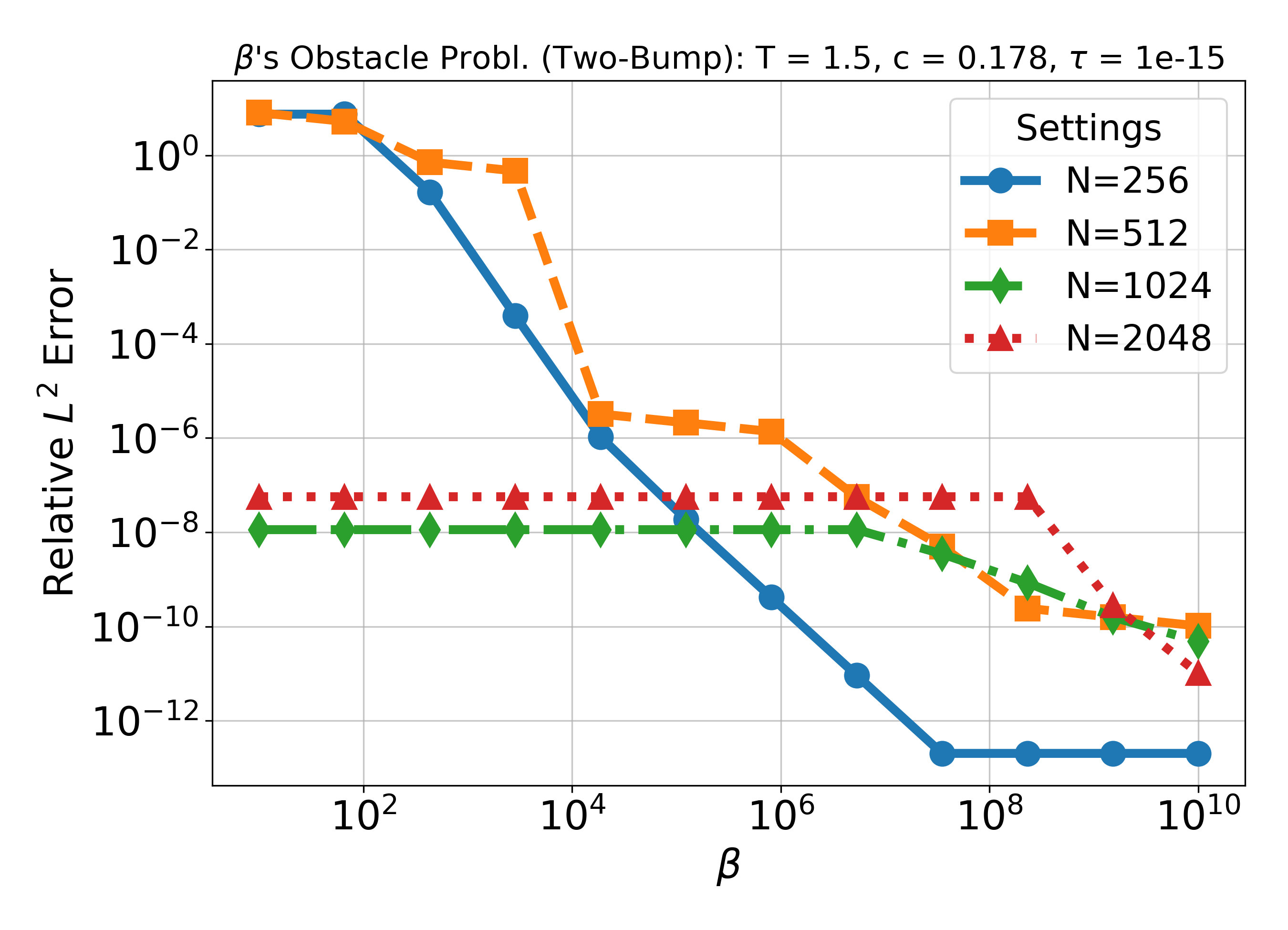}
    \caption{$T = 1.5, \tau = 10^{-15}, c = 0.178$}
    \end{subfigure}\hfill 
    \begin{subfigure}[t]{0.5\textwidth}
    \centering\includegraphics[width=1\linewidth]{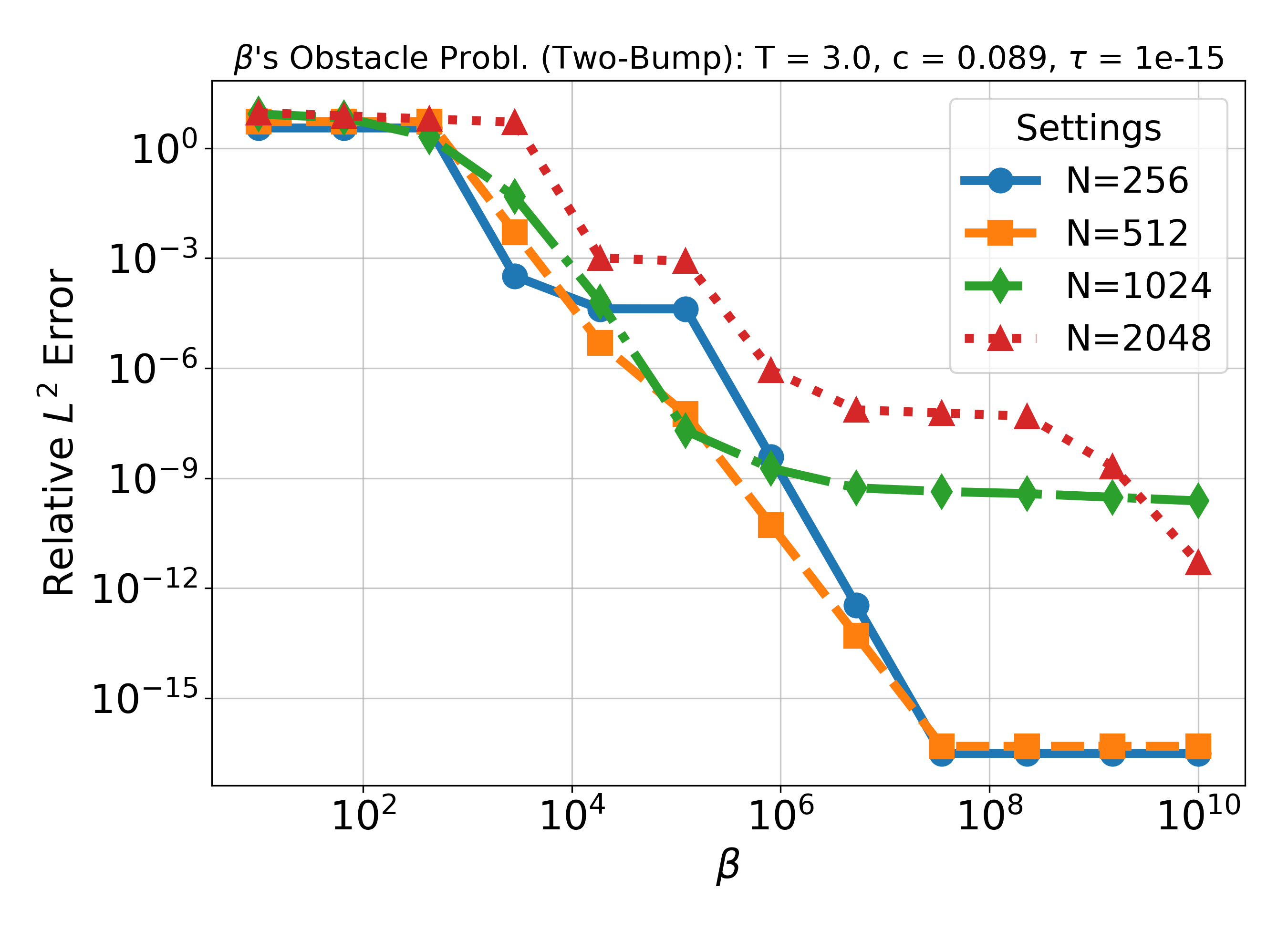}
    \caption{$T = 3.0, \tau = 10^{-15}, c = 0.134$}
    \end{subfigure}

\caption{Relative $\ell^2$ error versus penalty parameter $\beta$ for the obstacle problem using Gaussian RBF with truncated SVD stabilization. Results are shown for varying numbers of centers $N$ and different expansion factors $T$, with fixed tolerance $\tau = 10^{-15}$ and shape parameter $c$ chosen from the formula $c(\tau, T)$. The plots highlight the effect of the penalty parameter \(\beta\) for boundary conditions.}
\label{fig2_1d_obstacle_problem_benchmark_beta_parameter}
\end{figure}

Figure~\ref{fig2_1d_obstacle_problem_benchmark_beta_parameter} reports the corresponding results for the obstacle problem, including both one-bump and two-bump cases, for several values of the expansion factor \(T\), with \(\tau=10^{-15}\) fixed and \(c = c(T,\tau)\). As in the Poisson case, small \(\beta\) leads to large errors, while increasing \(\beta\) significantly improves the approximation. However, the obstacle problem exhibits greater sensitivity to \(\beta\), particularly for smaller \(T\) or smaller \(N\). In several cases, the error decreases sharply over an intermediate range of \(\beta\), then either stabilizes or slightly deteriorates when \(\beta\) becomes excessively large. This behavior suggests that overly strong boundary penalization may imbalance the variational problem in more challenging obstacle configurations. Hence, for the one-bump obstacle problem we let \(\beta =  10^{6}\) while \(\beta =  10^{8}\) in the two-bump problem.

These experiments show that \(\beta\) should be chosen large enough to enforce the boundary condition effectively, but not so large that it degrades accuracy. Across the examples considered here, the most effective choices lie in an intermediate-to-large range of \(\beta\), with the optimal regime depending on the problem type, the expansion factor \(T\), and the number of centers \(N\).

\subsection{Experiments on optimal value \(c = c(T, \tau)\)}\label{appendix:experiment_optimal_value_c}

\begin{figure}[p]
  \centering
  \begin{subfigure}[t]{0.5\textwidth}
    \centering\includegraphics[width=1\linewidth]{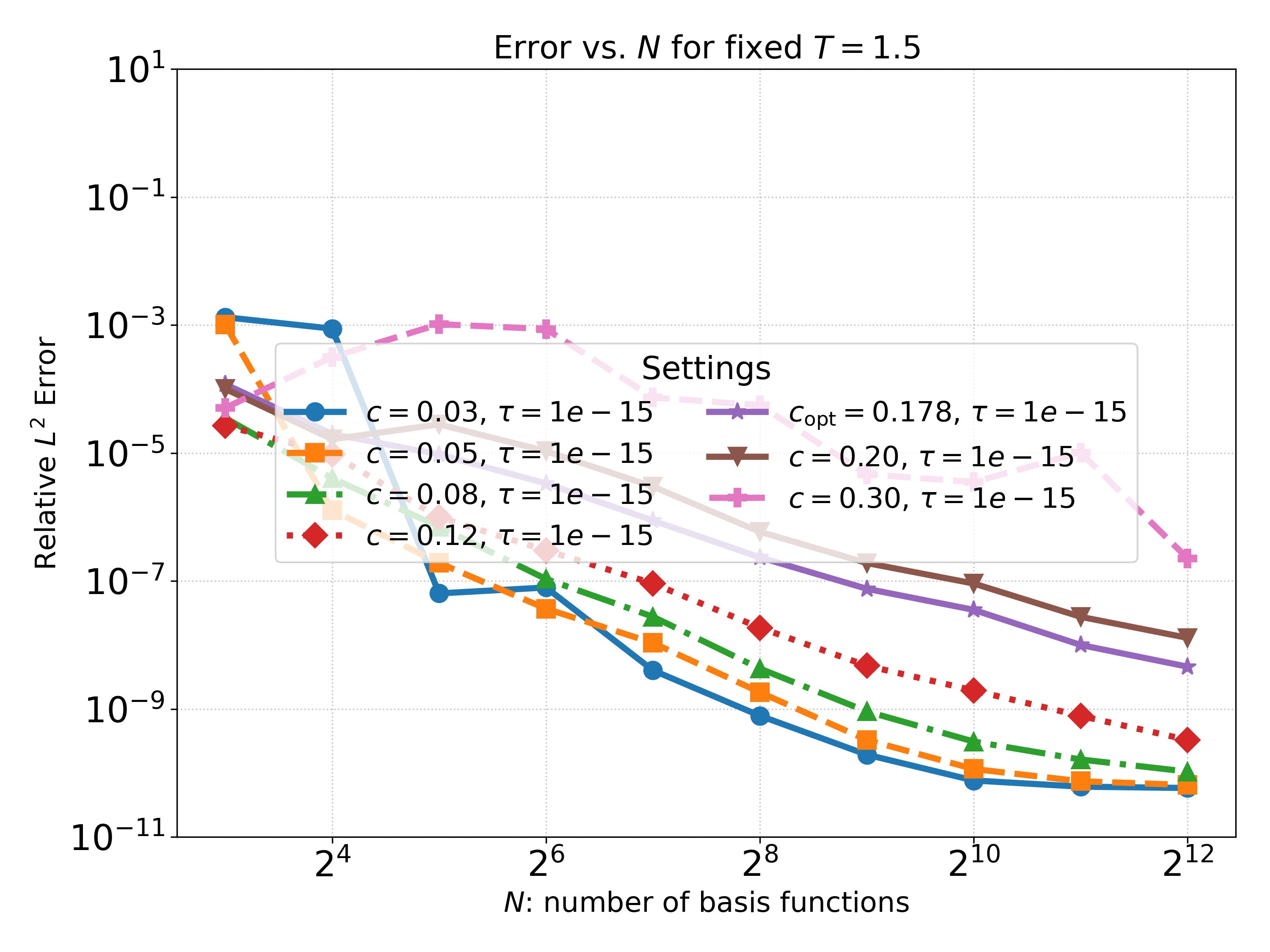}
  \end{subfigure}\hfill
  \begin{subfigure}[t]{0.5\textwidth}
    \centering\includegraphics[width=1\linewidth]{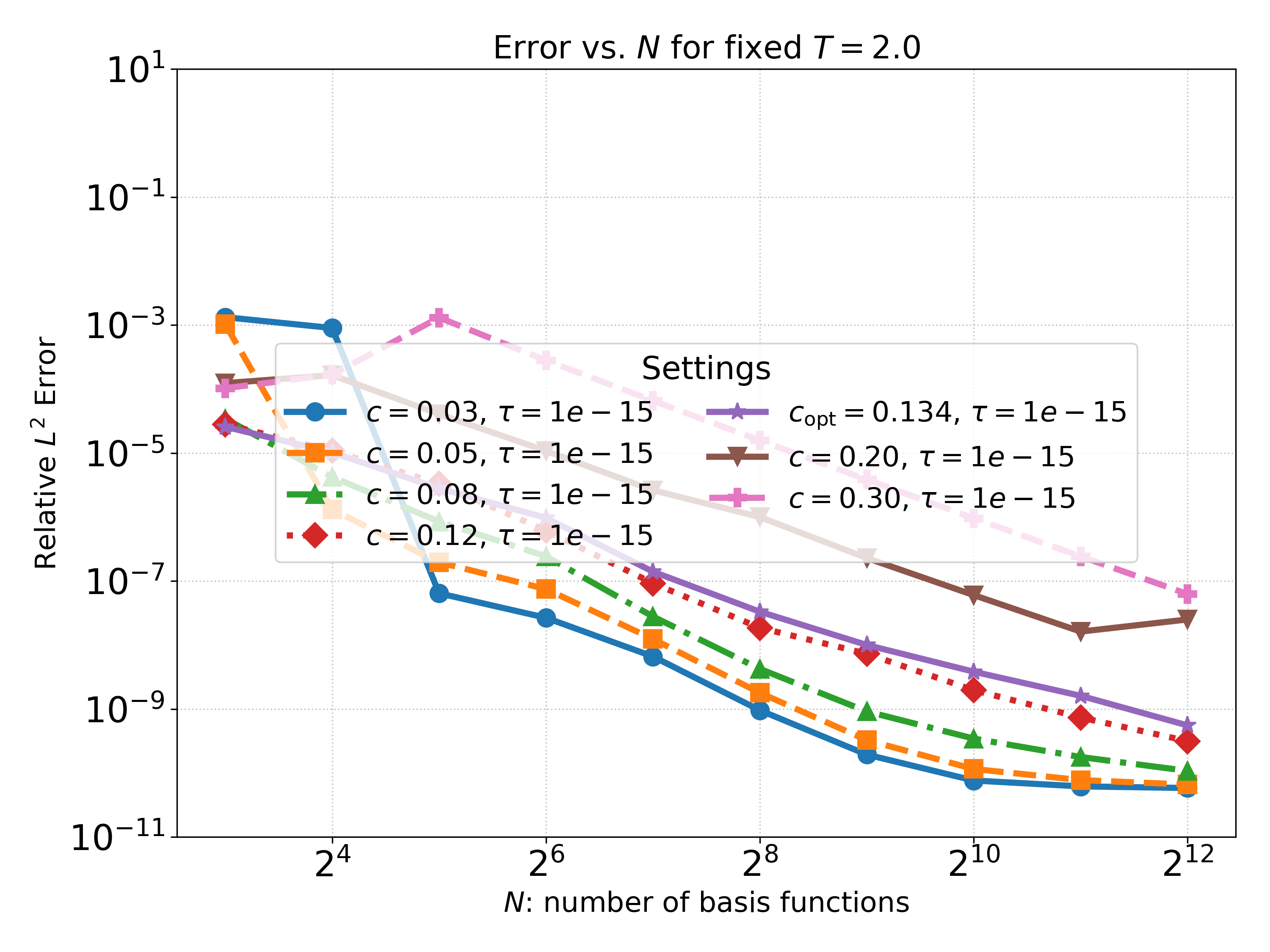}
  \end{subfigure}

  \vspace{0.1ex}

  \begin{subfigure}[t]{0.5\textwidth}
    \centering\includegraphics[width=1\linewidth]{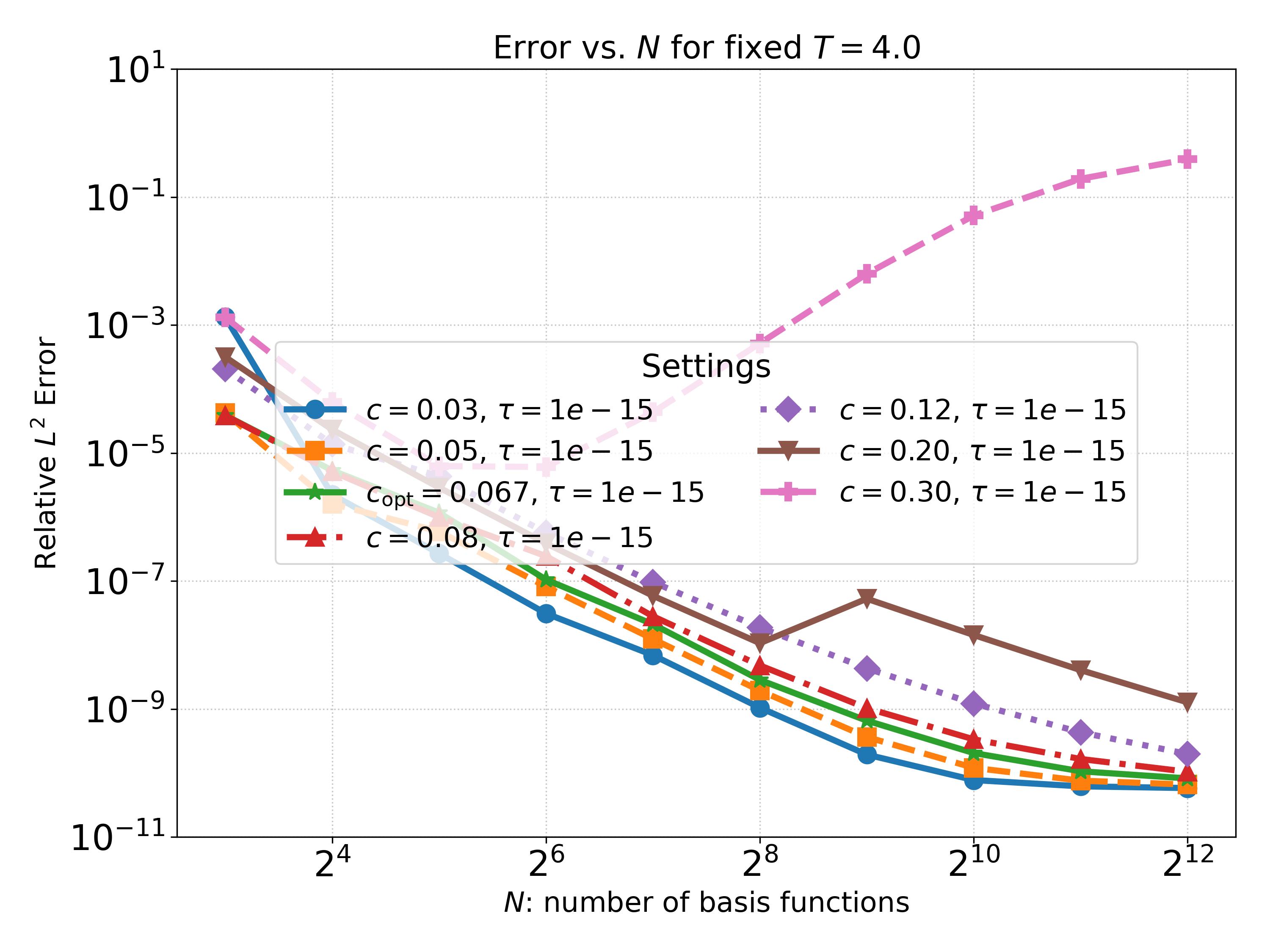}
  \end{subfigure}\hfill
  \begin{subfigure}[t]{0.5\textwidth}
    \centering\includegraphics[width=1\linewidth]{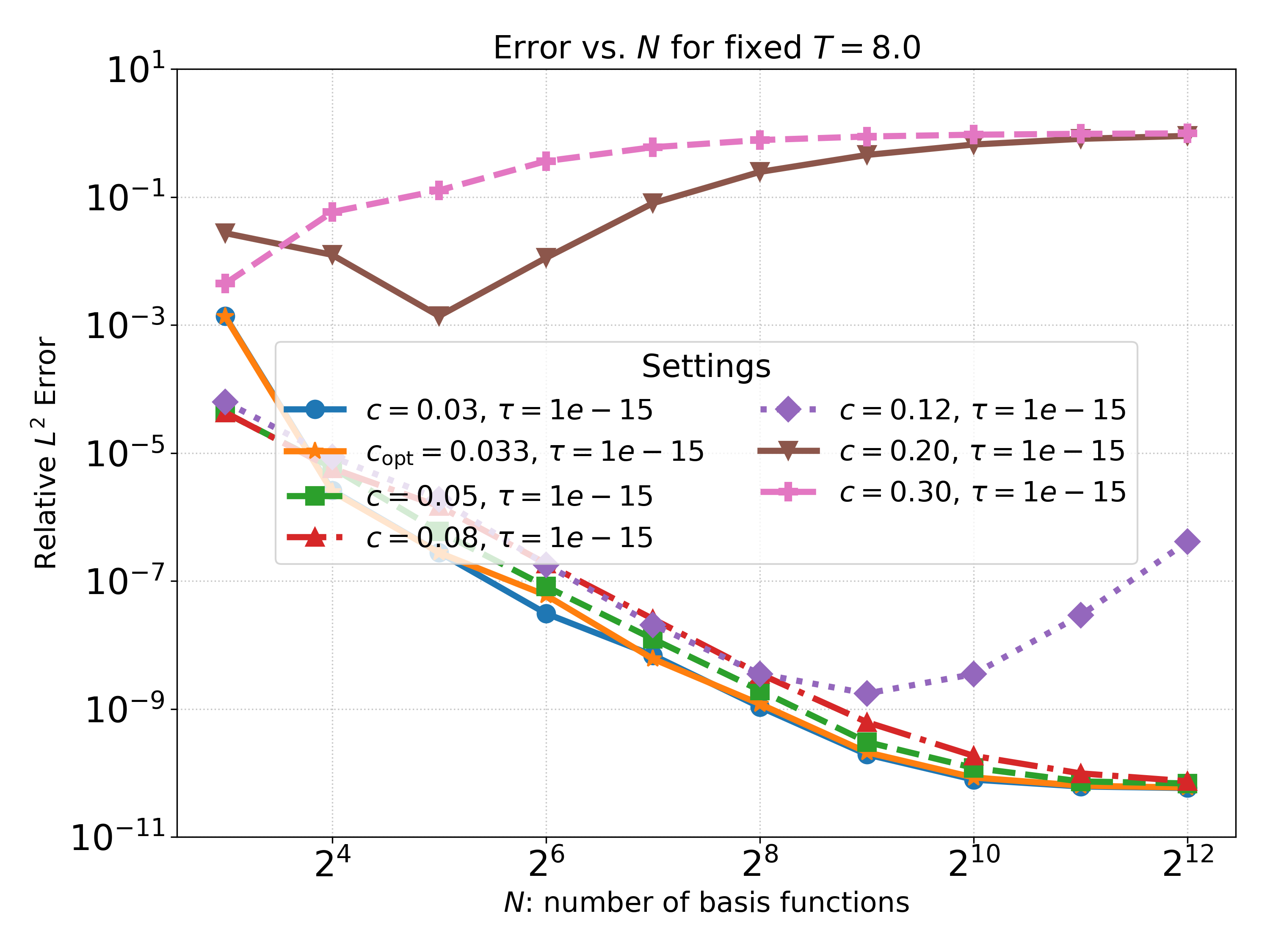}
  \end{subfigure}

\caption{\textbf{1D Poisson problem:} Relative \(\ell^2\) error versus the number of basis functions \(N\) for fixed expansion factor \(T\), with oversampling ratio \(\zeta=2\) and truncation threshold \(\tau=10^{-15}\). Each panel compares several values of the proportionality constant \(c\) in the linear shape-parameter scaling, including the selected value \(c_{\mathrm{opt}}(T,\tau)\). For moderate values of \(T\), the error decreases steadily with \(N\) for a range of suitable \(c\)-values, while overly large \(c\) leads to a loss of monotonicity in convergence.}
  \label{fig5_poisson_pde}
\end{figure}

\begin{figure}[p]
  \centering
  \begin{subfigure}[t]{0.5\textwidth}
    \centering\includegraphics[width=1\linewidth]{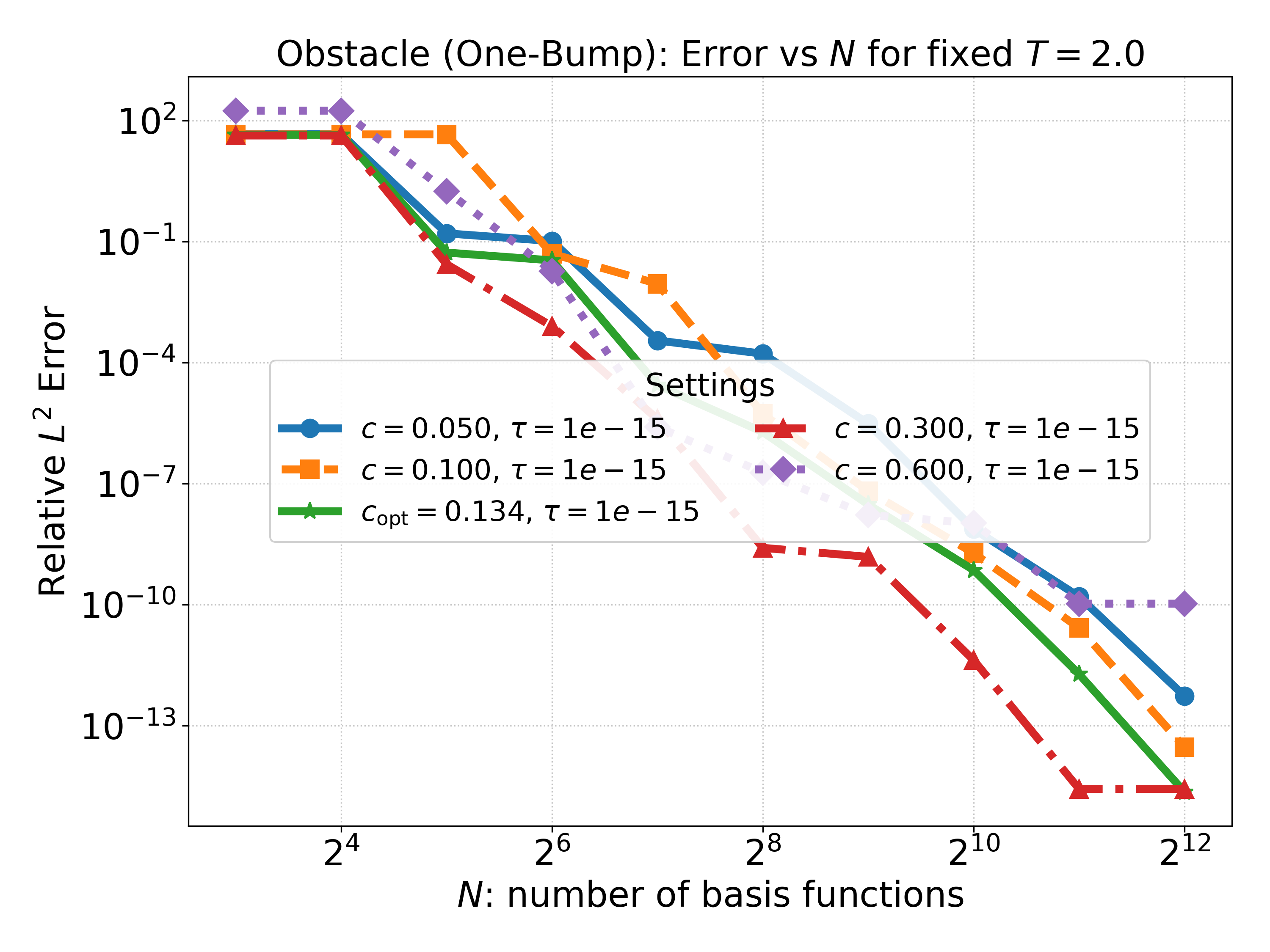}
  \end{subfigure}\hfill
  \begin{subfigure}[t]{0.5\textwidth}
    \centering\includegraphics[width=1\linewidth]{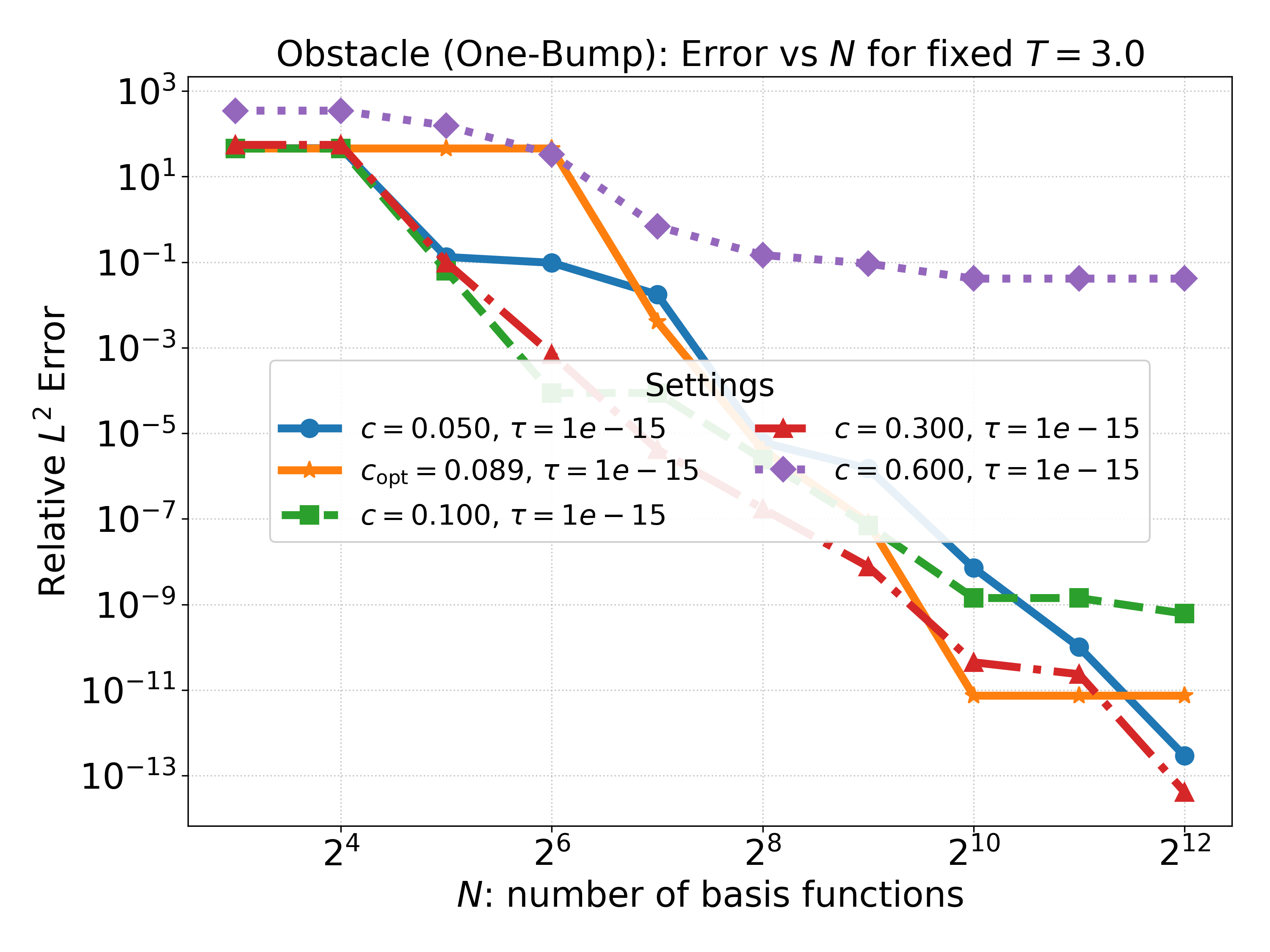}
  \end{subfigure}

  \vspace{0.1ex}

  \begin{subfigure}[t]{0.5\textwidth}
    \centering\includegraphics[width=1\linewidth]{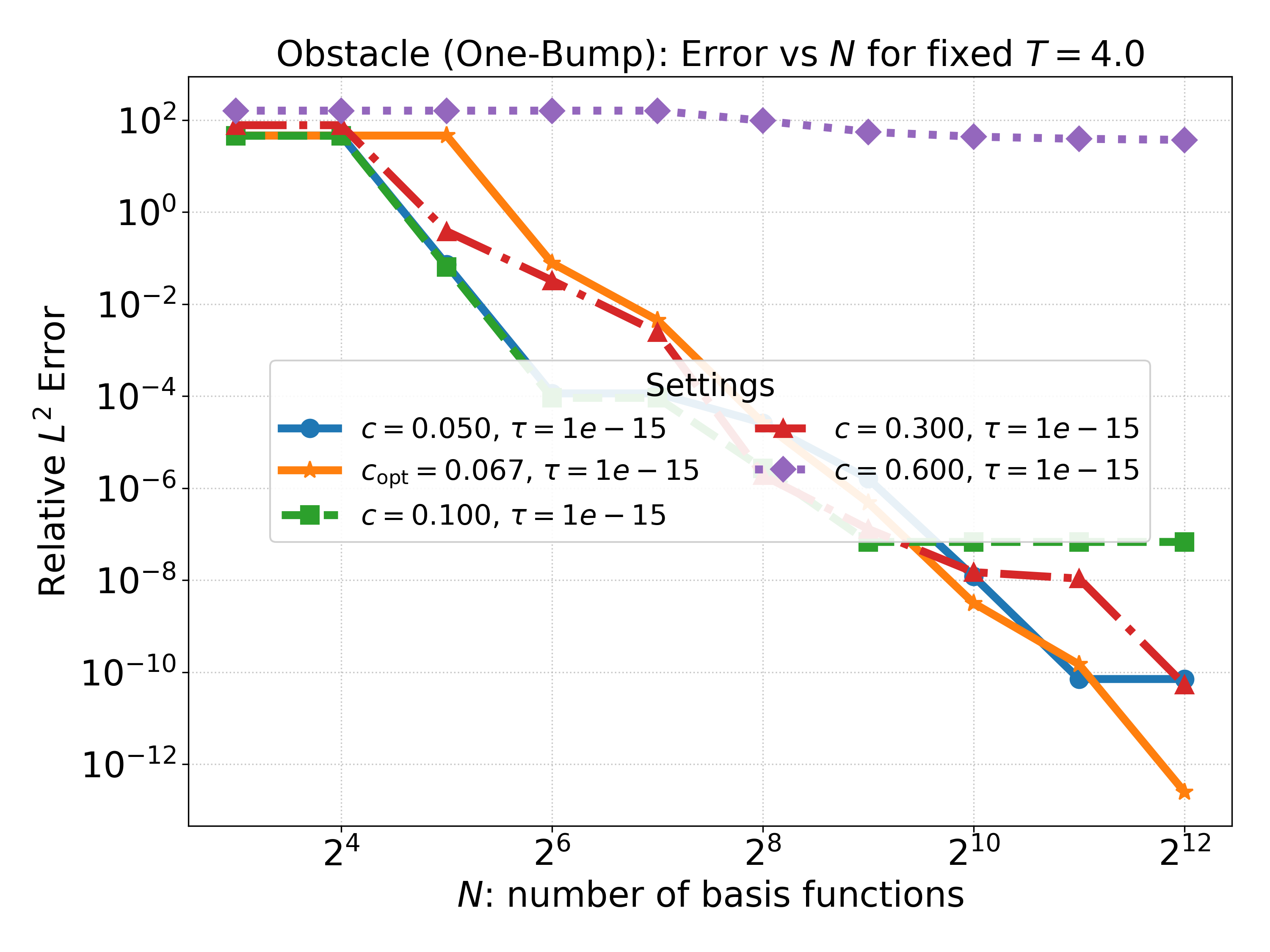}
  \end{subfigure}\hfill
  \begin{subfigure}[t]{0.5\textwidth}
    \centering\includegraphics[width=1\linewidth]{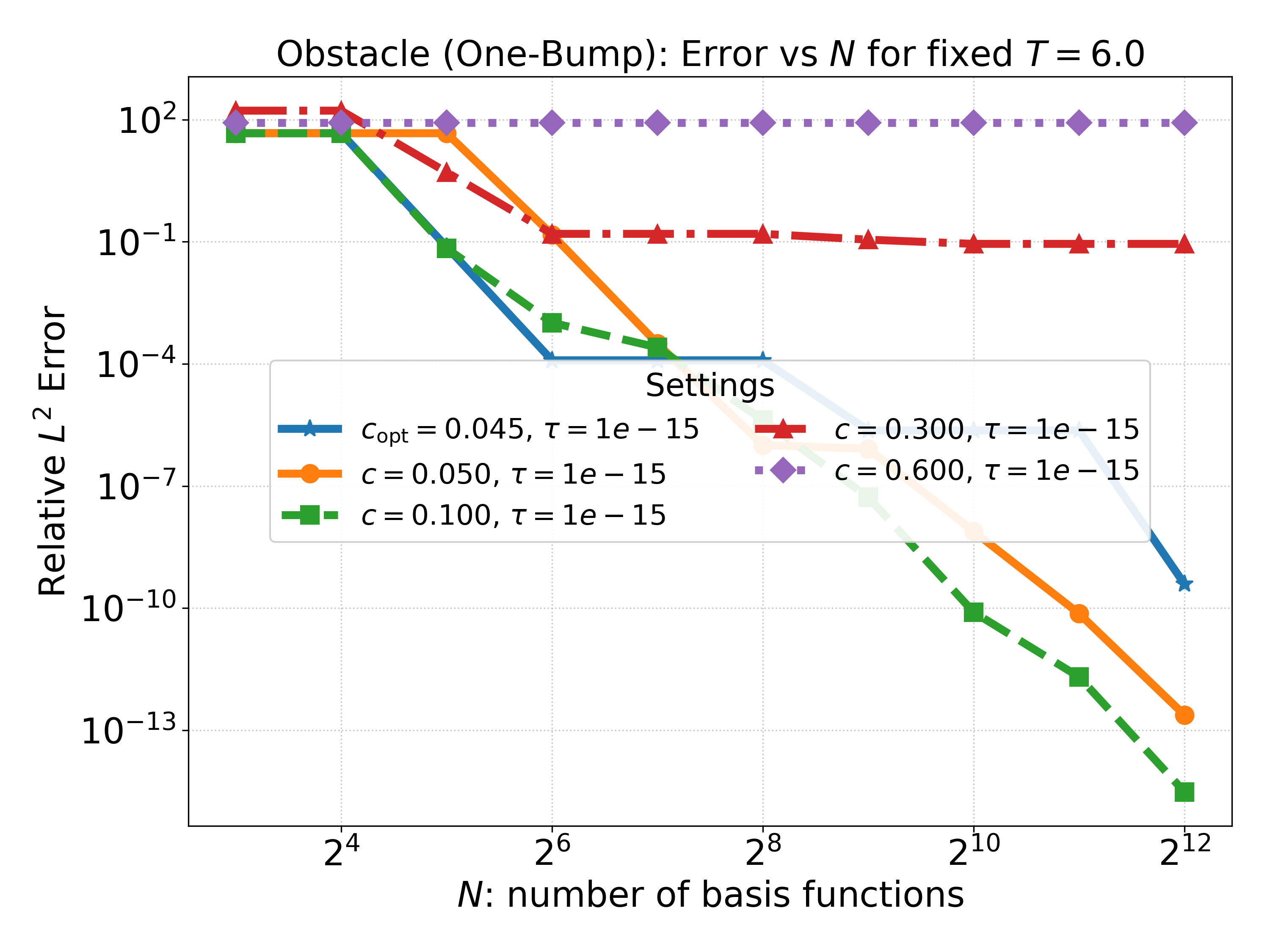}
  \end{subfigure}
  \caption{Relative \(\ell^2\) error versus the number of basis functions \(N\) for the one-bump obstacle problem, shown for fixed expansion factors \(T\in\{2,3,4,6\}\) in each panel. The oversampling ratio is fixed at \(\zeta=2\), the TSVD truncation threshold is \(\tau=10^{-15}\), and each curve corresponds to a different value of the proportionality constant \(c\) in the linear shape-parameter scaling, including the selected value \(c_{\mathrm{opt}}(T,\tau)\). The results show that the choice of \(c\) strongly affects monotonicity of convergence and asymptotic accuracy: moderate values of \(c\), in particular \(c_{\mathrm{opt}}(T,\tau)\) or nearby choices, yield sustained error decay as \(N\) increases, whereas overly large values may stagnate or become unstable.}

  \label{fig2_benchmark_obstacle1d_one_bump}
\end{figure}

\begin{figure}[p]
  \centering
  \begin{subfigure}[t]{0.5\textwidth}
    \centering\includegraphics[width=1\linewidth]{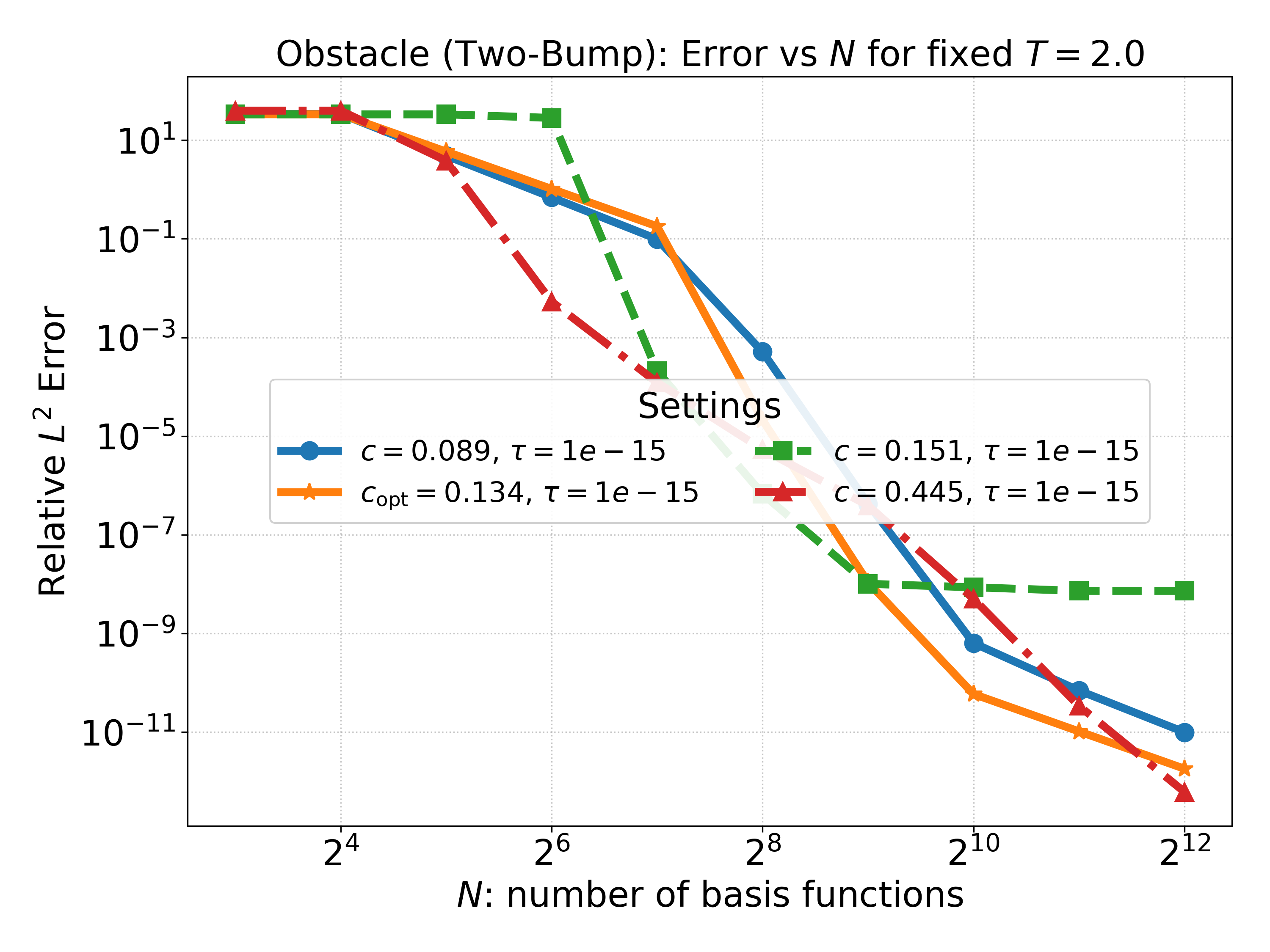}
  \end{subfigure}\hfill
  \begin{subfigure}[t]{0.5\textwidth}
    \centering\includegraphics[width=1\linewidth]{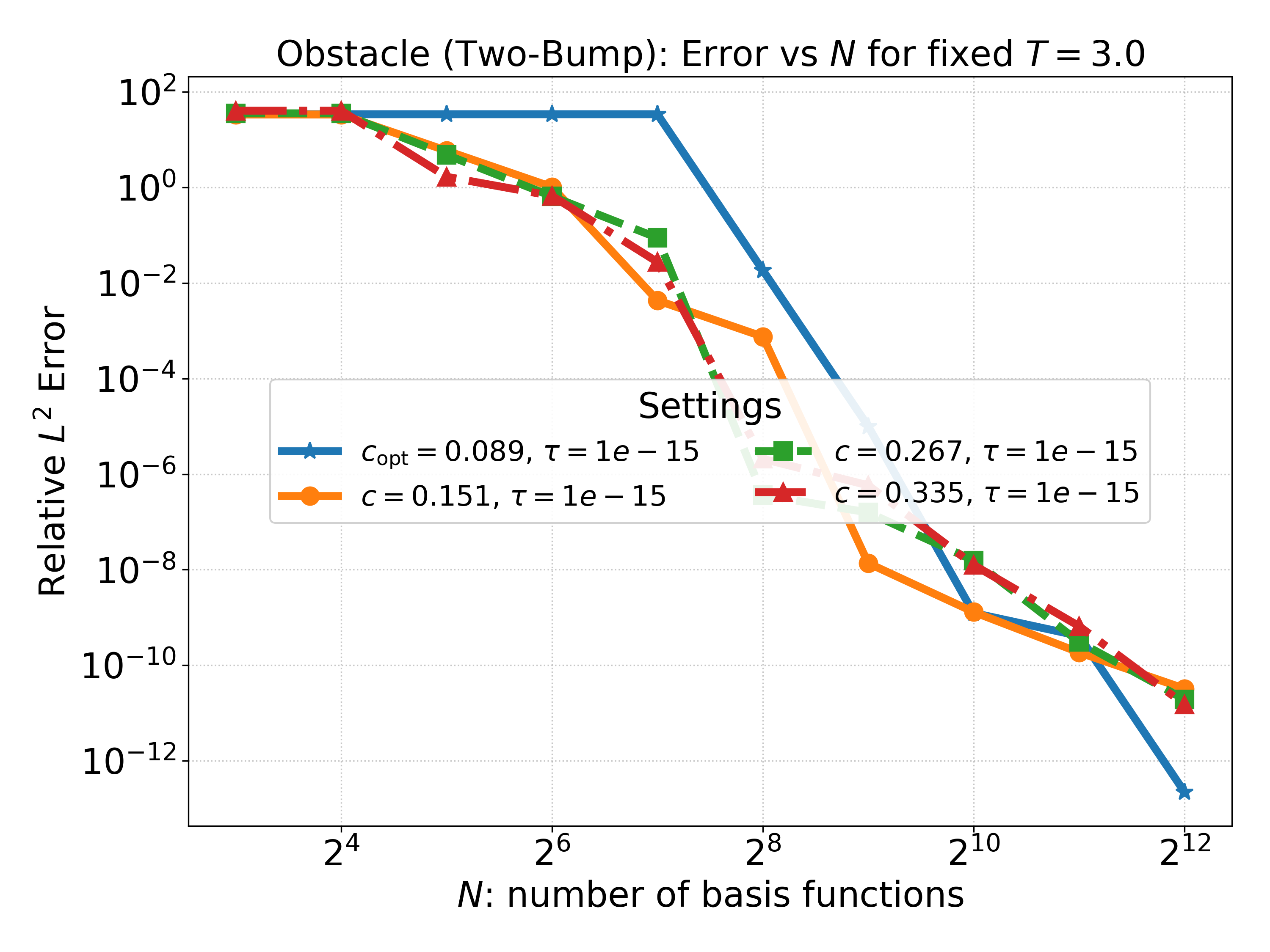}
  \end{subfigure}

  \vspace{0.1ex}

  \begin{subfigure}[t]{0.5\textwidth}
    \centering\includegraphics[width=1\linewidth]{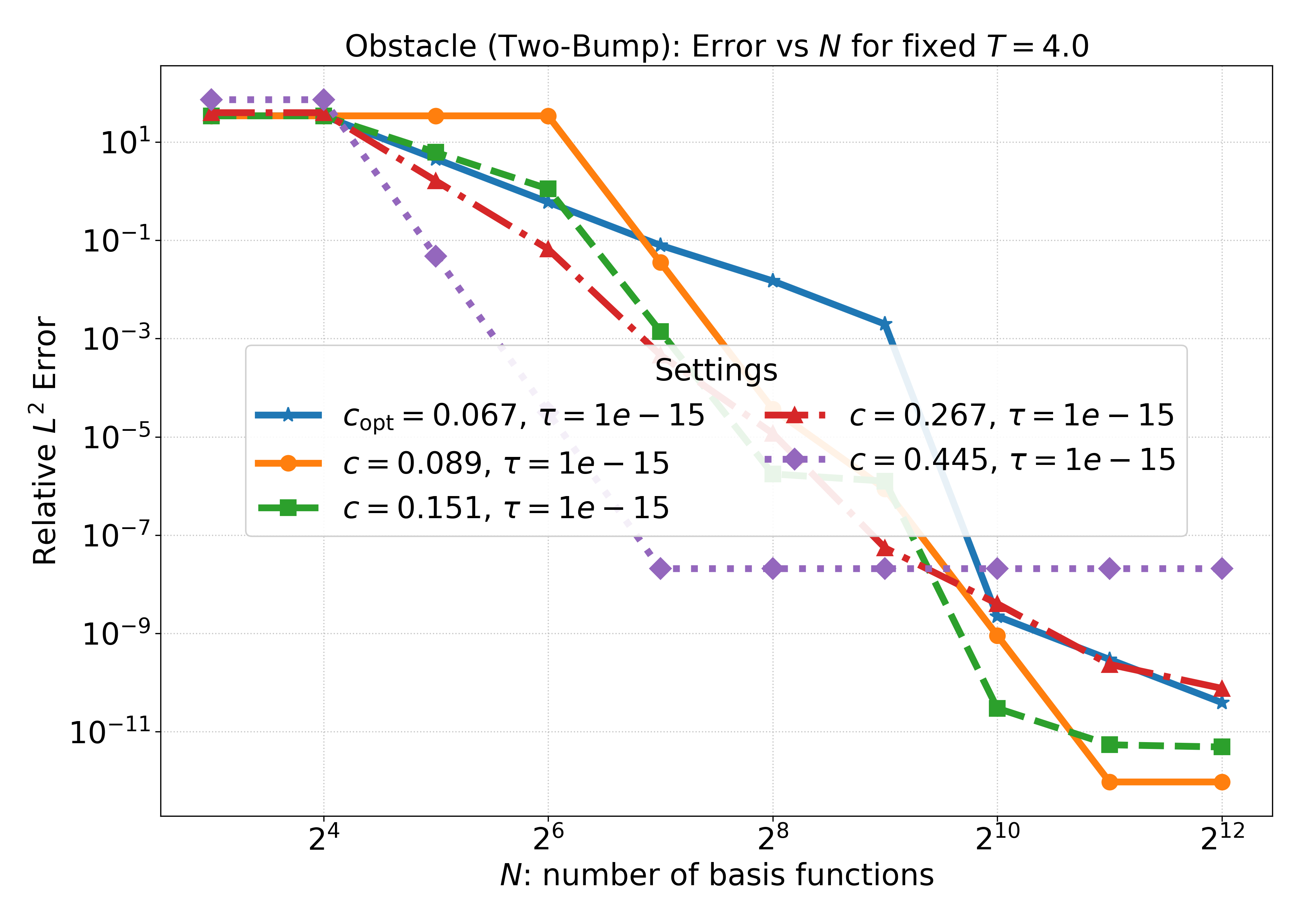}
  \end{subfigure}\hfill
  \begin{subfigure}[t]{0.5\textwidth}
    \centering\includegraphics[width=1\linewidth]{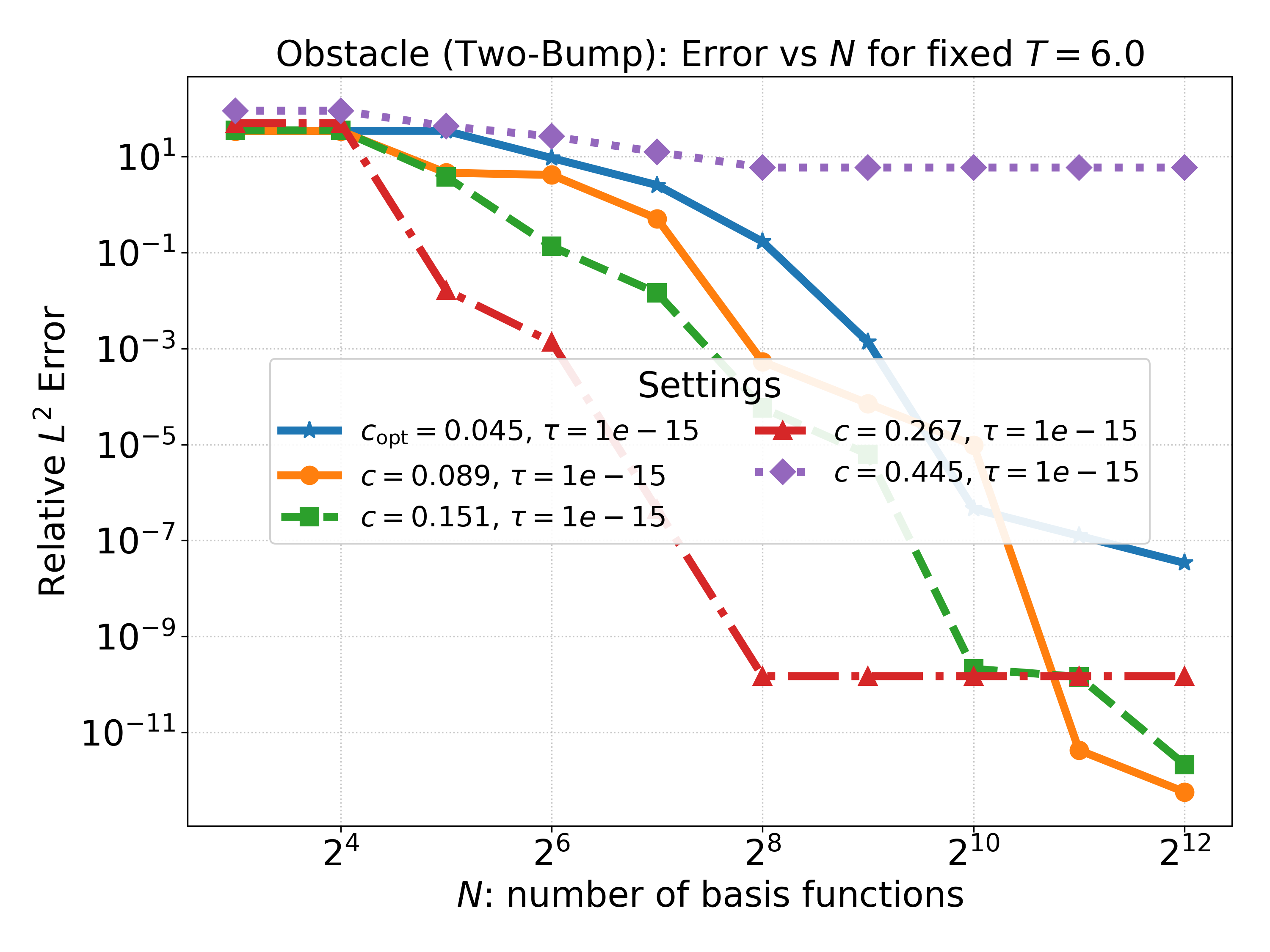}
  \end{subfigure}
  \caption{Relative \(\ell^2\) error versus the number of basis functions \(N\) for the two-bump obstacle problem, with fixed expansion factor \(T\) in each panel. The oversampling ratio is \(\zeta=2\), the TSVD truncation threshold is fixed at \(\tau=10^{-15}\), and the curves compare different values of the proportionality constant \(c\) in the linear shape-parameter scaling, including the selected choice \(c_{\mathrm{opt}}(T,\tau)\). As in the one-bump case, the approximation quality is sensitive to \(c\): suitable moderate values produce monotonous convergence and low final errors, while excessively large \(c\) can lead to stagnation or marked loss of accuracy, especially for larger \(T\).}
\label{fig2_benchmark_obstacle1d_two_bump}
\end{figure}

\begin{figure}[p]
  \centering
  \begin{subfigure}[t]{0.5\textwidth}
    \centering\includegraphics[width=1\linewidth]{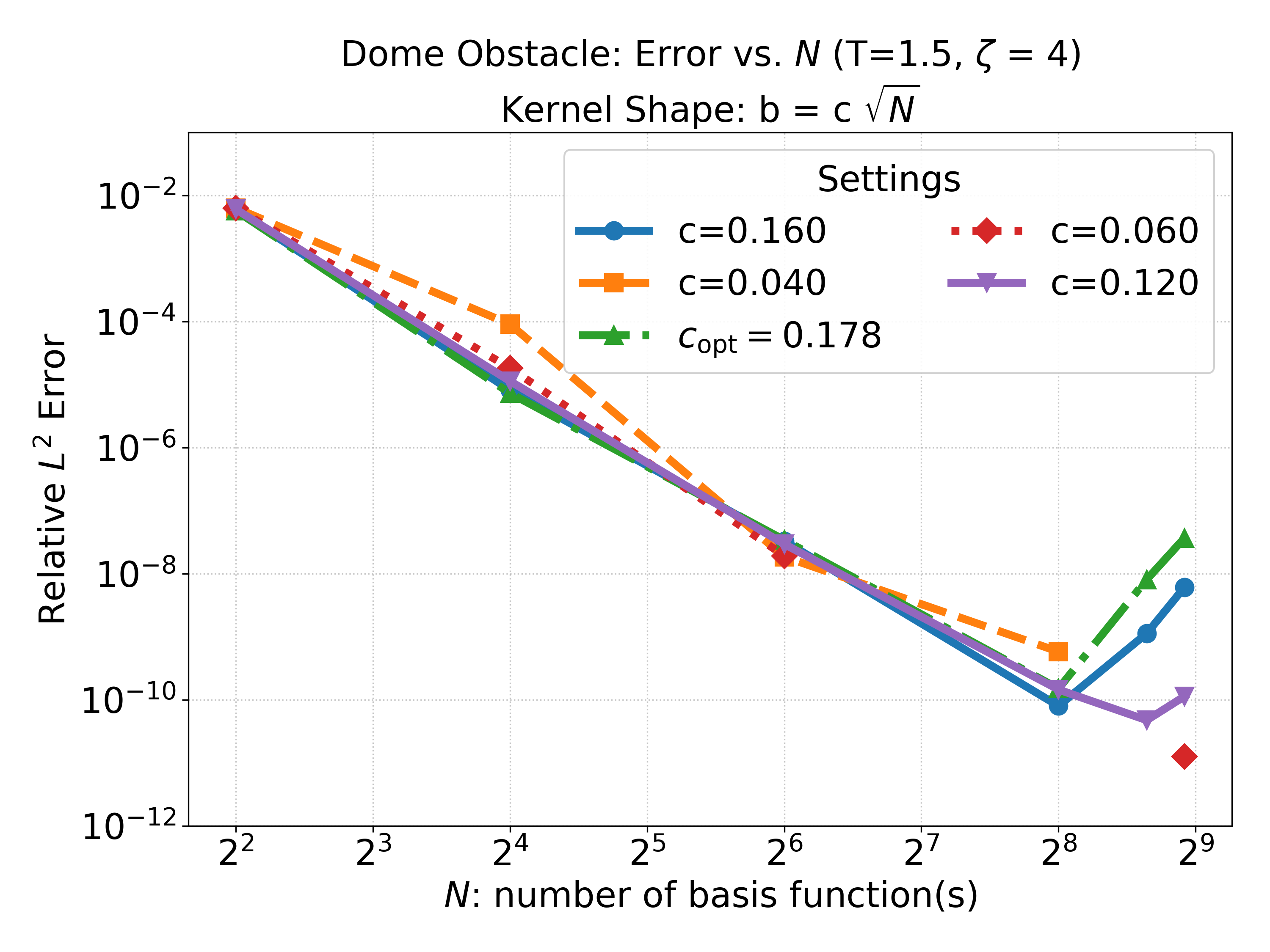}
  \end{subfigure}\hfill
  \begin{subfigure}[t]{0.5\textwidth}
    \centering\includegraphics[width=1\linewidth]{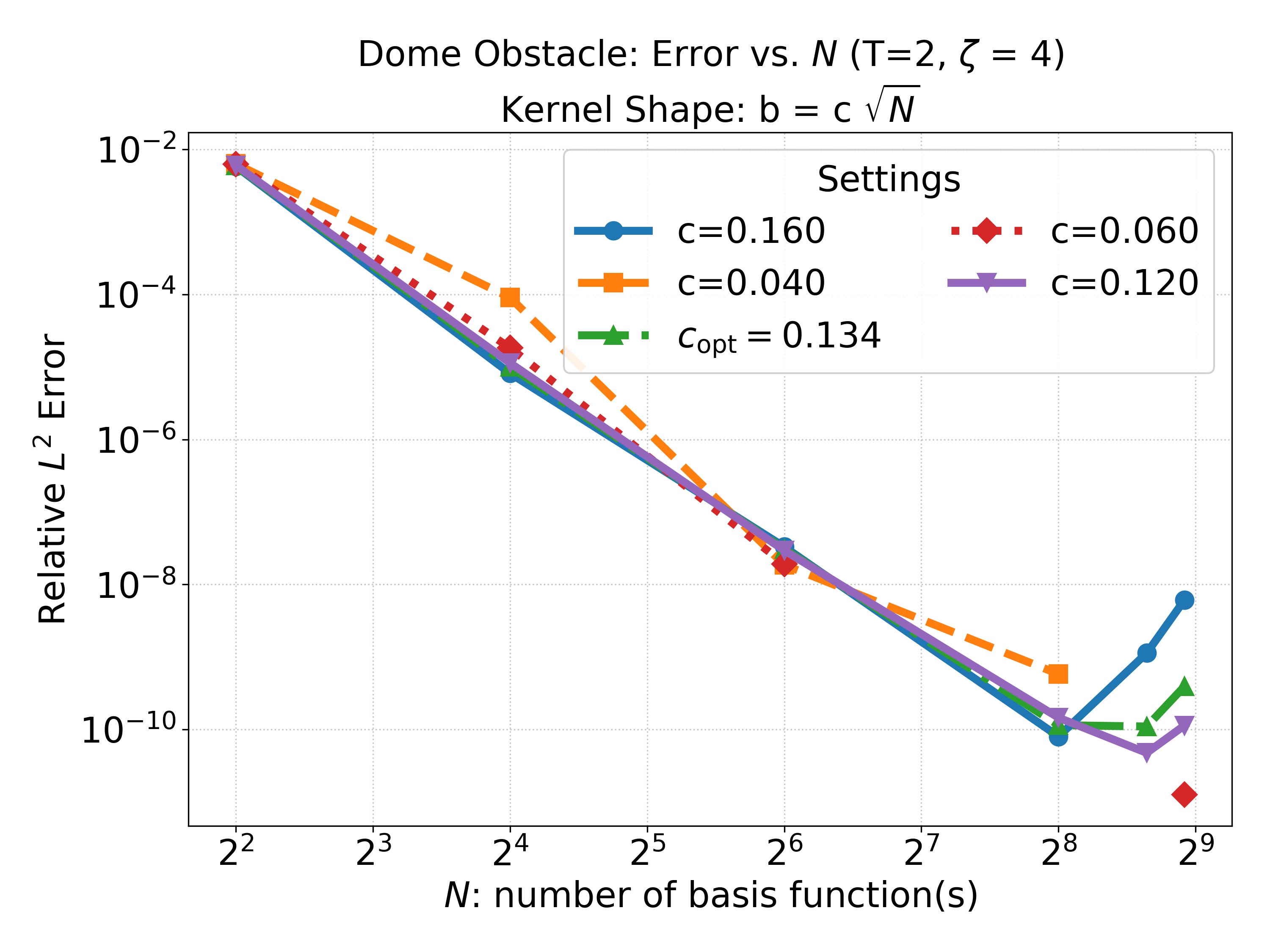}
  \end{subfigure}

  \vspace{0.1ex}

  \begin{subfigure}[t]{0.5\textwidth}
    \centering\includegraphics[width=1\linewidth]{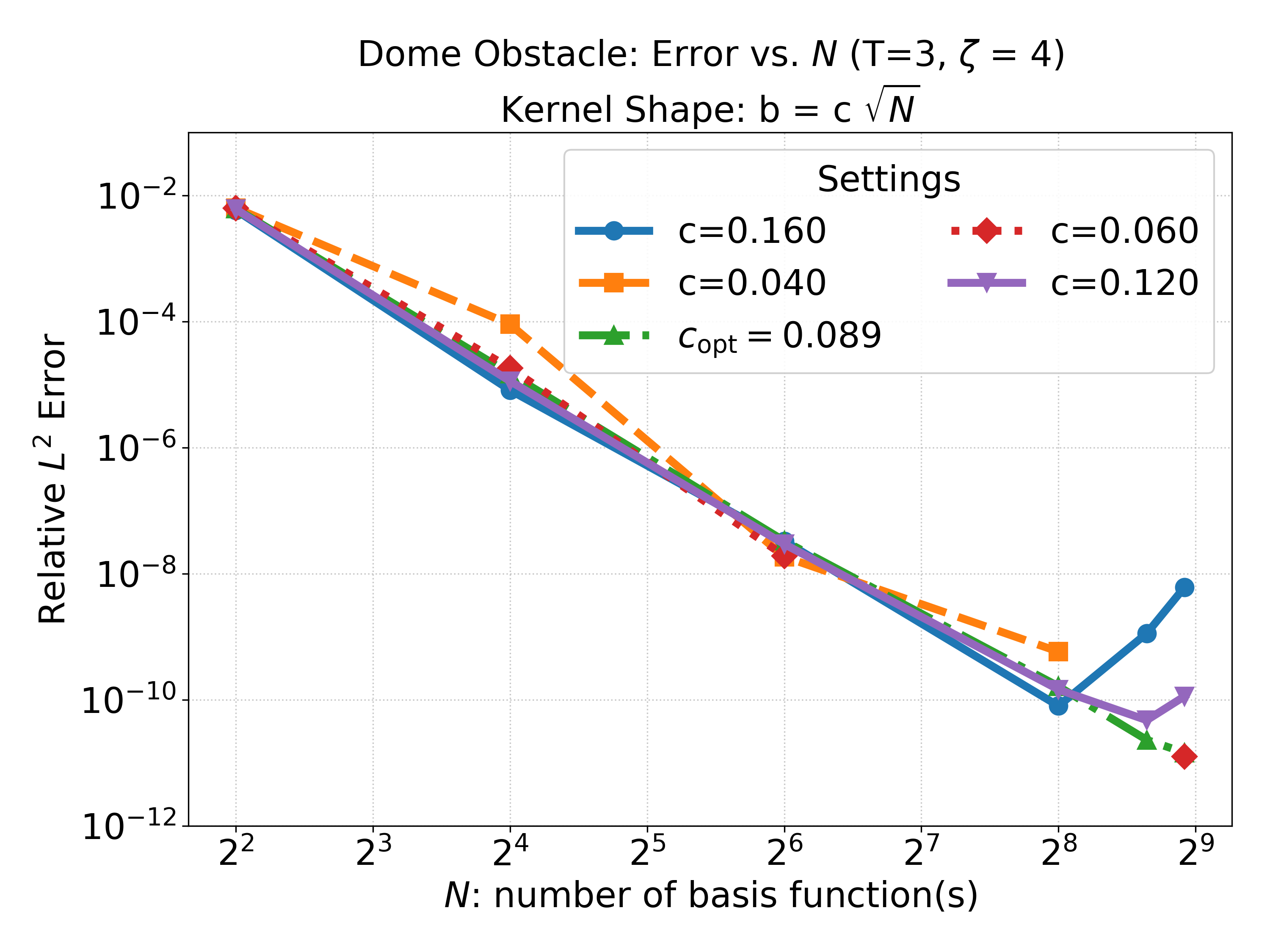}
  \end{subfigure}\hfill
  \begin{subfigure}[t]{0.5\textwidth}
    \centering\includegraphics[width=1\linewidth]{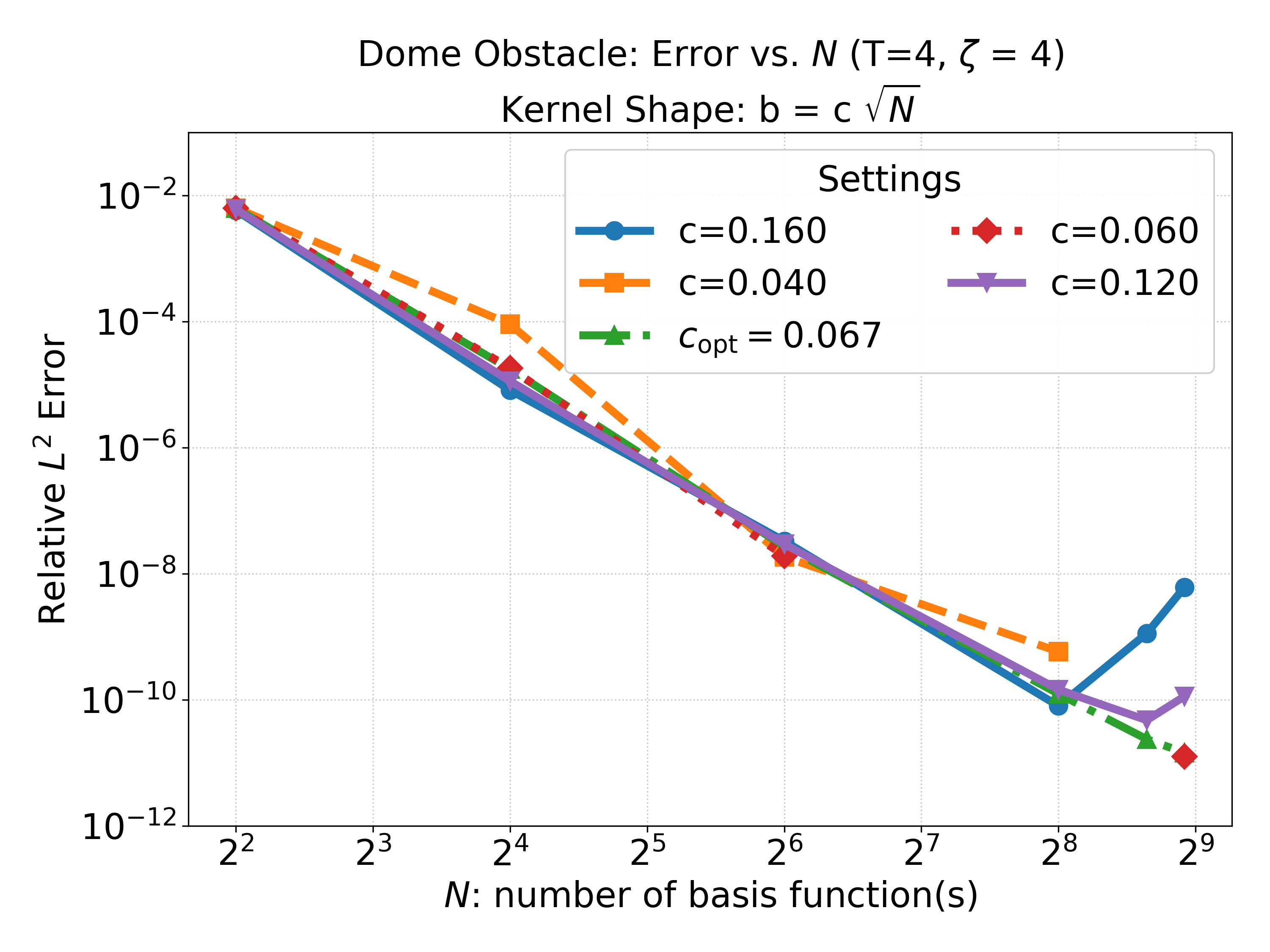}
  \end{subfigure}
\caption{Relative $\ell^2$ error vs.\ number of basis functions $N$ for the Dome obstacle with domain expansion $T\in\{1.5, 2, 3, 4\}$. Each panel compares shape parameters $c$, including $c_{\mathrm{opt}}(T, \tau)$; accordingly,  $c_{\mathrm{opt}}$ yields monotonous convergence for relative error to around $10^{-10}$, whereas overly large $c$ (e.g., $0.3$) becomes unstable for large $N$.}
\label{fig:2d_obstacle_problem_dome_benchmark_T_c_formula}
\end{figure}

We next examine the effect of the proportionality constant \(c\) in the shape-parameter scaling, with \(\tau=10^{-15}\) fixed and \(c\) varying across a range of values for each expansion factor \(T\). The results in Figures~\ref{fig5_poisson_pde}--\ref{fig:2d_obstacle_problem_dome_benchmark_T_c_formula} show that the accuracy is sensitive to the choice of \(c\), especially for larger \(N\). In general, values of \(c\) that are too large lead to loss of monotonicity in convergence as \(N\) increases, whereas values that are small provide monotonous convergence, but may not always give the lowest error.

For the 1D Poisson problem in Figure~16, the error decreases steadily with \(N\) for a range of moderate values of \(c\), while overly large values such as \(c=0.2\) or \(c=0.3\) lose monotonicity, particularly for larger expansion factors \(T\). The value \(c_{\mathrm{opt}}(T,\tau)\) consistently lies near the transition between monotonous decay and instability. A similar trend is observed for the one-bump and two-bump obstacle problems in Figures~\ref{fig2_benchmark_obstacle1d_one_bump} and~\ref{fig2_benchmark_obstacle1d_two_bump}: the proposed choice \(c_{\mathrm{opt}}(T,\tau)\) either gives the best performance or remains very close to the best-performing choice among all tested values, while overly large \(c\) frequently leads to stagnation or sharp deterioration.

The same behavior is also seen in the two-dimensional dome obstacle problem in Figure~\ref{fig:2d_obstacle_problem_dome_benchmark_T_c_formula}. Across all tested values of \(T\), the choice \(c_{\mathrm{opt}}(T,\tau)\) yields  monotonous decay and avoids the instability observed for larger values of \(c\). These results support the use of \(c=c(T,\tau)\) as a practical parameter-selection rule: \(c\) is sufficiently small to preserve monotonicity, yet large enough to maintain high approximation accuracy across different problems, dimensions, and expansion factors.

\subsection{Experiments on truncation threshold value \(\tau\)}
\label{appendix:experiment_truncation_threshold}

\begin{figure}[p]
  \centering
  \begin{subfigure}[t]{0.5\textwidth}
    \centering\includegraphics[width=1\linewidth]{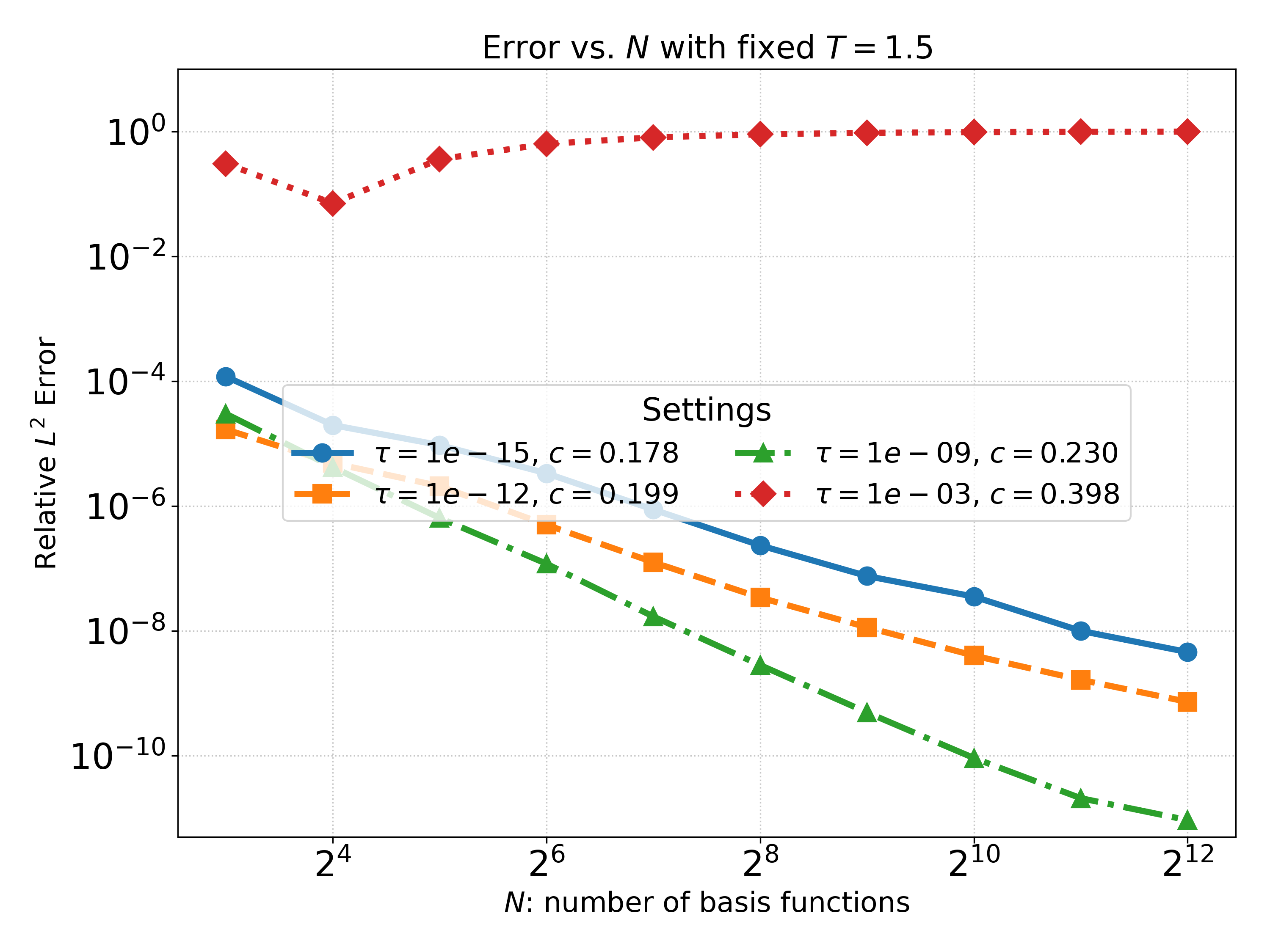}
  \end{subfigure}\hfill
  \begin{subfigure}[t]{0.5\textwidth}
    \centering\includegraphics[width=1\linewidth]{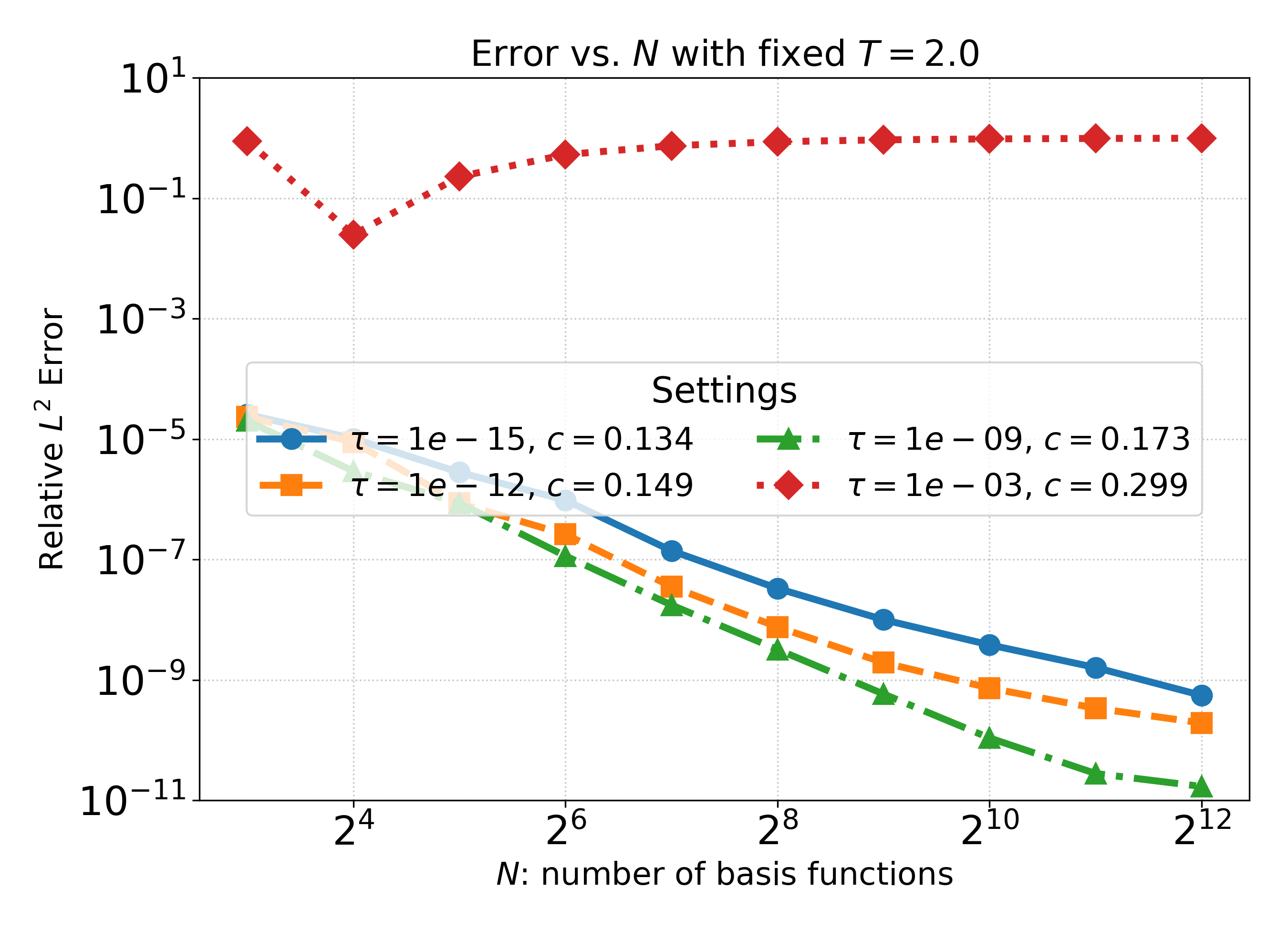}
  \end{subfigure}

  \vspace{0.1ex}

  \begin{subfigure}[t]{0.5\textwidth}
    \centering\includegraphics[width=1\linewidth]{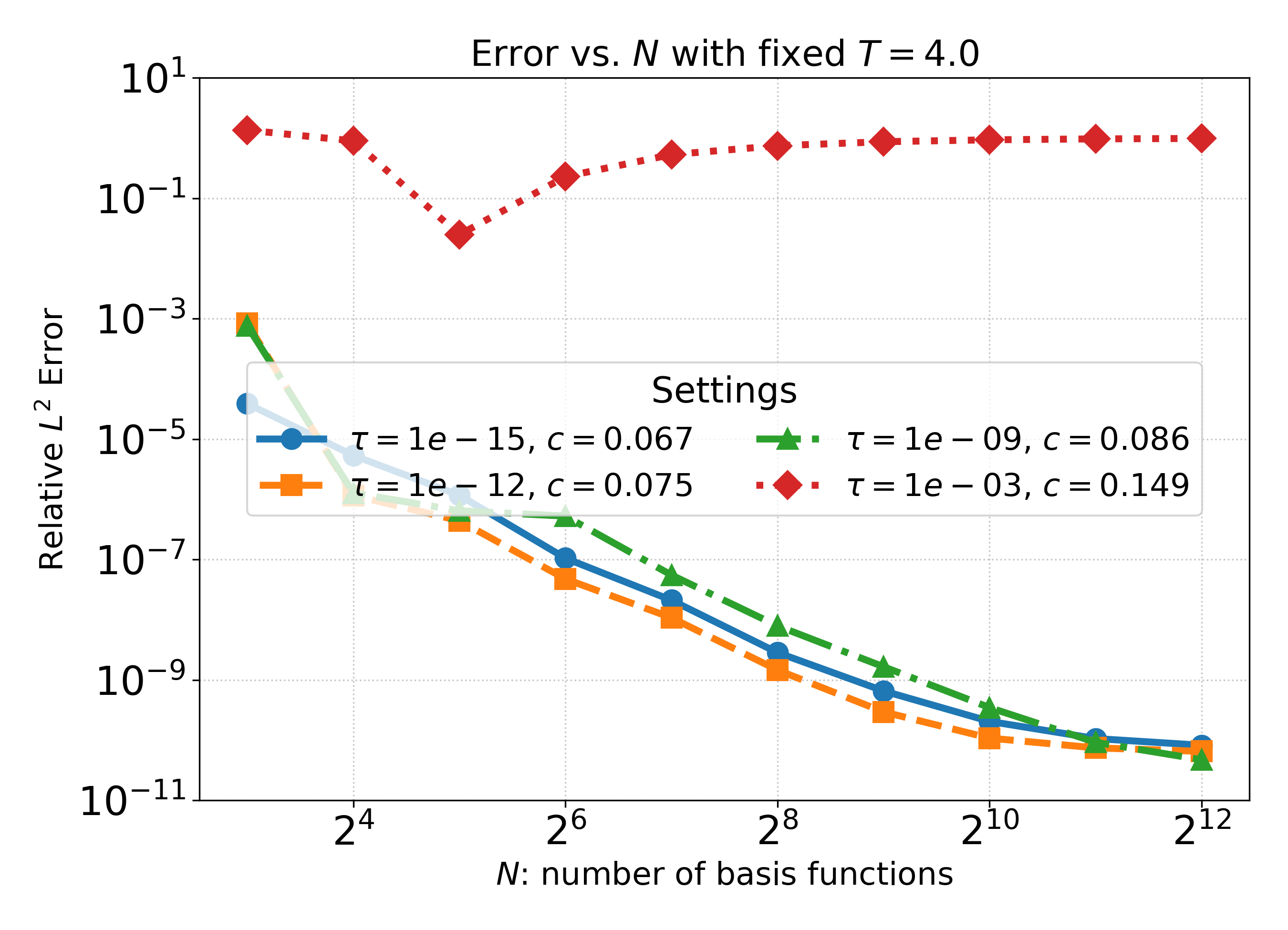}
  \end{subfigure}\hfill
  \begin{subfigure}[t]{0.5\textwidth}
    \centering\includegraphics[width=1\linewidth]{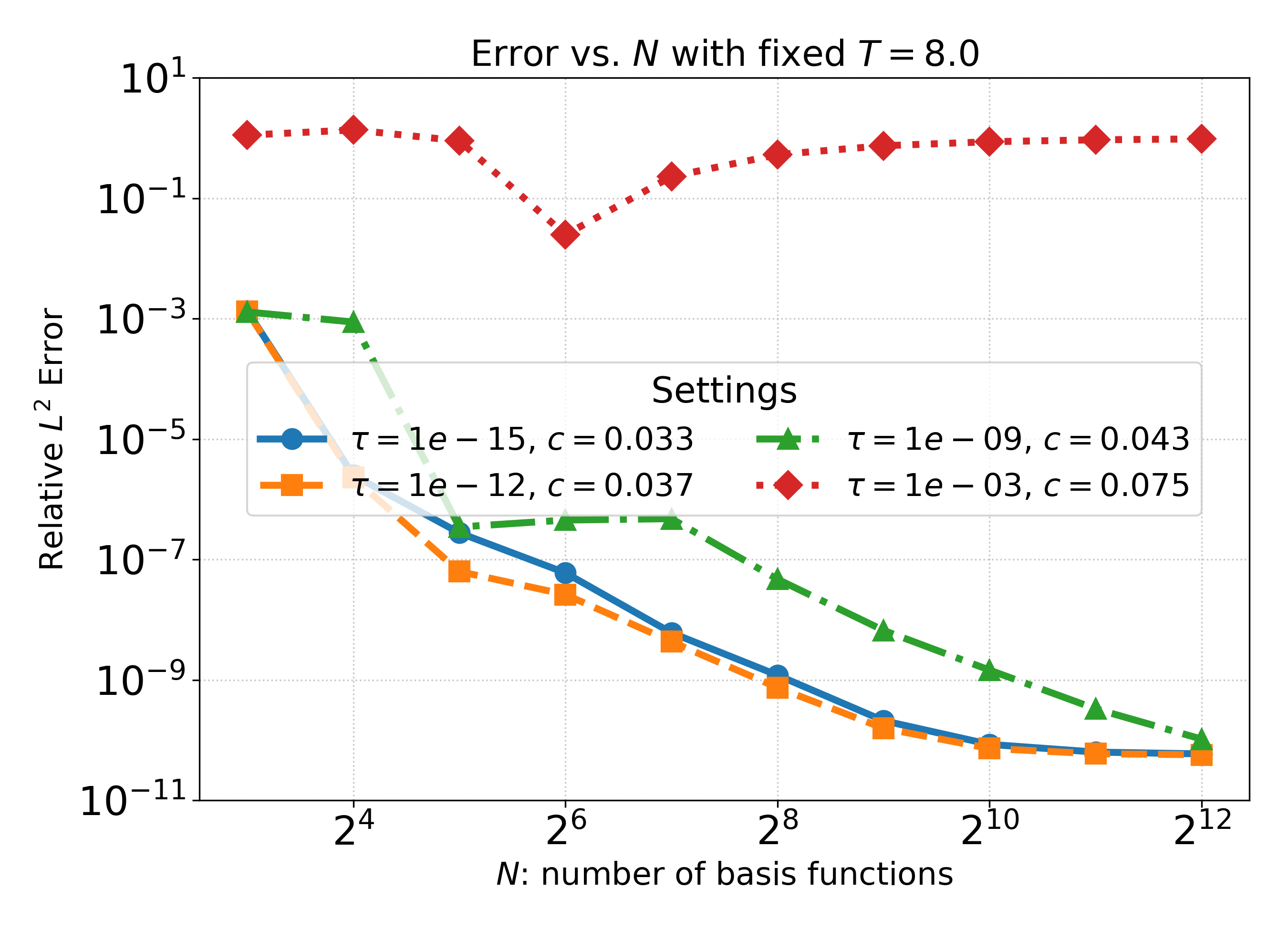}
  \end{subfigure}

\caption{Relative \(\ell^2\) error versus the number of basis functions \(N\) for the 1D Poisson problem at fixed expansion factor \(T\), with oversampling ratio \(\zeta=2\). In each panel, we vary the truncation threshold \(\tau \in \{10^{-3},10^{-9},10^{-12},10^{-15}\}\), and for each \(\tau\) the proportionality constant \(c\) is chosen according to the rule \(c=c(T,\tau)\). The results show that the truncation threshold has a pronounced effect on the large-\(N\) behavior: coarse truncation, such as \(\tau=10^{-3}\), leads to poor accuracy and early stagnation, whereas smaller thresholds yield stable error decay and substantially lower error floors. Overall, \(\tau=10^{-12}\) and \(\tau=10^{-15}\) provide the most accurate and robust performance across the tested values of \(T\).}
  \label{fig6_poisson_pde}
\end{figure}

\begin{figure}[p]
  \centering
  \begin{subfigure}[t]{0.5\textwidth}
    \centering\includegraphics[width=1\linewidth]{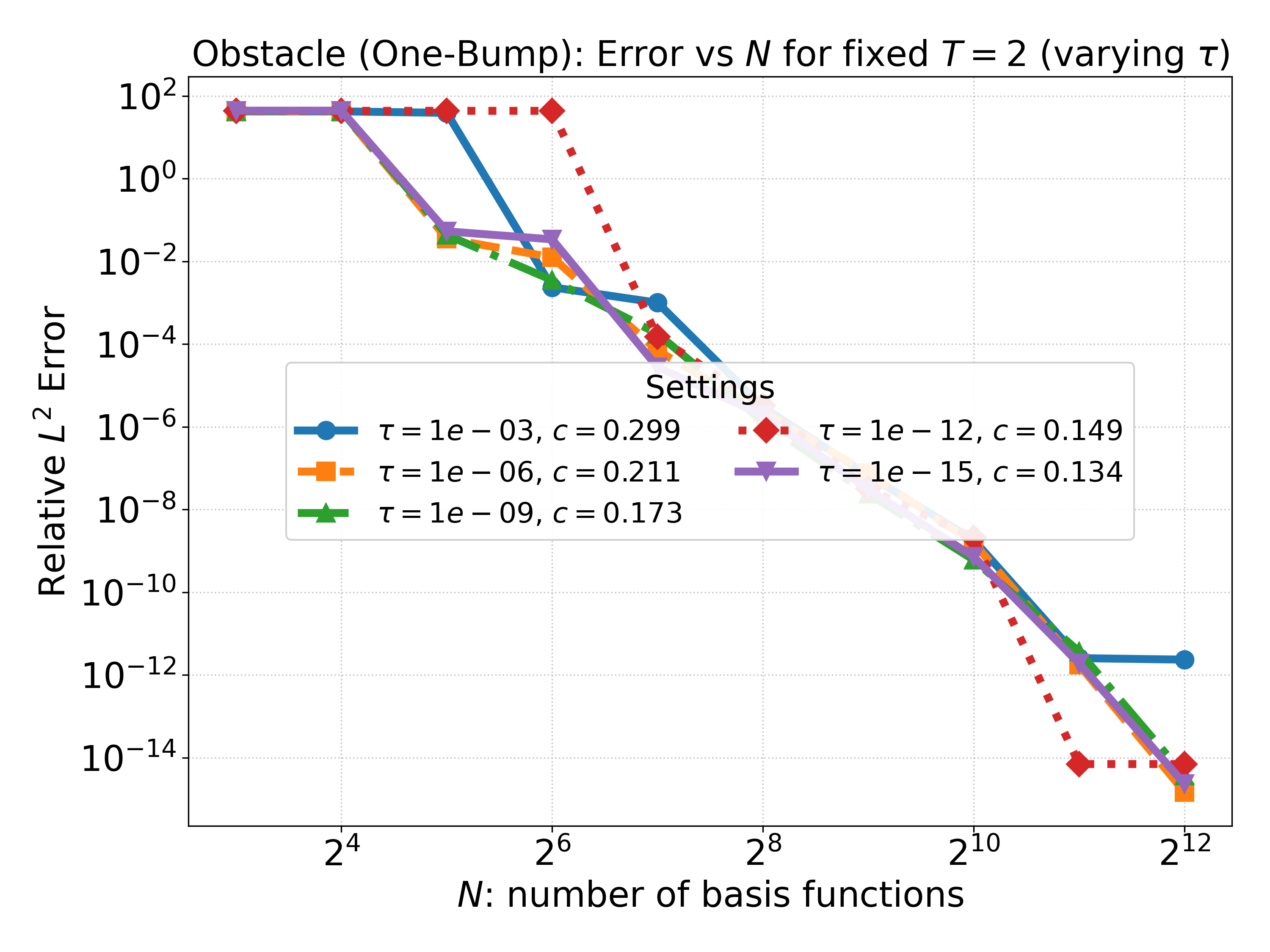}
  \end{subfigure}\hfill
  \begin{subfigure}[t]{0.5\textwidth}
    \centering\includegraphics[width=1\linewidth]{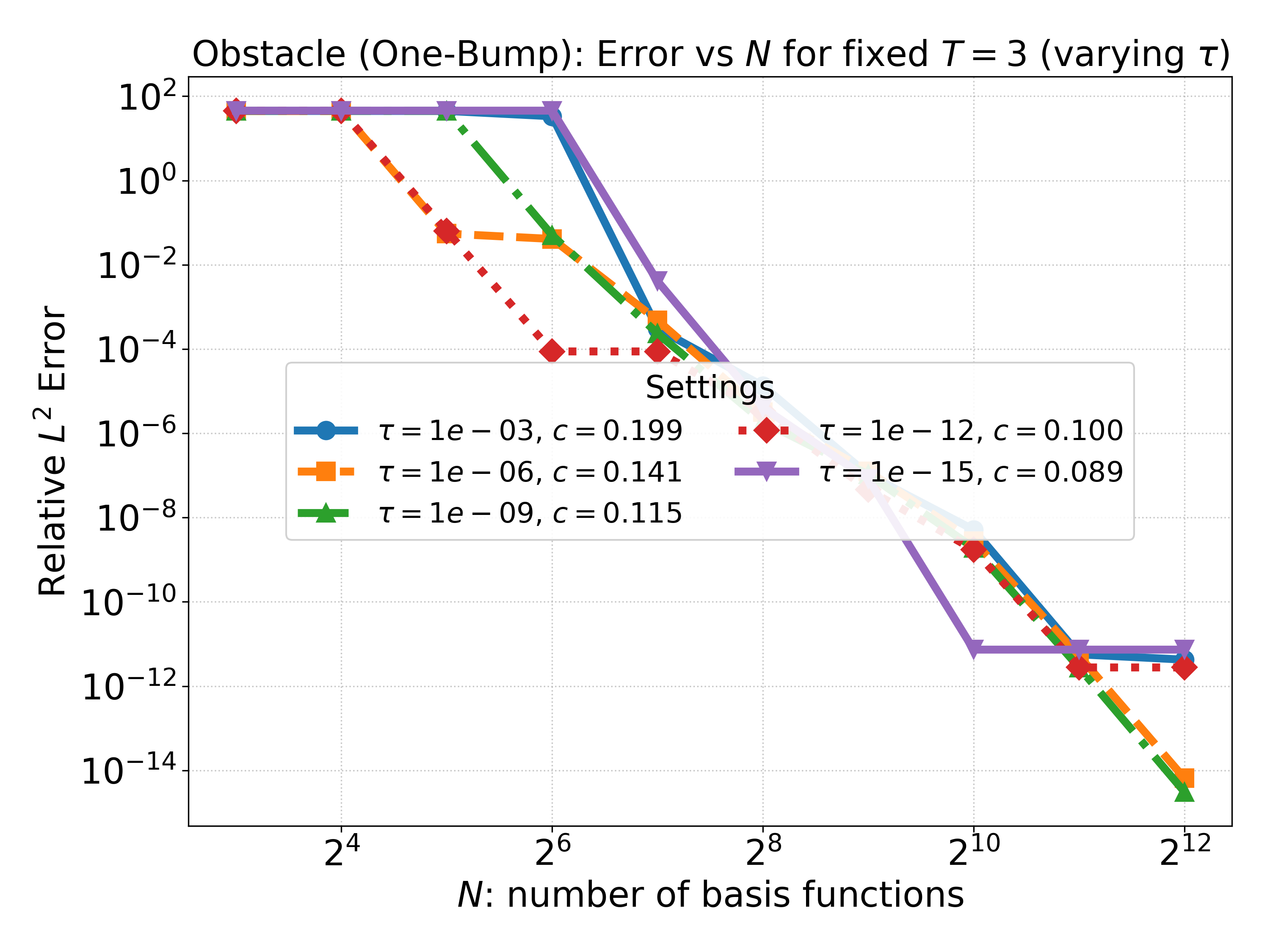}
  \end{subfigure}

  \vspace{0.1ex}

  \begin{subfigure}[t]{0.5\textwidth}
    \centering\includegraphics[width=1\linewidth]{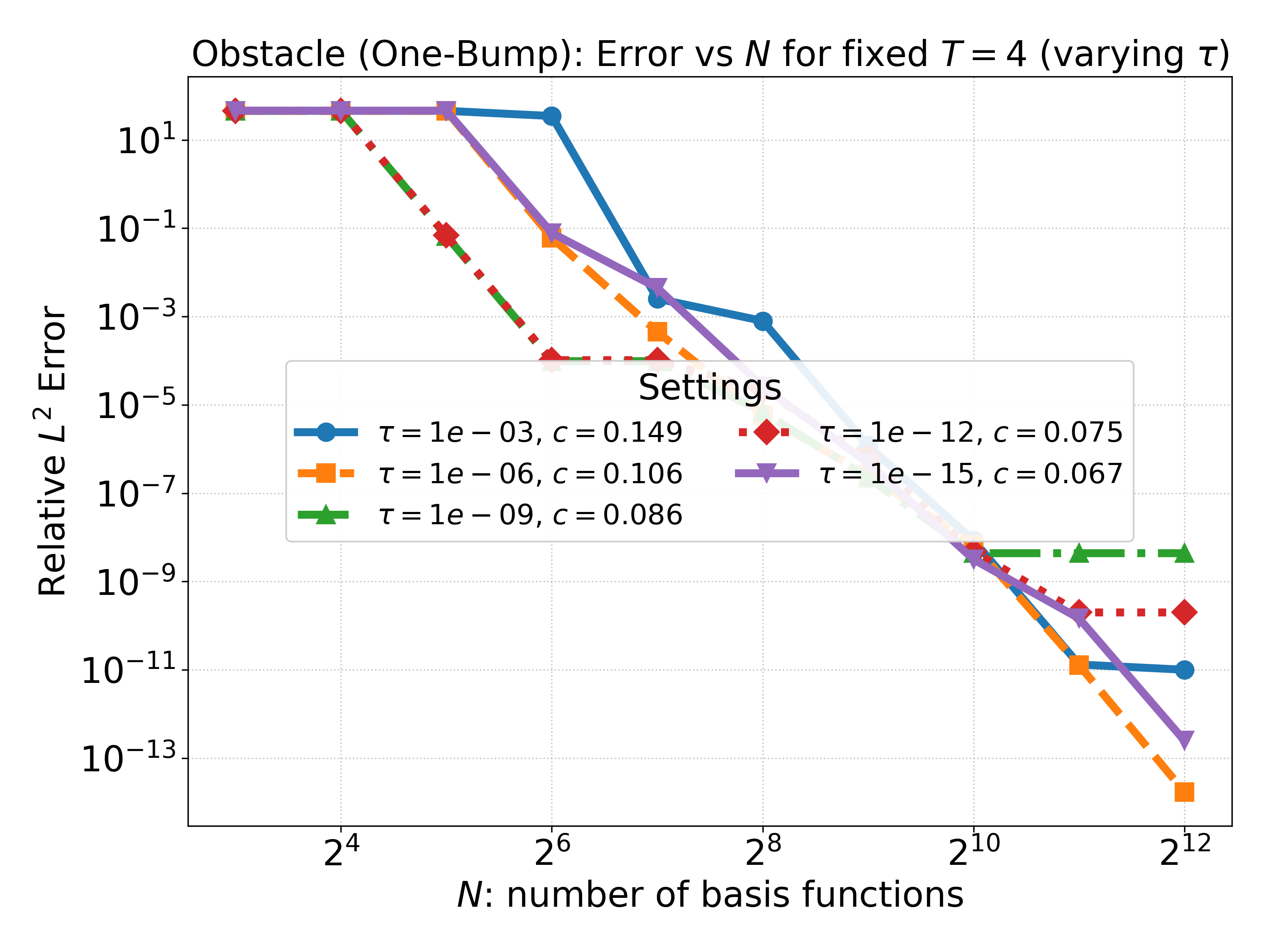}
  \end{subfigure}\hfill
  \begin{subfigure}[t]{0.5\textwidth}
    \centering\includegraphics[width=1\linewidth]{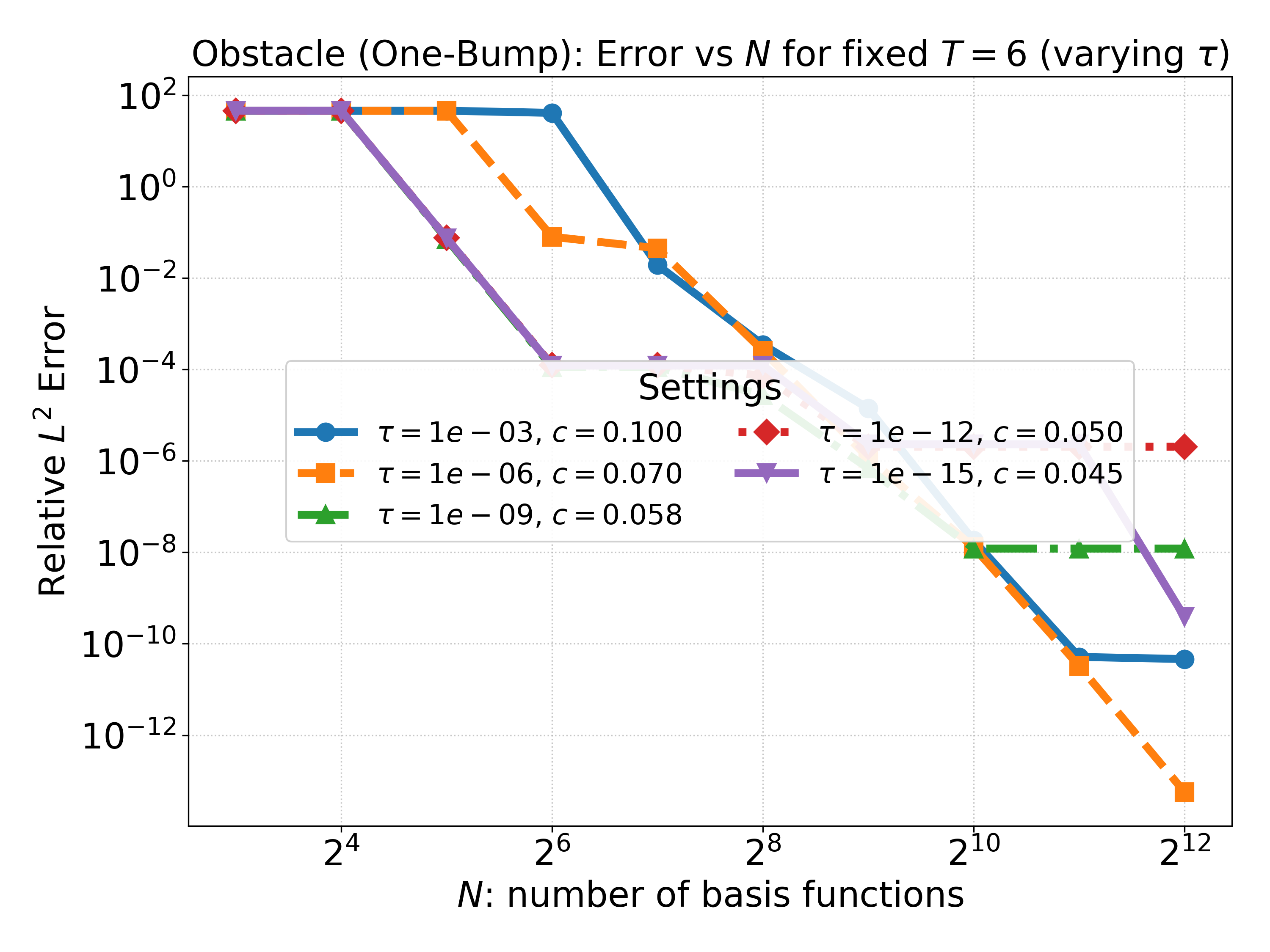}
  \end{subfigure}

\caption{Relative \(\ell^2\) error versus the number of basis functions \(N\) for the one-bump obstacle problem, shown for fixed expansion factors \(T \in \{2,3,4,6\}\), with oversampling ratio \(\zeta=2\). In each panel, the truncation threshold \(\tau\) is varied, and the proportionality constant in the linear shape-parameter scaling is selected as \(c=c(T,\tau)\). For all values of \(T\), the approximation improves substantially as \(N\) increases once \(\tau\) is chosen sufficiently small. In contrast, larger truncation thresholds may produce premature stagnation or a higher error floor. The results indicate that the method is most reliable for smaller values of \(\tau\), especially in the large-\(N\) regime.}
  
  \label{fig2_benchmark_obstacle1d_one_bump_varyingtau}
\end{figure}

\begin{figure}[p]
  \centering
  \begin{subfigure}[t]{0.5\textwidth}
    \centering\includegraphics[width=1\linewidth]{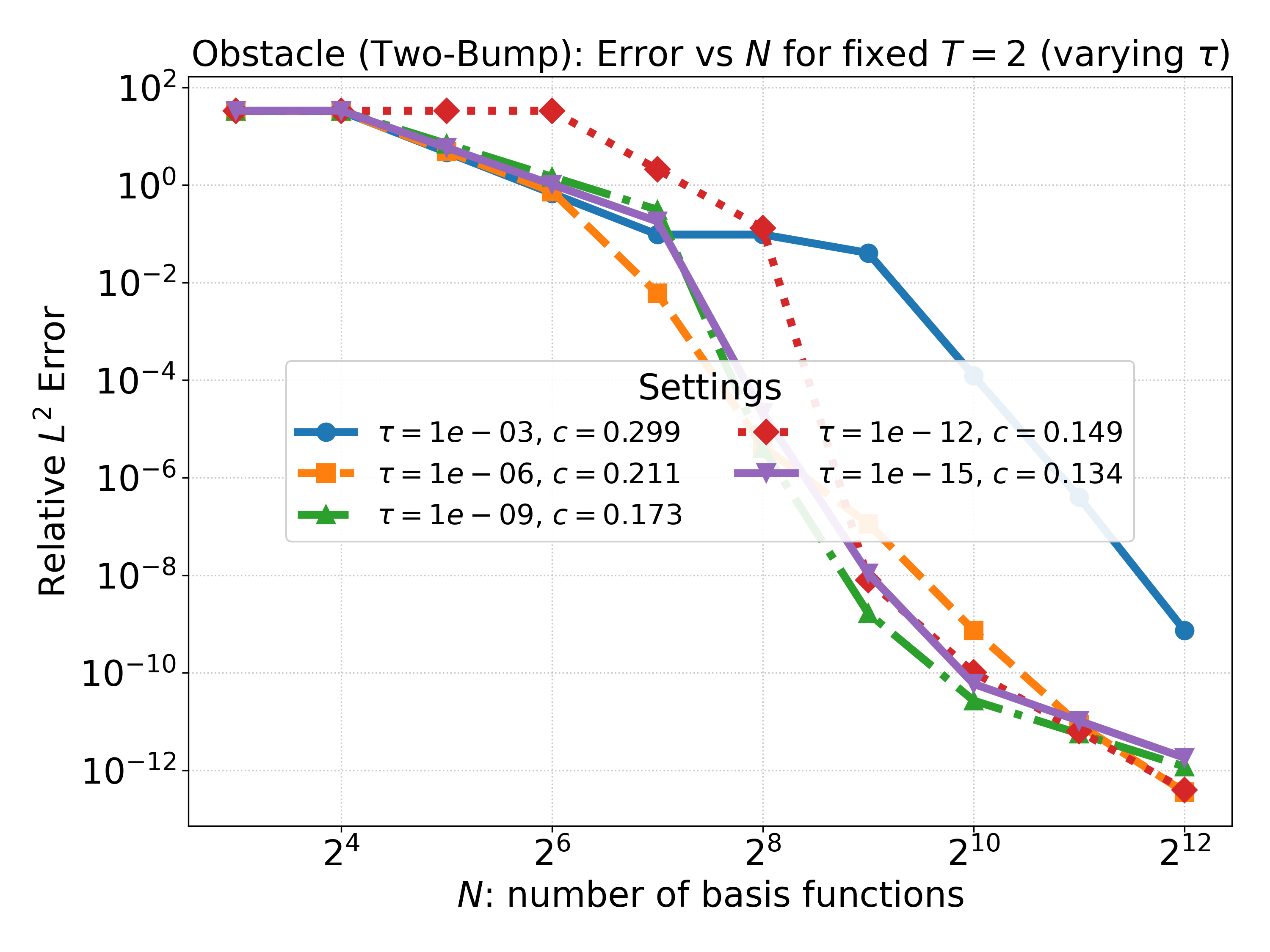}
  \end{subfigure}\hfill
  \begin{subfigure}[t]{0.5\textwidth}
    \centering\includegraphics[width=1\linewidth]{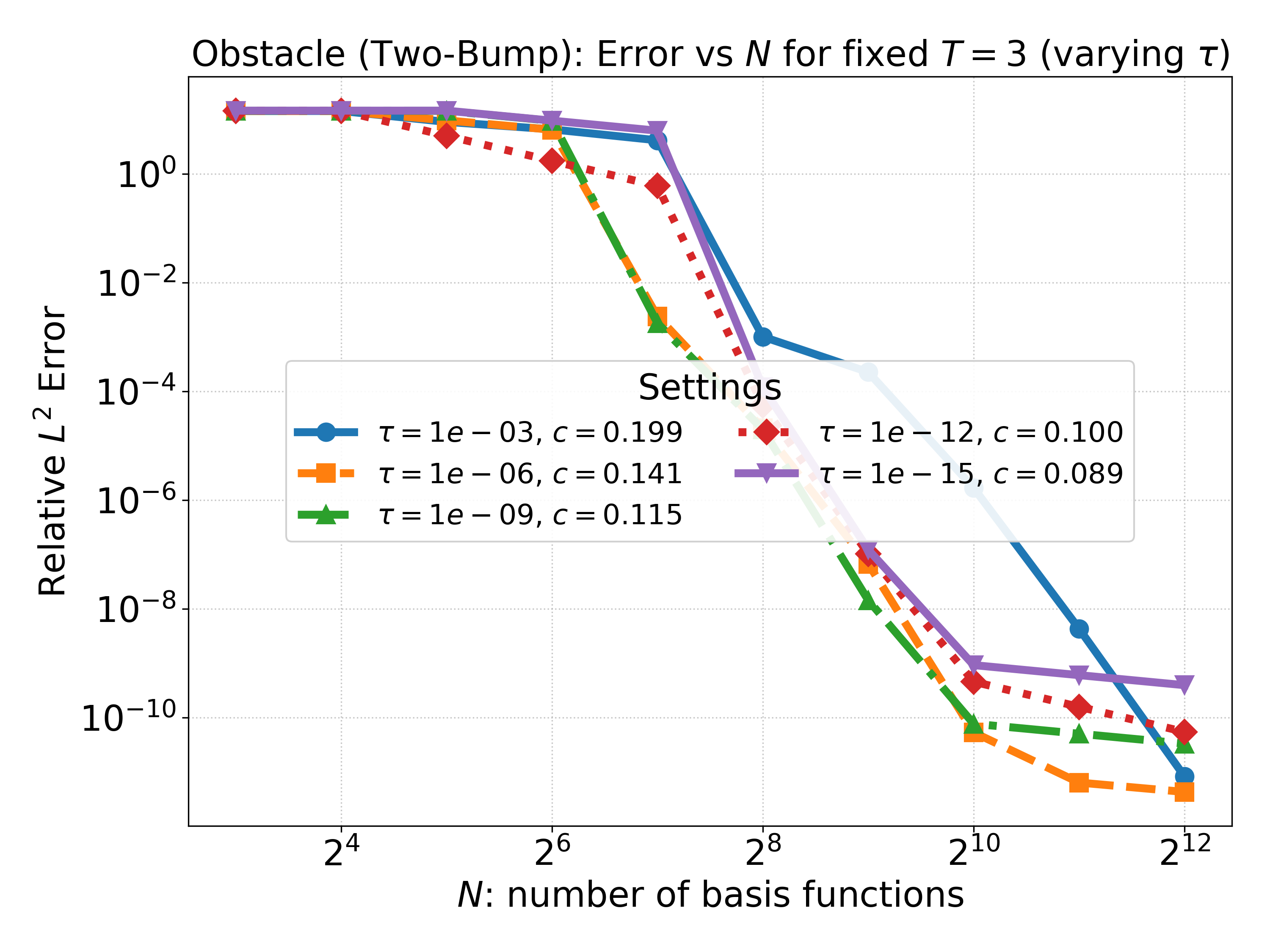}
  \end{subfigure}

  \vspace{0.1ex}

  \begin{subfigure}[t]{0.5\textwidth}
    \centering\includegraphics[width=1\linewidth]{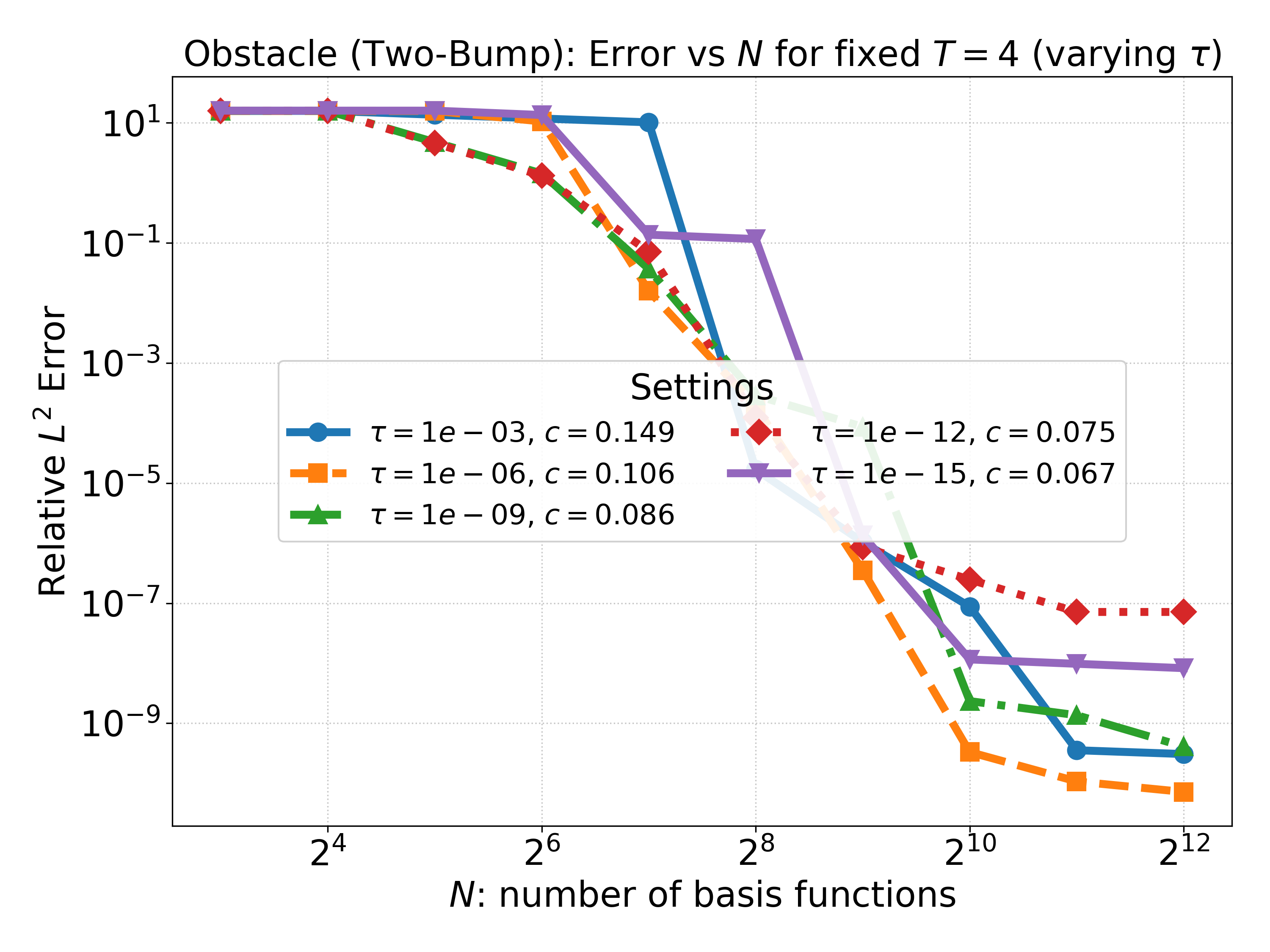}
  \end{subfigure}\hfill
  \begin{subfigure}[t]{0.5\textwidth}
    \centering\includegraphics[width=1\linewidth]{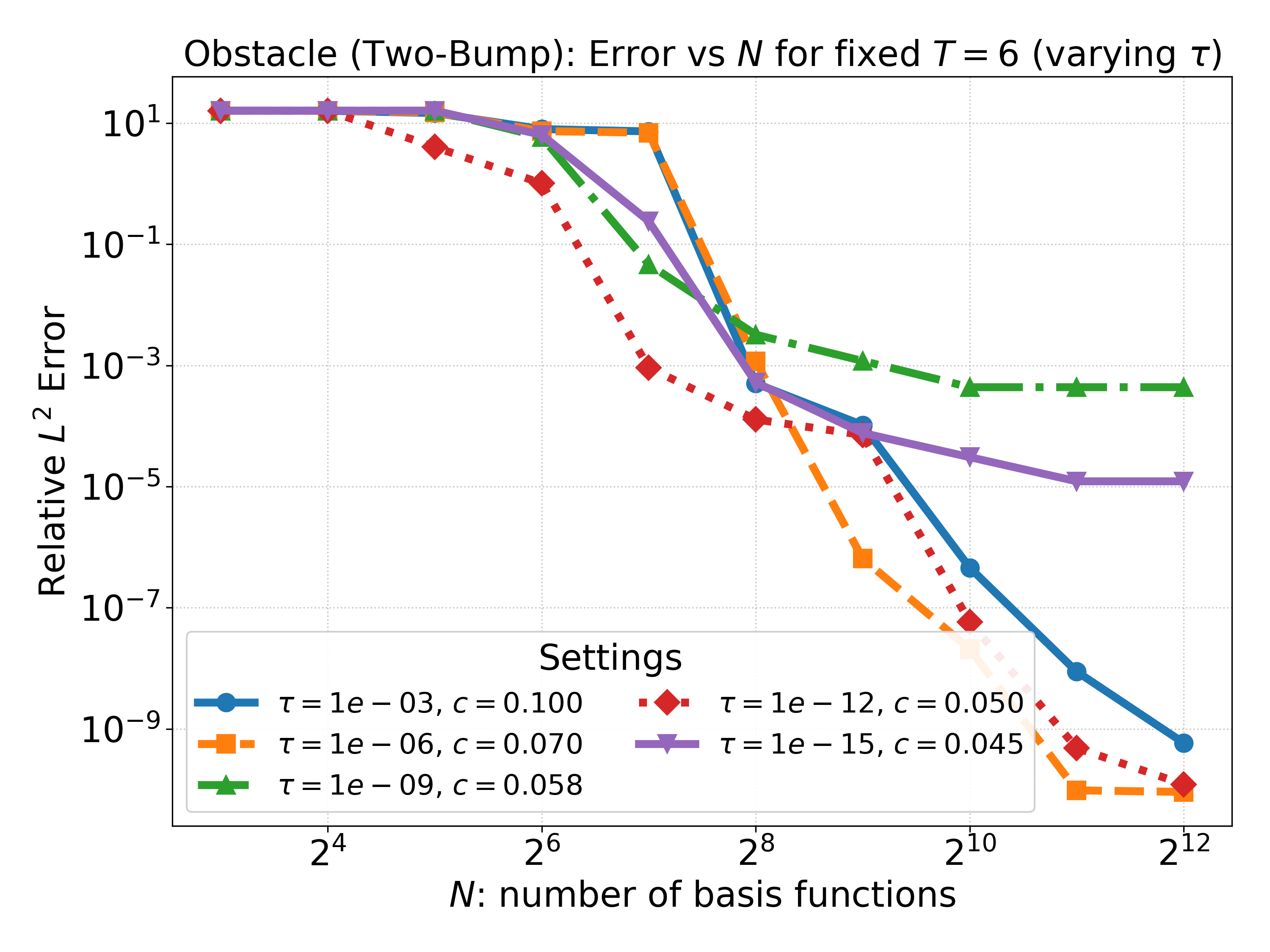}
  \end{subfigure}

\caption{Relative \(\ell^2\) error versus the number of basis functions \(N\) for the two-bump obstacle problem, with fixed expansion factor \(T\) in each panel, oversampling ratio \(\zeta=2\), and shape-parameter constant chosen by \(c=c(T,\tau)\). The figure compares several truncation thresholds \(\tau\). As in the one-bump case, smaller truncation thresholds generally lead to better asymptotic accuracy, whereas larger values of \(\tau\) may cause earlier stagnation or reduced stability. The sensitivity to \(\tau\) becomes more visible for larger \(T\), where aggressive truncation discards too much spectral information and limits the attainable accuracy.}

  \label{fig2_benchmark_obstacle1d_two_bump_varyingtau}
\end{figure}

\begin{figure}[p]
  \centering
  \begin{subfigure}[t]{0.5\textwidth}
    \centering\includegraphics[width=1\linewidth]{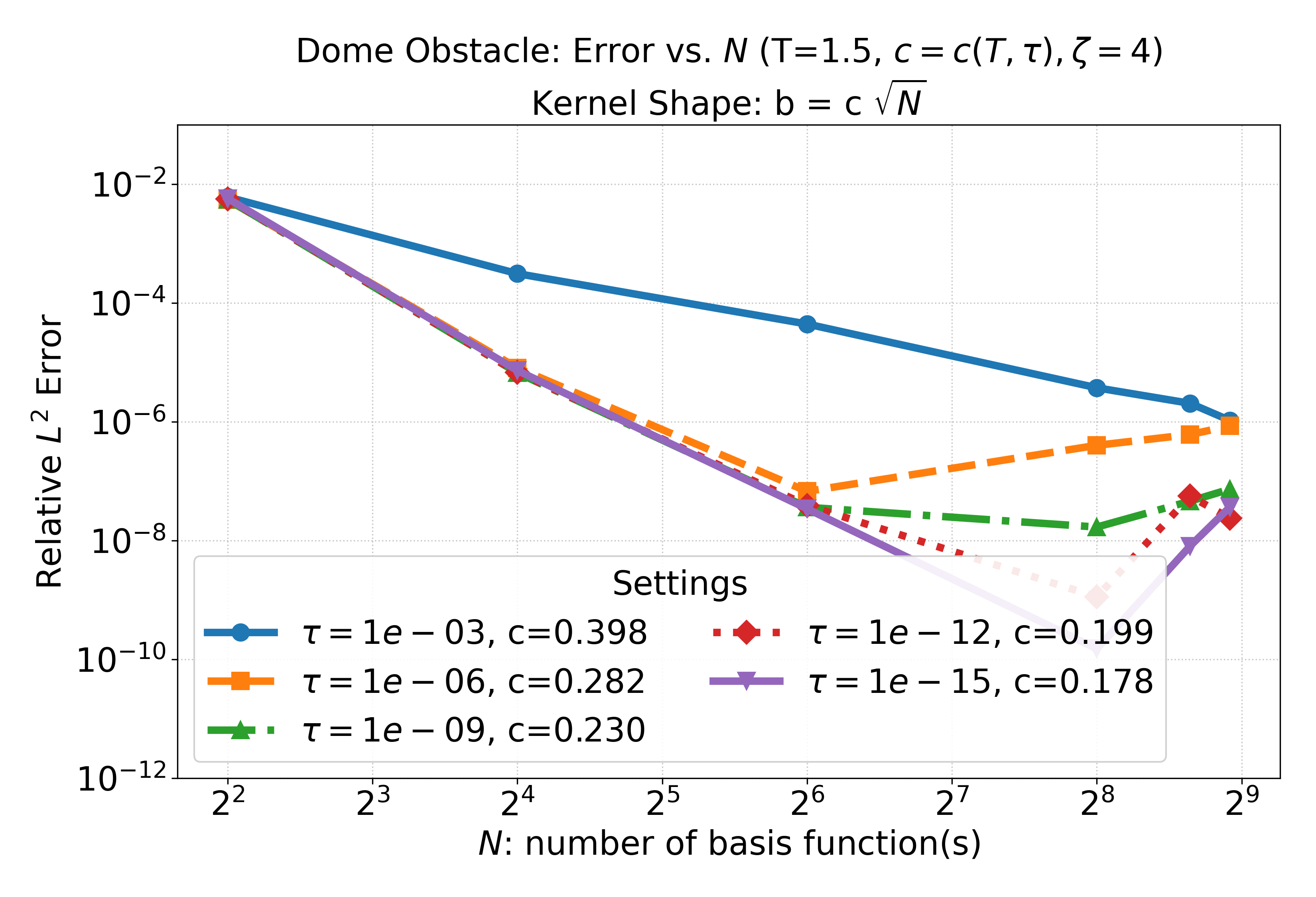}
  \end{subfigure}\hfill
  \begin{subfigure}[t]{0.5\textwidth}
    \centering\includegraphics[width=1\linewidth]{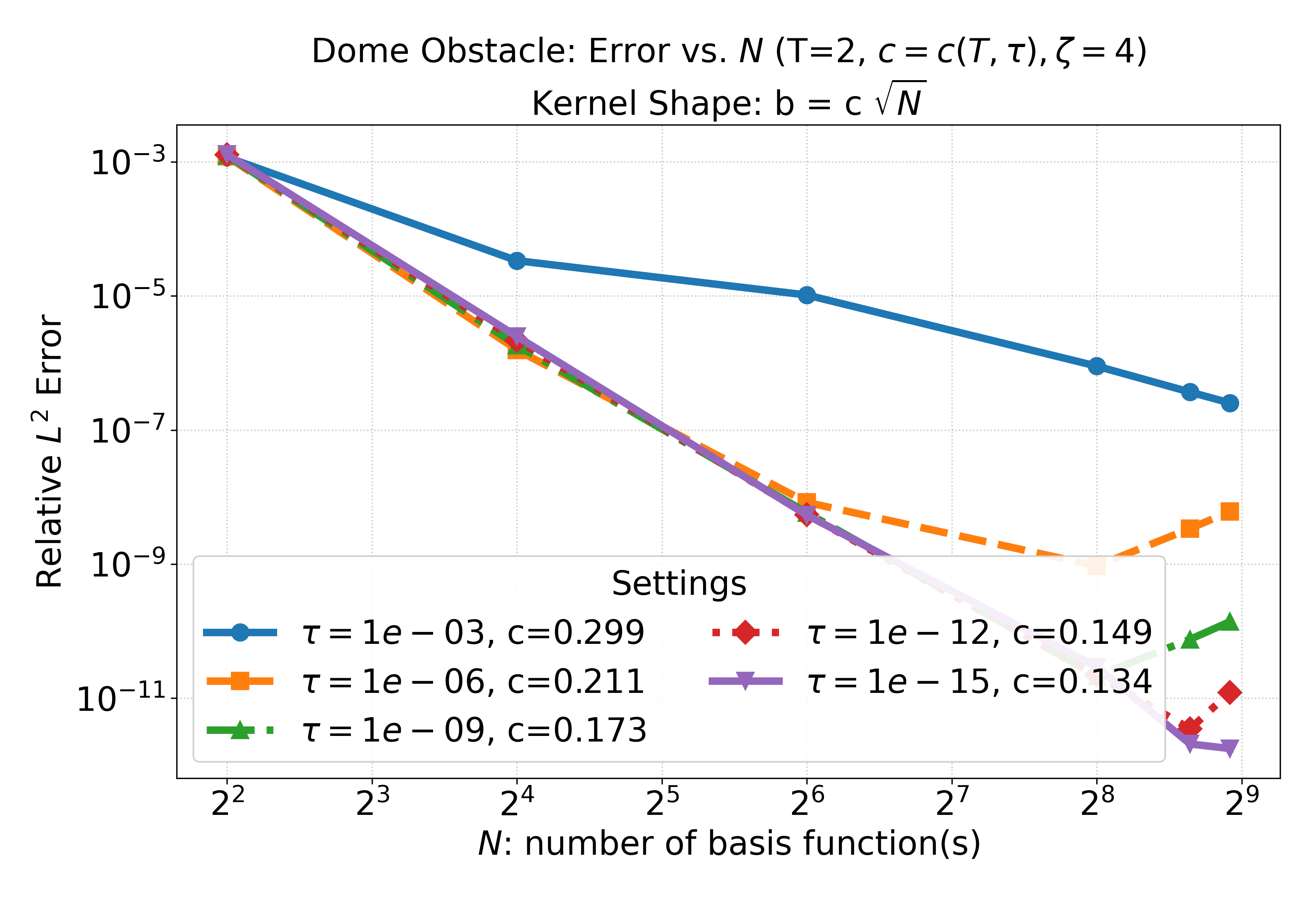}
  \end{subfigure}

  \vspace{0.1ex}

  \begin{subfigure}[t]{0.5\textwidth}
    \centering\includegraphics[width=1\linewidth]{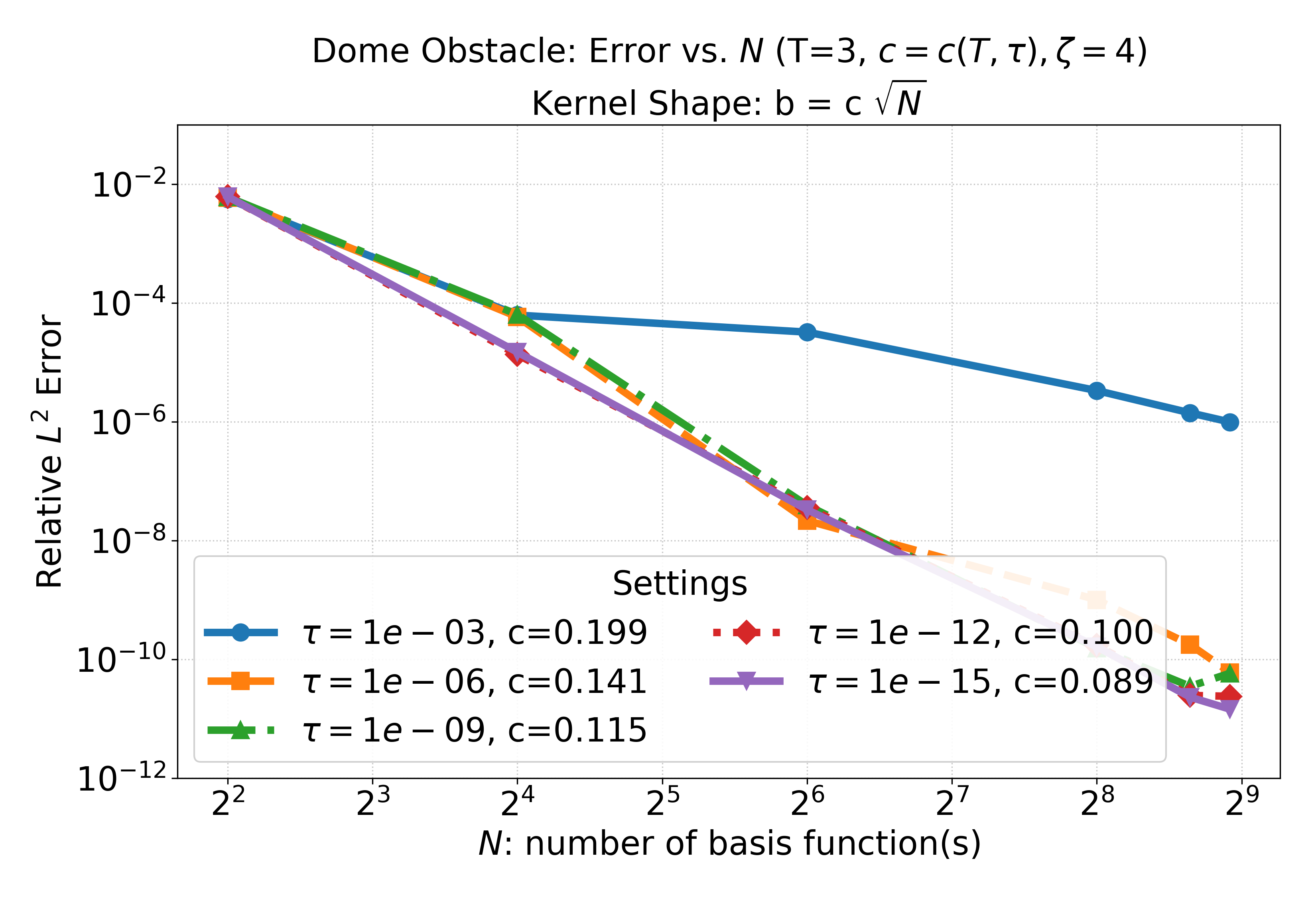}
  \end{subfigure}\hfill
  \begin{subfigure}[t]{0.5\textwidth}
    \centering\includegraphics[width=1\linewidth]{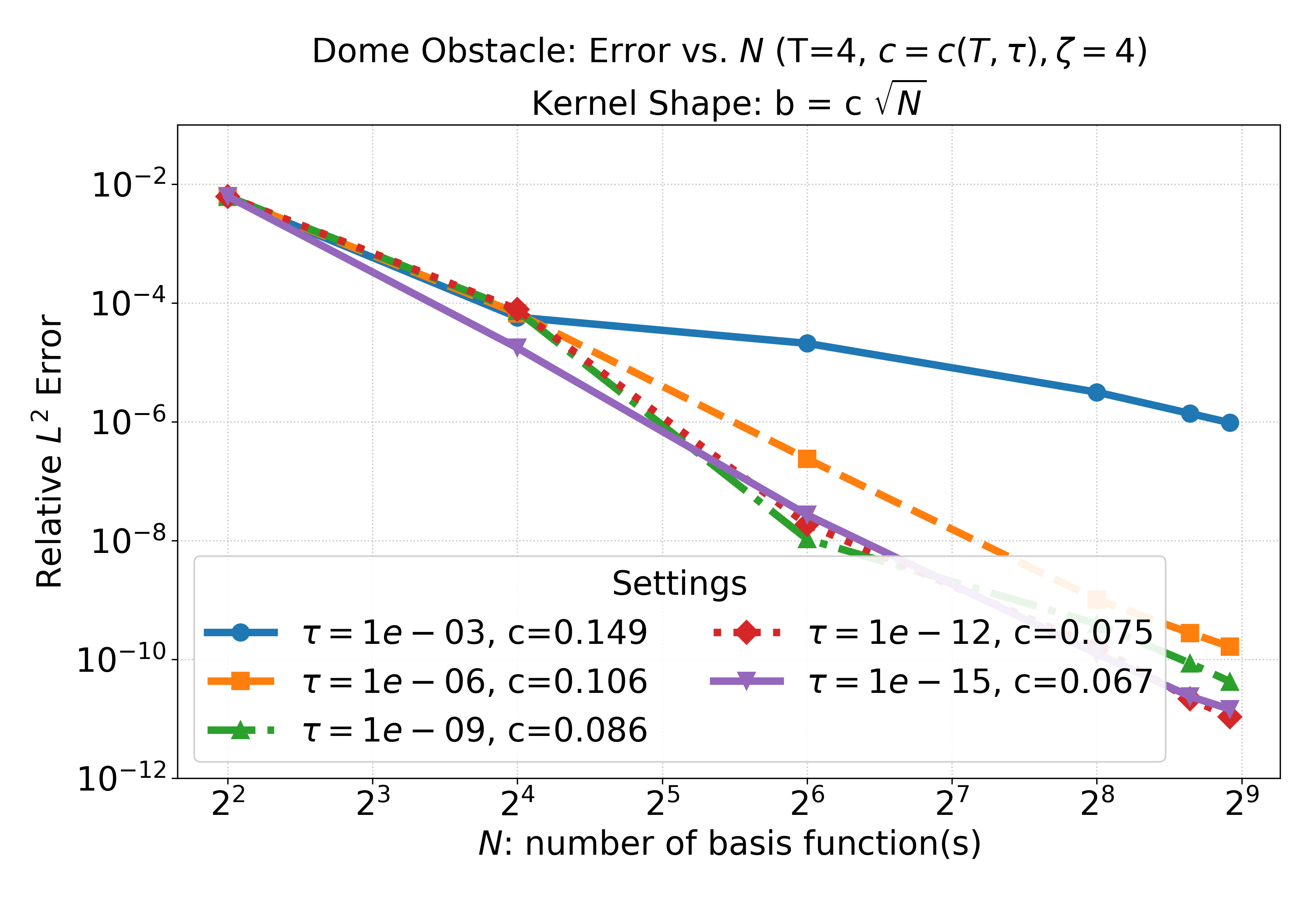}
  \end{subfigure}
\caption{Relative \(\ell^2\) error versus the number of basis functions \(N\) for the two-dimensional dome obstacle problem, shown for several expansion factors \(T\). In each panel, we vary the truncation threshold \(\tau\), while the shape parameter is selected through the rule \(c=c(T,\tau)\). The results show that smaller truncation thresholds produce lower error floors and more sustained decay as \(N\) increases, while larger thresholds lead to earlier saturation. Across all tested values of \(T\), the choices \(\tau=10^{-12}\) and \(\tau=10^{-15}\) yield the most accurate overall performance, with the difference between thresholds becoming more pronounced for larger expansion factors.}
\label{fig:2d_obstacle_problem_dome_benchmark_tau_parameter}
\end{figure}

We next investigate the effect of the truncation threshold \(\tau\) used in the TSVD stabilization. Figures~\ref{fig6_poisson_pde}--\ref{fig:2d_obstacle_problem_dome_benchmark_tau_parameter} compare several values of \(\tau\), where for each pair \((T,\tau)\) the proportionality constant is selected by the rule \(c=c(T,\tau)\) of Eq.~\eqref{eq:optimal-value-c-formula-from-adcock}. Overall, the results show that the truncation threshold has a strong influence on the large-\(N\) regime. If \(\tau\) is too large, then too many small singular components are removed, which leads to premature stagnation and a noticeably higher error floor. By contrast, smaller thresholds retain more spectral information and usually produce more sustained error decay as \(N\) increases.

For the 1D Poisson problem in Figure~\ref{fig6_poisson_pde}, \(\tau=10^{-12}\) and \(\tau=10^{-15}\) give the most accurate results across the tested values of \(T\), while \(\tau=10^{-3}\) performs poorly and \(\tau=10^{-9}\) is typically less accurate in the large-\(N\) regime. A similar trend is observed in the one-bump and two-bump obstacle problems in Figures~\ref{fig2_benchmark_obstacle1d_one_bump_varyingtau} and~\ref{fig2_benchmark_obstacle1d_two_bump_varyingtau}, although the sensitivity to \(\tau\) depends on the expansion factor \(T\). In particular, smaller values of \(\tau\) generally lead to lower final errors, whereas larger values may cause early plateauing. The same behavior is also visible in the dome obstacle problem in Figure~\ref{fig:2d_obstacle_problem_dome_benchmark_tau_parameter}, where smaller truncation thresholds produce lower error floors and more stable convergence.

Based on these experiments, we use \(\tau=10^{-15}\) in the main experiments. Although \(\tau=10^{-12}\) is sometimes comparable and can be slightly better in isolated cases, \(\tau=10^{-15}\) is the most consistently reliable choice across the full set of test problems. It avoids the premature truncation associated with larger thresholds and provides stable, high-accuracy performance in all examples considered here.

\pagebreak

\pagebreak

\pagebreak

\pagebreak

\pagebreak

\pagebreak

\pagebreak

\pagebreak

\pagebreak

\pagebreak

\pagebreak

\section{Reaction-Diffusion Boundary Value Problem}\label{appendix:sect_reaction-diffusion-boundary-problem}
In this section, we apply Algorithm~\ref{algorithm1-rbf-solver-for-elliptic-pde} to reaction-diffusion boundary value problems. We recall the general strong formulation for reaction-diffusion as follows: 
\begin{equation}
\begin{cases}
- \Delta u + u = f(\vec{x}),& \forall \vec{x} \in \Omega, \\ 
u(\vec{x}) = g(\vec{x}),& \forall \vec{x} \in \partial \Omega.
\end{cases}
\label{reaction-diffusion-bvp}
\end{equation}
The corresponding modified penalty method functional for the reaction-diffusion boundary value problem in Eq.~\eqref{reaction-diffusion-bvp} yields as follows: 
\begin{equation}
\mathcal{J}_{\beta, r.d.} (u) = \frac{1}{2} \int_{\Omega} \Big(  \| \nabla u \|_{2}^{2} + u^{2} \Big) \, d \vec{x} - \int_{\Omega} f \, u \, d\vec{x} + \frac{\beta}{2} \int_{\partial \Omega} ( u - g)^{2} \, d \vec{s}.
\label{weak-form-reaction-diffusion-problem}
\end{equation}
For the next section~\ref{appendix:sect_rd-bvp-variational-functional}, we consider the discretized variational functional.

\subsection{Discretized Functional for Reaction-Diffusion Boundary Value Problem}\label{appendix:sect_rd-bvp-variational-functional}

In Section~\ref{sec:methodology}, we selected $m$ collocation points from the closed domain $\overline{\Omega}$ which includes $m_{\Omega}$ interior points and $m_{\partial \Omega}$ boundary points. Using the RBF approximation in Eq.~\eqref{eq:main-approximation-form-rbf}, we convert the functional in Eq.~\eqref{weak-form-reaction-diffusion-problem} into a discretized functional. The minimization problem for the reaction-diffusion problem is given by
\begin{equation}
\min \limits_{u \in W^{1,2}} \mathcal{J}_{\beta, \mathrm{r.d.}}(u) \quad \text{where}\quad \mathcal{J}_{\beta, \mathrm{r.d.}}(u) =  \frac{1}{2} \int_{\Omega} \big( \| \nabla u \|_{2}^{2} + u^{2} \big) \, d \vec{x} - \int_{\Omega} f \, u \, d  \vec{x} + \frac{\beta}{2} \int_{\partial \Omega} (u - g)^{2} \, d  \vec{s}. 
\nonumber
\end{equation}
The discrete form of Eq.~\eqref{weak-form-reaction-diffusion-problem} is 
\begin{equation}
\widehat{\mathcal{J}}_{\beta,\mathrm{r.d.}}(\vec{x})  = \frac{| \Omega |}{m_{\Omega}} \sum_{j = 1}^{m_{\Omega}} \left( 
    \frac{1}{2} \big ( \| \nabla \hat{u}(\vec{x}^{(\Omega)}_{j}) \|_{2}^{2} 
    + \hat{u}^{2} (\vec{x}^{(\Omega)}_{j}) \big)
    - f(\vec{x}^{(\Omega)}_{j}) \hat{u}(\vec{x}^{(\Omega)}_{j}) 
\right)  + \frac{\beta \, | \partial \Omega| }{2 \, m_{\partial \Omega}} 
\sum_{n = 1}^{m_{\partial \Omega}} \left( 
    \hat{u}(\vec{x}^{(\partial \Omega)}_{n}) - g(\vec{x}^{(\partial \Omega)}_{n}) 
\right)^{2}.  
\label{eq:discretized-weak-form-reaction-diffusion} 
\end{equation}
More specifically, we solve 
\begin{align}
\vec{w}_{\mathrm{r.d.}}&= 
\operatorname*{\,min}_{\vec{w}\in\mathbb{R}^{N}}
\widehat{\mathcal{J}}_{\beta,\text{r.d.}}(\vec{w}),
\label{eq:reaction-diffusion_opt} \\[4pt]
\widehat{\mathcal{J}}_{\beta,\text{r.d.}}(\vec{w}) &=
\frac{1}{2}\,\vec{w}^{\top}A^{\small \small (\mathrm{r.d.})}_{1} \, \vec{w}
-\vec{f}^{\top}A_{2}\vec{w}
+\frac{\beta}{2}\,\lVert A_{3}\vec{w}-\vec{g}\rVert_{2}^{2}.
\label{eq:reaction-diffusion_loss}
\end{align}
The stiffness matrix in Eq.~\eqref{eq:reaction-diffusion_loss}, $A^{(\small \mathrm{r.d.})}_{1} \in \mathbb{R}^{N \times N}$ is given by 
\begin{equation}
(A^{(\small \mathrm{r.d.)}}_{1})_{i,\ell} = \frac{|\Omega|}{m_{\Omega} \, N }  \sum_{j} \nabla \phi^{\top} ( \| \vec{x}^{( \Omega)}_{j} - \vec{c}_{i}  \|_{2} ) \, \nabla \phi ( \| \vec{x}^{(\Omega)}_{j}  - \vec{c}_{\ell} \|_{2} ) + \phi(\| \vec{x}^{(\Omega)}_{j} - \vec{c}_{i} \|_{2} ) \, \phi^{\top} (\| \vec{x}^{(\Omega)}_{j} - \vec{c}_{\ell} \|_{2}, 
\label{eq:A1-matrix-term-entry-wise-reaction-diffusion}
\end{equation}
where \(1 \leq i, \ell \leq N\). 
This constitutes the only difference between the minimization problems in Eqs.~\eqref{eq:poisson_loss} and \eqref{eq:reaction-diffusion_loss} is the stiffness matrices according to Eqs.~\eqref{eq:A1-matrix-term-entry-wise-poisson} and ~\eqref{eq:A1-matrix-term-entry-wise-reaction-diffusion}; the other terms for the source function and boundary data are the same. 
\subsection{Numerical Experiments for Reaction-Diffusion Boundary Value Problem}\label{appendix:sec_numerical-experiment-rd-bvp}
In this section, we consider the following reaction-diffusion problem:  
\begin{equation}
\begin{cases}
- \Delta u(x) + u(x) = f(x), \quad \forall x \in \Omega = (0, L), \\
u(0) = u(L) = 0,
\end{cases}
\nonumber
\end{equation}
where $L \coloneq 5$ is the length of the domain, and $\partial \Omega = \{ 0, L\}$ denotes the boundary points. 
In this experiment, with
\[
f(x) = \mathcal{A} \, \sin(\pi x), \quad \forall x \in (0, L), 
\]
the solution is given by 
\[
u_{exact}(x) = \frac{\mathcal{A}}{1 + \pi^2} \, \sin(\pi x), \quad \forall x \in [0, L],
\]
where $\mathcal{A} \coloneq 50$ controls the driving force of the governing equation. 

We evaluate the accuracy, stability, and convergence behavior of our method. We compare the finite-difference method with a central second order scheme, a Galerkin method, the RBF-variational mesh-free method approximation of Eq.~\eqref{eq:approximation-solution-u-hat-tau-elliptic} and Algorithm~\ref{algorithm1-rbf-solver-for-elliptic-pde}, and a two-layer shallow neural network with \(N\) neurons. Because the governing equation is an inhomogeneous reaction-diffusion equation according to Eq.~\eqref{reaction-diffusion-bvp}, the presence of both diffusion and reaction terms introduces numerical challenges compared to the linear elliptic experiment in Section~\ref{sec:elliptic-boundary-value-problem-numerical-experiment}. 

\begin{figure}[p]
    \centering
\includegraphics[width=1.0\linewidth]{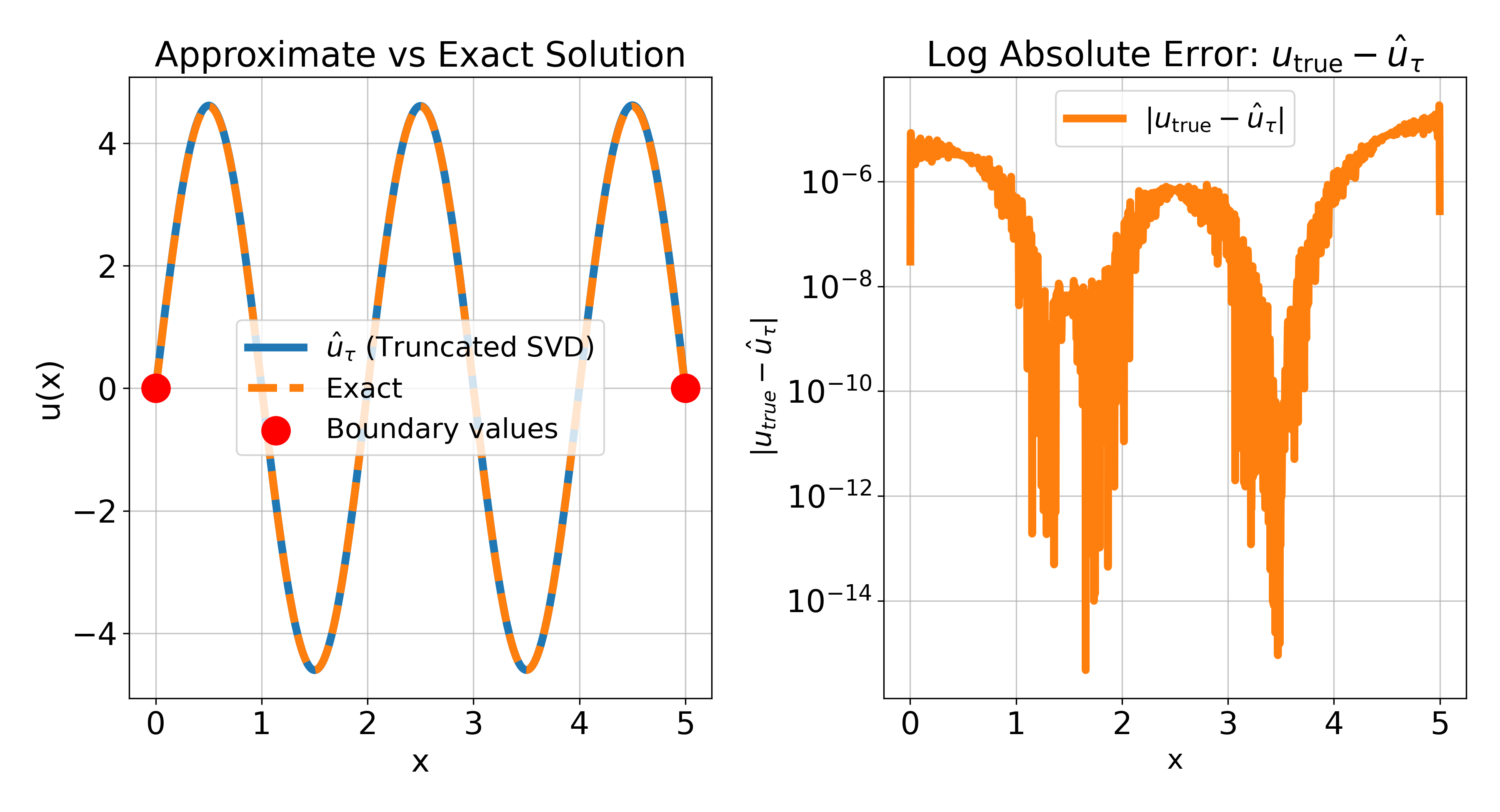}
    \caption{Numerical RBF approximation of $-\Delta u +u=f$ vs.\ analytical solution. The relative error is $E(N) = 7.697 \times 10^{-9}$. Setup: $N=2048$, $m=4096$, $\zeta=2$, $T=2, \tau =  10^{-15}, c = 0.134$.}
    \label{fig1_reactiondiffusion_pde}
\end{figure}

\begin{figure}[H]
    \centering
    \begin{subfigure}[t]{0.5\textwidth}
    \centering\includegraphics[width = 1.0\linewidth]{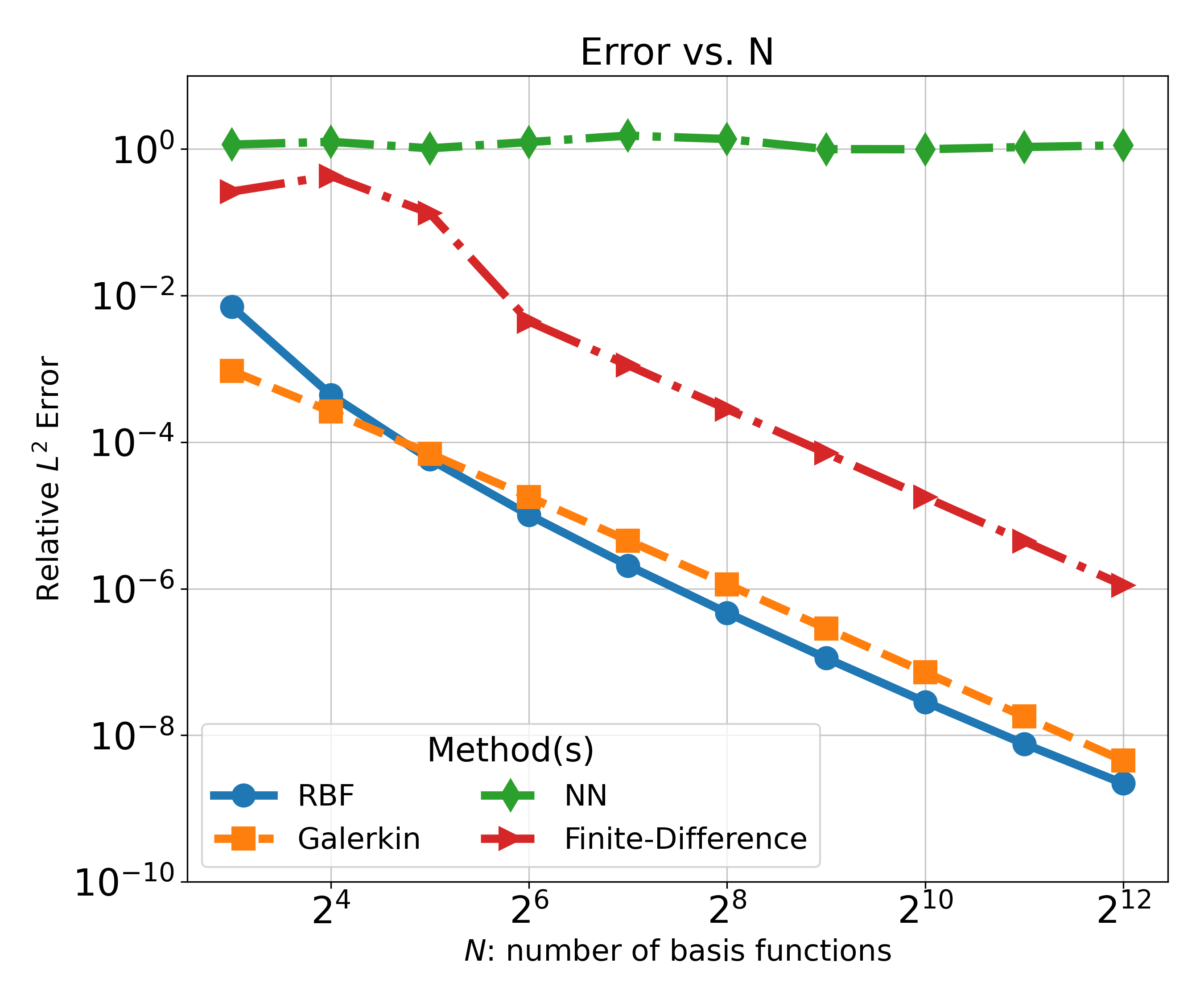}
    \caption{Convergence: relative $\ell^2$ error vs. $N$.}
    \end{subfigure}\hfill 
    \begin{subfigure}[t]{0.5\textwidth}
    \centering\includegraphics[width = 1.0\linewidth]{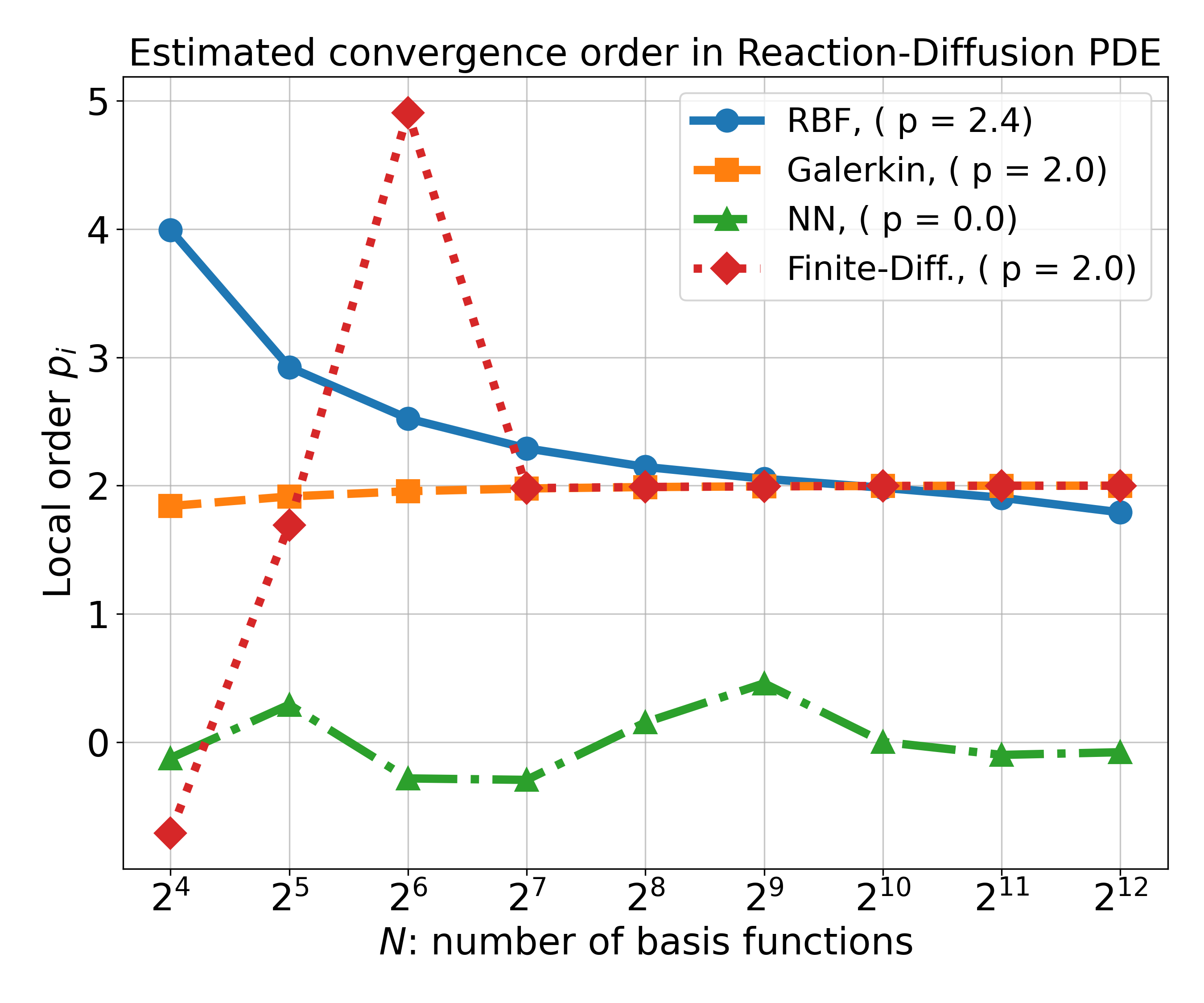}
    \caption{Estimation of convergence order $p_{i}$.}
    \end{subfigure}
    \caption{Convergence for Reaction-Diffusion $-\Delta u  + u  = f$: accuracy and convergence order.}
    \label{fig2_reactiondiffusion_pde}
\end{figure}

In Fig.~\ref{fig1_reactiondiffusion_pde}, we verify the TSVD solver with RBF representation using a resolution of $(N,m) = (2000, 4000)$. With a relative $\ell^2$ error of $7.935 \times 10^{-9}$, the RBF approximation is well-aligned with the analytical solution. The log-absolute error plot confirms the accuracy, with errors remaining uniformly small across the domain. Similarly, in Fig.~\ref{fig2_reactiondiffusion_pde} we study the convergence of the RBF method. It consistently achieves low relative $\ell^2$ errors compared to Galerkin, finite-difference, and neural network baselines. The estimated convergence order stabilizes around two for the RBF method, in agreement with the convergence rate of Galerkin, while other methods display larger fluctuations. 

\begin{figure}[H]
\begin{subfigure}[t]{0.5\textwidth}
\centering\includegraphics[width = 1.0\linewidth]{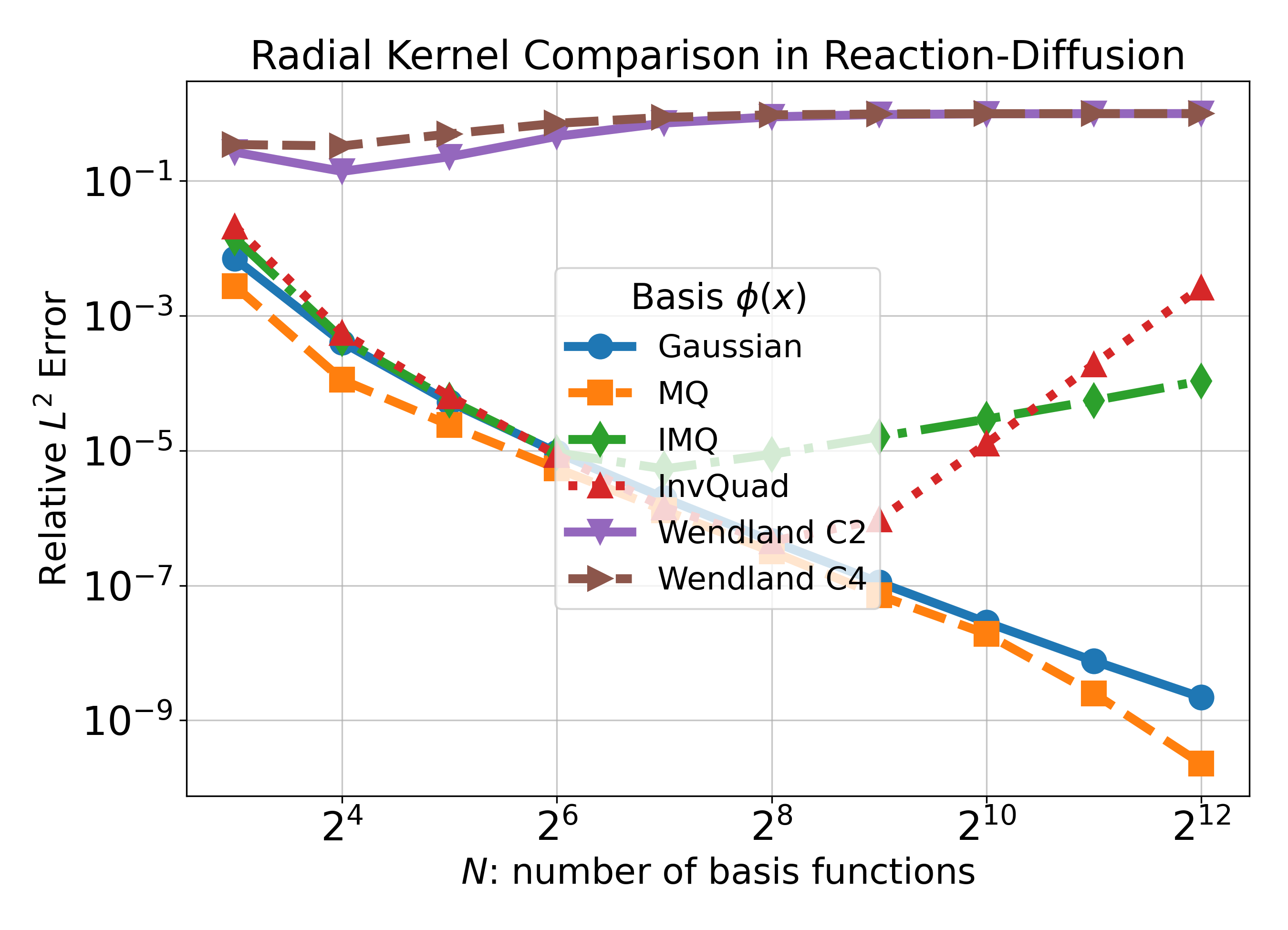}
\subcaption{Convergence vs.\ $N$ across RBF kernels.}
\end{subfigure}\hfill
\begin{subfigure}[t]{0.5\textwidth}
\centering\includegraphics[width = 1.0\linewidth]{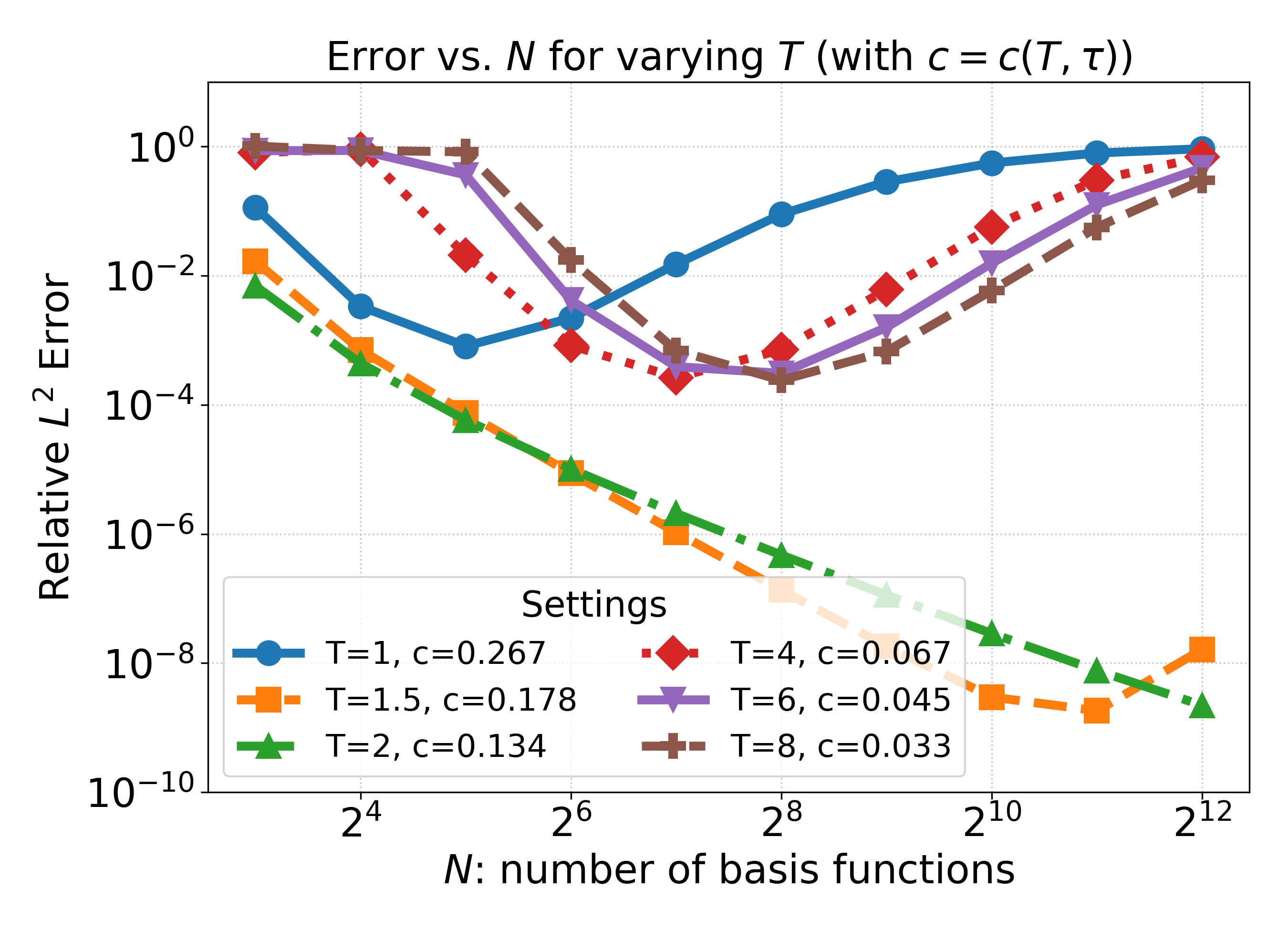}
\subcaption{Convergence vs.\ $N$ for varying $T$ and $c(T, \tau = 10^{-15})$.}
\end{subfigure}
\caption{Relative $\ell^2$ error vs.\ $N$ for Reaction–Diffusion: (a) Comparison across basis kernels using $T = 2.0, \tau = 10^{-15}$. (b). Effect of expansion factor $T$.}
\label{fig3_reactiondiffusion_pde}
\end{figure}

Fig.~\ref{fig3_reactiondiffusion_pde}, panel (a) compares the kernels in Table~\ref{tab:rbf-types}. It indicates that Gaussian, inverse quadratic, and multiquadric RBFs achieve the best accuracy, while compactly supported functions converge slowly. In panel (b), the effect of the domain expansion factor $T$ reveals that a moderate choice $T = 2$ provides good accuracy, whereas smaller or larger $T$ values may deteriorate accuracy. Additionally, the choice of the penalty parameter weight $\beta = 3 \times 10^{5}$ is justified as shown in panels (c)-(d) of Fig.~\ref{fig_beta_weight_poisson_reaction_diffusion_pde}. 

\begin{figure}[p]
  \centering
  \begin{subfigure}[t]{0.5\textwidth}
    \centering\includegraphics[width=1.0\linewidth]{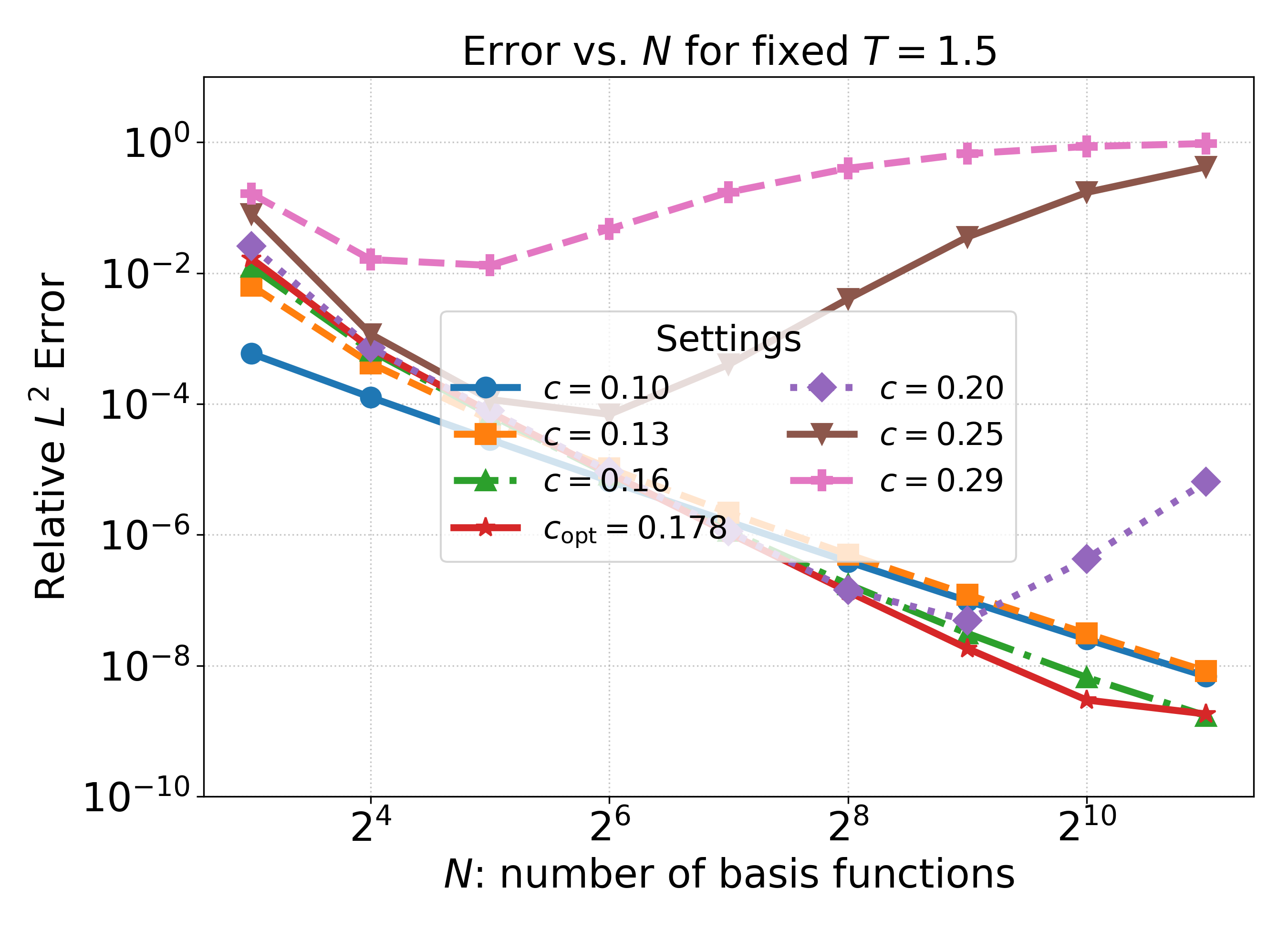}
  \end{subfigure}\hfill
  \begin{subfigure}[t]{0.5\textwidth}
    \centering\includegraphics[width=1.0\linewidth]{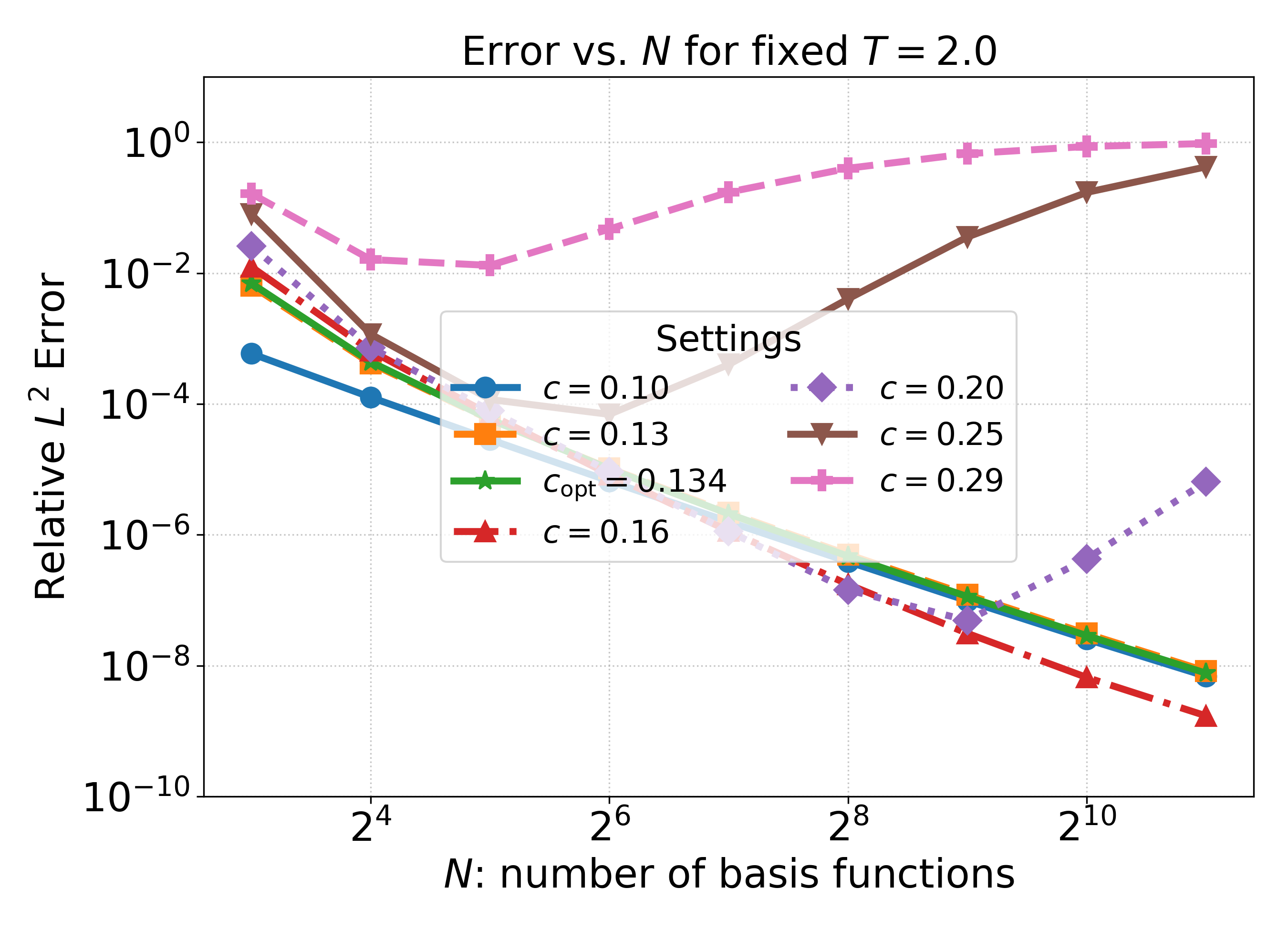}
  \end{subfigure}

  \vspace{0.1ex}

  \begin{subfigure}[t]{0.5\textwidth}
    \centering\includegraphics[width=1.0\linewidth]{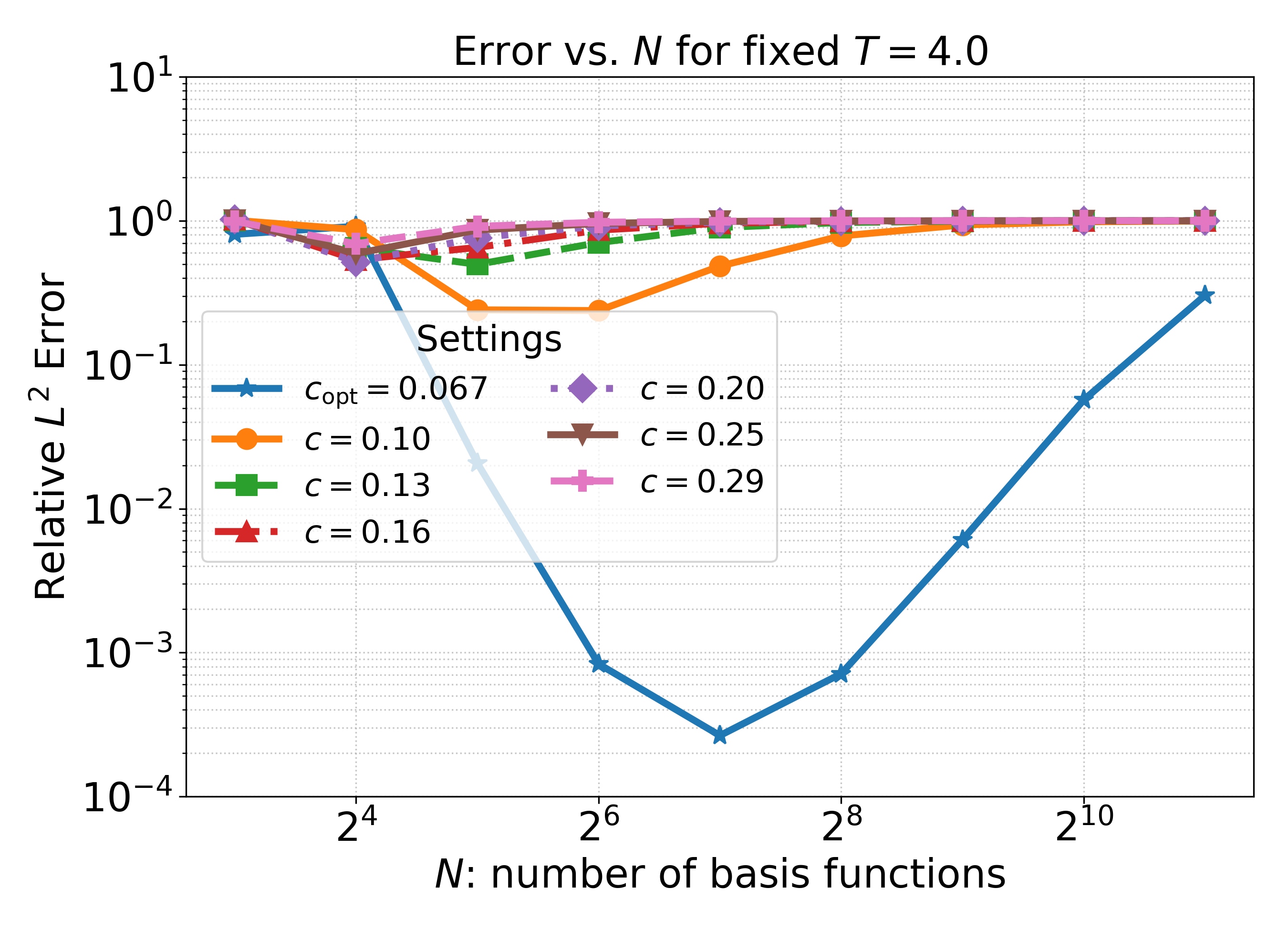}
  \end{subfigure}\hfill
  \begin{subfigure}[t]{0.5\textwidth}
    \centering\includegraphics[width=1.0\linewidth]{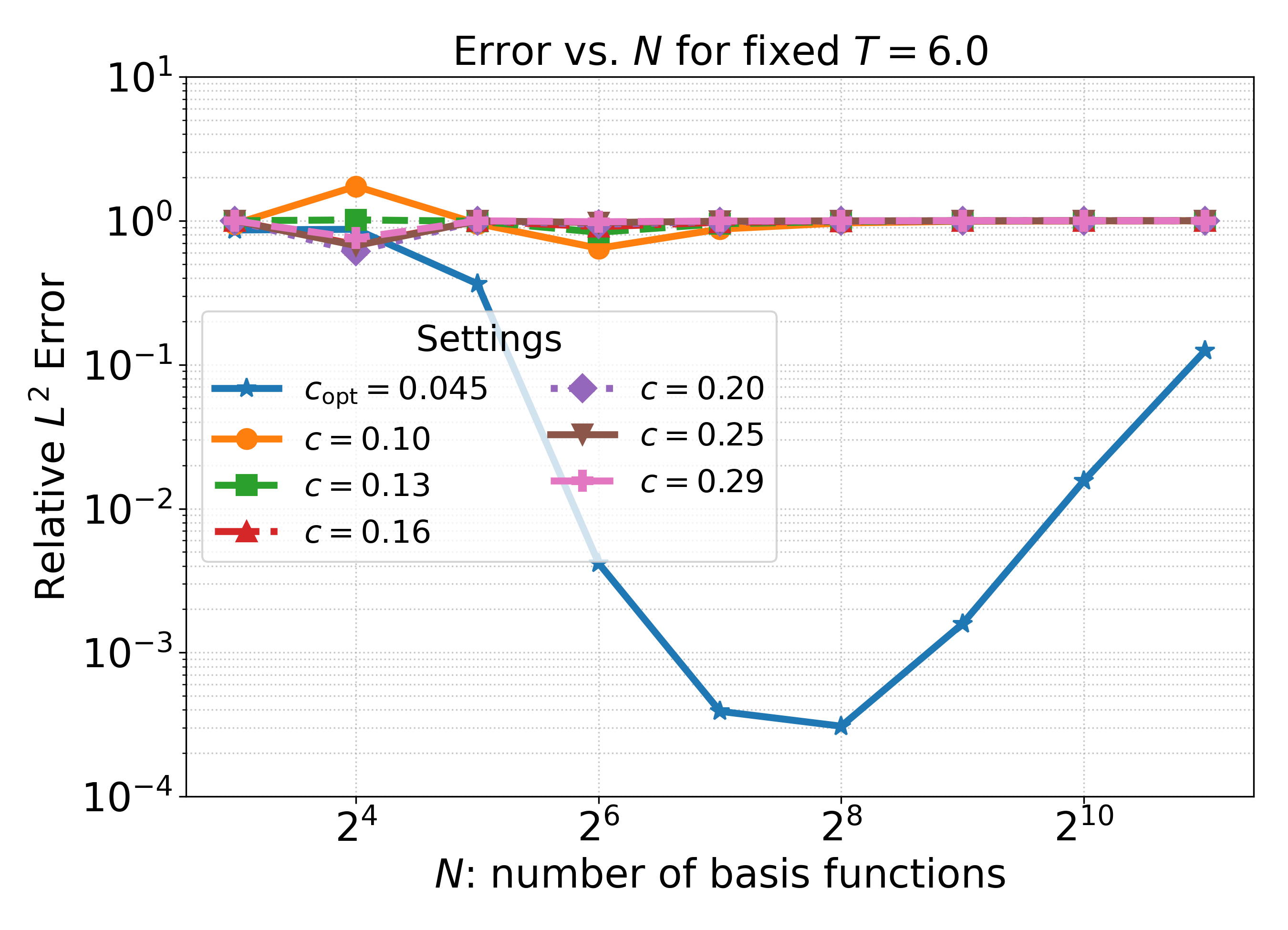}
  \end{subfigure}

  \caption{Reaction-diffusion problem. Relative $\ell^2$ error vs.\ $N$ for varying expansion factor $T$. Setup: $\tau = 10^{-15}$, $\zeta = 2$.}
  \label{fig4_reactiondiffusion_pde}
\end{figure}

According to Fig.~\ref{fig4_reactiondiffusion_pde}, for fixed expansion factors $T$ the error behavior versus $N$ depends strongly on the induced shape parameter $b = c N$. Overly large $T$ suppresses the benefits of increasing $N$, so the choice of a moderate value for $T$ maintains a balance between resolution and accuracy. 

\begin{figure}[p]
  \centering
  \begin{subfigure}[t]{0.5\textwidth}
    \centering\includegraphics[width=1.0\linewidth]{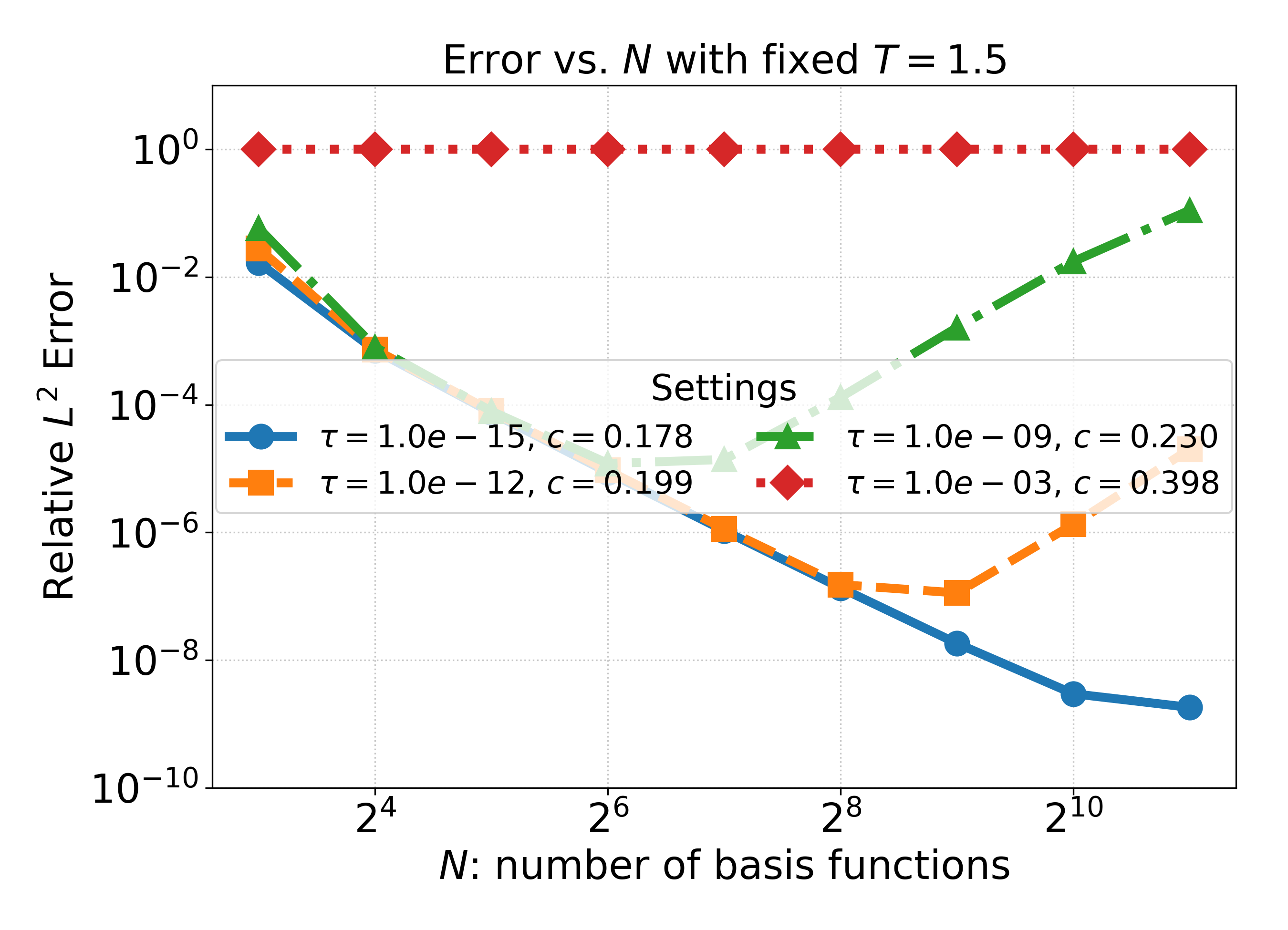}
  \end{subfigure}\hfill
  \begin{subfigure}[t]{0.5\textwidth}
    \centering\includegraphics[width=1.0\linewidth]{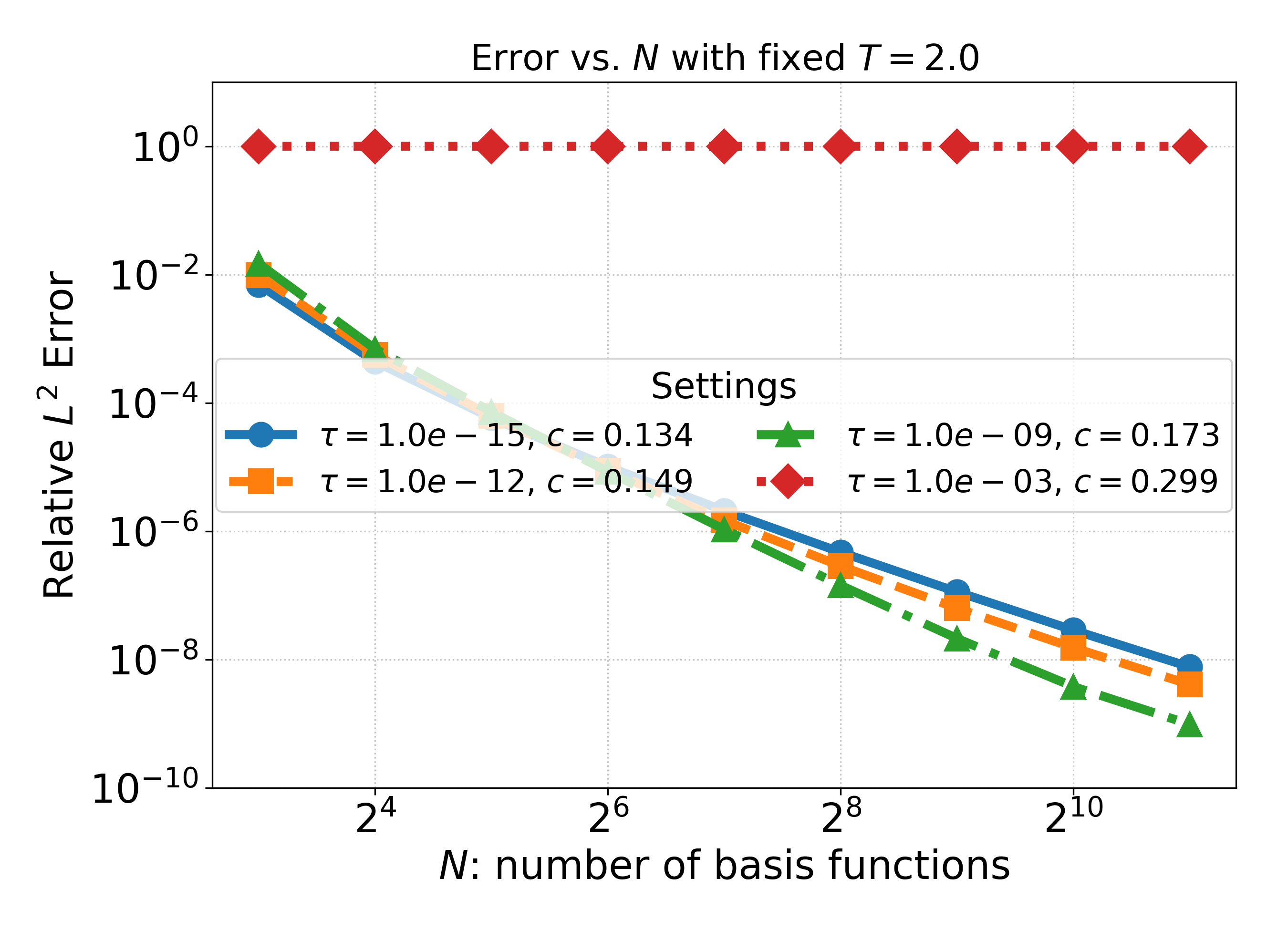}
  \end{subfigure}

  \vspace{0.1ex}

  \begin{subfigure}[t]{0.5\textwidth}
    \centering\includegraphics[width=1.0\linewidth]{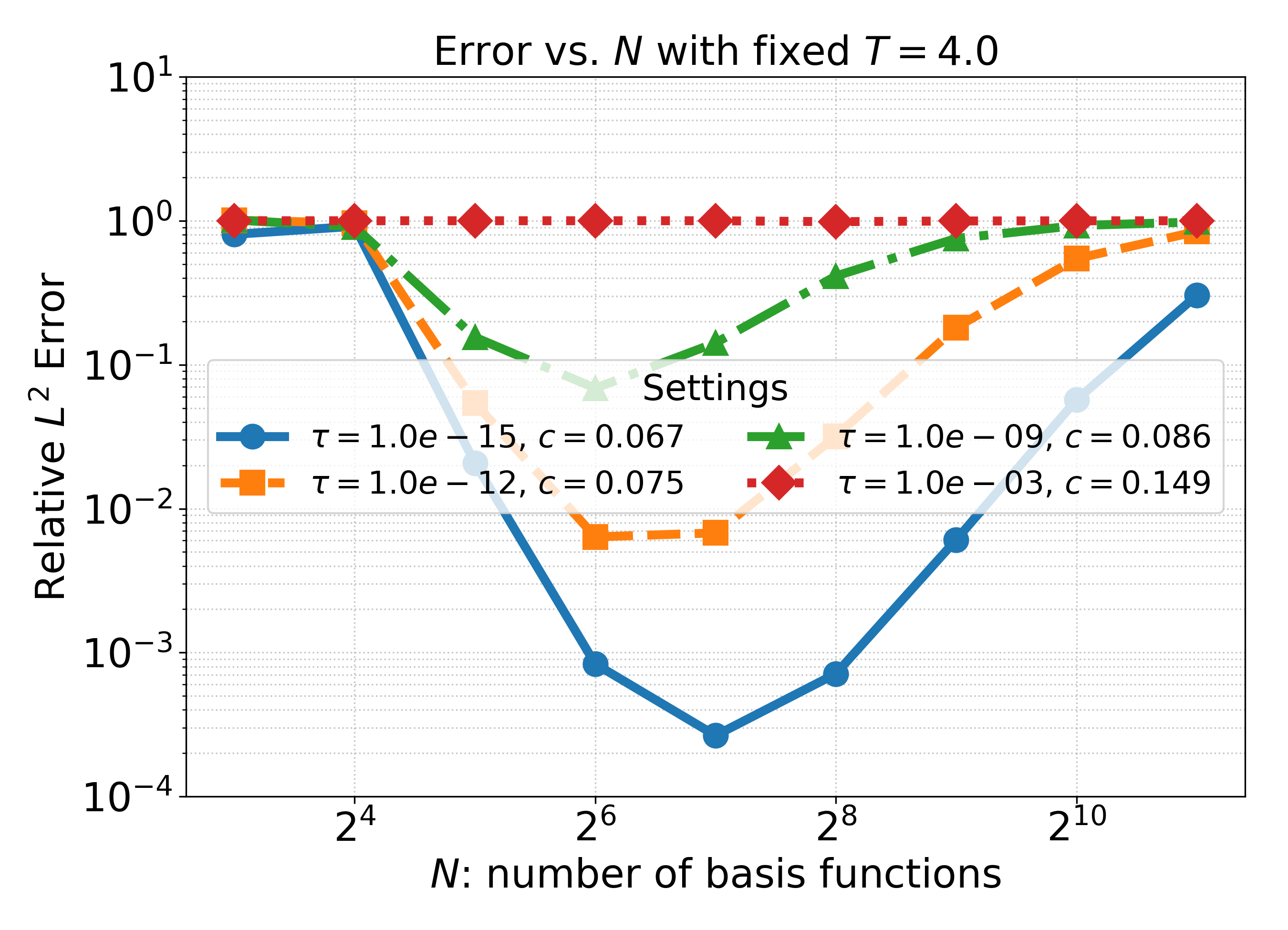}
  \end{subfigure}\hfill
  \begin{subfigure}[t]{0.5\textwidth}
    \centering\includegraphics[width=1.0\linewidth]{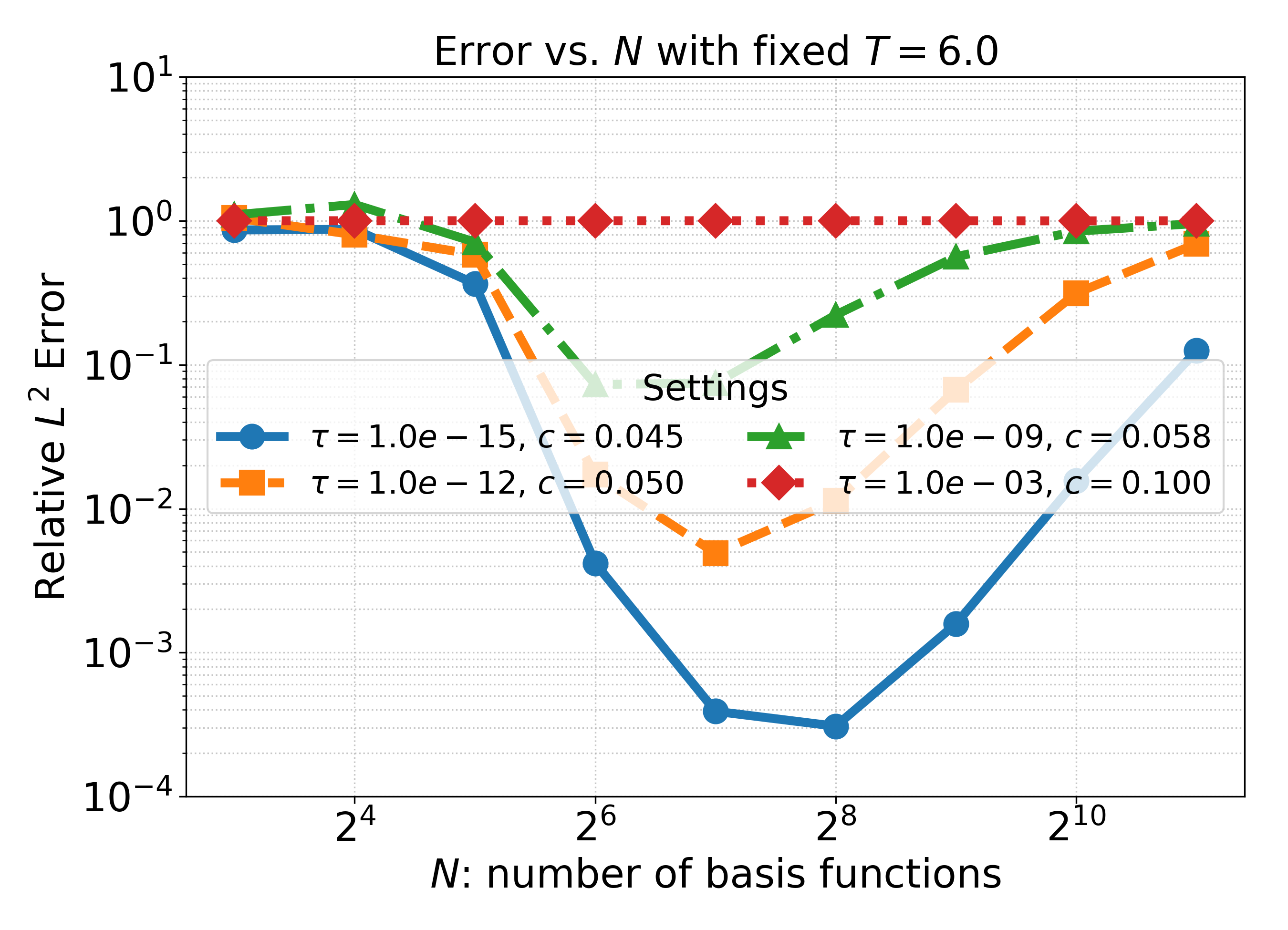}
  \end{subfigure}

  \caption{This figure is similar to Fig.~\ref{fig4_reactiondiffusion_pde}. Relative $\ell^2$ error vs.\ $N$ basis functions: we vary different threshold values $\tau$ and for each expansion factor $T$ we have the corresponding $c = c(T, \tau)$.}
  \label{fig5_reactiondiffusion_pde}
\end{figure}

\begin{figure}[p]
    \centering
    \begin{subfigure}[t]{0.5\textwidth}
    \centering\includegraphics[width = 1.0\linewidth]{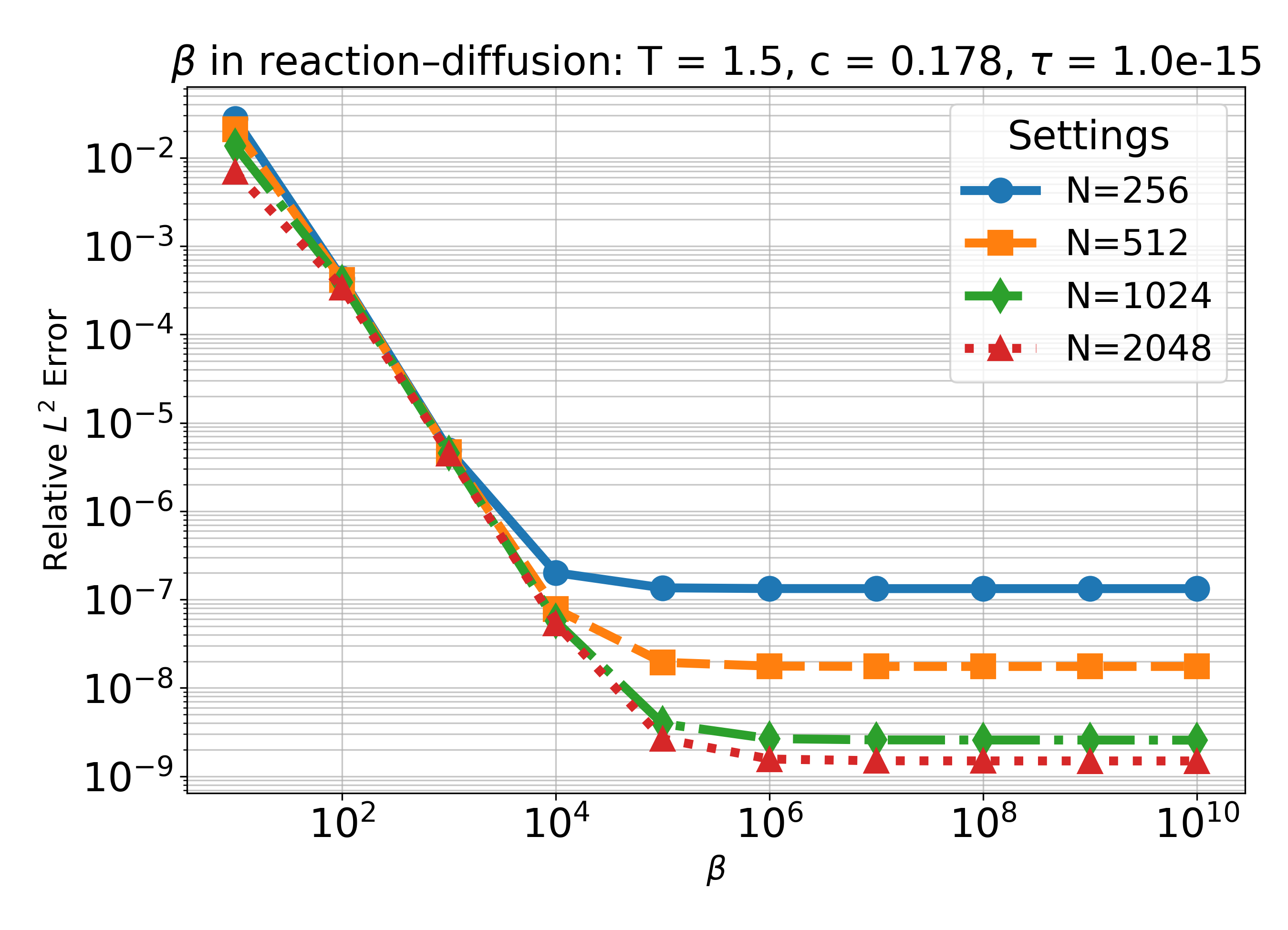}
    \caption{$T = 1.5, \tau =  10^{-15}, c = 0.178$}
    \end{subfigure}\hfill 
    \begin{subfigure}[t]{0.5\textwidth}
    \centering\includegraphics[width = 1.0\linewidth]{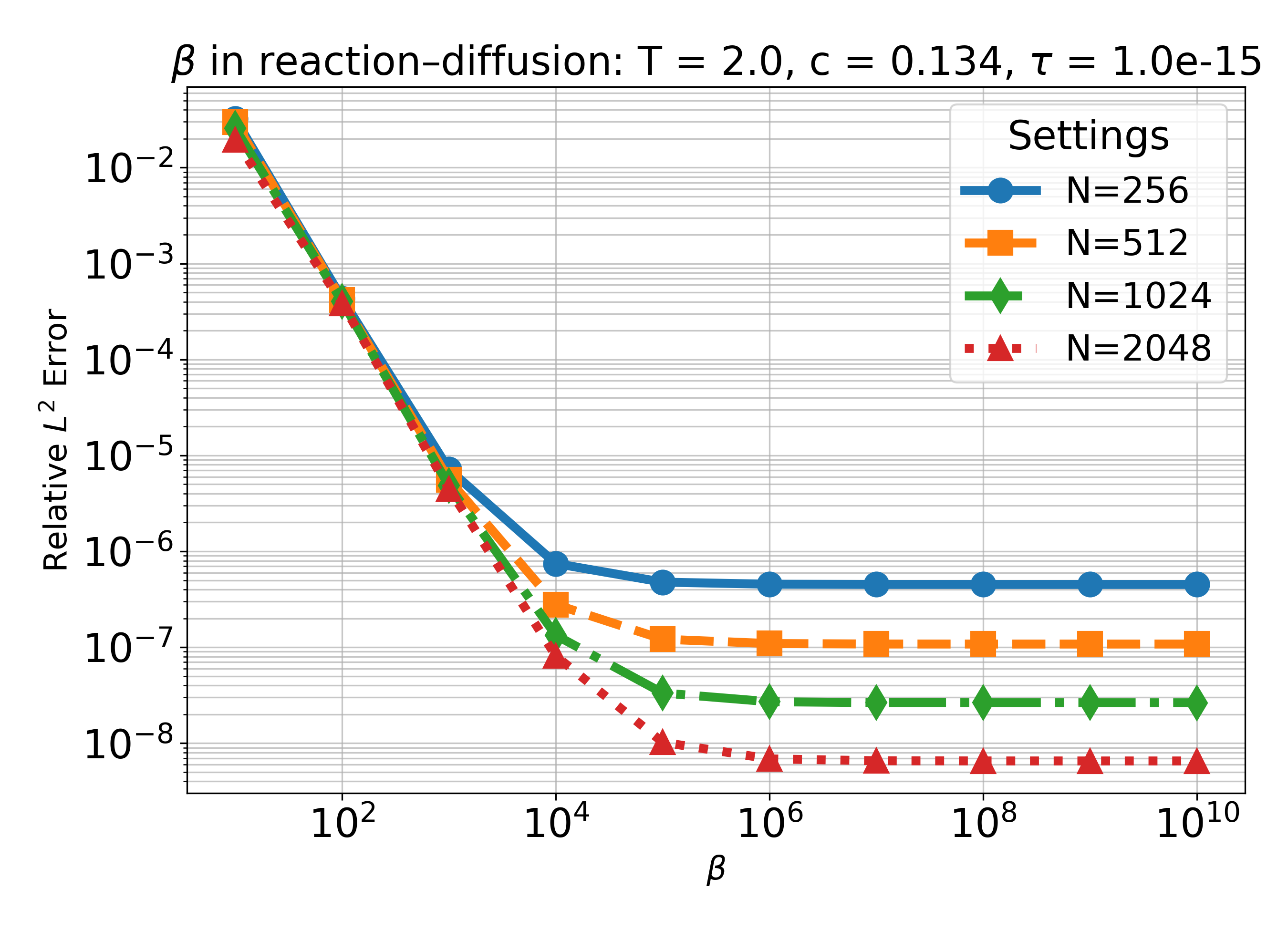}
    \caption{$T = 2.0, \tau = 1\times 10^{-15}, c = 0.134$.}
    \end{subfigure}
    \caption{\textbf{$\beta$ penalty parameter for boundary condition for Reaction-Diffusion PDE $-\Delta u + u = f$.}
    (a-b) Relative $\ell^{2}$ error vs.\ $\beta$ for different values of $N$ at $T = 1.5, 2.0$.
    }
\label{fig_beta_weight_poisson_reaction_diffusion_pde}
\end{figure}

In Fig.~\ref{fig5_reactiondiffusion_pde}, we vary the threshold $\tau$. Small thresholds ($\tau = 10^{-15}$) allow larger $b$ and yield lower errors at moderate $T$. However, larger thresholds ($\tau = 10^{-9}$) truncate early and flatten convergence. 

\pagebreak

\pagebreak

\end{document}